\documentclass[11pt]{elsarticle}
\usepackage{lineno,hyperref, booktabs}
\usepackage{mathrsfs, amsfonts,amssymb, amsmath,  amsthm, bm, bbm}
\usepackage{comment, xcolor, setspace, lipsum, geometry}
\usepackage{threeparttable, longtable, tabularx, multirow, array}
\usepackage{algorithmic, algorithm}
\usepackage{subcaption, float}

 
\newlength\myindent
\theoremstyle{plain}

\newtheorem{proposition}[]{Proposition}
\hypersetup{pdfauthor=author}
\usepackage[font=small,labelfont=bf,labelsep=period, singlelinecheck=false]{caption}

\usepackage{titlesec}

\modulolinenumbers[5]
\biboptions{authoryear}

    \journal{arXiv}

\begin{document}

\begin{frontmatter}
\title{Integrated Home Care Staffing and Capacity Planning: Stochastic Optimization Approaches}
\author[address_1]{Ridong Wang\corref{firstauthor}}
\address[address_1]{Department of Industrial Engineering, Tsinghua University, Beijing, China}

\author[address_2]{Karmel S. Shehadeh\corref{secondgauthor}}
\address[address_2]{Department of Industrial and Systems Engineering, Lehigh University, Bethlehem, PA, USA}

\author[address_1]{Xiaolei Xie\corref{thirdauthor}}

\author[address_1]{Lefei Li\corref{fourthauthor}}
\begin{abstract}

\noindent We propose stochastic optimization methodologies for a staffing and capacity planning problem arising from home care practice. Specifically, we consider the perspective of a home care agency that must decide the number of caregivers to hire (staffing) and the allocation of hired caregivers to different types of services (capacity planning) in each day within a specified planning horizon. The objective is to minimize the total cost associated with staffing (i.e., employment), capacity allocation, over-staffing, and under-staffing. We propose two-stage stochastic programming (SP) and distributionally robust optimization (DRO) approaches to model and solve this problem considering two types of decision-makers, namely an everything in advance decision-maker (EA) and a flexible adjustment decision-maker (FA). In the EA models, we determine the staffing and capacity allocation decisions in the first stage before observing the demand. In the FA models, we decide the staffing decisions in the first stage. Then, we determine the capacity allocation decisions based on demand realizations in the second stage. We derive equivalent mixed-integer linear programming (MILP) reformulations of the proposed nonlinear DRO model for the EA decision-maker that can be implemented and efficiently solved using off-the-shelf optimization software. We propose a computationally efficient column-and-constraint generation algorithm with valid inequalities to solve the proposed DRO model for the FA decision-maker. Finally, we conduct extensive numerical experiments comparing the operational and computational performance of the proposed approaches and discuss insights and implications for home care staffing and capacity planning.

\end{abstract}


\end{frontmatter}

\section{Introduction}\label{introduction}
Home care is the provision of providing services to people at their homes,  such as nursing, wound treatment, personal care assistance, and lifestyle support for the elderly. Home care can improve the quality of life for customers, especially elderly customers, as caregivers meet their needs at their homes. Moreover, it leads to cost-saving for the entire health system by mitigating the need for hospitalization, among others \citep{lanzarone2014robust, restrepo2020home}. The home care service industry has been rapidly growing worldwide due to the aging population, outspread of chronic and infectious diseases, emergent changes in family structures, work obligations, and extended work hours. In particular, COVID-19 pandemic has led to unprecedented increases in the demand for home healthcare, in particular to support vulnerable and high-risk populations. In Canada, it's estimated that the number of Canadians with home care needs will increase by 615,479 to just over 1.7 million by 2031 \citep{canadia}. As pointed out by \cite{zhan2021home}, the United States spends about \$82 billion on various forms of home healthcare services and \$300 billion on home repairs and maintenance, and these numbers are rapidly growing. In 2018, the global home service market was valued at around \$282 billion and is expected to reach \$1133.4 billion by 2026 \citep{tsang2021distributionally}. 

The demand for home care is associated with increased demand for cross-trained caregivers who can provide multiple service types. \cite{NDWA} shows that the demand for home care will create more than one million new home care jobs. Unfortunately, however, the home care industry faces a significant shortage of caregivers. In addition, \cite{Wash} recently reported that labor shortages would worsen in the next twenty years in Washington State. The state estimates that nearly 77,000 additional home care aides will be needed to by 2030. When the high turnover rate is taken into account, the number of additional home care aides needed surges to 125,000. 

The increase in demand for home care is also associated with an increase in service costs.  According to the newly released Genworth cost of care survey \citep{genworth}, the cost of home health aides, who provide personal assistance with activities such as bathing, dressing and eating, grew by 12.5 \% in 2020. Accordingly, the development of computationally efficient optimization tools for homecare staffing and capacity allocation is essential to support decision-making in all areas of the home care industry to meet the growing demand, improve service quality, and reduce costs. 

In this paper, we propose stochastic optimization methodologies for a staffing and capacity planning problem arising from home care practice. Specifically, we consider the perspective of a home care agency that must decide the number of (\textit{cross-trained}) caregivers to hire (\textit{staffing})  within a specified planning horizon and the allocation of hired caregivers to different types of services in each day in the planning horizon (\textit{capacity allocation}). Customers demand and durations for different types of services are random. The distribution of these random parameters may be unknown. The objective is to find employment/staffing and capacity allocation decisions that minimize total employment cost, capacity allocation cost, and penalty cost for daily over-staffing and under-staffing. We call this problem the \textit{\underline{h}ome care \underline{s}taffing and \underline{c}apacity planning (HSCP}) problem. 

Unfortunately, the HSCP problem is a challenging stochastic optimization problem for several reasons. First, it is a complex combinatorial optimization problem that requires deciding the number of different (cross-trained) caregivers to be hired and allocating (i.e., assigning) the capacity of hired caregivers to different types of services simultaneously. Second, there is significant variability in demand for home care services and service duration. In addition, our analysis of real-life data from a home care service provider suggest that there is a wide range of possible probability distributions for modeling the variability in the demand, indicating distributional \textit{ambiguity}, i.e., uncertainty in probability distribution (see this analysis in Section~\ref{motivation}). Such uncertainty and distributional ambiguity increase the complexity of modeling and solving the HSCP problem. However, ignoring uncertainty may lead to sub-optimal decisions and the inability to meet customer demand and high costs, among other disappointing consequences of adopting deterministic solutions. Failure to meet customer demand, in particular, may lead to adverse outcomes, especially in healthcare, as it impacts population health. It also impacts customers' satisfaction and thus the reputation of the home service providers. 

Motivated by these important issues and motivated by our collaboration with a home service provider in Beijing, in this paper, we propose new stochastic optimization approaches for modeling and solving the HSCP problem under uncertainty.  Specifically, given sets of service types, days in the planning horizon, and caregivers and their types (i.e., training), our models solves the following decisions problems simultaneously (a) a \textit{staffing problem} that determines the number of caregivers to hire within the planning horizon, and (b) an \textit{allocation problem} that determines the allocation of caregivers capacities (i.e., daily service hours) to service types during the planning horizon. We consider two types of decision-makers:\textit{ \underline{E}verything in \underline{A}dvance} (EA) and \textit{\underline{F}lexible \underline{A}djustment (FA)} decision-makers. The (EA) decision-maker decides the number of caregivers to hire and their daily capacity allocation to service types at the start of the planning horizon before the demand and service time are realized. In contrast,  the (FA) decision-maker decides the number of caregivers to hire at the beginning of the planning horizon. Then, s/he allocates the hired caregivers to service types once the demand is realized. 

To address uncertainty, we propose and analyze two-stage stochastic programming  (SP) and distributionally robust optimization (DRO) models for both decision-makers, assuming known and ambiguous distribution, respectively.  We obtain near-optimal solutions of our SP models using sample average approximation (SAA). We derive a mixed-integer linear programming (MILP) reformulation of our proposed DRO model for the (EA) decision-maker that can be implemented and efficiently solved using off-the-shelf optimization software. We propose a computationally efficient column-and-constraint generation (C\&CG) method to solve our proposed DRO model for the (FA) decision-maker. In addition, we derive valid inequalities to accelerate the C\&CG convergence.

We conduct extensive computational experiments comparing the proposed methodologies' computational and operational performances demonstrating where significant performance improvements can be gained.  Additionally, we derive several managerial insights relevant to practice. To the best of our knowledge, and according to our literature review in Section \ref{LR}, our paper is the first to propose stochastic optimization models and methods to model and solve the HSCP problem under uncertainty considering both (EA) and (FA) decision-makers. 

The remainder of this paper is organized as follows. In Section \ref{motivation}, we present a motivating example for modeling uncertainty and distributional ambiguity. In Section \ref{LR}, we review the relevant literature. In Section~\ref{First_model} we formally define the HSCP problem. In Sections~\ref{EA}-\ref{wait-and-see}, we present our proposed SP and DRO models.  In Section \ref{solution_approach}, we present our C\&CG algorithm and strategies to speed convergence. In section \ref{CE}, we present our computational results and managerial implications. Finally, we draw conclusions and discuss future directions in Section \ref{conclusion}.

\section{Motivations}\label{motivation}

This paper is based on our research with Pinetree Care, one of China's leading home care agencies. Pinetree Care provides more than 5 million care sessions to over 600,000 families through a network of more than 10,000 caregivers. One of the main challenges that Pinetree Care and other home care companies face when making staffing and capacity allocation decisions is the uncertainty in service demand and duration. To analyze this uncertainty, we use data from Pinetree Care collected electronically from June 2017 through December 2019, representing $\sim$205 K care sessions in Beijing. 

\begin{figure}[t!]
    \centering
    \subcaptionbox{Rehabilitation nursing\label{pdf_rn}}{
        \includegraphics[scale=0.37]{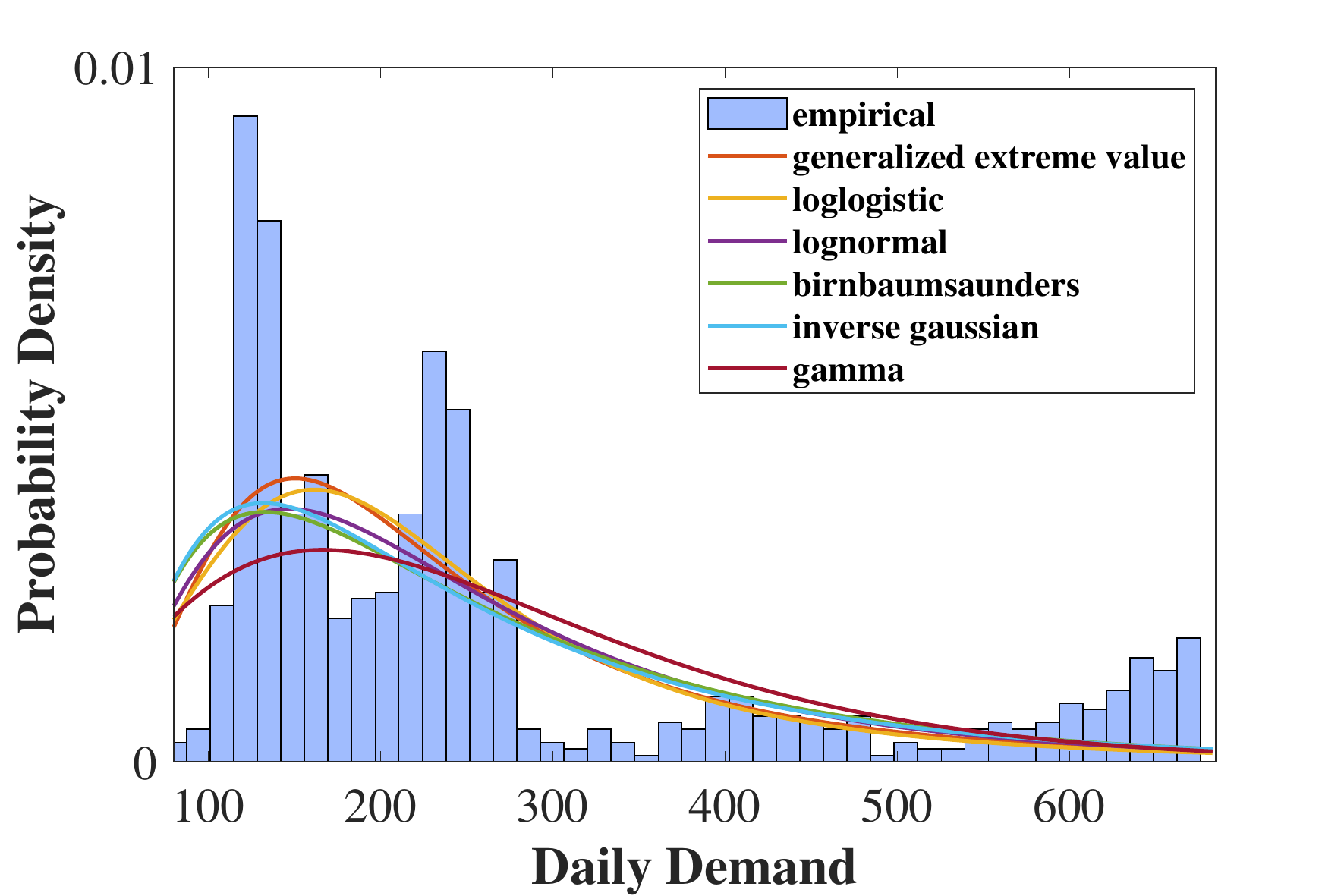}}
    \subcaptionbox{Health status assessment\label{pdf_ha}}{
        \includegraphics[scale=0.37]{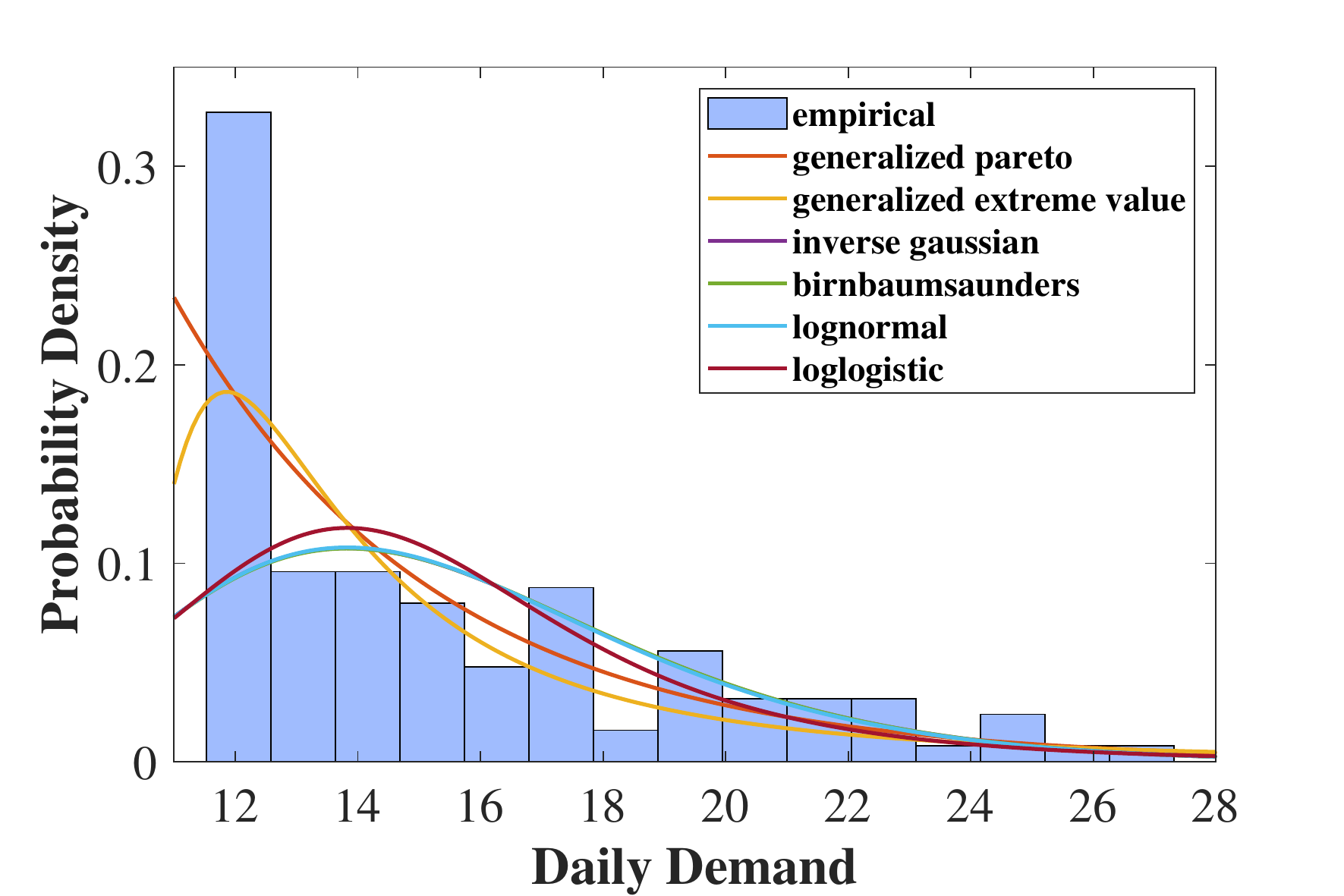}}
    \subcaptionbox{Rehabilitation nursing (fourth quarter of 2017)\label{pdf_rn_2017}}{
        \includegraphics[scale=0.37]{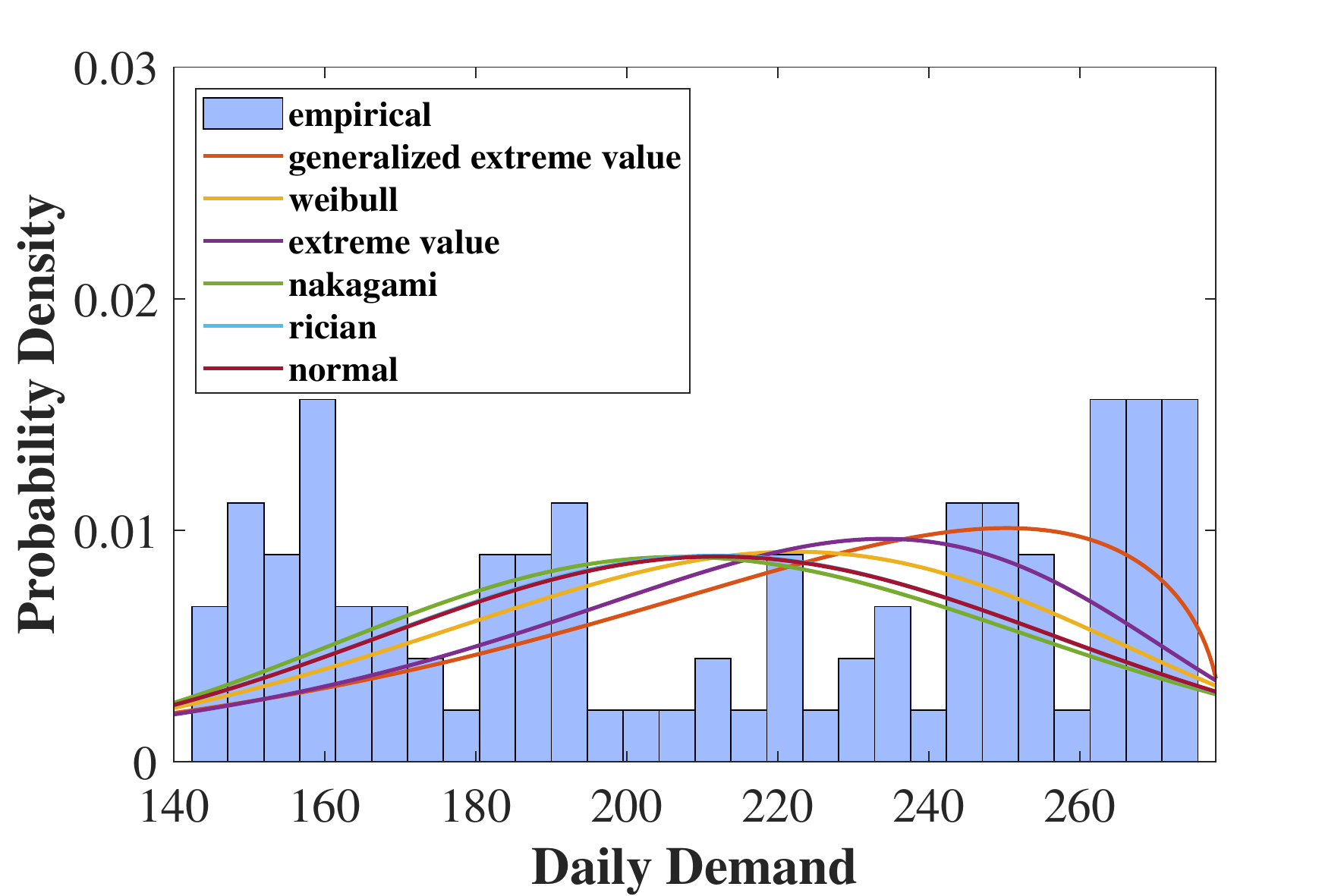}}
    \subcaptionbox{Rehabilitation nursing (fourth quarter of 2018)\label{pdf_rn_2018}}{
        \includegraphics[scale=0.37]{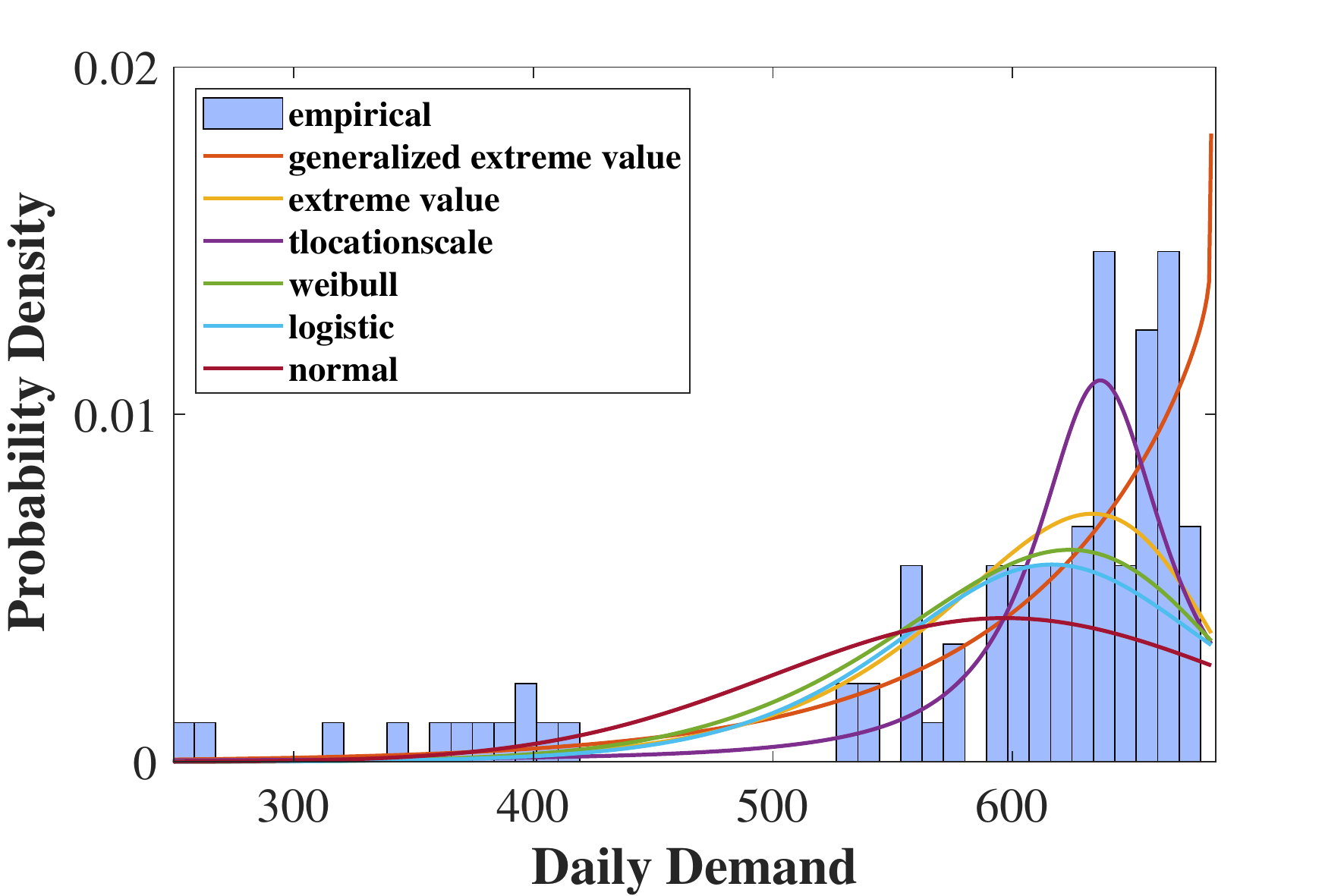}}
    \caption{The empirical and fitted probability distributions for actual demands.}
    \label{PDF}
\end{figure}

In the top panel of Figure~\ref{PDF}, we present the empirical and fitted probability distributions for actual demands for rehabilitation nursing and health status assessment, and in \ref{CPD} we present the associated the goodness-of-fit metrics (e.g., Negative of the Log-Likelihood (NLogL), Akaike Information Criterion (AIC), and Bayesian Information Criterion (BIC)). In the bottom panel of Figure ~\ref{PDF}, we present demand distributions during the fourth quarter of 2017 and the fourth quarter of 2018. 

We make the following observations from these results. First, Figures~\ref{pdf_rn}-\ref{pdf_ha} show significant variability in demand for the same service type and across service types.  Furthermore, these figures, along with goodness-of-fit metrics, demonstrate that a wide range of distributions can represent uncertainty in demand, suggesting distributional ambiguity (i.e,. uncertainty in distribution). Figures \ref{pdf_rn_2017} --\ref{pdf_rn_2018} show that demand variability and distributions are quite different between the fourth quarters of 2017 and 2018 (see also Figure \ref{dr_ts} in \ref{CPD}). Indeed, future uncertainty is often not distributed in the past.  These results support prior numerical studies that show that different distributions can typically explain raw data of uncertain parameters, indicating distributional ambiguity. Finally, these results motivate us to adopt a DRO approach to obtain robust decisions that could hedge against distributional ambiguity.

\section{Literature Review}\label{LR}
Operations research in home care service can be divide into three decision levels: strategical level, tactical level and operational level \citep{restrepo2020home}. In this section, we focus primarily on the literature in tactical level which our problem lies in, i.e., papers that apply stochastic optimization to address HSCP problems similar to ours. For comprehensive recent surveys of operations research methods applied to decisions in home health care, we refer readers to \cite{gutierrez2013home}, \cite{fikar2017home}, \cite{grieco2021operational} and \cite{di2021routing}. We also refer to \cite{michael2018scheduling} for a comprehensive survey of scheduling theory, applications and methods.

Existing stochastic home care service staffing and capacity planning models often consider the demand uncertainty only. SP approaches for HSCP include \cite{rodriguez2015staff}, \cite{restrepo2020home} and \cite{zheng2021stochastic}. \cite{rodriguez2015staff} proposed a two-stage stochastic
programming (SP) approach to solve the staff dimensioning problem
in home care services considering random demand. In their two-stage approach, the first stage calculates (near-)optimal levels of resources for possible demand scenarios, while the second stage computes the optimal number of caregiver for each profession. Notably, \cite{restrepo2020home} is the first to study integrated staffing and scheduling problem in the context of home healthcare. They considered staffing and scheduling decisions in the first stage and the temporary reallocation of caregivers to neighboring districts in the second stage. Their study show that the total cost presents an increase when there is less flexibility associated with the allocation of schedules.

The recent work of \cite{zheng2021stochastic} is closely related to ours. Specifically, \cite{zheng2021stochastic} considered a multi-resource and multi-demand network matching problem in which they assume deterministic service duration and that the demands of services follow a fully known probability distribution. Accordingly, \cite{zheng2021stochastic} proposed a two-stage stochastic MILP which seeks the optimal service authorization and capacity planning decision in first stage, then the service resource allocation in the second stage. In their model, the goal is to maximize the total expected revenue netted by the authorization cost, payment and the penalty cost of service shortage. However, they do not incorporate the penalty cost of over-staffing in the objective, which is an important source of operational expenses \citep{restrepo2020home}. In addition, ignoring random service duration may lead to sub-optimal solutions.

These SP models assume that the decision-makers are risk-neutral and have complete knowledge about the underlying uncertainty through a known probability distribution \citep{shapiro2014lectures}. Hence, the applicability of the SP approach is limited to the case in which the decision-makers have sufficient data to model the distribution of random parameters which is implausible in practice. Given the distributional ambiguity (see our analysis in Section \ref{motivation}), if we calibrate a model to a misspecified distribution, the resulting optimal SP decisions may have a disappointing out-of-sample performance under the true distribution which is termed the \emph{optimizer’s curse} \citep{smith2006optimizer}.

Robust optimization (RO) is alternative technique to model, analyze and optimize decisions under uncertainty and ambiguity (where the distributions are unknown). It assumes that the decision maker is risk-averse and has no distributional knowledge about the underlying uncertainty, except for its support (partial information), and the model minimizes the worst-case cost over an uncertainty set \citep{ben2009robust}. In home care service area, several studies propose robust approaches to address the uncertainty \citep{carello2014cardinality, cappanera2018demand, shi2019robust, cappanera2021addressing}. Notably, \cite{cappanera2018demand} addressed uncertainty of patient demand over a multiple-day time horizon and jointly studied the assignment, scheduling and routing decisions via a non-standard cardinality-constrained robust approach. \cite{shi2019robust} proposed a robust model for both uncertain service and travel time and discuss the heuristic solution approaches to find the optimal route for each caregiver and appointment time for each patient.

Although RO is a powerful technique in dealing with uncertainty in optimization, its solutions can be too conservative \citep{thiele2010note,roos2020reducing}. Additionally, RO does not fully utilize any distributional information of uncertainty. Hence, RO may yield poor expected performance and sub-optimal decisions for other more-likely scenarios. For comprehensive recent surveys of robust optimization, we refer to \cite{bertsimas2011theory}, \cite{gabrel2014recent} and \cite{gorissen2015practical}.

Distributionally robust optimization (DRO) is another approach for modeling uncertainty that bridges the gap between the conservatism of robust optimization and the specificity of stochastic programming \citep{goh2010distributionally}. In DRO, we assume only partial distributional information is known. Then, we model the distribution of uncertainty as a decision variable that belongs to an ambiguity set (i.e., a family of distributions characterized through certain known properties of the unknown data-generating distribution, \citealp{esfahani2018data}). Optimization is based on the worst-case distribution within the ambiguity set. And we can use easy to compute information such as the mean and range of random parameters to construct the ambiguity sets and build DRO models. We refer to \cite{rahimian2019distributionally} for a comprehensive survey  on DRO techniques.

DRO has the following striking benefits. First, it adopts a worst-case approach regularizing the optimization problem and thereby mitigating the optimizer’s curse characteristic for SP. Second, DRO avoids the well-known over-conservatism and the poor expected performance of RO. Finally, DRO models are often more computationally tractable than SP and RO \citep{esfahani2018data}. We observe that despite the potential advantages, only one study, \cite{tsang2021distributionally}, using DRO approach to hedge against uncertainty in home care service area. However, their study focus on the routing and appointment scheduling problem for a single-provider. To the best of our knowledge, there is no DRO approach for the specific HSCP problem that we address in this paper.

In this work, we develop two DRO models to consider the realistic lack of distributional information and compare these models to SP models. We compare our work with \cite{zheng2021stochastic} and \cite{restrepo2020home}, which are recent studies relevant to our work. First, the three papers both address stochastic demand. However, \cite{zheng2021stochastic} and \cite{restrepo2020home} assumed that the service time is deterministic and do not address the potential distributional ambiguity of demand which is captured by real world data (see our analysis in Section \ref{motivation}). As mentioned earlier, this may lead to sub-optimal and unrealistic solutions. In contrast to them, we extend the considered problem by incorporating uncertainty and distributional ambiguity for demand and service duration. Second, we propose different SP and DRO models for two types of decision-makers (EA and FA). Third, \cite{zheng2021stochastic} and \cite{restrepo2020home} assumed that the caregivers have single skill. In contrast to these papers, we consider cross-trained caregivers representing flexible capacity. To the best of our knowledge, our work is the first to consider cross-training in the HSCP problem.

\section{Problem Description}\label{First_model}

Let us now start to introduce our HSCP problem. Suppose that there is a home care company that provides $L$ different types of home care services (e.g., rehabilitation nursing, health status assessment, etc.). The decision-maker of this company wants to determine the number of caregivers to hire within an arbitrary planning horizon of $T$ days to provide home care services to customers.  Suppose that there are $K$ types of caregivers to hire. Depending on their training, skills, and qualifications, each caregiver of type $k \in [K]:=\{1, 2, \ldots, K\}$, can provide a subset $R_{k} \subset [L]$ of services. In addition, the daily service capacity for type $k \in [K]$ caregivers is $h_{k}$. The fixed cost of hiring type $k$ caregiver is $c_{k}$.  

The demand, $d_{l,t}$, for service type $l\in [L]$ in each day $t\in [T]$ is random. Service time, $s_{l,t}$, for each service type $l\in [L]$ in each day $t\in [T]$ is also random. Due to the uncertainty of demand and service time, one or multiple of the following scenario can happen: (1) over-staffing, or and, (2) under-staffing (i.e., failure to meet customers demand). The over-staffing and under-staffing unitary penalty costs are respectively $c_{l,t}^{o}$ and $c_{l,t}^{u}$,  for all $l\in [L]$  and $t \in [T]$. A complete listing of the problem parameters can be found in Table~\ref{notation}.

Given $T$, $K$, and $L$, our models determine the number of caregivers to hire and the capacity of type $k$ caregivers allocated to type $l\in [L]$ in each day $t \in [T]$. In particular, we consider two types of decision-makers. The first decides the number of caregivers to hire of each type and their daily capacity allocation to service types at the beginning of the planning horizon before the demand and service time are realized. We call this decision-maker \underline{E}verything in \underline{A}dvance (EA) decision-maker. The second decision-maker decides the number of caregivers to hire before realizing the demand (i.e., at the beginning of the planning horizon). S/he then allocates these caregivers to service types once the demand is realized. We call this decision maker \underline{F}lexible \underline{A}djustment (FA) decision-maker.

In the next sections, we propose two-stage SP and DRO models for (EA) and (FA), considering the case of known and unknown distributions of demand and service time, respectively.  In the first stage of the (EA) model, we decide the number of caregivers to hire from each type $k \in [K]$ and their capacity allocation to service type $l \in [L]$ in each day $t \in [T]$. And in the second stage, after the demand and service time are realized, we compute over-staffing and under-staffing. The objective is to minimize the total costs: fixed staffing (or hiring) and allocation (assignment) costs and a weighted sum of random costs related to over-staffing and under-staffing. In the first stage of the (FA) model, we decide the number of caregivers to hire, and in the second stage we decide capacity allocation and compute over-staffing and under-staffing. The objective of this model is to minimize the total costs: fixed hiring cost and a weighted sum of the random cost related to capacity  allocation (assignment) cost, over-staffing cost, and under staffing cost.  

\noindent \textbf{Additional notation}: For $a\in \mathbb{Z}$, we define $[a] := \{1,2,...,a\}$. For $a, b\in \mathbb{Z}$, $U[a; b]$ refers to uniform distribution over the interval $[a; b]$. We use boldface notation to denote vectors, e.g., $\bm{d}:=[d_{1}, d_{2},...,d_{N}]^{T}$. Table \ref{notation} summarizes other notations.

\begin{table}[t!]  
\small
\center
\renewcommand{\arraystretch}{0.8}
\caption{Notation.}
\begin{tabular}{ll}
\hline
\multicolumn{2}{l}{\textbf{Indices}}\\
    $l$ & index of service types, $l \in [L]$\\
    $t$ & index of days, $t \in [T]$\\
    $k$ & index of caregivers types, $k \in [K]$\\
    \multicolumn{2}{l}{\textbf{Parameters and sets}}\\
    $L$ & number of service types\\
    $T$ & number of days in the planning horizon\\
    $K$ & number of caregiver types\\
    $R_{k}$ & the skills set of type $k$ caregivers\\
    $c_{k}$ & the fixed cost for hire type $k$ full-time caregivers \\
    $c_{k,l,t}$ & the unitary service cost of type $k$ full-time caregivers for serve type $l$ customers in day $t$\\
    $c_{l,t}^{u}$ & the unitary penalty cost due to shortage of service capacity of type $l$ service\\
    $c_{l,t}^{o}$ & the unitary penalty cost due to over-staffing of service capacity of type $l$ service\\
    $\overline{w}$ & the upper bound number of full-time caregivers\\
    $\underline{w}$ & the lower bound number of full-time caregivers\\
    $h_{k}$ & the service capacity (service hours) of one type $k$ full-time caregiver in one day\\
    \multicolumn{2}{l}{\textbf{Scenario-dependent Parameters}} \\
    $d_{l,t}$ & the demand quantity for service type $l$ in day $t$\\
    $s_{l,t}$ & the service time for service type $l$ of caregiver in day $t$\\
    \hline
\end{tabular}\label{notation}
\end{table}

\section{Models for the (EA) Decision-maker}\label{EA}
\subsection{Two-stage SP for (EA)}
In this section, we present our proposed two-stage formulation of the HSCP problem considering the (EA) decision-maker. This model (E-SP), assumes that we know the distribution $\mathbb{P}$ of uncertain parameters $\bm{\xi} := [\bm{d}, \bm{s}]\top$. First, let us define the variables defining our E-SP model. Let integer decision variables $x_k$ represent the number of type $k$ caregivers hired, for all $k \in [K]$. Let continuous decision variables $y_{k,l,t}$ represent the capacity of type $k$ caregivers allocated for type $l$ service in day $t$, for all $k\in [K], l\in [L], t\in [T]$. Let continuous decision variables $o_{l,t}$ and $u_{l,t}$ represent over-staffing and under-staffing for type $l$ service in each day t, for all $l\in [L], t\in [T]$. Our E-SP model can now be stated as follows:
\allowdisplaybreaks
{
\setlength{\abovedisplayskip}{4pt}
\setlength{\belowdisplayskip}{4pt}
\begin{subequations}
\label{eq:optim2}
    \begin{alignat}{3}
       \text{(E-SP)}\ \min  &\quad \Bigg \{ \sum_{k=1}^{K} c_{k}x_{k} + \sum_{k=1}^{K} \sum_{l=1}^{L} \sum_{t=1}^{T} c_{k,l,t} y_{k,l,t} + \mathbb{E}_{\mathbb{P}} \ [Q(\textbf{\emph{x}}, \textbf{\emph{y}}, && \bm{\xi})] \Bigg\} \quad \label{eq:cont21}\\
        \mbox{s.t.} &\quad \underline{w} \leq \sum_{k=1}^{K} x_{k} \leq \overline{w}, && \label{eq:cont22}\\
        &\quad  \sum_{l\in R_{k}} y_{k,l,t} \leq x_{k}h_{k}, && \forall k \in [K], t \in [T], \label{eq:cont23}\\
        &\quad x_{k}\in \mathbb{Z}_{+}, y_{k,l,t} \in \mathbb{R}_{+},  &&\forall k\in [K], l\in [L], t\in [T], \label{eq:cont24}
    \end{alignat}
\end{subequations}}
\noindent where for each feasible decision $\bm{x}$, $\bm{y}$ and a joint realization of uncertain parameters $\bm{\xi} := [\bm{d}, \bm{s}]^\top$, our second-stage (recourse) formulation is as follows:
\allowdisplaybreaks
\begin{subequations}
\setlength{\abovedisplayskip}{3pt}
\label{eq:optim3}
\begin{alignat}{3}
    Q(\textbf{\emph{x}}, \textbf{\emph{y}}, \bm{\xi}):= &\quad \mbox{min}&&\quad \sum_{t=1}^{T}\sum_{l=1}^{L} (c_{l,t}^{o}o_{l,t} + c_{l,t}^{u}u_{l,t}) &&\quad \label{eq:const31} \\
    &\quad\mbox{s.t.} &&\quad \sum_{k: l\in R_{k}} y_{k,l,t} + u_{l,t} - o_{l,t} = d_{l,t}s_{l,t}, &&\quad \forall l \in [L], t \in [T], \label{eq:const32}\\
    &\quad &&\quad o_{l,t}, u_{l,t} \in \mathbb{R}_{+},  &&\quad \forall l\in [L], t\in [T]. \label{eq:const33}
\end{alignat}
\end{subequations}
\indent Formulation \eqref{eq:optim2} seeks to find staffing and allocation decisions ($\pmb{x,y}$) that minimize the sum of (1) fixed  employment and assignment costs (first two terms), and (2) expectation of the second stage (recourse) random cost $Q(\textbf{\emph{x}}, \textbf{\emph{y}}, \bm{\xi})$ (third term) subject to uncertainty $\bm{\xi} \sim \mathbb{P}$. The second stage recourse function minimize the random over-staffing (first-term) and under-staffing (second-term) costs.  Constraints \eqref{eq:cont22} ensure that the total number of staff is at least $\underline{w}$ and at most $\overline{w}$ (these can be adjusted by the decision-maker depending on relevant practical constraints).  Constraints \eqref{eq:cont23} ensure that capacity allocations to service types does not exceed the available capacity.  Constraints \eqref{eq:cont24} specifies feasible ranges of the decision variables. Second stage constraint \eqref{eq:const32} computes the over-staffing and under-staffing. 

\subsection{DRO model for (EA)}
In this section, we present our proposed two-stage DRO formulation of the HSCP with (EA) decision-maker, which does not assume that the probability distributions of demand and service durations are known. We, however, assume that we know the mean and support of these random parameters. (Recall that these parameters can be estimated based on, e.g., expert knowledge, where support could represent the dispersion or range of variability that we seek protection against). First, let us define some additional parameters, sets, and notations defining our ambiguity set and DRO model. We let $\mu_{l,t}^d$ and $\mu_{l,t}^{s}$ respectively represent the mean value of $d_{l,t}$ and $s_{l,t}$, for all $l \in [L]$ and $t \in [T]$.  The random vector $\bm{\xi}$ is a measurable function $\bm{\xi}:\Omega\rightarrow \mathcal{U}$ with a measurable space $(\Omega, \mathcal{F})$, where  $\mathcal{U} = \mathcal{U}^{d}\times \mathcal{U}^{s}$ is the bounded support of $\bm{\xi}$.  $\mathcal{U}^{d}$ and $\mathcal{U}^{s}$ are respectively the support set of $\pmb{d}$ and $\pmb{s}$ defined as follows:
\allowdisplaybreaks
{
\setlength{\abovedisplayskip}{4pt}
\setlength{\belowdisplayskip}{4pt}
\begin{align}
\setlength{\abovedisplayskip}{3pt}
    &\mathcal{U}^{d}:=\{\pmb{d}\geq  0: \underline{d}_{l,t} \leq d_{l,t} \leq \overline{d}_{l,t}, \forall l\in [L], t\in [T]\}\label{support_set_1},\\
    &\mathcal{U}^{s}:=\{\pmb{s}\geq 0: \underline{s}_{l,t} \leq s_{l,t} \leq \overline{s}_{l,t}, \forall l\in [L], t\in [T]\}\label{support_set_2}.
\end{align}}
\indent We define $\mathcal{P}(\mathcal{U})$ as the set of probability distributions supported on $\mathcal{U}$, and $\mathbb{E}_{\mathbb{P}}$ as the expectation under $\mathbb{P}$. Finally, we donate $\bm \mu:=\mathbb{E}_\mathbb{P}[\bm{\xi}]=[\bm{\mu}^{d}, \bm{\mu}^{s}]^\top$. Using this notation, we construct the following mean-support ambiguity set:
\allowdisplaybreaks
{
\setlength{\abovedisplayskip}{4pt}
\setlength{\belowdisplayskip}{4pt}
\begin{align}\label{ambiguity_set}
      \mathcal{F}(\mathcal{U}, \bm{\mu}):=\left\{ \mathbb{P}\in \mathcal{P}(\mathcal{U})
    \Bigg|\begin{array}{lr}
    \int_{\mathcal{U}} d\mathbb{P} = 1\\
    \mathbb{E}_{\mathbb{P}}[\bm{\xi}] = \bm{\mu}
    \end{array} 
    \right\}.
\end{align}}
\indent Using the ambiguity set $\mathcal{F}(\mathcal{U}, \bm{\mu})$, we formulate our DRO model of the HSCP with (EA) decision-maker (denoted as E-DHSCP) as follow:
{
\setlength{\abovedisplayskip}{4pt}
\setlength{\belowdisplayskip}{4pt}
\begin{equation}
    \text{(E-DHSCP)} \ \    \underset{(\bm{x},\bm{y})\in \mathcal{X}}{\mbox{min}}  \Bigg \{\sum_{k=1}^{K} c_{k}x_{k} + \sum_{k=1}^{K} \sum_{l=1}^{L} \sum_{t=1}^{T} c_{k,l,t}y_{k,l,t} + 
    \underset{\mathbb{P}\in \mathcal{F}(\mathcal{U}, \bm{\mu})}{\mbox{sup}}
    \mathbb{E}_{\mathbb{P}} [Q(\textbf{\emph{x}}, \textbf{\emph{y}}, \bm{\xi})] \Bigg\} \label{min-max},
\end{equation}}
where $\mathcal{X}:=\{(\bm{x},\bm{y})$: \eqref{eq:cont21}-\eqref{eq:cont24}\} is the feasible region of first stage decision variables. Formulation \eqref{min-max} seeks first stage decisions $\bm{x}$ and $\bm{y}$ that minimizes the summation of fixed hiring cost (first term), the capacity allocation cost (second term) and the worst-case (maximum) expectations of the second-stage operational cost (third term), where the expectation is taken over all distributions residing in $\mathcal{F}(\mathcal{U}, \bm{\mu})$.

\subsubsection{Reformulation}
Recall that $Q(\cdot)$ is defined by a minimization problem; hence, in \eqref{min-max}, we have an inner max-min problem. As such, it is not straightforward to solve \eqref{min-max} in its presented form. In this section, we derive an equivalent formulation of the min-max model \eqref{min-max} that is solvable. First, in Proposition~\ref{pro_deter}, we present an equivalent reformulation of the inner maximization problem $\underset{\mathbb{P}\in \mathcal{F}(\mathcal{U}, \bm{\mu})}{\mbox{sup}} \mathbb{E}_{\mathbb{P}} [Q(\textbf{\emph{x}}, \textbf{\emph{y}}, \bm{\xi})]$ in \eqref{min-max} (see \ref{proof_pos_d} for a detailed proof).

\begin{proposition}
\label{pro_deter}
For any $\bm{(x,y)}\in \mathcal{X}$, problem  $\underset{\mathbb{P}\in \mathcal{F}(\mathcal{U}, \bm{\mu})}{\mbox{sup}}
    \mathbb{E}_{\mathbb{P}} [Q(\bm{x}, \bm{y}, \bm{\xi})]$ in \eqref{min-max} is equivalent to.
{
\setlength{\abovedisplayskip}{4pt}
\setlength{\belowdisplayskip}{4pt}
\begin{equation}\label{deter}
   \begin{split}
    {\rm(E-DHSCP)}\ \min_{\bm{\alpha , \beta}} \Bigg\{
    \sum_{l=1}^{L}\sum_{t=1}^{T} \mu_{l,t}^d\alpha_{l,t} & + \sum_{l=1}^{L}\sum_{t=1}^{T}\mu_{l,t}^s\beta_{l,t} \\
    &+ \max_{(\bm{d,s})\in \mathcal{U}} \Big\{Q(\bm{x}, \bm{y}, \bm{\xi}) - \sum_{l=1}^L \sum_{t =1}^T(d_{l,t}\alpha_{l,t} + s_{l,t}\beta_{l,t})\Big\}
    \Bigg\}.
   \end{split}
\end{equation}}
\end{proposition}

Again, the problem in \eqref{deter} involves an inner max-min problem that is not straightforward to solve in its presented form.  Next, we use the structural properties of the recourse $Q(\cdot)$ to derive an equivalent linear programming  reformulation of the inner maximization problem in \eqref{deter}. Note that for fixed $\bm{x, y}, \bm{\xi}$,  $Q(\textbf{\emph{x}}, \textbf{\emph{y}}, \bm{\xi})$ is a linear program (LP). The dual formulation of $Q(\textbf{\emph{x}}, \textbf{\emph{y}}, \bm{\xi})$ is as follows:
\allowdisplaybreaks
{
\setlength{\abovedisplayskip}{4pt}
\setlength{\belowdisplayskip}{4pt}
\begin{subequations}
\label{eq:optim8}
\begin{alignat}{3}
    Q(\textbf{\emph{x}}, \textbf{\emph{y}}, \bm{\xi})=\quad \max_{\bm{\rho}}  &\quad \sum_{l=1}^{L} \sum_{t=1}^{T} (d_{l,t}s_{l,t} - \sum_{k: l\in R_{k}} y_{k,l,t})\rho_{l,t}\quad \label{eq:const81} \\
    \mbox{s.t.} &\quad - c_{l,t}^{o} \leq \rho_{l,t} \leq c_{l,t}^{u}, &&\quad \forall l\in [L], t\in [T]. \label{eq:const82}
\end{alignat}
\end{subequations}}
where $\rho_{l,t}$ is the dual variables associated with constraints \eqref{eq:const32}. In view of the dual formulation in \eqref{eq:optim8}, the inner problem $\max_{(\bm{d,s})\in \mathcal{U}} \left\{Q(\bm{x}, \bm{y}, \bm{\xi}) - \sum_{l=1}^L \sum_{t =1}^T(d_{l,t}\alpha_{l,t} + s_{l,t}\beta_{l,t})\right\}$ in \eqref{deter} is equivalent to problem \eqref{inner_max} in Proposition~\ref{min_inner} (see \ref{Appx:min_inner} for a detailed proof).

\begin{proposition}\label{min_inner}
For fixed $\bm{x}, \bm{y}, \bm{\alpha}$, and $\bm{\beta}$, problem $\underset{(\bm{d,s})\in \mathcal{U}}{\max} \Big\{Q(\bm{x}, \bm{y}, \bm{\xi}) - \sum\limits_{l=1}^L \sum\limits_{t =1}^T(d_{l,t}\alpha_{l,t} + s_{l,t}\beta_{l,t})\Big\}$ is equivalent to.
{
\setlength{\abovedisplayskip}{4pt}
\setlength{\belowdisplayskip}{4pt}
\begin{subequations}
\label{inner_max}
    \begin{alignat}{3}
       \min_{\pmb{\eta}}
        &\quad \sum_{l=1}^L \sum_{t=1}^T \eta_{l,t} &&\quad \label{inner_max_1}\\
        \mbox{s.t.}
        &\quad \eta_{l,t} \geq \overline{d}_{l,t}\overline{s}_{l,t}c_{l,t}^{u} - \sum_{k: l\in R_{k}} y_{k,l,t}c_{l,t}^{u} - \overline{d}_{l,t}\alpha_{l,t} -\overline{s}_{l,t}\beta_{l,t},&&\quad \forall l\in[L], t\in[T],\label{inner_max_first}\\
        &\quad \eta_{l,t} \geq -\underline{d}_{l,t}\underline{s}_{l,t}c_{l,t}^{o} + \sum_{k: l\in R_{k}} y_{k,l,t}c_{l,t}^{o} - \underline{d}_{l,t}\alpha_{l,t} - \underline{s}_{l,t}\beta_{l,t},&&\quad \forall l\in[L], t\in[T],\\
        &\quad \eta_{l,t} \geq -\overline{d}_{l,t}\overline{s}_{l,t}c_{l,t}^{o} + \sum_{k: l\in R_{k}} y_{k,l,t}c_{l,t}^{o} - \overline{d}_{l,t}\alpha_{l,t} - \overline{s}_{l,t}\beta_{l,t},&&\quad \forall l\in[L], t\in[T],\\
        &\quad \eta_{l,t} \geq \underline{d}_{l,t}\overline{s}_{l,t}c_{l,t}^{u} - \sum_{k: l\in R_{k}} y_{k,l,t}c_{l,t}^{u} - \underline{d}_{l,t}\alpha_{l,t} - \overline{s}_{l,t}\beta_{l,t},&&\quad \forall l\in[L], t\in[T],\\
        &\quad \eta_{l,t} \geq \overline{d}_{l,t}\underline{s}_{l,t}c_{l,t}^{u} - \sum_{k: l\in R_{k}} y_{k,l,t}c_{l,t}^{u} - \overline{d}_{l,t}\alpha_{l,t} - \underline{s}_{l,t}\beta_{l,t},&&\quad \forall l\in[L], t\in[T],\\
        &\quad \eta_{l,t} \geq \underline{d}_{l,t}\underline{s}_{l,t}c_{l,t}^{u} - \sum_{k: l\in R_{k}} y_{k,l,t}c_{l,t}^{u} - \underline{d}_{l,t}\alpha_{l,t} - \underline{s}_{l,t}\beta_{l,t},&&\quad \forall l\in[L], t\in[T],\\
        &\quad \eta_{l,t} \geq -\overline{d}_{l,t}\underline{s}_{l,t}c_{l,t}^{o} + \sum_{k: l\in R_{k}} y_{k,l,t}c_{l,t}^{o} - \overline{d}_{l,t}\alpha_{l,t} - \underline{s}_{l,t}\beta_{l,t},&&\quad \forall l\in[L], t\in[T],\\
        &\quad \eta_{l,t} \geq -\underline{d}_{l,t}\overline{s}_{l,t}c_{l,t}^{o} + \sum_{k: l\in R_{k}} y_{k,l,t}c_{l,t}^{o} - \underline{d}_{l,t}\alpha_{l,t} - \overline{s}_{l,t}\beta_{l,t},&&\quad \forall l\in[L], t\in[T].\label{inner_max_last}
    \end{alignat}
\end{subequations}
}
\end{proposition}

Replacing the inner problem  $\max \limits_{(\bm{d,s})\in \mathcal{U}} \{ \cdot\}$ in \eqref{deter} with its equivalent reformulation in \eqref{inner_max}, then combining with the outer minimization problems  in \eqref{deter} and \eqref{min-max}, we derive the following equivalent MILP reformulation of the E-DHSCP model.
\begin{subequations}
\setlength{\abovedisplayskip}{4pt}
\label{final_1}
    \begin{alignat}{3}
        \mbox{min} &\ \Big\{\sum_{k=1}^{K} c_{k}x_{k} + \sum_{k=1}^{K} \sum_{l=1}^{L} \sum_{t=1}^{T} c_{k,l,t}y_{k,l,t} + \sum_{l=1}^{L}\sum_{t=1}^{T} \mu_{l,t}^d\alpha_{l,t} + \sum_{l=1}^{L}\sum_{t=1}^{T}\mu_{l,t}^s\beta_{l,t} + \sum_{l=1}^{L}\sum_{t=1}^{T}\eta_{l,t} \Big\} &&\quad \label{final_1_1}\\
        \mbox{s.t.} &\quad (\bm{x,y}) \in \mathcal{X}, \bm{\alpha}\in \mathbb{R}, \bm{\beta}\in \mathbb{R},&&\label{final_1_2}\\
        &\quad \eqref{inner_max_first} - \eqref{inner_max_last}. &&\label{final_1_3}
    \end{alignat}
\end{subequations}

\section{Models for the (FA) Decision-maker}\label{wait-and-see}
\subsection{Two-stage SP for (FA)}

In this section, we present our proposed two-stage SP formulation of the HSCP problem considering the (FA) decision-maker. Specifically, we keep $x_{k},\ \forall k\in [K]$ as first stage decisions, and move $y_{k,l,t}$ to the second-stage, for all $k \in [K], l \in [L], t \in [T]$. In addition, we define non-negative continuous variables $o_{k,t}$ representing the capacity surplus of type $k \in [K]$ caregiver in day $t\in [T]$. This is to account scenarios where the total capacity of hired caregivers is larger than the total realized demand. Accordingly, our F-SP model can now be stated as follows:
\allowdisplaybreaks
{
\setlength{\abovedisplayskip}{4pt}
\setlength{\belowdisplayskip}{4pt}
\begin{subequations}
    \label{FA}
    \begin{alignat}{3}
        \mbox{(F-SP)}\ {\text{min}} &\quad \sum_{k=1}^{K} c_{k}x_{k} + \mathbb{E} \left[ Q^{A}(\textbf{\emph{x}}, \bm{\xi})\right]  &&\quad \label{FA_1} \\
        \mbox{s.t.} &\quad  \underline{w} \leq \sum_{k=1}^{K} x_{k} \leq \overline{w}, &&\quad \label{FA_2} \\
        &\quad x_{k}\in \mathbb{Z}_{+},  &&\forall k\in [K], \label{FA_3}
    \end{alignat}
\end{subequations}
}
where for each feasible decision $\bm{x}$ and a joint realization of uncertain parameters $\bm{\xi} := [\bm{d}, \bm{s}]^\top$, our second-stage (recourse) formulation is as follows:
{
\setlength{\abovedisplayskip}{4pt}
\setlength{\belowdisplayskip}{4pt}
\begin{subequations}
    \label{FA_inner}
    \begin{alignat}{4}
        Q^{A}(\textbf{\emph{x}}, \bm{\xi})= &\mbox{min} &&\ \sum_{k=1}^{K}\sum_{l=1}^{L}\sum_{t=1}^{T} c_{k,l,t}y_{k,l,t} + \sum_{k=1}^{K}\sum_{t=1}^{T} c_{k,t}^{o}o_{k,t} + && \sum_{l=1}^{L}\sum_{t=1}^{T} c_{l,t}^{u}u_{l,t} \quad \label{FA_inner_1} \\
        &\mbox{s.t.} &&\  \sum_{l\in R_{k}} y_{k,l,t} + o_{k,t} = x_{k}h_{k}, && \forall k \in [K], t\in [T] \label{FA_inner_2},\\
        &\quad  &&\ \sum_{k: l\in R_{k}} y_{k,l,t} + u_{l,t} = d_{l,t}s_{l,t}, && \forall l \in [L], t\in [T] \label{FA_inner_3},\\
        &\quad &&\ y_{k,l,t}, o_{k,t}, u_{l,t} \in \mathbb{R}_{+},  &&\forall k\in [K], l\in [L], t\in [T]. \label{FA_inner_4}
    \end{alignat}
\end{subequations}
}
\indent Formulation \eqref{FA} seeks to find staffing decisions ($\pmb{x}$) that minimize the sum of (1) fixed  employment cost (first term), and (2) expectation of the second stage (recourse) random cost $Q(\bm{x}, \bm{\xi})$ (second term) subject to uncertainty $\bm{\xi} \sim \mathbb{P}$. The second stage recourse function minimize the random assignment cost (first-term), over-staffing (second-term) and under-staffing (third term) costs.  Constraints \eqref{FA_inner_2} compute the over-staffing. Constraints \eqref{FA_inner_3} compute the under-staffing.

\subsection{DRO model for (FA)}
In this section, we present our proposed two-stage DRO formulation of the HSCP with (FA) decision-maker. Using the ambiguity set $\mathcal{F}(\mathcal{U}, \bm{\mu})$ defined in \eqref{ambiguity_set}, we formulate our DRO model of the HSCP with (FA) decision-maker (denoted as F-DHSCP) as follow:
{
\setlength{\abovedisplayskip}{3pt}
\setlength{\belowdisplayskip}{3pt}
\begin{equation}
    \text{(F-DHSCP)}  \ \  \underset{\bm{x}\in \mathcal{X}_{2}}{\mbox{min}}  \Bigg \{\sum_{k=1}^{K} c_{k}x_{k} +
    \underset{\mathbb{P}\in \mathcal{F}(\mathcal{U}, \bm{\mu})}{\mbox{sup}}
    \mathbb{E}_{\mathbb{P}} [Q^{A}(\textbf{\emph{x}}, \bm{\xi})] \Bigg\}, \label{min-max_A}
\end{equation}
}
where $\mathcal{X}_{2}:=\{\eqref{FA_2}-\eqref{FA_3}\}$ is the feasible region of first stage decision variables. Formulation \eqref{min-max_A} seeks the first stage decisions $\bm{x}$ that minimize the sum of fixed employment cost (first term) and the worst-case (maximum) expectations of the second-stage operational cost (second term), where the expectation is taken over all distributions residing in $\mathcal{F}(\mathcal{U}, \bm{\mu})$.

Using the same techniques in Proposition \ref{pro_deter}, we derive the following equivalent reformulation of the inner maximization problem $\mbox{sup}_{\mathbb{P}\in \mathcal{F}(\mathcal{U}, \bm{\mu})} \mathbb{E}_{\mathbb{P}} [Q^{A}(\textbf{\emph{x}}, \bm{\xi})]$ in \eqref{min-max_A}.
\begin{equation}\label{A_deter}
\setlength{\abovedisplayskip}{6pt}
    \text{(F-DHSCP)}\ \underset{\bm{\alpha}, \bm{\beta}}{\mbox{min}} \left\{
    \sum_{t=1}^{T}\sum_{l=1}^{L} (\mu_{l,t}^d\alpha_{l,t} + \mu_{l,t}^s\beta_{l,t}) + \underset{(\bm{d,s})\in \mathcal{U}}{\mbox{max}} \Big\{Q^{A}(\bm{x}, \bm{\xi}) - \sum_{l=1}^{L}\sum_{t=1}^{T}(d_{l,t}\alpha_{l,t} + s_{l,t}\beta_{l,t})\Big\} \right\}.
\end{equation}
\indent Note that it is not straightforward to solve formulation \eqref{A_deter} in its presented form due to the inner max-min problem. Therefore, our goal is to derive an equivalent formulation of \eqref{A_deter} that is solvable. As a first step, we re-write the recourse $Q^{A}(\bm{x}, \bm{\xi})$ using its equivalent dual formulation as (note that for fixed $\bm{x}$ and $\bm{\xi}$, $Q^{A}(\bm{x}, \bm{\xi})$ is a bounded LP):
\allowdisplaybreaks
{
\setlength{\abovedisplayskip}{4pt}
\setlength{\belowdisplayskip}{4pt}
\begin{subequations}
\label{A_dual}
\begin{alignat}{3}
    Q^{A}(\textbf{\emph{x}}, \bm{\xi})=\quad \max_{\bm{\rho}}  &\sum_{t=1}^{T}\Bigg \{ \sum_{l=1}^{L} d_{l,t}s_{l,t}\rho_{l,t} + \sum_{k\in K} \lambda_{k,t}x_{k}h_{k} \Bigg \}\quad \label{A_dual_1} \\
    \mbox{s.t.} &\quad \rho_{l,t} + \lambda_{k,t} \leq c_{k,l,t}, && \forall t \in [T], k \in [K], l\in R_{k}, \label{A_dual_2}\\
    &\quad \rho_{l,t} \leq c_{l,t}^{u}, && \forall t \in [T], l \in [L], \label{A_dual_3}\\
    &\quad \lambda_{k,t} \leq c_{k,t}^{o}, && \forall t \in [T], k \in [K], \label{A_dual_4}
\end{alignat}
\end{subequations}}
where $\rho_{l,t}$ and $\lambda_{k,t}$ is the dual variables associated with constraints \eqref{FA_inner_2} and \eqref{FA_inner_3}. Replacing $Q^{A}(\cdot)$ in \eqref{A_deter} with its dual formulation in \eqref{A_dual}, we derive the following equivalent reformulation of \eqref{A_deter}.
{
\setlength{\abovedisplayskip}{4pt}
\setlength{\belowdisplayskip}{4pt}
\begin{equation}\label{A_inner}
\setlength{\abovedisplayskip}{6pt}
    \text{(F-DHSCP)}\ \underset{\bm{\alpha}, \bm{\beta}}{\mbox{min}} \left\{
    \sum_{t=1}^{T}\sum_{l=1}^{L} (\mu_{l,t}^d\alpha_{l,t} + \mu_{l,t}^s\beta_{l,t}) + \underset{(\bm{d,s})\in \mathcal{U}}{\mbox{max}} h(\bm{x},\bm{d}, \bm{s}, \bm{\alpha}, \bm{\beta})
    \right\},
\end{equation}}
where
{
\setlength{\abovedisplayskip}{4pt}
\setlength{\belowdisplayskip}{4pt}
\begin{equation}\label{pro_deter_A_sub}
\setlength{\abovedisplayskip}{6pt}
    h(\bm{x}, \bm{d}, \bm{s}, \bm{\alpha}, \bm{\beta}) := \underset{(\bm{\rho, \lambda})\in P}{\mbox{max}} \sum_{t=1}^{T} \left\{\sum_{l=1}^{L} d_{l,t}s_{l,t}\rho_{l,t} +\sum_{k\in K} \lambda_{k,t}x_{k}h_{k} -\sum_{l=1}^{L}(d_{l,t}\alpha_{l,t} + s_{l,t}\beta_{l,t})\right\}.
\end{equation}}
$P:=\{\eqref{A_dual_2} - \eqref{A_dual_4}\}$ is the feasible region of variables $\bm{\rho}$ and $\bm{\lambda}$. Again, the problem in \eqref{A_inner} is not straightforward to solve in its presented form due to the inner maximization $\underset{(\bm{d,s})\in \mathcal{U}}{\mbox{max}} h(\bm{x}, \bm{d}, \bm{s}, \bm{\alpha}, \bm{\beta})$. However, next, in Proposition~\ref{LinearH}, we present an equivalent MILP of $\underset{(\bm{d,s})\in \mathcal{U}}{\mbox{max}} h(\bm{x}, \bm{d}, \bm{s}, \bm{\alpha}, \bm{\beta})$ that is solvable (see  \ref{Appx:H_Reform} for a detailed proof).

\begin{proposition}\label{LinearH}
Let $\Delta d_{l,t} = \overline{d}_{l,t} - \underline{d}_{l,t}$ and $\Delta s_{l,t} = \overline{s}_{l,t} - \underline{s}_{l,t}$, for all $l \in [L]$ and $t \in [T]$. Then, fixed $\bm{x}$, $\bm{\alpha}$, and $\bm{\beta}$, solving  $\underset{(\bm{d,s})\in \mathcal{U}}{\mbox{max}} h(\bm{x}, \bm{d}, \bm{s}, \bm{\alpha}, \bm{\beta})$ is equivalent to solving the following MILP:
{
\setlength{\abovedisplayskip}{0pt}
\setlength{\belowdisplayskip}{4pt}
\begin{subequations}
    \label{h_final}
    \begin{alignat}{3}
        \underset{\substack{\bm{\rho}, \bm{\lambda}, \bm{g}, \bm{z},\\ \bm{v}, \bm{r}, \bm{q}, \bm{e}}}{\mbox{max}} & \sum_{t=1}^{T} \Bigg\{ \sum_{k=1}^{K} x_{k}h_{k}\lambda_{k,t} -\sum_{l=1}^{L}[\underline{d}_{l,t}\alpha_{l,t} + \Delta d_{l,t}g_{l,t}\alpha_{l,t} + \underline{s}_{l,t}\beta_{l,t} + \Delta s_{l,t} && z_{l,t} \beta_{l,t}] \quad \notag\\
        &\quad + \sum_{l=1}^{L} \left[\underline{d}_{l,t}\underline{s}_{l,t}\rho_{l,t} +  \Delta s_{l,t}\underline{d}_{l,t}r_{l,t} + \Delta d_{l,t}\underline{s}_{l,t}q_{l,t} + \Delta d_{l,t}\Delta s_{l,t}e_{l,t} \right] &&\Bigg \} \quad \label{h_final_1}\\
        \mbox{s.t.}&\quad \eqref{A_dual_2}-\eqref{A_dual_4}, &&\label{h_final_2}\\
        &\quad \rho_{l,t} \geq \underline{\rho}_{l,t}, \quad v_{l,t} - g_{l,t} \leq 0,&& \forall l\in [L], t\in [T],\\
        &\quad v_{l,t} \geq 0, \quad v_{l,t} - z_{l,t} \leq 0, &&\forall l\in [L], t\in [T],\\
        &\quad v_{l,t}- g_{l,t} - z_{l,t} \geq -1, && \forall l\in [L], t\in [T],\\
        &\quad e_{l,t} - \underline{\rho}_{l,t} v_{l,t}\geq 0, \quad e_{l,t}  - c_{l,t}^{u}v_{l,t} \leq 0, && \forall l\in [L], t\in [T], \label{Mac1}\\
        &\quad e_{l,t} - \rho_{l,t} + \underline{\rho}_{l,t}(1-v_{l,t})\leq 0,\quad e_{l,t}- \rho_{l,t} +c_{l,t}^u (1-v_{l,t}) \geq 0,&& \forall l\in [L], t\in [T],\\
        &\quad q_{l,t} - \underline{\rho}_{l,t} g_{l,t}\geq 0,\quad q_{l,t}  - c_{l,t}^{u}g_{l,t} \leq 0, && \forall l\in [L], t\in [T],\\
        &\quad q_{l,t} - \rho_{l,t} + \underline{\rho}_{l,t}(1-g_{l,t})\leq 0,\quad q_{l,t}- \rho_{l,t} +c_{l,t}^u (1-g_{l,t}) \geq 0, && \forall l\in [L], t\in [T], \\
        &\quad r_{l,t} - \underline{\rho}_{l,t} v_{l,t}\geq 0,\quad r_{l,t}  - c_{l,t}^{u}v_{l,t} \leq 0, && \forall l\in [L], t\in [T],\\
        &\quad r_{l,t} - \rho_{l,t} + \underline{\rho}_{l,t}(1-v_{l,t})\leq 0,\quad r_{l,t}- \rho_{l,t} +c_{l,t}^u (1-v_{l,t}) \geq 0, &&\forall l\in [L], t\in [T], \label{Maclast}\\
        &\quad g_{l,t}, z_{l,t}\in \{0,1\},v_{l,t}, r_{l,t}, q_{l,t}, e_{l,t}, \rho_{l,t} \in \mathbb{R} && \forall l\in [L], t\in [T]. \label{h_final_5}
    \end{alignat}
\end{subequations}}
\end{proposition}
Note that  McCormick inequalities \eqref{Mac1}--\eqref{Maclast} in formulation \eqref{h_final} involve big-M coefficients $\underline{\rho}_{l,t}$, for all $l \in [L]$ and $t \in [T]$, that can undermine computational efficiency if they are set too large. Therefore, in Proposition~\ref{pro_rho}, we derive tight lower bounds (i.e., $\underline{\rho}_{l,t}$) on variables $\rho_{l,t}$ to strengthen the MILP reformulation (see  \ref{Appx:H_Reform} for a detailed proof).
\begin{proposition}
\label{pro_rho}
$\underline{\rho}_{l,t} = \underset{k:l\in R_{k}}{\mbox{min}} (c_{k,l,t} - c_{k,t}^o)$ is a valid lower bounds on  $\rho_{l,t}$, for all $l\in [L]$ and $t\in [T]$.
\end{proposition}

Replacing $\underset{(\bm{d,s})\in \mathcal{U}}{\mbox{max}} h(\bm{x}, \bm{d}, \bm{s}, \bm{\alpha}, \bm{\beta})$ in \eqref{A_inner} with \eqref{h_final}, then combining the inner problem in the form of \eqref{pro_deter_A_sub} with the outer minimization problem in the F-DHSCP model in \eqref{min-max_A}, we derive the following equivalent reformulation of the F-DHSCP model:
\allowdisplaybreaks
{
\setlength{\abovedisplayskip}{4pt}
\setlength{\belowdisplayskip}{4pt}
\begin{subequations}
\label{A_final}
    \begin{alignat}{4}
        \text{(F-DHSCP)}\quad \underset{\bm{x}\in \mathcal{X}_{2}, \bm{\alpha}, \bm{\beta}}{\mbox{min}} &\quad \Bigg \{ \sum_{k=1}^{K} c_{k}x_{k} + \sum_{t=1}^{T}\sum_{l=1}^{L} (\mu_{l,t}^d\alpha_{l,t} + \mu_{l,t}^s\beta_{l,t}) + \sum_{t=1}^{T}\delta_{t} \Bigg\} &&\quad \label{A_final_1}\\
        \mbox{s.t.} &\quad \delta_{t} \geq \underset{\substack{\bm{\rho}_t, \bm{\lambda}_t, \bm{g}_t, \bm{z}_t,\\ \bm{v}_t, \bm{r}_t, \bm{q}_t, \bm{e}_t}}{\mbox{max}} \Bigg \{ \sum_{k=1}^{K} x_{k}h_{k}\lambda_{k,t} -\sum_{l=1}^{L}\Big[\underline{d}_{l,t}\alpha_{l,t} + \Delta d_{l,t}g_{l,t}\alpha_{l,t} \nonumber\\
        & \qquad \qquad \qquad + \underline{s}_{l,t}\beta_{l,t} + \Delta s_{l,t}z_{l,t}\beta_{l,t} + \underline{d}_{l,t}\underline{s}_{l,t}\rho_{l,t} +  \Delta s_{l,t}\underline{d}_{l,t}r_{l,t} \nonumber\\
        & \qquad \qquad \qquad + \Delta d_{l,t}\underline{s}_{l,t}q_{l,t}\ + \Delta d_{l,t}\Delta s_{l,t}e_{l,t}\Big]: \eqref{h_final_2}-\eqref{h_final_5} \Bigg\}, &&t\in [T]. \label{C_max_A_Final}
    \end{alignat}
\end{subequations}
}

\section{Solution approach} \label{solution_approach}
\noindent In this section, we present a Monte Carlo Optimization algorithm to obtain near-optimal solutions to the SP models (Section \ref{MCO}), and a column-and-constraint generation (C\&CG) algorithm to solve the F-DHSCP model (Section \ref{CAG}). Note that the E-DHSCP model can be implemented and solved directly using off-the-shelf optimization software.

\subsection{Monte Carlo Optimization}\label{MCO}
Note that it is difficult to obtain an exact optimal solution to the two-stage SP models in \eqref{eq:optim2} and \eqref{FA}. Hence, we first present a Monte Carlo approach base on \cite{homem2014monte}, \cite{kleywegt2002sample} and \cite{shehadeh2021using} to obtain near-optimal solutions. We replace the distribution of random parameters with a empirical distribution $N$ independent and identically distributed (i.i.d) samples of services demands and durations at first and then we solve the sample average approximation (SAA) formulations of \eqref{eq:optim2} and \eqref{FA} in \ref{SAA}.

Algorithm \ref{algorithm1} in \ref{SAA} summarizes the Monte Carlo Optimization (MCO) algorithm. We start with a small initial sample size. In first step, we will solve the SAA formulation with designated sample size $N$ and evaluate the solution $x_{N}^{k}$ by Monte Carlo simulation with the generated $N'$ random samples and record the objective value $v_{N}^{k}$ of SAA and $v_{N'}^{k}$ of Monte Carlo simulation. We will run step 1 with $K$ replicates. In second step, we calculate the average of $v_{N}^{k}$ and $v_{N'}^{k}$ as $\overline{v}_{N}$ and $\overline{v}_{N'}$. According to \cite{mak1999monte}, $\overline{v}_{N}$ and $\overline{v}_{N'}$ are statistical lower and upper bound of the optimal value of original formulation. In step 3, we calculate the approximate optimality index as $| AOI_{N} = \frac{\overline{v}_{N'} - \overline{v}_{N}}{\overline{v}_{N'}} |$. $AOI_{N}$ approximately estimate the optimality gap of lower bound and upper bound of original formulation. Therefore, if $AOI_{N}$ is smaller than the very small predetermined termination tolerance $\epsilon$ (note that $\epsilon$ can be zero), the algorithm will terminate and output the termination sample size $N$ , objective value $\overline{v}_{N}$, $\overline{v}_{N'}$ and approximate optimality index $AOI_{N}$. Otherwise, we will increase the sample size ($N\leftarrow 2N$) then go to step 1.

\subsection{Column-and-constraint generation algorithm}\label{CAG}

The F-DHSCP formulation \eqref{A_final} involves a maximization problem in constraints \eqref{C_max_A_Final} and thus we cannot solve it directly using standard techniques. In this section, we develop a C\&CG algorithm to solve the F-DHSCP model \eqref{A_final}. The motivation of this algorithm is as follows. Considering the inner maximization problem $\underset{(\bm{d,s})\in \mathcal{U}}{\mbox{max}} h(\bm{x}, \bm{d}, \bm{s}, \bm{\alpha}, \bm{\beta})$, the optimal $\bm{d}^*\in\mathcal{U}^{d}$ would only take value $\pmb{ \underline{d}}$ or $\pmb{\overline{d}}$. Similarly, the optimal $\bm{s}^* \in \mathcal{U}^{s}$ would only take value $\pmb{\underline{s}}$ or   $\pmb{\overline{s}}$. Thus, the inner maximization over $(\bm{d,s})\in \mathcal{U}$ could be taken over all the combination of $\pmb{d}$ and $\pmb{s}$. However, only a few critical scenarios of $\bm{d}$ and $\bm{s}$ play a role in deriving an optimal solution to the F-DHSCP model. As such, instead of solving a problem with exponentially many constraints, we use C\&CG to identify these scenarios and obtain an optimal solution.

Algorithm \ref{Alg1:CAG} presents our C\&CG algorithm which is based on the original C\&CG of \cite{Zeng_Zhao:2013} and is implemented in a master-sub problem framework. The master problem \eqref{MP} employs scenario-based constraints \eqref{MP_third}-\eqref{MP_last} and the sub-problem generates demand and service time scenarios and optimality cuts from the worst-case distribution. The algorithm proceeds as follows. At each iteration, we first solve the master problem, which includes only a subset of demand and service time scenarios, and obtain the optimal solution (i.e., hiring decisions). Thus, by considering a subset of demand and service time scenarios, the master problem is a relaxation of the original problem and thus provides a lower bound on the optimal value of F-DHSCP. Then, given the optimal solution from the master problem, we identify demand $\pmb{d^*}$ and service time $\pmb{s^*}$ scenarios by solving the sub-problem. Note that solutions of the master problem are feasible to the original problem, then by solving the sub-problem using these solutions, we obtain an upper bound. Third, we pass the identified demand and service time scenarios to the master problem on an as-needed basis, along with introducing relevant second-stage variables and constraints. Then, we solve the master problem again with the new information (a larger set of scenarios $\mathbb{\textbf{O}}$) from the sub-problems. This process continues until the gap between the lower and upper bound on the optimal value of the problem obtained in each iteration satisfies a predetermined termination tolerance (note that we set this tolerance to 0.02).
\begin{algorithm}[!]
\small
\caption{Column-and-constraint generation (C\&CG).}
\label{Alg1:CAG}
\noindent \textbf{Step 1. Initialization.} Set $LB=0$, $UB=+\infty$, $N=0$, $\mathbb{\textbf{O}}=\emptyset$

\textbf{Step 2. Master Problem.} Solve the following master problem:
\begin{subequations}
    \label{MP}
        \begin{alignat}{3}
            Z = \underset{\bm{x}, \bm{\alpha}, \bm{\beta}, \bm{\delta}}{\mbox{min}} &\quad \Bigg \{ \sum_{k=1}^{K} c_{k}x_{k} + \sum_{t=1}^{T}\sum_{l=1}^{L} (\mu_{l,t}^d\alpha_{l,t} + \mu_{l,t}^s\beta_{l,t}) + \sum_{t=1}^{T} \delta_{t} \Bigg\} &&\quad  \label{MP_first}\\
            \mbox{s.t.} &\quad \bm{x}\in \mathcal{X}_{2}, &&\\
            &\quad \delta_{t} \geq  \sum_{k=1}^{K}\sum_{l=1}^{L} c_{k,l,t} y_{k,l,t}^{n} + \sum_{k=1}^{K} c_{k,t}^{o}o_{k,t}^{n} + \sum_{l=1}^{L} c_{l,t}^{u}u_{l,t}^{n}-\sum_{l=1}^{L}(d_{l}^{t,n}\alpha_{l,t} + s_{l}^{t,n}\beta_{l,t}) ,&&\forall t\in [T], \forall n \in \mathbb{\textbf{O}} \label{MP_second},\\
            &\quad  \sum_{l\in R_{k}} y_{k,l,t}^{n} + o_{k,t}^{n} = x_{k}h_{k}, \qquad \forall k \in [K], t\in [T], n\leq N \label{MP_third},\\
            &\quad\sum_{k: l\in R_{k}} y_{k,l,t}^{n} + u_{l,t}^{n} = d_{l,t}^{n}s_{l,t}^{n}, \qquad  \forall l \in [L], t\in [T], n\leq N,\\
            &\quad y_{k,l,t}^{n}, o_{k,t}^{n}, u_{l,t}^{n} \in \mathbb{R}_{+}, \qquad \forall k\in [K], l\in [L], t\in [T], n\leq N. \label{MP_last}
        \end{alignat}
\end{subequations}

Record an optimal solution ($\bm{x}^*, \bm{\alpha}^*, \bm{\beta}^*, \bm{\delta}^*$) and set $LB=Z^{*}$.

\textbf{Step 3: Solve sub-problem.} 

\begin{itemize}
\item[3.1] With ($\bm{x}, \bm{\alpha}, \bm{\beta}$) fixed to be ($\bm{x}^*, \bm{\alpha}^*, \bm{\beta}^*$), solve the following problem \eqref{sub} for all $t\in [T]$.
\begin{equation}\label{sub}
        \begin{split}
            W_{t} &=\underset{\substack{\bm{\rho}_t, \bm{\lambda}_t, \bm{g}_t, \bm{z}_t,\\ \bm{v}_t, \bm{r}_t, \bm{q}_t, \bm{e}_t}}{\mbox{max}} \Bigg \{ \sum_{k=1}^{K} x_{k}h_{k}\lambda_{k,t} -\sum_{l=1}^{L}\Big[\underline{d}_{l,t}\alpha_{l,t} + \Delta d_{l,t}g_{l,t}\alpha_{l,t} + \underline{s}_{l,t}\beta_{l,t} + \Delta s_{l,t}z_{l,t}\beta_{l,t}\Big]\\
            & \qquad \qquad + \sum_{l=1}^{L} \Big[\underline{d}_{l,t}\underline{s}_{l,t}\rho_{l,t} +  \Delta s_{l,t}\underline{d}_{l,t}r_{l,t} + \Delta d_{l,t}\underline{s}_{l,t}q_{l,t} + \Delta d_{l,t}\Delta s_{l,t}e_{l,t}\Big]: \eqref{h_final_2}-\eqref{h_final_5} \Bigg\}.
        \end{split}
    \end{equation}
    \item[3.2] Record the optimal solution ($\bm{\rho}_t^*, \bm{\lambda}_t^*, \bm{g}_t^*, \bm{z}_t^*$) and optimal value $W_{t}^*$. 
    \item[3.3] Update $UB=\mbox{min}\{UB, \sum_{t=1}^{T}W_{t}^{*}+(LB-\sum_{t=1}^{T}\delta_{t}^{*})\}$. 
    \item[3.4] \textbf{If}$\frac{UB-LB}{UB} \leq \varepsilon$ stop and return $\bm{x}^*$ as the optimal solution of \eqref{min-max_A}; \textbf{else} go to step 4.
\end{itemize}

\textbf{Step 4. Column-and-Constraint Generation Routine}
\begin{itemize}\itemsep0em
\item[4.1] Use $g_{l,t}^{*}$ to compute $d_{l,t}^{N+1}=\underline{d}_{l,t}+( \overline{d}_{l,t}- \underline{d}_{l,t}) g_{l,t}^{*}$ and $z_{l,t}^{*}$ to compute  $s_{l,t}^{N+1} = \underline{s}_{l,t}+(\overline{s}_{l,t}- \underline{s}_{l,t})z_{l,t}^{*}$, for all $l\in [L]$ and $t\in [T]$.
  \item[4.2] Add variables $y_{k,l,t}^{N+1}$, $u_{l,t}^{N+1}$, $o_{k,t}^{N+1}, \forall k\in [K], l\in [L], t\in [T]$, and the following constraints to the master problem:
    \begin{align*}
        &\delta_{t} \geq \sum_{k=1}^{K}\sum_{l=1}^{L} c_{k,l,t} y_{k,l,t}^{N+1} + \sum_{k=1}^{K} c_{k,t}^{o}o_{k,t}^{N+1} + \sum_{l=1}^{L} c_{l,t}^{u}u_{l,t}^{N+1} & -\sum_{l=1}^{L}(&d_{l,t}^{N+1}\alpha_{l,t} + s_{l,t}^{ N+1}\beta_{l,t}), \quad \forall t\in [T],\\
        &\sum_{l\in R_{k}} y_{kl}^{t, N+1} + o_{k}^{t, N+1} = x_{k} h_{k},& &\quad \forall k \in [K], t\in [T],\\
        &\sum_{k: l\in R_{k}} y_{k,l,t}^{N+1} + u_{l,t}^{N+1} = d_{l,t}^{N+1} s_{l,t}^{N+1},& &\quad \forall l\in [L], t\in [T],\\
        &y_{k,l,t}^{N+1}, o_{k,t}^{N+1}, u_{l,t}^{N+1} \in \mathbb{R}_{+},& &\quad \forall k\in [K], l\in [L], t\in [T].
    \end{align*}
    update $N\leftarrow N+1$, $\mathbb{\textbf{O}}\leftarrow \mathbb{\textbf{O}}\cup \lbrace N+1 \rbrace$ and go to step 2.
\end{itemize}
\end{algorithm}
Recall that the F-DHSCP formulation's recourse (second-stage) problem is feasible for every feasible first-stage decision and ($(\bm{d}, \bm{s}) \in \mathcal{U}$). Thus, we do not need to add any feasibility cuts \citep{Zeng_Zhao:2013}. Moreover, we only have a finite number of scenarios. Thus, the algorithm terminates in a finite number of iterations (see \citealp{Zeng_Zhao:2013} for detailed proofs).

\subsubsection{Valid inequalities}\label{enhancing}
Although C\&CG theoretically provides a good way to solve our problem, it sometimes does not perform well from a computational viewpoint. In this section, we aim to incorporate more second-stage information into the master problem by exploiting the properties of the recourse problem.

\noindent (1) Valid lower bound (LB) inequalities

Observe that the Dirac measure on $\pmb{\mu}$ lies in $\mathcal{F}(\mathcal{U}, \pmb{\mu})$. Therefore, we have
\begin{equation} \label{eqn:global_LB}
\underset{\mathbb{P}\in \mathcal{F}(\mathcal{U}, \bm{\mu})}{\mbox{sup}}
    \mathbb{E}_{\mathbb{P}} [Q^{A}(\textbf{\emph{x}}, \bm{\xi})] \geq Q^{A}(\textbf{\emph{x}}, \bm{\mu}) =: L,
\end{equation}
where the lower bound is a deterministic problem with a single scenario $\pmb{\xi}=\bm{\mu}$. This means we can impose the following constraint 
{
\setlength{\abovedisplayskip}{4pt}
\setlength{\belowdisplayskip}{4pt}
\begin{equation}\label{LB_inquality}
    \delta_{t} + \sum_{l=1}^{L}(d_{l,t}\alpha_{l,t} + s_{l,t}\beta_{l,t}) \geq L, 
\end{equation}}
which serves as a global lower bound on the objective of the DRO model. Second, from \eqref{eqn:global_LB}, for any given first-stage decision, the recourse value with $\pmb{\xi}=\bm{\mu}$ provides a lower bound. Therefore, in the initialization step of the C\&CG method, we could include the scenario $\pmb{\xi}=\bm{\mu}$.

\noindent (2) Valid upper and lower bound of variables $\bm{\alpha}$ and $\bm{\beta}$

According to the structure of \eqref{A_inner}, we propose valid upper and lower bounds on variables $\bm{\alpha}$ and $\bm{\beta}$ in Proposition \ref{pos_alpha} (see \ref{proof_alpha} for a detailed proof). In Section \ref{E-of-init}, we demonstrate the computational advantages gained by incorporating these valid inequalities.
\begin{proposition}\label{pos_alpha}
$\overline{\alpha}_{l,t}$ and $\overline{\beta}_{l,t}$ are respectively valid upper bound on variables $\alpha_{l,t}$ and $\beta_{l,t}$. $\underline{\alpha}_{l,t}$ and $\underline{\beta}_{l,t}$, defined in \eqref{a_b_lower}, are respectively valid lower bound on variables $\alpha_{l,t}$ and $\beta_{l,t}$.
{
\setlength{\abovedisplayskip}{4pt}
\setlength{\belowdisplayskip}{4pt}
\begin{equation}\label{a_b_lower}
\setlength{\abovedisplayskip}{6pt}
    \overline{\alpha}_{l,t} = \overline{s}_{l,t} c_{l,t}^u,
    \quad
    \underline{\alpha}_{l,t} = - \overline{s}_{l,t} |\underline{\rho}_{l,t}|,
    \quad
    \overline{\beta}_{l,t} = \overline{d}_{l,t} c_{l,t}^u,
    \quad
    \underline{\beta}_{l,t} = - \overline{d}_{l,t} |\underline{\rho}_{l,t}|.
\end{equation}}
\end{proposition}

\section{Computational experiments}\label{CE}
\noindent In this section, we compare the computational and operational performance of the proposed models for each decision maker and derive insights into the HSCP problem. In Section \ref{setup}, we describe the set of test instances that we constructed and discuss other experimental setups. In Section \ref{ONS-SP}, we obtain near-optimal solutions to the E-SP and F-SP via their SAA formulations. In Section \ref{CPU_analysis}, we compare solution times of the proposed models. In Section \ref{E-of-init}, we demonstrate efficiency of the proposed valid inequalities for the F-DHSCP model. In Section \ref{AOR-EA}, we compare optimal solutions of the E-SP and E-DHSCP models and their out-of-sample performance. Finally, we compare the optimal solutions of the F-DHSCP and F-SP models and their out-of-sample performance in Section \ref{AOR-FA}.

\subsection{Description of Experiments}\label{setup}
To the best of our knowledge, there are no standard benchmark instances in the literature for the HSCP problem. Therefore, we construct 15 HSCP instances based on the same parameters settings and assumptions made in recent related literature \citep[see e.g.][]{lanzarone2014robust, rodriguez2015staff, zhan2018vehicle, zhan2021home, zheng2021stochastic}. Each of these instances is characterized by the number of service types ($L\in \{4,6\}$ as in \cite{mankowska2014home, rodriguez2015staff, rodriguez2018home, zheng2021stochastic}), caregivers types ($K \in \{4,6,8\}$), days in the planning horizon $T\in \{30, 90, 180\}$. We summarize these instances in Table \ref{instance} in \ref{instance_a}.

We set the caregivers' standard work hours  $h_{k}$ to 480 minutes (i.e., 8 hours, see, e.g., \cite{zhan2018vehicle}). We generate the monthly fixed cost for employing one full-time caregivers, $c_{k}$, from $U[3000, 6000]$ \citep{rodriguez2018home}. We generate the quarterly and semi-annual cost,  $c_{k}$, from $U[9000, 18000]$ and $U[18000, 36000]$, respectively. Given that multi-skills (cross-trained) caregivers are more costly, we generated their additional monthly cost from $U[500, 1000]$.  To ensure that the (FA) models will allocate caregivers capacity to different service types, we have $c_{l,t}^{u} + c_{k,t}^{o} > c_{k,l,t}$ for all $k, l$ and $t$. This is because it is always optimal to set $y_{k,l,t} = 0$ when $c_{l,t}^{u} + c_{k,t}^{o} < c_{k,l,t}$.

We followed the same procedures in the stochastic scheduling, home care service, and DRO literature to generate random parameters as follows. We generate the mean values of the demand, $\mu^{d}_{l,t}$, from $U[40; 60]$ and set the standard deviation $\sigma^d_{l,t}=\gamma^d_{l,t}\mu_{l,t}^d$, where $\gamma^{d}_{l,t}$ is generated from $U[0.5, 1]$, for each $l \in [L]$ and $t \in [T]$. Similarly, we generate the mean service time $\mu^{s}_{l,t}$ from $U[40; 60]$ and set $\sigma^{s}_{l,t} = \gamma^{s}_{l,t}\mu^{s}_{l,t}$, for each $l \in [L]$ and $t \in [T]$.

For each instance, we sample $N$ scenarios of $d_{1,1}^{n}, \ldots, d_{L,T}^{n}$, for $n=1, \ldots, N$, by following lognormal (LogN) distributions with the generated  $\pmb{\mu}^d$ and $\pmb{\sigma}^d$. We truncated these scenarios on basic demand range $[\underline{d}, \overline{d}$]=[40, 60] \citep{rodriguez2015staff, rodriguez2018home}. Similarly, for each instance, we generated $N$ scenarios of $s_{1,1}^{n}, \ldots, s_{L,T}^{n}$, for $n=1, \ldots, N$, by following LogN distributions with the generated $\pmb{\mu}^s$ and $\pmb{\sigma}^s$. We truncated these scenarios on $[\underline{s}, \overline{s}]$=U[20, 80].  We round each generated parameter to the nearest integer \citep{tsang2021distributionally}. We chose the LogN distribution to generate the in-sample data because it is typically used to model service time and demand in the home care, scheduling, and other literature \citep{cayirli2006designing, jiang2019data, tsang2021distributionally}. 

According to staffing standards for various types of elderly care institutions in China, the minimum staffs number for serving less than 100 customers is three or five. Hence, we set the minimum number of caregivers, $\underline{w}$, to 3. We implement all models and solution methods using the Python 3.8.8 programming language calling Gurobi V9.1.2 as a solver with default settings. We ran all experiments on a computer with Intel(R) Xeon(R) E3-1231-v3 processor with two 3.40GHz CPUs and 32 GB shared RAM memory.

\subsection{Obtaining near-optimal solutions to the SP models}\label{ONS-SP}

For each instance in Table \ref{instance}, we implemented the MCO algorithm in Section \ref{MCO} with $N_o=10$, $N'=1000$ and $K=10$. Table \ref{MCO_result} in \ref{CR_1} and Table \ref{MCO_result_2} in \ref{CR_2} present the approximate optimality index ($AOI_N=\frac{\overline{v}_{N^\prime}-\overline{v}_{N}}{\overline{v}_{N^\prime}}$) and the 95\% confidence intervals (95\%CI) of the statistical lower and upper bounds ($\overline{v}_{N}$ and $\overline{v}_{N^\prime}$, receptively) on the objective values of the E-SP and F-SP models  for each instance under $(c_{l,t}^{u}, c_{l,t}^{o} (c_{k,t}^{o})) = (20, 2)$ and truncated LogN distribution. We observe similar results under different choices of these cost parameters.  In the next section, we present and analyze the ranges of solution time of the SAA formulations. Herein we summarize the key findings of the MCO algorithm. Clearly, $N=320$ and $N=100$ scenarios of the demand and service duration are sufficient to obtain a near-optimal solution to the E-SP and F-SP models via their SAA formulations, respectively. First, $AOI_N$ at $N=320$ and $N=100$ ranges from 0 to 0.0098 and from 0 to 0.0055 for the E-SP and F-SP models, respectively. Second, the 95\%CI of $\overline{v}_{N}^{E}$ ($\overline{v}_{N}^{F}$) and $\overline{v}_{N^\prime}^{E}$ ($\overline{v}_{N^\prime}^{F}$) with $N=320$ (100) are very tight (i.e., have a small variance). These results qualify $\overline{v}_{N}^{E}$ and $\overline{v}_{N}^{F}$ as tight statistical estimates for the lower and upper bounds on the optimal value of each instance, respectively. In addition, we also use truncated normal distribution and uniform distribution to generate the sample data. We obtain similar conclusion for these distributions (see \ref{CR_1} and \ref{CR_2}).

It is worth mentioning that solving these models with larger $N$ resulted in longer solution times without consistent and significant improvements in the relative approximation gaps. However, as shown in \ref{CR_1}--\ref{CR_2}, solution times are still reasonable under $N=500$ scenarios. And the $AOI_N$ is close to the value at termination sample size.

\subsection{CPU time and computational details}\label{CPU_analysis}
In this section, we analyze the solution times of the proposed SP and DRO models. In addition to the based range of the demand described in Section~\ref{setup}, we consider another range $[40, 100]$, which  corresponding to a higher amound of demand and variation level.  For each of the 15 instances in Table \ref{instance} and demand range, we generated and solved 10 instances as described in Section \ref{setup} for a total of 300 random instances. Then, we solved these instances using the proposed models for each decision-maker.

\begin{table}[t!]
  \centering
  \renewcommand{\arraystretch}{0.6}
  \small
  \caption{CPU times in seconds for EA model}
    \begin{tabular}{cccccccccccccc}
    \toprule
    \multirow{3}[6]{*}{Inst} & \multicolumn{6}{c}{Demand range: [40, 60]} &   & \multicolumn{6}{c}{Demand range: [40, 100]} \\
\cmidrule{2-7}\cmidrule{9-14}      & \multicolumn{3}{c}{E-DHSCP} & \multicolumn{3}{c}{E-SP} &   & \multicolumn{3}{c}{E-DHSCP} & \multicolumn{3}{c}{E-SP} \\
\cmidrule{2-7}\cmidrule{9-14}      & Min & Avg & Max & Min & Avg & Max &   & Min & Avg & Max & Min & Avg & Max \\
    \midrule
    1 & 0.1  & 0.1  & 0.1  & 1.8  & 2.1  & 2.4  &   & 0.1  & 0.1  & 0.1  & 2.0  & 2.2  & 2.4  \\
    2 & 0.3  & 0.3  & 0.4  & 5.9  & 6.2  & 6.8  &   & 0.3  & 0.3  & 0.4  & 6.0  & 6.3  & 7.3  \\
    3 & 0.7  & 0.8  & 0.9  & 12.1  & 12.8  & 13.5  &   & 0.6  & 0.7  & 0.8  & 12.5  & 13.4  & 14.7  \\
    4 & 0.1  & 0.1  & 0.2  & 1.9  & 2.2  & 2.5  &   & 0.1  & 0.1  & 0.1  & 2.0  & 2.1  & 2.3  \\
    5 & 3.7  & 4.0  & 4.4  & 5.9  & 6.3  & 6.8  &   & 0.2  & 0.3  & 0.3  & 5.7  & 6.0  & 6.4  \\
    6 & 2.5  & 2.8  & 3.2  & 13.3  & 14.4  & 16.0  &   & 0.5  & 0.6  & 0.7  & 12.9  & 13.6  & 14.1  \\
    7 & 0.1  & 0.1  & 0.1  & 2.1  & 2.1  & 2.3  &   & 0.1  & 0.1  & 0.1  & 2.1  & 2.3  & 2.5  \\
    8 & 0.3  & 0.3  & 0.4  & 6.4  & 6.6  & 6.8  &   & 0.4  & 0.4  & 0.4  & 6.5  & 6.7  & 6.9  \\
    9 & 0.7  & 0.7  & 0.9  & 12.5  & 13.0  & 14.0  &   & 0.8  & 0.8  & 0.9  & 13.7  & 14.0  & 14.5  \\
    10 & 0.1  & 0.1  & 0.2  & 2.8  & 2.9  & 3.0  &   & 0.1  & 0.1  & 0.1  & 2.8  & 2.8  & 3.0  \\
    11 & 0.5  & 0.5  & 0.5  & 8.7  & 8.8  & 8.9  &   & 0.5  & 0.5  & 0.5  & 8.8  & 8.9  & 9.1  \\
    12 & 0.9  & 1.0  & 1.0  & 18.1  & 18.6  & 19.4  &   & 1.0  & 1.0  & 1.2  & 16.7  & 17.5  & 19.2  \\
    13 & 0.1  & 0.1  & 0.2  & 2.7  & 2.7  & 2.8  &   & 0.1  & 0.1  & 0.1  & 2.7  & 2.7  & 2.8  \\
    14 & 0.4  & 0.4  & 0.5  & 8.5  & 9.1  & 11.4  &   & 0.3  & 0.4  & 0.4  & 8.6  & 8.9  & 9.1  \\
    15 & 0.9  & 1.0  & 1.0  & 17.8  & 18.1  & 18.6  &   & 0.8  & 0.9  & 1.0  & 17.4  & 18.7  & 20.4  \\
    \bottomrule
    \end{tabular}%
  \label{CPU_time_EA}%
\end{table}%

\begin{table}[t!]
  \renewcommand{\arraystretch}{0.6}
  \centering
  \small
  \caption{CPU times in seconds for FA model}
    \begin{tabular}{ccccccccccccccccc}
    \toprule
    \multirow{3}[6]{*}{Inst} & \multicolumn{6}{c}{Demand range: [40, 60]} &   & \multicolumn{6}{c}{Demand range: [40, 100]} \\
\cmidrule{2-7}\cmidrule{9-14}      & \multicolumn{3}{c}{F-DHSCP} & \multicolumn{3}{c}{F-SP} &   & \multicolumn{3}{c}{F-DHSCP} & \multicolumn{3}{c}{F-SP} \\
\cmidrule{2-7}\cmidrule{9-14}      & Min & Avg & Max & Min & Avg & Max &   & Min & Avg & Max & Min & Avg & Max \\
    \midrule
    1 & 0.5  & 0.5  & 0.6  & 2.8  & 3.0  & 3.3  &   & 0.9  & 0.9  & 1.0  & 2.8  & 3.0  & 3.3  \\
    2 & 1.4  & 1.6  & 1.9  & 16.0  & 17.9  & 21.3  &   & 1.3  & 1.5  & 1.7  & 16.6  & 18.3  & 20.7  \\
    3 & 2.8  & 2.8  & 2.9  & 52.2  & 55.9  & 61.0  &   & 4.0  & 4.1  & 4.2  & 55.7  & 59.3  & 62.4  \\
    4 & 1.8  & 1.9  & 1.9  & 2.9  & 3.0  & 3.1  &   & 1.3  & 1.4  & 1.5  & 3.2  & 3.3  & 3.4  \\
    5 & 2.2  & 2.2  & 2.2  & 12.8  & 13.2  & 13.7  &   & 5.0  & 5.2  & 5.3  & 15.0  & 16.0  & 17.2  \\
    6 & 3.0  & 3.1  & 3.3  & 41.5  & 43.6  & 45.5  &   & 40.5  & 41.0  & 42.1  & 48.2  & 52.6  & 58.1  \\
    7 & 1.0  & 1.0  & 1.1  & 5.8  & 6.1  & 6.8  &   & 1.2  & 1.2  & 1.3  & 5.1  & 5.5  & 6.2  \\
    8 & 3.4  & 3.5  & 3.6  & 36.7  & 41.4  & 45.6  &   & 3.6  & 3.7  & 3.7  & 31.7  & 34.3  & 35.6  \\
    9 & 3.8  & 3.9  & 4.0  & 139.6  & 148.4  & 159.8  &   & 69.2  & 70.3  & 71.1  & 104.9  & 112.2  & 116.5  \\
    10 & 1.8  & 1.8  & 1.9  & 5.4  & 5.4  & 5.5  &   & 1.1  & 1.1  & 1.2  & 5.2  & 5.3  & 5.4  \\
    11 & 2.7  & 2.8  & 2.8  & 32.8  & 35.3  & 43.2  &   & 1.9  & 2.0  & 2.1  & 32.6  & 33.7  & 34.6  \\
    12 & 4.1  & 4.2  & 4.3  & 114.8  & 127.5  & 141.1  &   & 4.0  & 4.1  & 4.1  & 111.1  & 119.2  & 127.3  \\
    13 & 1.0  & 1.0  & 1.1  & 5.6  & 5.8  & 6.0  &   & 0.7  & 0.7  & 0.8  & 6.1  & 6.2  & 6.5  \\
    14 & 4.8  & 4.9  & 5.1  & 27.5  & 28.6  & 30.2  &   & 8.5  & 8.6  & 8.7  & 31.4  & 32.1  & 32.9  \\
    15 & 4.7  & 4.8  & 5.0  & 86.2  & 90.3  & 94.3  &   & 10.4  & 10.5  & 10.6  & 101.0  & 105.0  & 112.3  \\
    \bottomrule
    \end{tabular}%
  \label{CPU_time_FA}%
\end{table}%

Let us first compare the computational performance of the  E-SP and E-DHSCP models. Table \ref{CPU_time_EA} presents the minimum (Min), average (Avg), and maximum (Max) CPU time (in seconds) for these models. We first observe that the computational effort tend to increase with the scale ($L\times K\times T$). Second, both models can quickly solve all of the instances to optimality. However, the E-DHSCP model takes a shorter time than the E-SP model to solve each instance. These results are consistent under both ranges of the demand. 

Next, we compare the solution times of the F-SP and F-DHSCP models presented in Table~\ref{CPU_time_FA}. We make the following observations from this table. First, it is clear that solution times of these models are much longer than the E-SP and E-DHSCP models. For example, the average solution time of instance 15 (largest instance) using the E-SP and F-SP models under demand range [40,60] are 18.1 and 90.3 seconds, respectively. And the average solution time of this instance using the E-DHSCP and F-DHSCP models are 1 and 4.8 seconds, respectively. Since the (EA) and (FA) models have different structures (i.e., decision variables and objectives), these results make sense. Note that the F-SP model has a larger second stage than the E-SP model. In addition, we solve the E-DHSCP directly using Groubi, but we employ C\&CG to solve the E-DHSCP model. Second, again the DRO model (F-DHSCP) can solve each instance faster than the SP model (F-SP). Solution times of the F-SP and F-DHSCP ranges from 2.8-159.8 and 0.5-71.1 seconds, respectively. 

Overall, these results demonstrate the computational efficiency of the proposed approaches. To confirm our conclusion for the scenario-based SAA models, we solve the SAA formulations of E-SP and F-SP models with 500 samples. We solve 6 instances with $L=6$ under these settings and present solution times in  Tables \ref{CPU_500} (see \ref{CR_1}) and Table \ref{CPU_500_2} (see \ref{CR_2}). Clearly, the E-SP model can solve these instances within 30 seconds, and the F-SP model can solve all of them within 30 minutes.

\subsection{Efficiency of valid inequalities}\label{E-of-init}
\begin{figure}[t!]
    \centering
    \subcaptionbox{LB values\label{LB_m}}{
        \includegraphics[scale=0.3508]{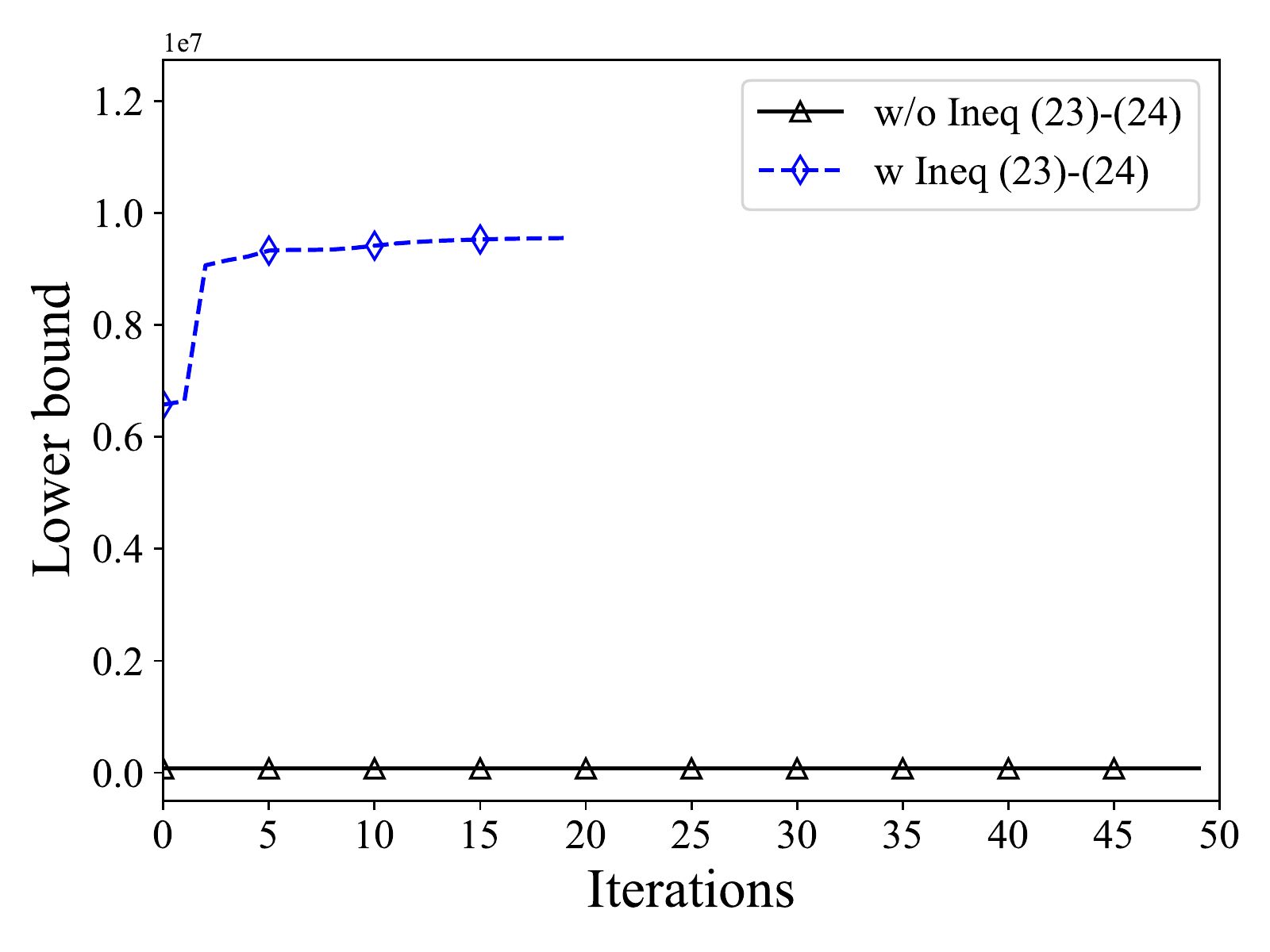}}
    \subcaptionbox{Gap values\label{gap_m}}{
        \includegraphics[scale=0.35]{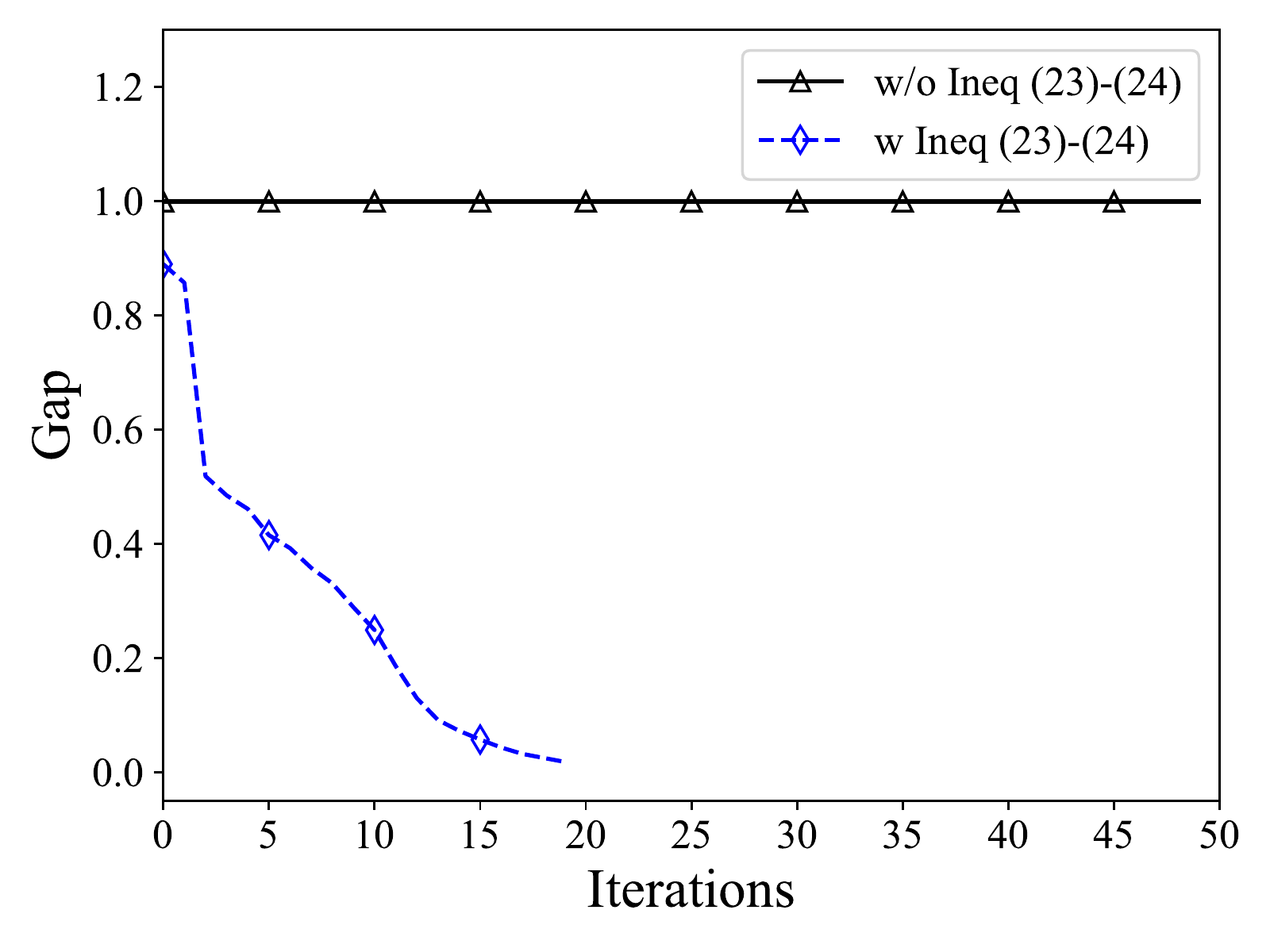}}
    \caption{Comparisons of lower bound and gap values with and without valid inequality \eqref{LB_inquality}--\eqref{a_b_lower}}
    \label{FA-LB-efficiency}
\end{figure}

In this section, we investigate the efficiency of the proposed valid inequalities \eqref{LB_inquality}--\eqref{a_b_lower} for the E-DHSCP and C\&CG in Section \ref{enhancing}. Given that we observe similar results for all instances, for brevity and illustrative purposes, we use Instance 15 ($L=6, K=8, T=180$) in this experiment. We separately solve the E-DHSCP model via C\&CG with (w/) and without (w/o) inequalities \eqref{LB_inquality}--\eqref{a_b_lower}. We first note the algorithm takes hundreds of iterations and a longer time until convergence w/o inequalities \eqref{LB_inquality}--\eqref{a_b_lower}. Therefore, in Figure \ref{FA-LB-efficiency}, we present the lower bound and gap values w/ and w/o inequalities \eqref{LB_inquality}--\eqref{a_b_lower} from the first 50 iterations. It is obvious that both the lower bound and gap values converge faster when we introduce inequalities \eqref{LB_inquality}--\eqref{a_b_lower} into the master problem. Moreover, because of the better bonding effect, the algorithm converges to the optimal solution in fewer iterations and shorter solution times. For example, the algorithm takes 531 seconds and 19 iterations to solve this instance w/ inequalities \eqref{LB_inquality}--\eqref{a_b_lower} and terminates at the specified time limit (i.e., 7200 seconds) with gap 99.95\% w/o these inequalities.

\subsection{Analysis of results of EA models}\label{AOR-EA}

In this section, we analyze the optimal staffing patterns of the E-SP and E-DHSCP models (Section~\ref{OSP-EA}) and their out-of-sample performances (Section~\ref{OOP-EA-M}). For illustrative purposes and brevity, we present results for Instance 7 ($L=4, K=8, T=30$) and Instance 10 ($L=6, K=6, T=30$), which corresponds to middle-sized instances. In addition, we consider the following types of caregivers. In Instance 7 (which consists of $K=8$ types of caregivers), we let caregivers types $k \in \{1,\ldots,4\}$ represent \textit{specialized} caregivers, meaning they can only provide one type of service. And we let type $k \in \{5,\ldots,8\}$ represents \textit{cross-trained} caregivers who are trained or qualified for two types of services. In Instance 10 (which consist of $K=6$ types of caregivers), all caregivers are cross-trained.

\subsubsection{Optimal staffing patterns of EA models}\label{OSP-EA}
\begin{figure}[t!]
    \centering
    \subcaptionbox{$c_{l,t}^o = 1$\label{EA-10-6-6-30}}{
        \includegraphics[scale=0.35]{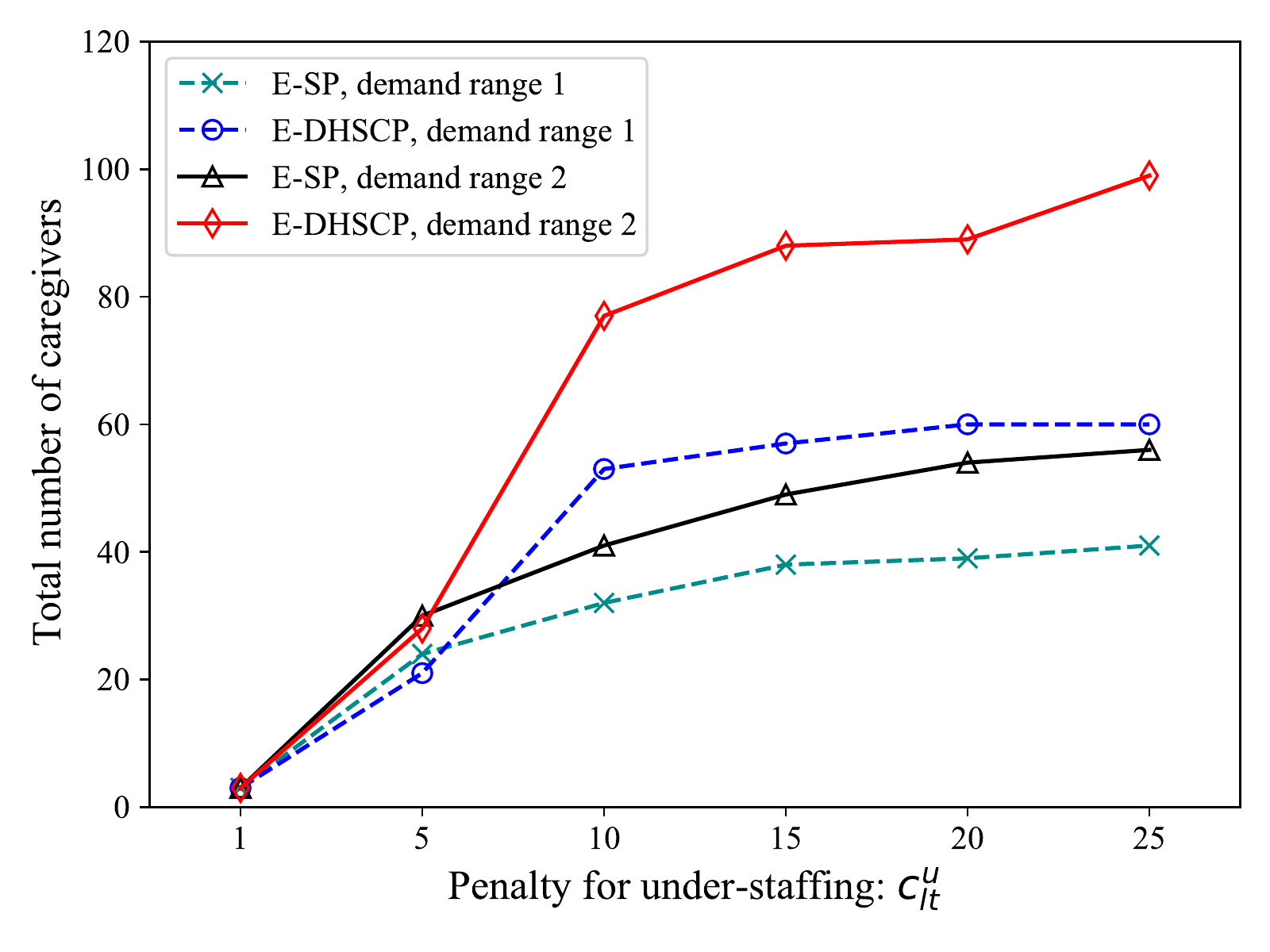}}
    \subcaptionbox{$c_{l,t}^u = 10$\label{EA-10-6-6-30-o}}{
        \includegraphics[scale=0.35]{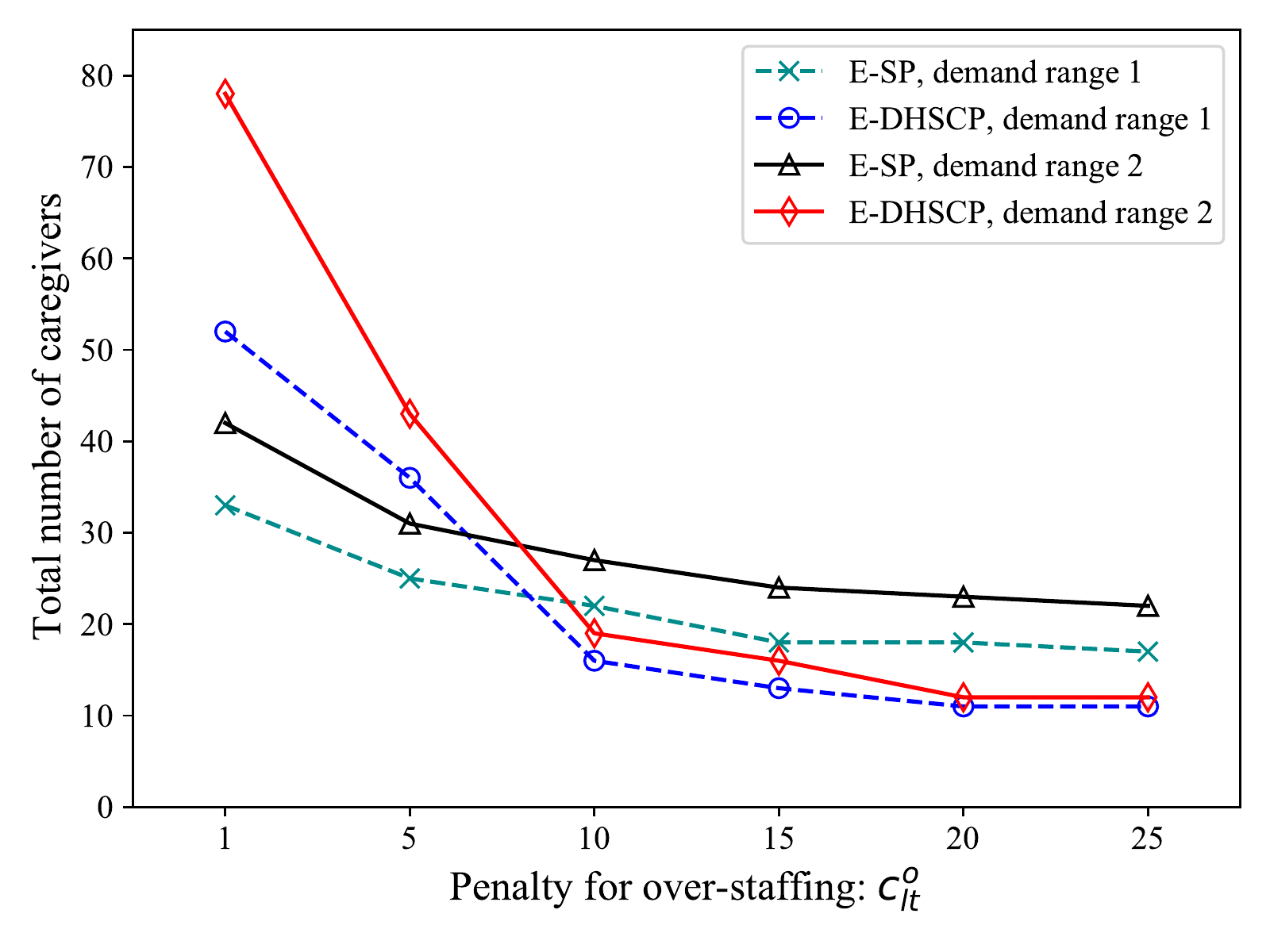}}
    \caption{Total number of caregivers hired by E-SP and E-DHSCP models for Instance 10.}
    \label{total_num_fig}
\end{figure}
In this section, we compare the optimal staffing pattern (i.e., the optimal number of hired caregivers) of the E-SP and E-DHSCP models under different values of the under-staffing $c_{l,t}^u$ and over-staffing $c_{l,t}^o$ penalty costs. First, in Figure~\ref{EA-10-6-6-30}, we present the optimal number of hired caregivers under demand range 1 ($[40, 60]$) and range 2 ($[40, 100]$) with $c_{l,t}^o=1$ and $c_{l,t}^u \in \{1, 5,  10, 15, 20, 25\} $ for Instance 10. We observe the same staffing patterns for Instance 7 (see Figure \ref{total_num_fig_7} in \ref{OSP_A}).

We observe the following from Figure~\ref{EA-10-6-6-30}. First, both models hire more caregivers under range 2 than range 1. This makes sense because, in range 2, both the volume and variability of the demand are higher. Second, it is clear that both models hire more caregivers as the under-staffing cost increases from 1 to 25. However, the number of hired caregivers stabilize after $c_{l,t}^u=10$. Third, when  $c_{l,t}^u \leq 5$ (i.e., low under-staffing cost), both models hire approximately the same number of caregivers. In contrast, when $c_{l,t}^u>5$, the E-DHSCP models hire more caregivers than the E-SP model to mitigate excessive under-staffing. These results indicate that $c_{l,t}^u = 10$ may be a suitable penalty cost of under-staffing for (EA) decision-makers. We observe a different pattern when we fix $c_{l,t}^u=10$ and increase the penalty for over-staffing  $c_{l,t}^o$ from 1 to $\{5,  10, 15, 20, 25\}$. Specifically, as shown in Figure~\ref{EA-10-6-6-30-o} both models hire less caregivers as increasing $c_{l,t}^o$. And the number of hired caregivers stabilize after $c_{l,t}^o=10$.

Next, we analyze the detailed staffing patterns under cost structure $(c_{l,t}^{u}, c_{l,t}^{o}) = (10, 1)$ presented in Figure \ref{EA-solution} (see \ref{OSP_A}). First, if there are only cross-trained caregivers (Instance 10), all models tend to hire all types of caregivers. Second, if there are both specialized and cross-trained caregivers (Instance 7), the E-SP model tend to hire specialized caregivers under both demand ranges while the E-DHSCP model hire a mixture of specialized and cross-trained caregivers. In addition, the E-DHSCP model hire more cross-trained caregivers under demand range 2 than under demand range 1. This could indicate that cross-trained caregivers are more favorable under a higher variation level.

\subsubsection{Out-of-sample performance of EA model}\label{OOP-EA-M}
In this section, we compare the operational performance of the optimal solutions of the  E-DHSCP and E-SP via out-of-sample simulation. We test the out-of-sample performance (i.e., the objective value obtained by simulating the optimal solution of a model under a larger unseen data) as follows. First, we solve each instance and obtain the optimal first-stage decisions. Second, we fix the obtained first-stage decisions, then re-optimize the second stage of the SP using the following sets of $N'=10,000$ out-of-sample data to compute the second stage cost (i.e., under-staffing plus over-staffing costs). We compute the out-of-sample total cost as first-stage cost plus the out-of-sample second-stage cost.
\begin{itemize}\itemsep0em
\item[Set 1.] We assume that we have perfect distributional information. That is, we use the same distributions (LogN) and parameter settings that we use for generating the $N$ data  in the optimization to generate $N^\prime$ data. This is to simulate the performance when the true distribution is the same as the one used in the optimization.

\item[Set 2.]  In this set, we assume that we have used misspecified distributions of the demand and service duration in the optimization (motivated by our analysis in Section~\ref{motivation}). We follow a similar out-of-sample simulation testing procedure described in \cite{shehadeh2020DMFRS} and \cite{wang2020distributionally} and perturb the support of the random demand and service durations by a parameter $\Delta$ and used parameterized uniform distribution $U[(1 - \Delta)$ lower bound, $(1 + \Delta)$ upper bound] to generate the $N^\prime=10,000$ data. We apply $\Delta \in \{0, 0.1, 0.25, 0.5\}$. A higher $\Delta$ corresponds to a higher variation level.

\end{itemize}

For brevity, we focus on Instance 10 and present results under $(c_{l,t}^{u}, c_{l,t}^{o}) = (10, 1)$. We observe similar results for Instance 7 (see \ref{oop-EA-7}). First, let us analyze simulation results under Set 1 (i.e., perfect distributional information case). Figure \ref{EA-oop-in-6} in \ref{oop-EA-10} presents histograms of the out-of-sample total cost, second-stage costs, and under-staffing. It is clear from this figure that the E-DHSCP model yields a higher total cost and lower second-stage cost than the E-SP model on average and at all quantiles. This makes sense because, as we discussed in Section \ref{OSP-EA}, the E-DHSCP tends to hire more caregivers to hedge against under-staffing (reflected by significantly lower under-staffing on average and all quantile than the E-SP model in Figure \ref{EA-inu-10-2-6-6-30-4060}). However, this results in a higher fixed cost associated with hiring more caregivers.

\begin{figure}[t!]
    \centering
    \subcaptionbox{$\Delta = 0$\label{EA-outc2-10-2-6-6-30-4060}}{
        \includegraphics[scale=0.3]{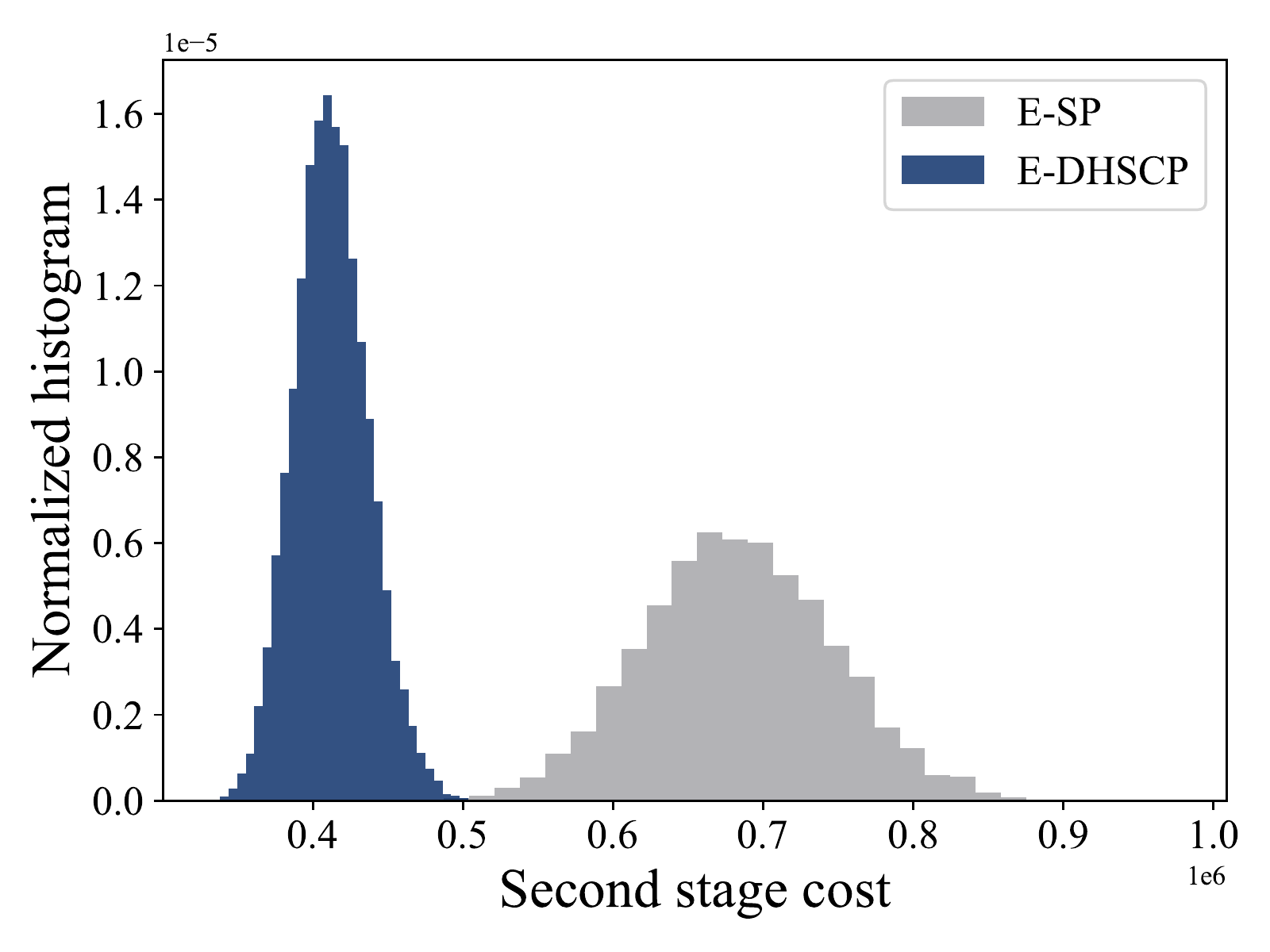}}
    \subcaptionbox{$\Delta = 0.25$\label{EA-out025c2-10-2-6-6-30-4060}}{
        \includegraphics[scale=0.3]{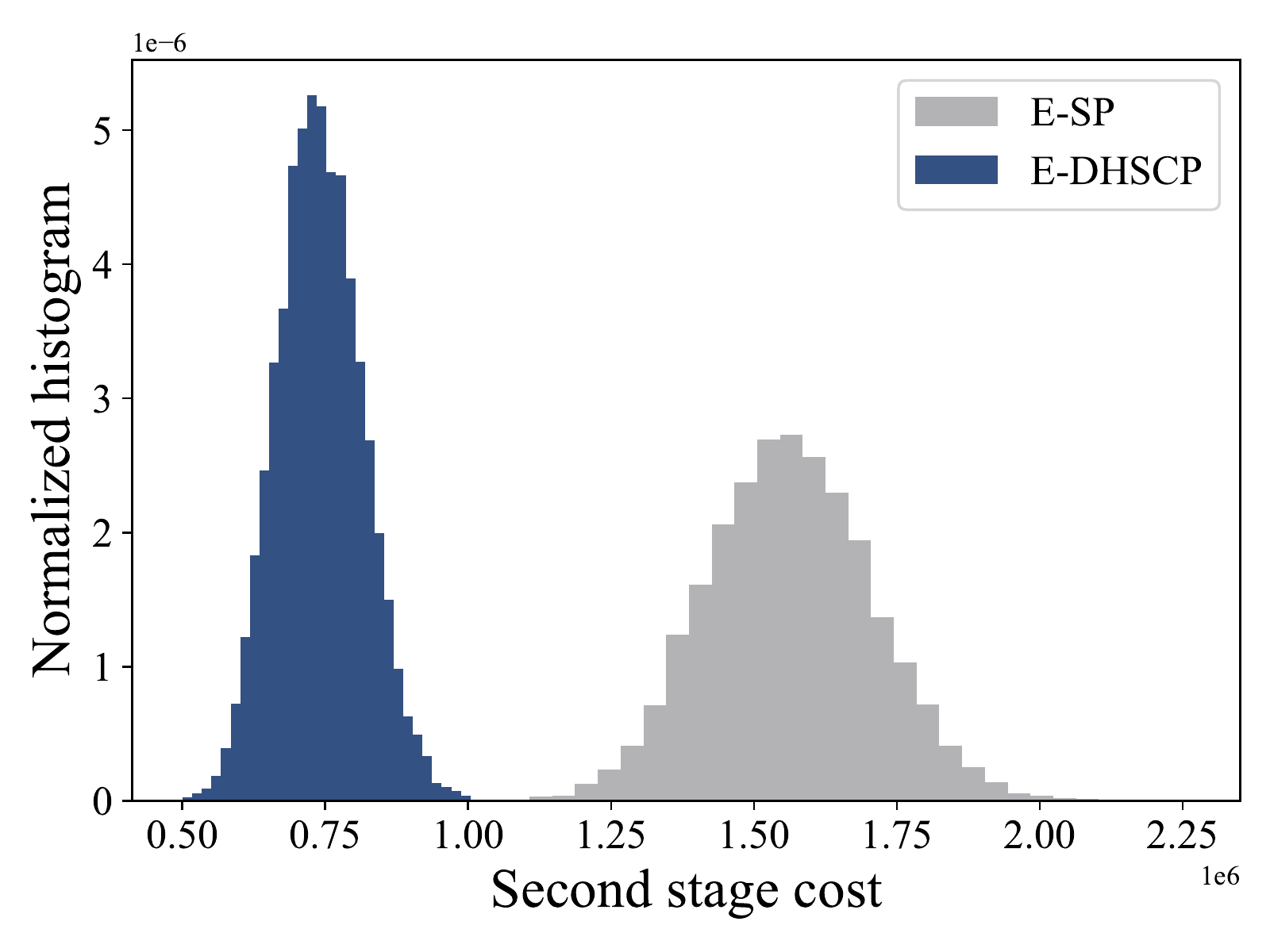}}
    \subcaptionbox{$\Delta = 0.5$\label{EA-out05c2-10-2-6-6-30-4060}}{
        \includegraphics[scale=0.3]{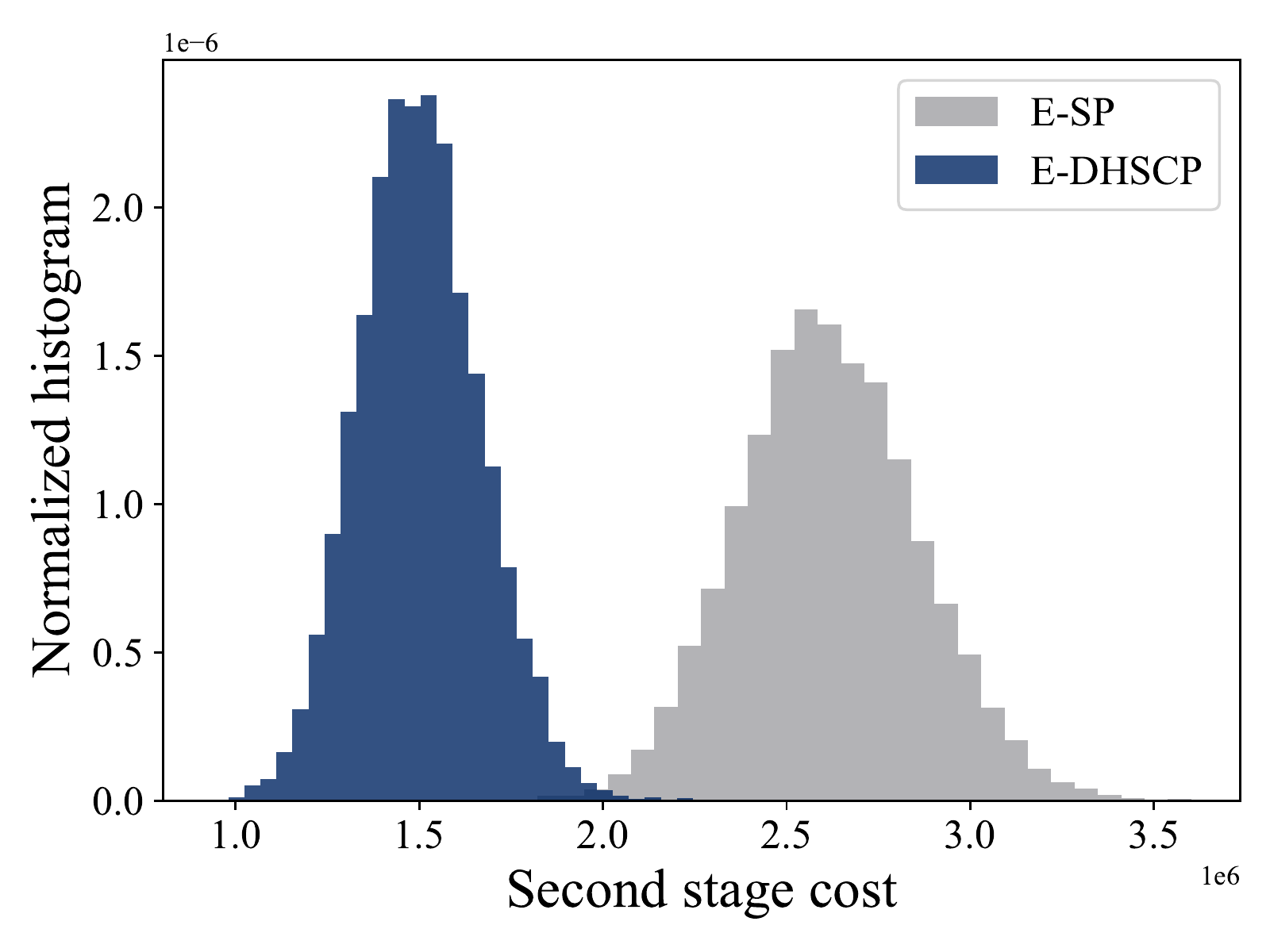}}
    \caption{Out-of-sample second stage cost for Instance 10, demand range 1 under Set 2}
    \label{EA-oop-6}
\end{figure}

We observe the following from simulation results under Set 2 (i.e., misspecified distributional information case) with demand range 1 presented in Figures \ref{EA-oop-6} (see Figure \ref{EA-oop-6-1080} in \ref{oop-EA-10} for results with range 2).  Notably, the E-DHSCP model consistently outperforms the E-SP model under all levels of variation ($\Delta$) and across the criteria of mean and all quantiles of the second-stage costs under the two ranges of the demand. In particular, as shown in Figure~\ref{EA-oop-u-6-4060} (see \ref{oop-EA-10}), the E-DHSCP yield significantly lower under-staffing cost than the E-SP model. Note that a lower second-stage and under-staffing cost indicates better operational performance (i.e., lower shortage and meeting greater demand) and thus has a significant practical impact, especially in healthcare.  Second, the E-DHSCP solutions appear to be more stable with a significantly smaller standard deviation  (i.e., variations) in the total and second-stage costs than the E-SP solutions. The superior performance of the DRO model reflects the value of modeling uncertainty and distributional ambiguity. 

Next, we investigate the value of distributional robustness from the perspective of out-of-sample disappointment, which measures the extent to which the out-of-sample cost exceeds the model's optimal value \citep{van2021data, wang2020distributionally}. Following \cite{shehadeh2020DMFRS} and \cite{wang2020distributionally}, we define the  out-of-sample disappointment as $\mbox{max} \Big\{0, \frac{V_{out} - V_{opt}}{V_{opt}}\Big\}\times 100\%$, where $V_{opt}$ and $V_{out}$ are as the model's optimal value and the out-of-sample objective value, respectively. A disappointment of zero indicates that the model's optimal value is equal to or larger than the out-of-sample (actual) cost, suggesting that the model is more conservative and avoids underestimating costs. In contrast, a larger disappointment implies a higher level of over-optimism because, in this case, the actual cost ($V_{out}$) of implementing the optimal solution of a model is larger than the estimated cost ($V_{opt}$). 

\begin{figure}[t!]
    \centering
    \subcaptionbox{$\Delta=0$\label{EA-outdis-10-2-6-6-30-4060}}{
        \includegraphics[scale=0.3]{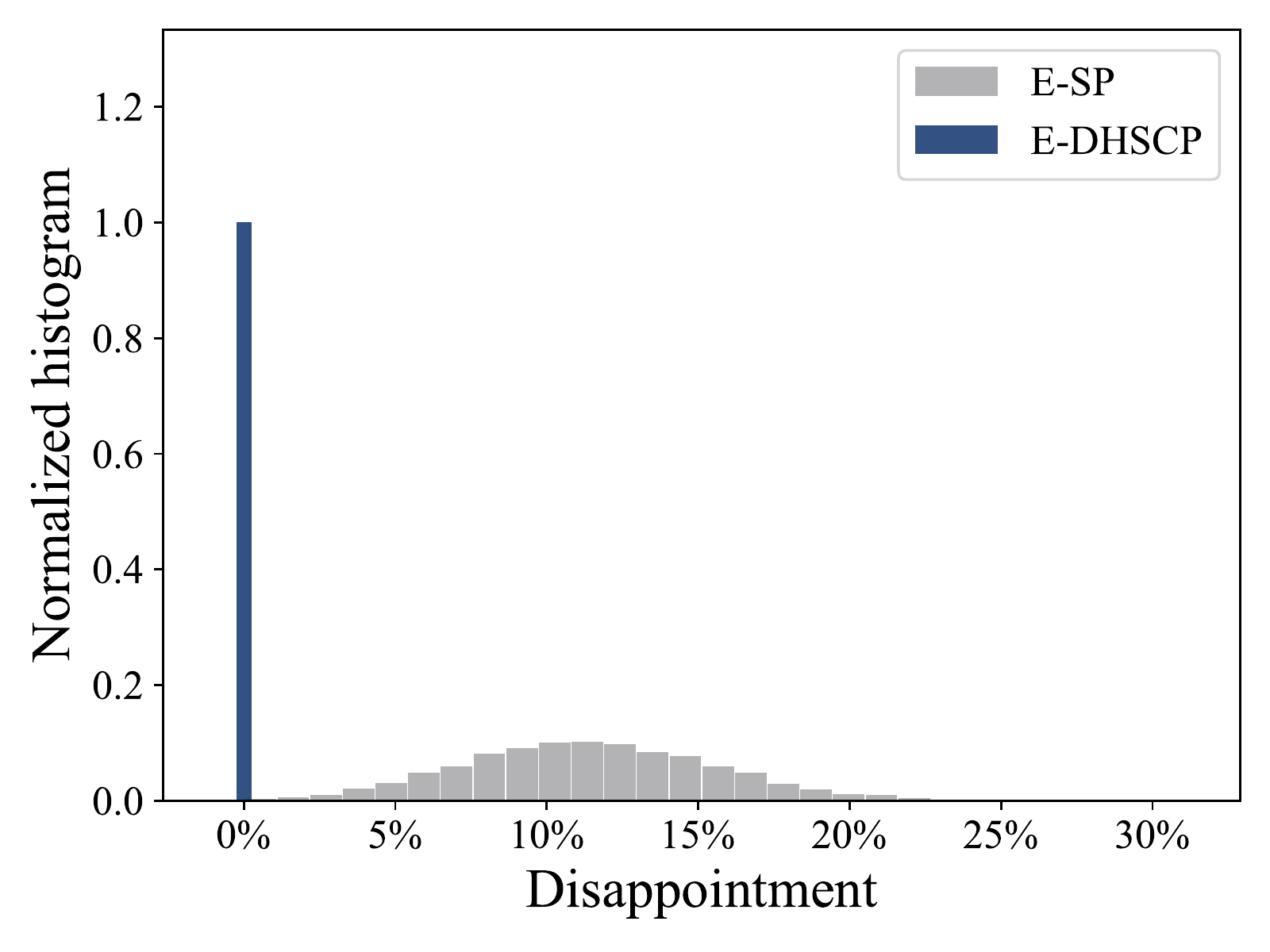}}
    \subcaptionbox{$\Delta=0.25$\label{EA-out025dis-10-2-6-6-30-4060}}{
        \includegraphics[scale=0.3]{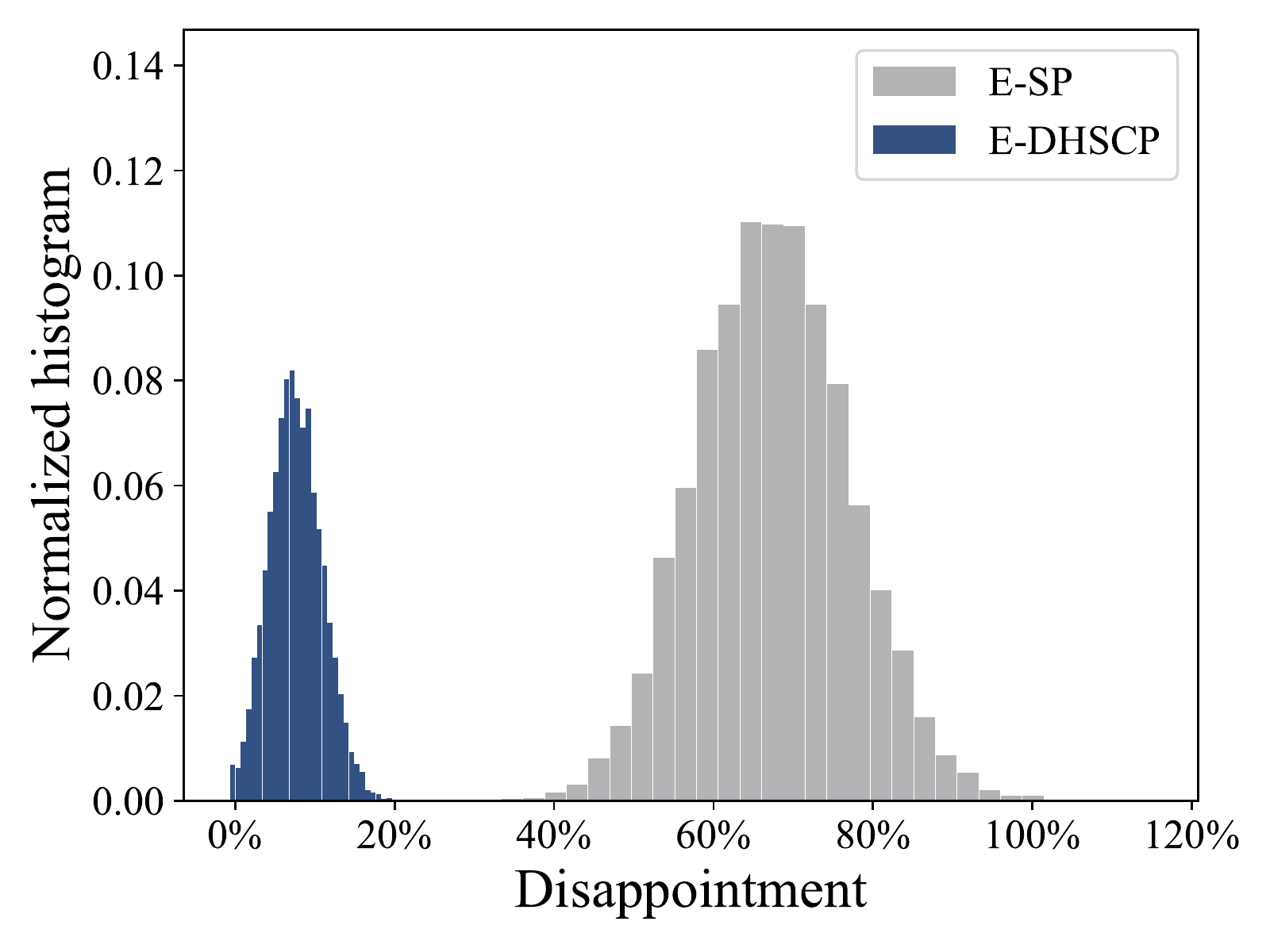}}
    \subcaptionbox{$\Delta=0.5$\label{EA-out05dis-10-2-6-6-30-4060}}{
        \includegraphics[scale=0.3]{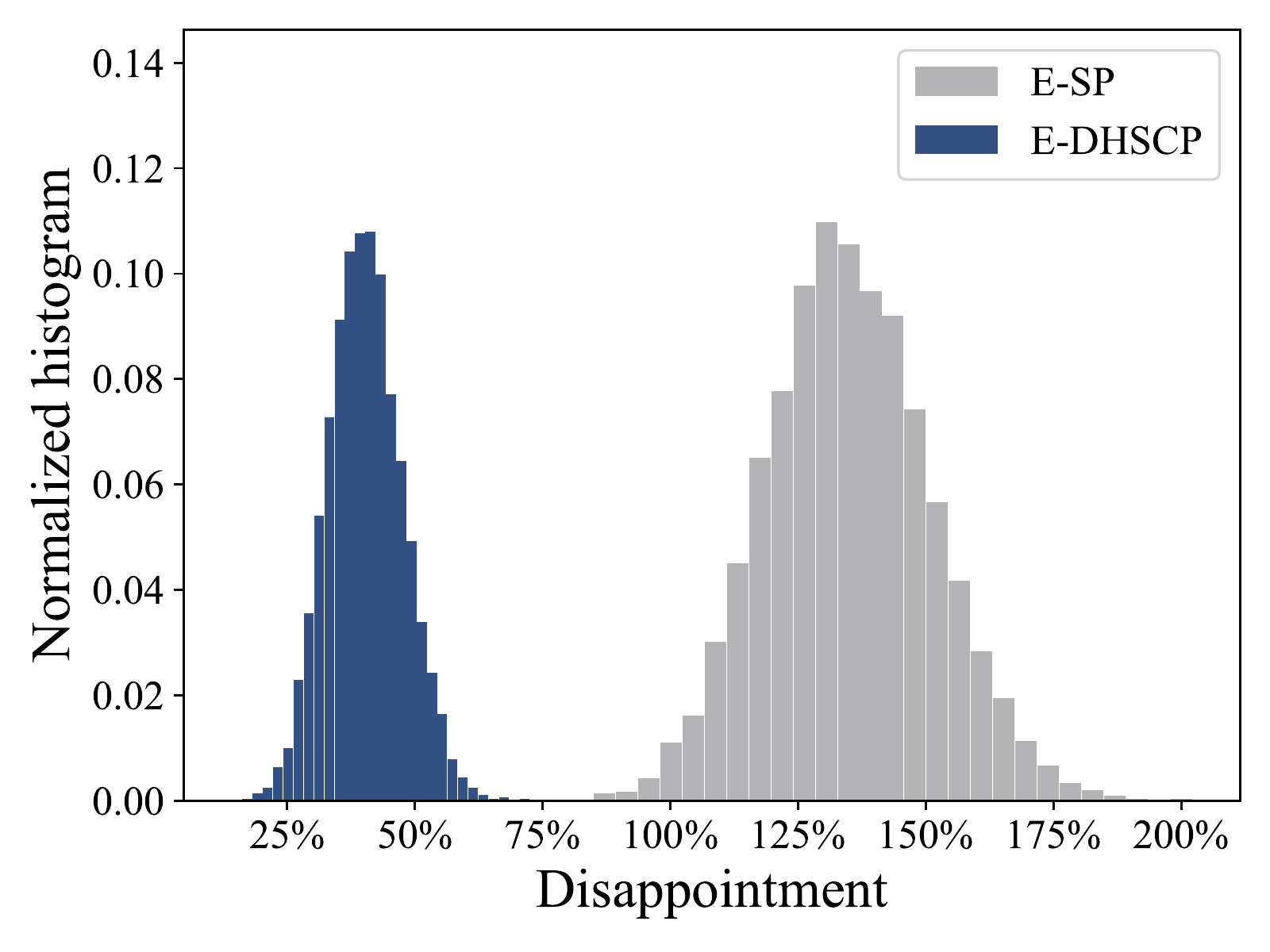}}
    \caption{Out-of-sample disappointment for Instance 10 with demand range 1 under Set 2.}
    \label{EA-oop-dis-6}
\end{figure}

Figure~\ref{EA-oop-dis-6} presents histograms of the out-of-sample disappointment. It is clear from this figure that the E-DHSCP  solutions yield substantially smaller out-of-sample disappointments on average and at all quantiles than the E-SP model. In particular, the disappointment of the E-SP model significantly increases under higher variation, exceeding 100\% under $\Delta=0.5$. In contrast, the E-DHSCP model can control disappointment with a smaller and stable range. These results confirm that the DRO model provides a more robust estimate of the actual cost that we will incur in practice. In a decision-making context where the goal is to minimize expenses, disappointments are more harmful than positive surprises (overestimated costs).

\subsection{Analysis of results of FA models}\label{AOR-FA}
In this section, we analyze the optimal staffing patterns of the F-SP and F-DHSCP models (Section \ref{OSP-FA}) and their out-of-sample performance (Section \ref{OOP-FA}). We present the results for Instance 7 and Instance 10 as defined in Section \ref{AOR-EA}.
\subsubsection{Optimal staffing patterns of FA models}\label{OSP-FA}

\begin{figure}[t!]
    \centering
    \subcaptionbox{$c_{k,t}^o = 1$\label{FA-10-6-6-30}}{
        \includegraphics[scale=0.35]{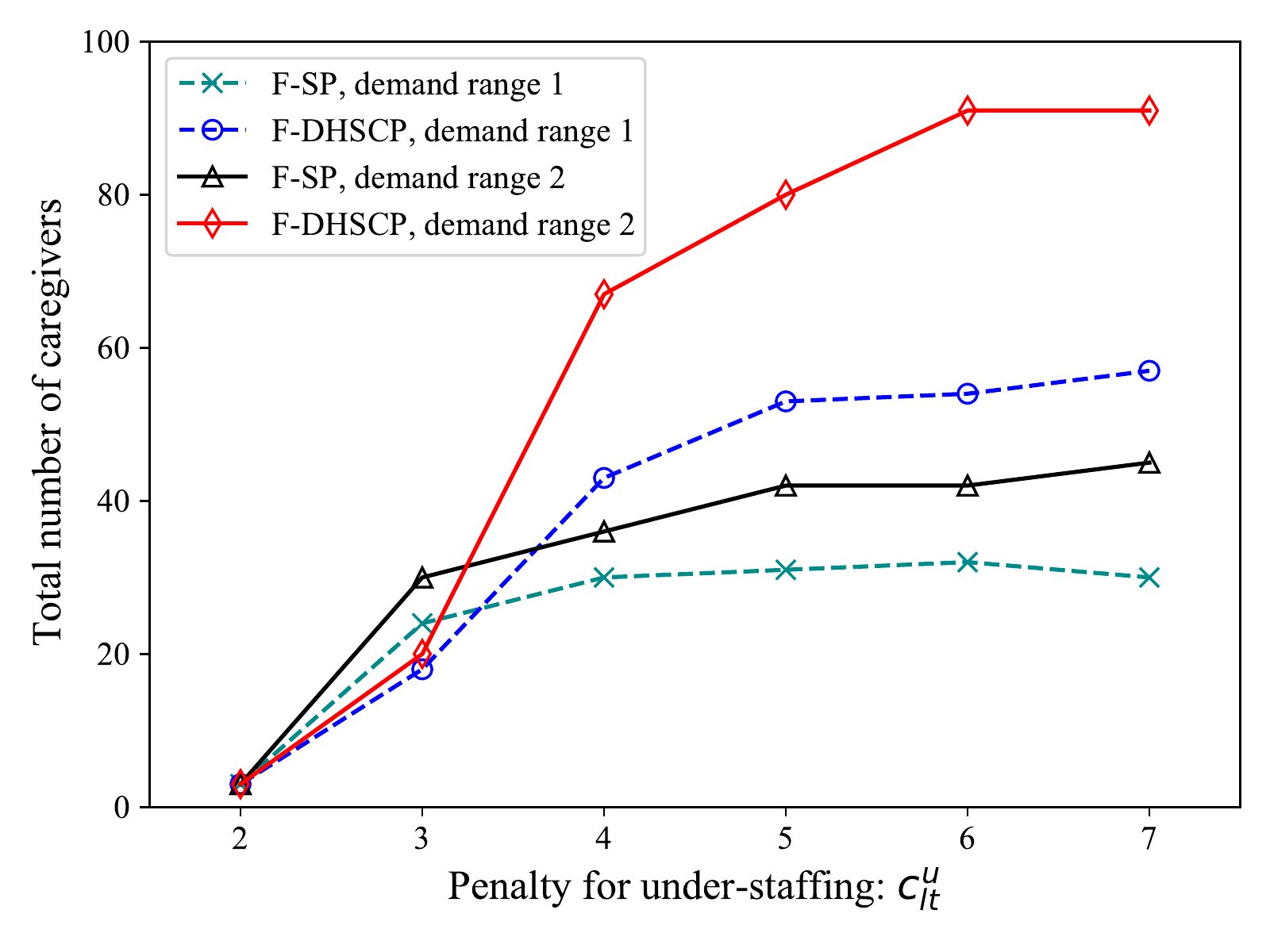}}
    \subcaptionbox{$c_{l,t}^u = 5$\label{FA-10-6-6-30-o}}{
        \includegraphics[scale=0.35]{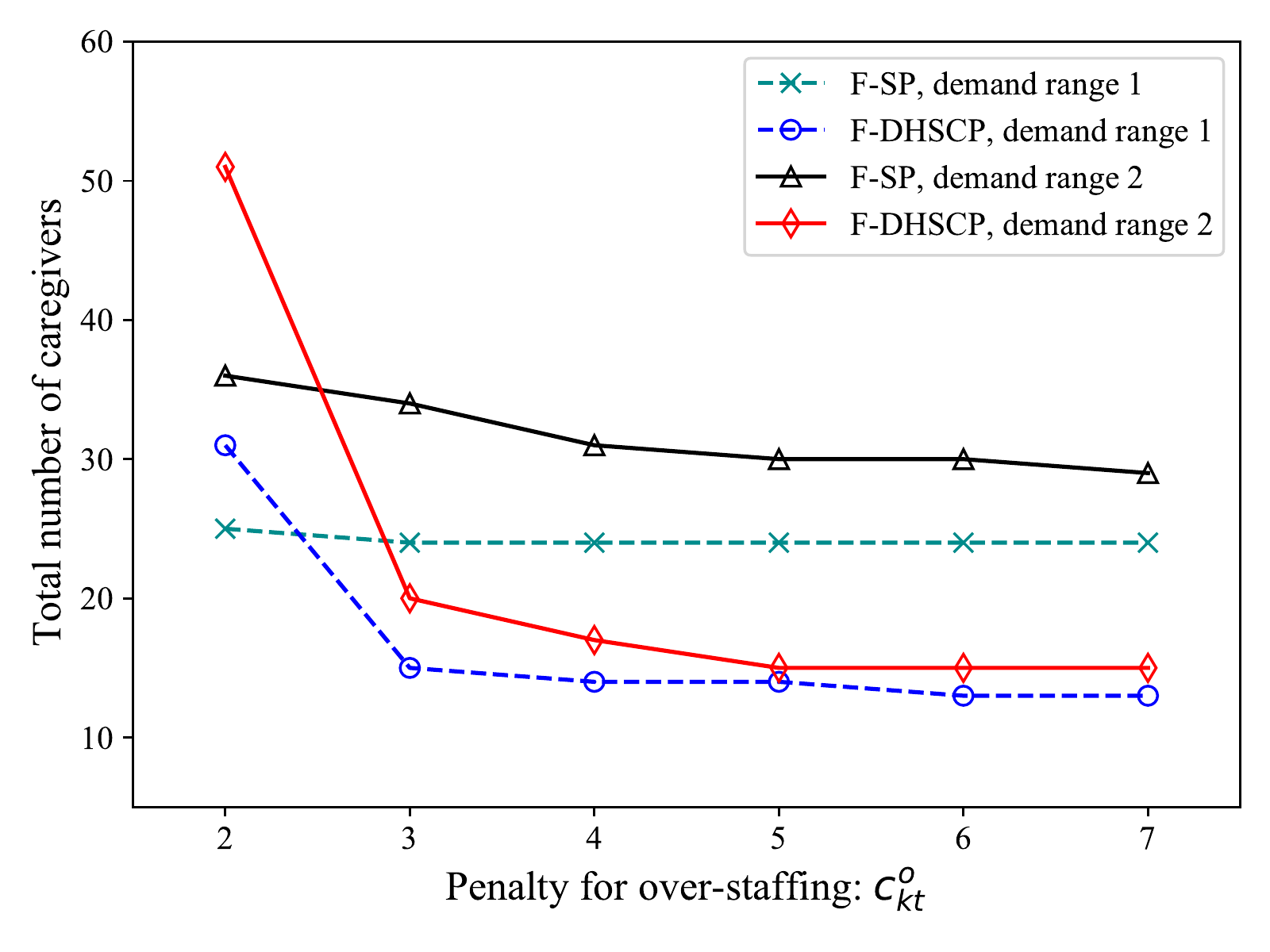}}
    \caption{Total number of caregivers hired by F-SP and F-DHSCP models for Instance 10.}
    \label{total_num_fig_f}
\end{figure}

In this section, we compare the optimal staffing patterns of the F-SP and F-DHSCP under different values of the under-staffing $c_{l,t}^u$ and over-staffing $c_{k,t}^o$
penalty costs. In Figure \ref{FA-10-6-6-30}, we present the optimal number of hired caregivers under demand
range 1 and range 2 with $c_{k,t}^o=1$ and $c_{l,t}^u \in \{2,3,4,5,6,7\}$ for Instance 10. We
observe the same staffing patterns for Instance 7 (see Figure \ref{total_num_fig_f_a} in \ref{OSP_F}).

We observe the following from Figure~\ref{total_num_fig_f}. First, both models hire more caregivers under demand range 2 than under range 1 and the number of hired caregivers stabilize after $c_{l,t}^u=5$. Second, when $c_{l,t}^u\leq 3$, both models hire approximately the same number of caregivers. In contrast, when $c_{l,t}^u > 3$ the F-DHSCP models hire more caregivers than the E-SP models. These results indicate that $c_{l,t}^u = 5$ may be a suitable penalty cost of under-staffing for the (FA) decision-makers. We observe a different pattern when we fix $c_{l,t}^u = 5$ and increase the penalty for over-staffing $c_{k,t}^o$ from 2 to $\{3,4,5,6,7\}$. Specifically, as shown in Figure \ref{FA-10-6-6-30-o} both models hire less caregivers as  $c_{k,t}^o$ increases. And the number of hired caregivers stabilize after $c_{k,t}^o=3$.

Next, we analyze the detailed staffing patterns under cost structure $(c_{l,t}^{u}, c_{k,t}^{o}) = (5, 2)$ presented in Figure \ref{FA-solution} (see \ref{OSP_F}). Compared with EA models, both F-SP and F-DHSCP models prefer cross-trained caregivers under two demand ranges for Instance 7.

\subsubsection{Out-of-sample performance of FA model}\label{OOP-FA}
In this section, we compare the operational performance of the optimal solutions of the F-DHSCP and F-SP via out-of-sample simulation. Using the same methods in Section \ref{OOP-EA-M}, we generate two sets of $N^{'} = 10,000$ out-of-sample data (Set 1-2). For brevity, we focus on Instance 10 and present results under $(c_{l,t}^{u}, c_{k,t}^{o}) = (5, 2)$ (see \ref{oop-FA-7} for results related to Instance 7). 

\begin{figure}[t!]
    \centering
    \subcaptionbox{\label{FA-1outc-10-2-6-6-30-4060}}{
        \includegraphics[scale=0.3]{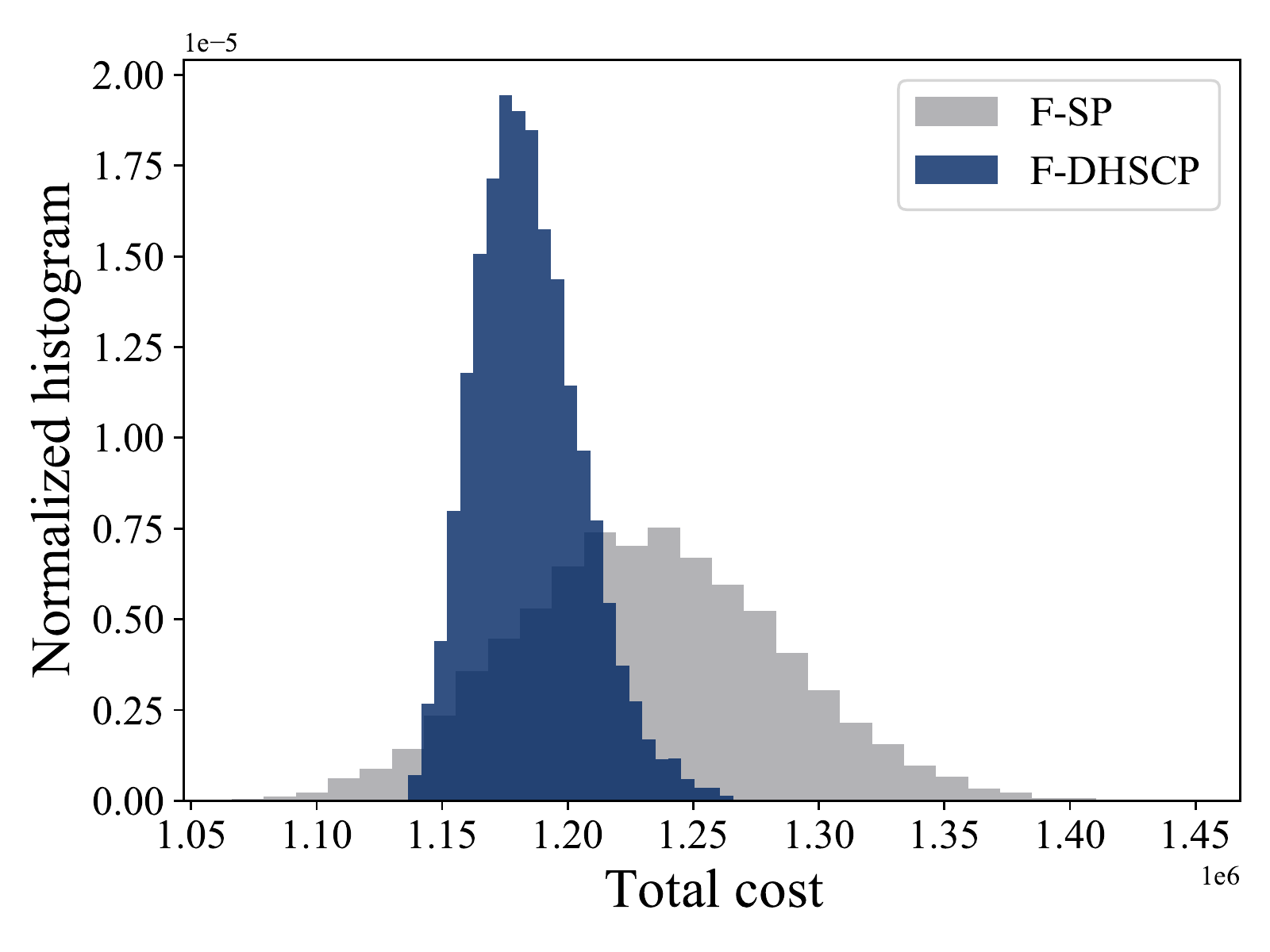}}
    \subcaptionbox{\label{FA-1outc2-10-2-6-6-30-4060}}{
        \includegraphics[scale=0.3]{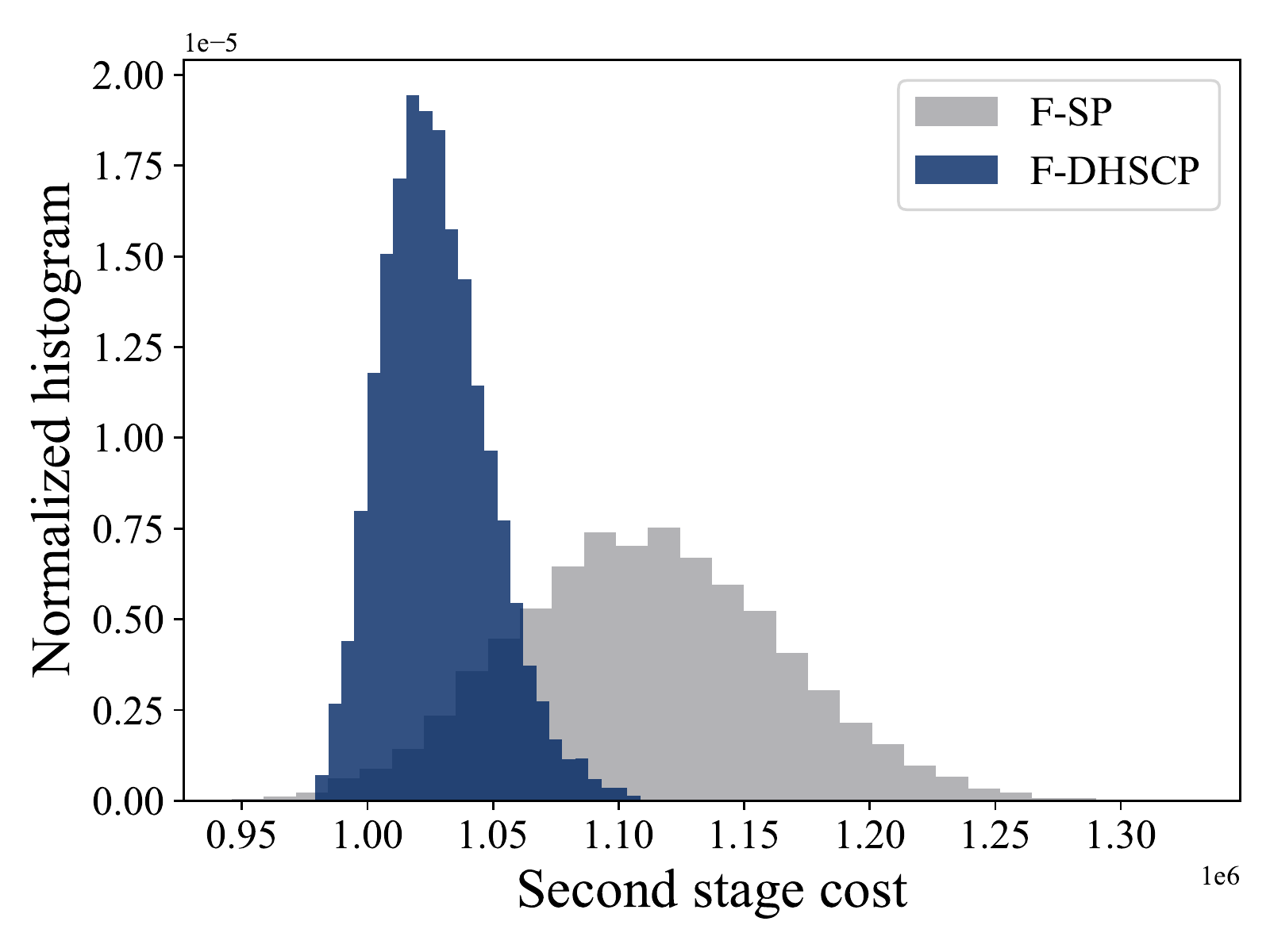}}
    \subcaptionbox{\label{FA-1outu-10-2-6-6-30-4060}}{
        \includegraphics[scale=0.3]{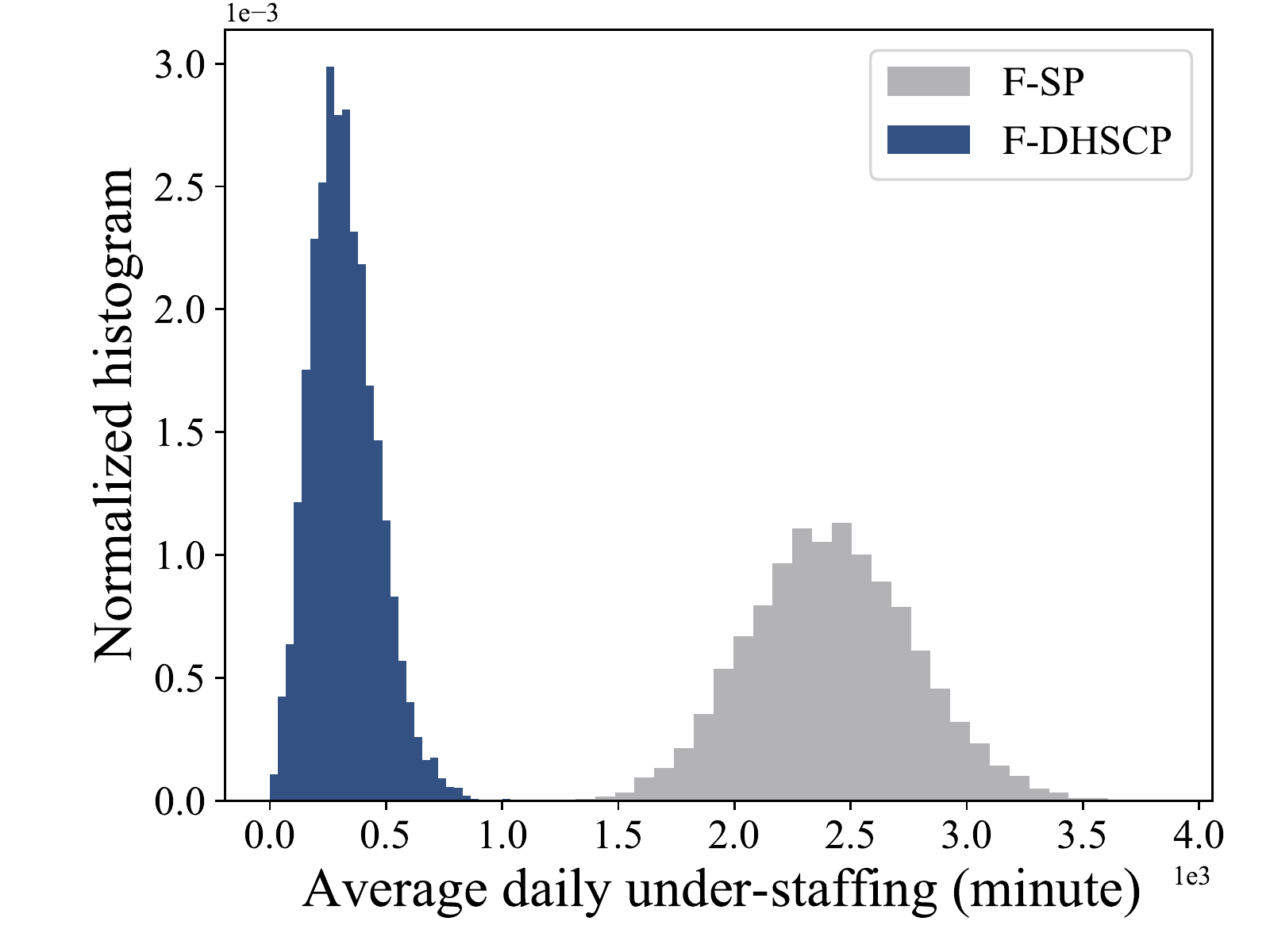}}
    \caption{Out-of-sample performance of FA models for Instance 10, demand range 1 under Set 2, $\Delta = 0$}
    \label{FA-oop-o-6}
\end{figure}

We observe the following from the simulation results under Set 2, $\Delta=0$ (i.e., misspecified distributional information case) with demand range 1 presented in Figure \ref{FA-oop-o-6} (see Figure \ref{FA-oop-o-6-2} in \ref{oop-FA-10} for results with range 2). The F-DHSCP model outperforms the F-SP model under all levels of variation ($\Delta$) and across the criteria of mean and all quantiles of the second-stage costs and under-staffing under the two ranges of the demand (see Figure \ref{FA-oop-6}--\ref{FA-oop-u-6-4060} in \ref{oop-FA-10}). Specifically, as shown in Figure \ref{FA-1outc-10-2-6-6-30-4060}, the F-DHSCP model yields a lower total cost than the F-SP model on average and at all quantiles. This makes sense because F-DHSCP make more flexible capacity allocation, which reduce the number of hired caregivers then the fixed cost associated with hiring more caregivers.

\begin{figure}[t!]
    \centering
    \subcaptionbox{$\Delta=0$\label{FA-outdis-10-2-6-6-30-4060}}{
        \includegraphics[scale=0.3]{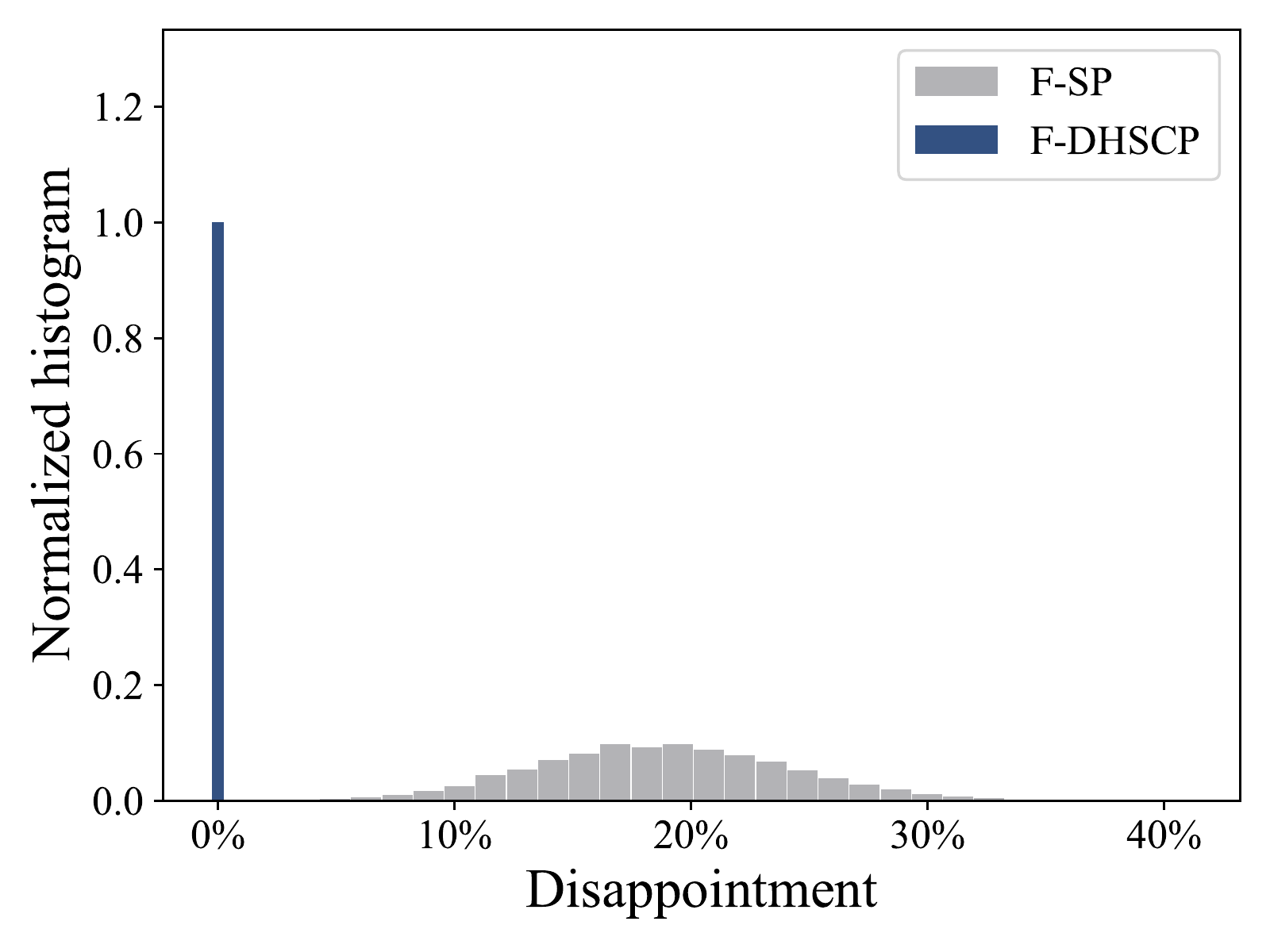}}
    \subcaptionbox{$\Delta=0.25$\label{FA-out025dis-10-2-6-6-30-4060}}{
        \includegraphics[scale=0.3]{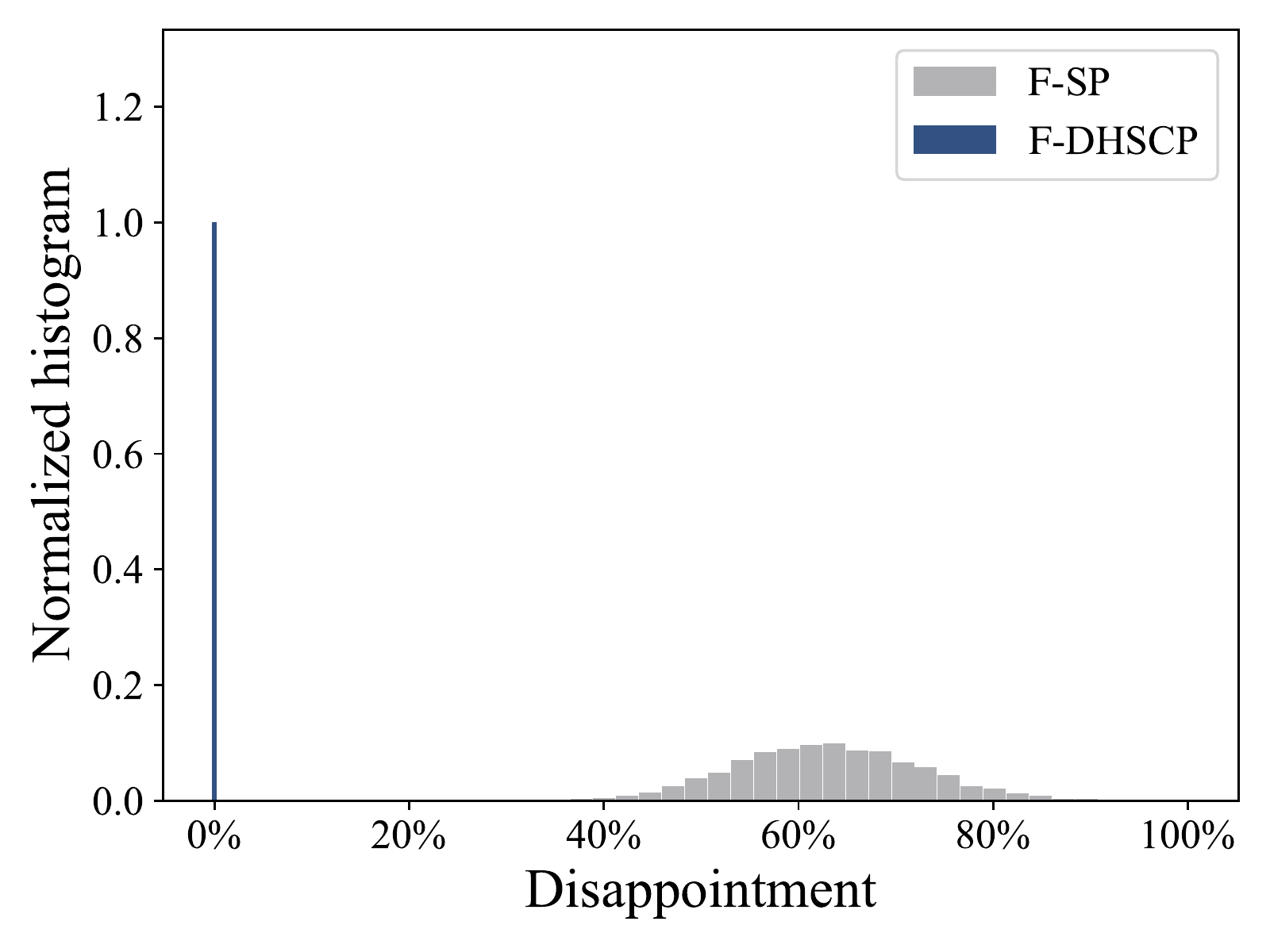}}
    \subcaptionbox{$\Delta=0.5$\label{FA-out05dis-10-2-6-6-30-4060}}{
        \includegraphics[scale=0.3]{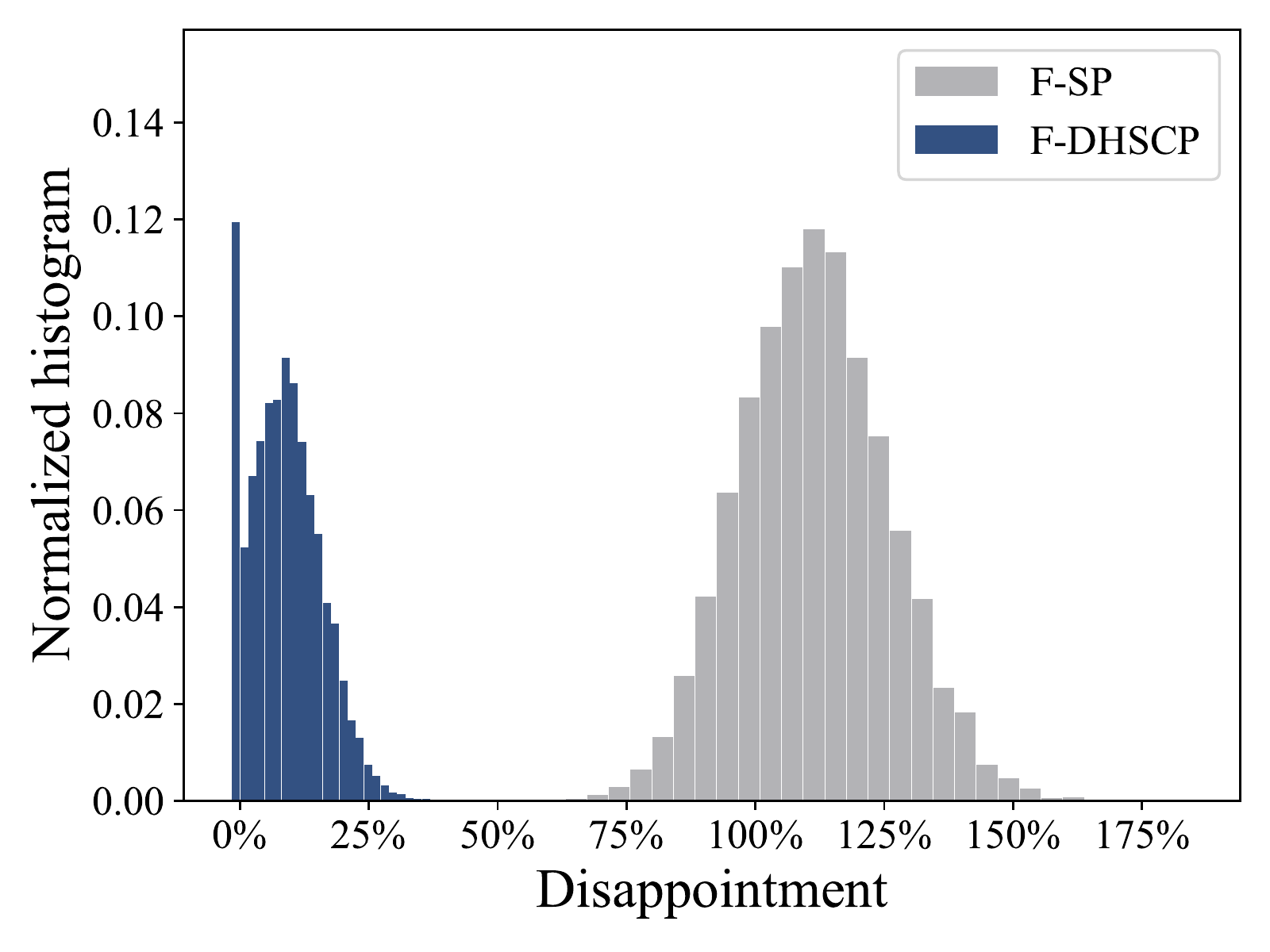}}
    \caption{Out-of-sample disappointment for Instance 10, demand range 1 under Set 2}
    \label{FA-oop-dis-6}
\end{figure}

Finally, it is clear from Figure~\ref{FA-oop-dis-6} that the F-DHSCP solutions yield smaller out-of-sample disappointments on average and at all quantiles than the F-SP solutions under all levels of variation ($\Delta$). Notably, the disappointment of the F-DHSCP solutions is smaller than 50\% under a higher variation $\Delta=0.5$. These results confirm that the F-DHSCP model provides a more robust estimate of the actual cost that we will incur in practice.

\section{Conclusion}\label{conclusion}

In this paper, we proposed stochastic optimization methodologies for the   HSCP problem, which seeks to find the number of caregivers to hire and their allocation decisions to different types of services in each day within a specified planning horizon. The objective is to minimize the total cost associated with staffing (i.e., employment), capacity allocation, over-staffing, and under-staffing. To address uncertainty in service duration and the demand for home service, we proposed and analyzed SP and DRO models considering two types of decision-makers, namely everything in advance decision-maker (for which we proposed E-SP and E-DHSCP models) and a flexible adjustment decision-maker (FA) (for which we proposed (F-SP and F-DHSCP). We derive equivalent MILP reformulations of the proposed nonlinear E-DHSCP model for the EA decision-maker that can be implemented and efficiently solved using off-the-shelf optimization software. We propose a computationally efficient column-and-constraint generation algorithm with valid inequalities to solve the proposed F-DHSCP model for the FA decision-maker.

Finally, we conduct extensive numerical experiments comparing the operational and computational performance of the proposed approaches and discuss insights and implications for home care staffing and capacity planning. The results demonstrate (1) the DRO models tend to hedge against uncertainty and ambiguity by hiring more caregivers and thus have a higher fixed cost but significantly lower operational (i.e., under-staffing) costs than the SP models, (2) the DRO solutions are more robust under various distributions with substantially lower disappointments and thus higher reliability than the SP models, (3) the computational efficiency of the proposed approaches, indicating their applicability in practice. We also derive insights into the HSCP under different settings. In particular, we show how one can use the proposed models to obtain staffing and capacity allocation decisions under different cost structures that the decision-maker can easily change in the models. These results address the primary objective of this paper. That is, to investigate the value of SP and DRO approaches for this specific HSCP problem. More broadly, our results motivate the need for considering different approaches when modeling uncertainty and draw attention to the need to model the distributional ambiguity of uncertain problem data in strategic real-world stochastic optimization problems such as the HSCP problem.

We suggest the following areas for future research. First, we want to extend our approach by incorporating other multi-modal sources of uncertainty (e.g., travel time). Second, our models can be considered as the first step toward building data-driven and robust home care service planning including the creating routing plans for caregivers and assigning appointment times for customers, among other important considerations. Finally, we aim to incorporate other aspects, such as the possibility of hiring part-time caregivers and the unexpected availability of the hired full-time caregivers on the day of service.

\noindent \textbf{Acknowledgment}

\noindent  We want to thank all of our colleagues who have contributed significantly to the related literature. Dr.~Karmel S.~Shehadeh dedicates her effort in this paper to every little dreamer in the whole world who has a dream so big and so exciting. Believe in your dreams and do whatever it takes to achieve them--the best is yet to come for you.

\vspace{4mm}

\noindent \textbf{References}

\addcontentsline{toc}{section}{References}\textbf{}
{
\onehalfspacing
\small
\setlength{\bibsep}{1ex}
\bibliography{HSCP}
}

\newpage
\appendix
\begin{center}
    \Large Integrated Home Care Staffing and Capacity Planning: Stochastic Optimization Approaches (Appendices)\\
    
    Ridong Wang, Karmel S. Shehadeh, Xiaolei Xie, Lefei Li
\end{center}
\section{Candidate Probability Distributions for demands}\label{CPD}
\setcounter{figure}{0}
\setcounter{table}{0}  
In this section, we first show the fitting distributions for demands and their goodness of
fit metrics; Negative of the Log Likelihood (NLogL), Akaike Information Criterion (AIC), and Bayesian Information Criterion (BIC). Then, in Figure \ref{dr_ts}, we show the rehabilitation service demand range change in different quarters.

\begin{table}[htbp]
  \centering
  \caption{Candidate probability distributions for demands of rehabilitation nursing service}
    \begin{tabular}{lccc}
    \hline
    Distribution name & NLogL & BIC & AIC \\
    \hline
    generalized extreme value & 4951.50  & 9923.02  & 9909.01  \\
    loglogistic & 4967.45  & 9948.24  & 9938.90  \\
    lognormal & 4968.68  & 9950.71  & 9941.37  \\
    birnbaumsaunders & 4989.83  & 9992.99  & 9983.65  \\
    inverse gaussian & 4991.17  & 9995.69  & 9986.35  \\
    gamma & 4994.32  & 10001.97  & 9992.63  \\
    weibull & 5029.16  & 10071.66  & 10062.32  \\
    nakagami & 5043.96  & 10101.26  & 10091.92  \\
    rayleigh & 5058.13  & 10122.93  & 10118.26  \\
    rician & 5058.13  & 10129.60  & 10120.26  \\
    generalized pareto & 5068.65  & 10157.31  & 10143.30  \\
    tlocationscale & 5104.24  & 10228.49  & 10214.49  \\
    logistic & 5153.31  & 10319.95  & 10310.61  \\
    exponential & 5169.56  & 10345.80  & 10341.13  \\
    normal & 5171.11  & 10355.56  & 10346.22  \\
    extreme value & 5336.91  & 10687.17  & 10677.83  \\
    \hline
    \end{tabular}%
  \label{CPD_t_1}%
\end{table}%

\begin{table}[htbp]
  \centering
  \caption{Candidate probability distributions for demands of health status assessment service}
    \begin{tabular}{lccc}
    \hline
    Distribution name & NLogL & BIC & AIC \\
    \hline
    generalized pareto & 297.65  & 609.69  & 601.31  \\
    generalized extreme value & 321.55  & 657.49  & 649.10  \\
    inverse gaussian & 333.13  & 675.85  & 670.25  \\
    birnbaumsaunders & 333.24  & 676.08  & 670.48  \\
    lognormal & 333.33  & 676.25  & 670.66  \\
    loglogistic & 335.70  & 680.99  & 675.40  \\
    gamma & 337.40  & 684.40  & 678.80  \\
    nakagami & 342.40  & 694.40  & 688.81  \\
    logistic & 346.23  & 702.06  & 696.47  \\
    rician & 348.16  & 705.91  & 700.32  \\
    tlocationscale & 345.89  & 706.18  & 697.79  \\
    normal & 348.79  & 707.17  & 701.58  \\
    weibull & 352.38  & 714.35  & 708.76  \\
    rayleigh & 380.76  & 766.32  & 763.52  \\
    extreme value & 379.28  & 768.14  & 762.55  \\
    exponential & 451.12  & 907.03  & 904.23  \\
    \hline
    \end{tabular}%
  \label{CPD_t_2}%
\end{table}%

\begin{figure}[!ht]
    \centering
    \includegraphics[scale=0.47]{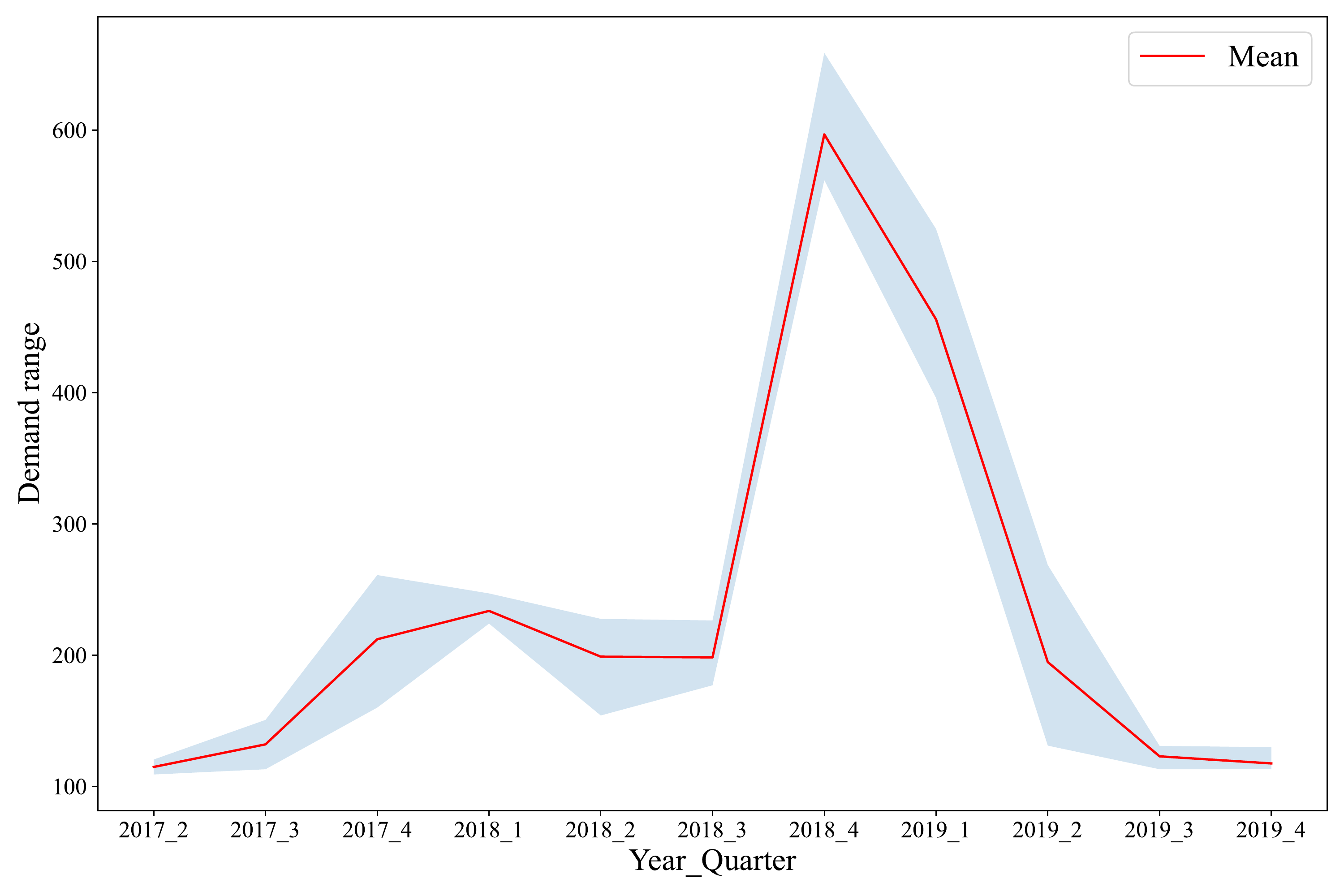}
    \caption{Rehabilitation service demand range in time series}
    \label{dr_ts}
\end{figure}

\clearpage
\newpage

\section{Proof of Proposition \ref{pro_deter}}\label{proof_pos_d}
\setcounter{figure}{0}
\setcounter{table}{0}  
\begin{proof}
\noindent For fixed $\bm{x}$ and $\bm{y}$, we can formulate problem $\mbox{sup}_{\mathbb{P}\in \mathcal{F}(\mathcal{U}, \bm{\mu})} \mathbb{E}_{\mathbb{P}} [Q(\textbf{\emph{x}}, \bm{y}, \bm{\xi})]$ as following linear optimization problem.
{
\setlength{\abovedisplayskip}{4pt}
\setlength{\belowdisplayskip}{4pt}
\begin{subequations}
    \label{sup2}
    \begin{alignat}{4}
        \underset{\mathbb{P}\geq 0}{\mbox{max}} &\quad \int_{\mathcal{U}}[ Q(\textbf{\emph{x}}, \bm{y},  \bm{\xi})]d\mathbb{P} && &&\label{sup_21}\\
        \mbox{s.t.} &\quad \int_{\mathcal{U}} d_{l,t}d\mathbb{P} = \mu_{l,t}^d, &&\forall l\in [L], t\in [T], &&\quad \alpha_{l,t}\label{sup_22}\\
        &\quad \int_{\mathcal{U}} s_{l,t}d\mathbb{P} = \mu_{l,t}^s, &&\forall l\in [L], t\in [T], &&\quad \beta_{l,t} \label{sup_23}\\
        &\quad \int_{\mathcal{U}} d\mathbb{P} = 1 && &&\quad \theta \label{sup_24}
    \end{alignat}
\end{subequations}}
\indent Note that $\mathcal{U}$ is compact, and the recourse is a continuous function in $d$ and $s$. Under standard assumptions that $\mu_{l,t}^d$ belongs to the interior of set \{ $\int_{\mathcal{U}}d_{l,t}d\mathbb{Q}: \mathbb{Q}$ is a probability distribution over $\mathcal{U}$ \} for all $l\in [L]$ and $t\in [T]$ and that $s_{l,t}$ belongs to the interior of set \{ $\int_{\mathcal{U}}s_{l,t}d\mathbb{Q}: \mathbb{Q}$ is a probability distribution over $\mathcal{U}$ \} for all $l\in [L]$ and $t\in [T]$ , slater condition will hold \citep{boyd2004convex}. Then, strong duality holds and the worst-case expectation equals:
{
\setlength{\abovedisplayskip}{4pt}
\setlength{\belowdisplayskip}{4pt}
\begin{subequations}
    \label{deter_proof}
    \begin{alignat}{4}
        \underset{\bm{\alpha}, \bm{\beta}, \theta}{\mbox{min}} &\quad \left\{
        \sum_{l=1}^{L}\sum_{t=1}^{T} \mu_{l,t}^d\alpha_{l,t} + \sum_{l=1}^{L}\sum_{t=1}^{T}\mu_{l,t}^s\beta_{l,t} + \theta \right\} && \label{deter_proof_1}\\
        \mbox{s.t.} &\quad \sum_{l=1}^{L}\sum_{t=1}^{T}(d_{l,t}\alpha_{l,t} + s_{l,t}\beta_{l,t}) + \theta \geq Q(\textbf{\emph{x, y}}, \bm{\xi}), &&\quad \forall (\bm{d,s})\in \mathcal{U}\label{deter_proof_2}
    \end{alignat}
\end{subequations}}
where $\alpha_{l,t}, \beta_{l,t}$ and $\theta$ are the dual variables associated with constraints \eqref{sup_22}, \eqref{sup_23}, and \eqref{sup_24}. Note that we can rewrite constraint \eqref{deter_proof_2} as follows:
{
\setlength{\abovedisplayskip}{4pt}
\setlength{\belowdisplayskip}{4pt}
\begin{equation*}
    \theta \geq \underset{(\bm{d,s})\in \mathcal{U}}{\mbox{max}} \left\{ Q(\textbf{\emph{x, y}}, \bm{\xi}) -\sum_{l=1}^{L}\sum_{t=1}^{T}(d_{l,t}\alpha_{l,t} + s_{l,t}\beta_{l,t})\right\}
\end{equation*}}
Given that we are minimizing $\theta$ in \eqref{deter_proof_1}, we can  rewrite formulation \eqref{deter_proof} as follows:
{
\setlength{\abovedisplayskip}{4pt}
\setlength{\belowdisplayskip}{4pt}
\begin{equation*}
    \underset{\bm{\alpha}, \bm{\beta}}{\mbox{min}}\quad \left\{
        \sum_{l=1}^{L}\sum_{t=1}^{T} \mu_{l,t}^d\alpha_{l,t} + \sum_{l=1}^{L}\sum_{t=1}^{T}\mu_{l,t}^s\beta_{l,t} + \underset{(\bm{d,s})\in \mathcal{U}}{\mbox{max}} \left\{Q(\textbf{\emph{x, y}}, \bm{\xi}) -\sum_{l=1}^{L}\sum_{t=1}^{T}(d_{l,t}\alpha_{l,t} + s_{l,t}\beta_{l,t})\right\} \right\}
\end{equation*}}
\end{proof}

\clearpage 
\newpage

\section{Proof of Proposition \ref{min_inner}}\label{Appx:min_inner}
\begin{proof}
\setcounter{figure}{0}
\setcounter{table}{0}  

Not that problem $Q(\textbf{\emph{x}}, \textbf{\emph{y}}, \bm{\xi})$  in \eqref{eq:optim8} is separable by each $t$ and $l$. Therefore, w.l.o., we can rewrite $Q(\textbf{\emph{x}}, \textbf{\emph{y}}, \bm{\xi})= \sum_{l} \sum_{t} Q_{l,t}(\textbf{\emph{x}}, \textbf{\emph{y}}, d_{l,t}, s_{l,t})$, where for fixed $\pmb{x}$, $\pmb{y}$, $l$, $t$, $d_{l,t} \in [\underline{d}_{l,t},  \overline{d}_{l,t}]$, and $s_{l,t} \in[\underline{s}_{l,t},  \overline{s}_{l,t}]$.
{
\setlength{\abovedisplayskip}{4pt}
\setlength{\belowdisplayskip}{4pt}
\begin{subequations}\label{Q_decom}
    \begin{align}
        Q_{l,t}(\textbf{\emph{x}}, \textbf{\emph{y}}, d_{l,t}, s_{l,t}) &= \max_{\rho_{l,t}} \Big \{(d_{l,t}s_{l} - \sum_{k: l\in R_{k}} y_{k,l,t})\rho_{l,t} \Big\} \label{Qdecomp:obj}\\
        &    \mbox{s.t.}  \ \ P_{l,t}:=\{- c_{l,t}^{o} \leq \rho_{l,t} \leq c_{l,t}^{u}\}.
    \end{align}
\end{subequations}}
\indent Accordingly, 
{
\setlength{\abovedisplayskip}{4pt}
\setlength{\belowdisplayskip}{4pt}
\begin{align}
    \max \limits_{(\bm{d,s})\in \mathcal{U}} \left\{Q(\bm{x}, \bm{y}, \bm{\xi}) - \sum \limits_{l=1}^L \sum \limits_{t =1}^T(d_{l,t}\alpha_{l,t} + s_{l,t}\beta_{l,t})\right\}& 
    \equiv \sum \limits_{l=1}^{L} \sum \limits_{t=1}^{T}\max\limits_{\substack{d_{l,t} \in [\underline{d}_{l,t},  \overline{d}_{l,t}]\\ s_{l,t} \in[\underline{s}_{l,t},  \overline{s}_{l,t}] }} \Big\{ Q_{l,t}(\textbf{\emph{x}}, \textbf{\emph{y}}, d_{l,t}, s_{l,t}) \nonumber\\
    &\qquad \qquad \qquad \qquad \qquad- d_{l,t}\alpha_{l,t} - s_{l,t}\beta_{l,t} \Big\}, \nonumber\\
    &\equiv \sum \limits_{l=1}^{L} \sum \limits_{t=1}^{T}\max\limits_{\substack{d_{l,t} \in [\underline{d}_{l,t}, \overline{d}_{l,t}]\\ s_{l,t} \in[\underline{s}_{l,t}, \overline{s}_{l,t}]}} h_{l,t}(d_{l,t}, s_{l,t}, \bm{y}_{l,t}, \alpha_{l,t}, \beta_{l,t})\label{maxH},
\end{align}}
where
{
\setlength{\abovedisplayskip}{4pt}
\setlength{\belowdisplayskip}{4pt}
\begin{equation}
    h_{l,t}(d_{l,t}, s_{l,t}, \bm{y}_{l,t}, \alpha_{l,t}, \beta_{l,t}) := \underset{\rho_{l,t}\in P_{l,t}}{\mbox{max}} \left\{d_{l,t}s_{l,t}\rho_{l,t} - \sum_{k: l\in R_{k}} y_{k,l,t}\rho_{l,t} - d_{l,t}\alpha_{l,t} - s_{l,t}\beta_{l,t}\right\}\label{func}.
\end{equation}}
\indent It is straightforward to verify that for fixed $\pmb{x}$, $\pmb{y}$,  $d_{l,t} \in [\underline{d}_{l,t},  \overline{d}_{l,t}]$, and $s_{l,t} \in[\underline{s}_{l,t},  \overline{s}_{l,t}]$ function \eqref{func} is convex in variable $\rho_{l,t}$. Hence, problem \eqref{func} is a convex maximization problem. It follows from the fundamental convex analysis (see \cite{boyd2004convex}) that here exists an optimal solution $\hat \rho_{l,t}$ at one of the extreme points of polyhedron $P_{l,t}:=\{- c_{l,t}^{o} \leq \rho_{l,t} \leq c_{l,t}^{u}\}$. In any extreme point, constraint $- c_{l,t}^{o} \leq \rho_{l,t} \leq c_{l,t}^{u}$ is binding at either the lower bound or upper bound. Thus, given that $d_{l,t} \in [\underline{d}_{l,t},  \overline{d}_{l,t}]$, and $s_{l,t} \in[\underline{s}_{l,t},  \overline{s}_{l,t}]$, it is easy to verify that the optimal value $\eta_{l,t}^{*}$ of problem $\max\limits_{\substack{d_{l,t} \in [\underline{d}_{l,t},  \overline{d}_{l,t}]\\ s_{l,t} \in[\underline{s}_{l,t},  \overline{s}_{l,t}] }} h_{l,t}(d_{l,t}, s_{l,t}, \bm{y}_{l,t}, \alpha_{l,t}, \beta_{l,t})$ in \eqref{maxH} equals to:
{
\setlength{\abovedisplayskip}{4pt}
\setlength{\belowdisplayskip}{4pt}
\begin{equation}
    \eta_{l,t}^{*} = \mbox{max} \left\{
    \begin{array}{rcl}
        &\overline{d}_{l,t}\overline{s}_{l,t}c_{l,t}^{u} - \sum_{k: l\in R_{k}} y_{k,l,t}c_{l,t}^{u} - \overline{d}_{l,t}\alpha_{l,t} -\overline{s}_{l,t}\beta_{l,t}\\
        &-\underline{d}_{l,t}\underline{s}_{l,t}c_{l,t}^{o} + \sum_{k: l\in R_{k}} y_{k,l,t}c_{l,t}^{o} - \underline{d}_{l,t}\alpha_{l,t} - \underline{s}_{l,t}\beta_{l,t}\\
        &-\overline{d}_{l,t}\overline{s}_{l,t}c_{l,t}^{o} + \sum_{k: l\in R_{k}} y_{k,l,t}c_{l,t}^{o} - \overline{d}_{l,t}\alpha_{l,t} - \overline{s}_{l,t}\beta_{l,t}\\
        &\underline{d}_{l,t}\overline{s}_{l,t}c_{l,t}^{u} - \sum_{k: l\in R_{k}} y_{k,l,t}c_{l,t}^{u} - \underline{d}_{l,t}\alpha_{l,t} - \overline{s}_{l,t}\beta_{l,t}\\
        &\overline{d}_{l,t}\underline{s}_{l,t}c_{l,t}^{u} - \sum_{k: l\in R_{k}} y_{k,l,t}c_{l,t}^{u} - \overline{d}_{l,t}\alpha_{l,t} - \underline{s}_{l,t}\beta_{l,t}\\
        &\underline{d}_{l,t}\underline{s}_{l,t}c_{l,t}^{u} - \sum_{k: l\in R_{k}} y_{k,l,t}c_{l,t}^{u} - \underline{d}_{l,t}\alpha_{l,t} - \underline{s}_{l,t}\beta_{l,t}\\
        &-\overline{d}_{l,t}\underline{s}_{l,t}c_{l,t}^{o} + \sum_{k: l\in R_{k}} y_{k,l,t}c_{l,t}^{o} - \overline{d}_{l,t}\alpha_{l,t} - \underline{s}_{l,t}\beta_{l,t}\\
        &-\underline{d}_{l,t}\overline{s}_{l,t}c_{l,t}^{o} + \sum_{k: l\in R_{k}} y_{k,l,t}c_{l,t}^{o} - \underline{d}_{l,t}\alpha_{l,t} - \overline{s}_{l,t}\beta_{l,t}
    \end{array} \right\}
\end{equation}}
Hence,  $\max\limits_{\substack{d_{l,t} \in [\underline{d}_{l,t},  \overline{d}_{l,t}]\\ s_{l,t} \in[\underline{s}_{l,t},  \overline{s}_{l,t}] }} h_{l,t}(d_{l,t}, s_{l,t}, y, \alpha_{l,t}, \beta_{l,t})$ is equivalent to the following minimization problem.
\begin{spacing}{0.9}
{
\setlength{\abovedisplayskip}{4pt}
\setlength{\belowdisplayskip}{4pt}
\begin{subequations}
\label{A:inner_max}
    \begin{alignat}{3}
        \mbox{min}
        &\quad \eta_{l,t} &&\quad \label{A:inner_max_1}\\
        \mbox{s.t.}
        &\quad \eta_{l,t} \geq \overline{d}_{l,t}\overline{s}_{l,t}c_{l,t}^{u} - \sum_{k: l\in R_{k}} y_{k,l,t}c_{l,t}^{u} - \overline{d}_{l,t}\alpha_{l,t} -\overline{s}_{l,t}\beta_{l,t},&&\label{A:inner_max_first}\\
        &\quad \eta_{l,t} \geq -\underline{d}_{l,t}\underline{s}_{l,t}c_{l,t}^{o} + \sum_{k: l\in R_{k}} y_{k,l,t}c_{l,t}^{o} - \underline{d}_{l,t}\alpha_{l,t} - \underline{s}_{l,t}\beta_{l,t},&&\\
        &\quad \eta_{l,t} \geq -\overline{d}_{l,t}\overline{s}_{l,t}c_{l,t}^{o} + \sum_{k: l\in R_{k}} y_{k,l,t}c_{l,t}^{o} - \overline{d}_{l,t}\alpha_{l,t} - \overline{s}_{l,t}\beta_{l,t},&&\\
        &\quad \eta_{l,t} \geq \underline{d}_{l,t}\overline{s}_{l,t}c_{l,t}^{u} - \sum_{k: l\in R_{k}} y_{k,l,t}c_{l,t}^{u} - \underline{d}_{l,t}\alpha_{l,t} - \overline{s}_{l,t}\beta_{l,t},&&\\
        &\quad \eta_{l,t} \geq \overline{d}_{l,t}\underline{s}_{l,t}c_{l,t}^{u} - \sum_{k: l\in R_{k}} y_{k,l,t}c_{l,t}^{u} - \overline{d}_{l,t}\alpha_{l,t} - \underline{s}_{l,t}\beta_{l,t},&&\\
        &\quad \eta_{l,t} \geq \underline{d}_{l,t}\underline{s}_{l,t}c_{l,t}^{u} - \sum_{k: l\in R_{k}} y_{k,l,t}c_{l,t}^{u} - \underline{d}_{l,t}\alpha_{l,t} - \underline{s}_{l,t}\beta_{l,t},&&\\
        &\quad \eta_{l,t} \geq -\overline{d}_{l,t}\underline{s}_{l,t}c_{l,t}^{o} + \sum_{k: l\in R_{k}} y_{k,l,t}c_{l,t}^{o} - \overline{d}_{l,t}\alpha_{l,t} - \underline{s}_{l,t}\beta_{l,t},&&\\
        &\quad \eta_{l,t} \geq -\underline{d}_{l,t}\overline{s}_{l,t}c_{l,t}^{o} + \sum_{k: l\in R_{k}} y_{k,l,t}c_{l,t}^{o} - \underline{d}_{l,t}\alpha_{l,t} - \overline{s}_{l,t}\beta_{l,t}.&& \label{A:inner_max_last}
    \end{alignat}
\end{subequations}
}
\end{spacing}

Replacing $\max h(\cdot)$ in \ref{maxH} with its equivalent reformulation in \eqref{A:inner_max} and summing over all $l \in [L]$ and $t \in[T]$, we derive the following equivalent reformulation.
\begin{spacing}{0.9}
{
\setlength{\abovedisplayskip}{4pt}
\setlength{\belowdisplayskip}{4pt}
\begin{subequations}
\label{inner_max_appendix}
    \begin{alignat}{3}
       \min_{\pmb{\eta}}
        &\quad \sum_{l=1}^L \sum_{t=1}^T \eta_{l,t} &&\quad \label{inner_max_appendix_1}\\
        \mbox{s.t.}
        &\quad \eta_{l,t} \geq \overline{d}_{l,t}\overline{s}_{l,t}c_{l,t}^{u} - \sum_{k: l\in R_{k}} y_{k,l,t}c_{l,t}^{u} - \overline{d}_{l,t}\alpha_{l,t} -\overline{s}_{l,t}\beta_{l,t},&&\quad \forall l\in[L], t\in[T],\label{inner_max_appendix_first}\\
        &\quad \eta_{l,t} \geq -\underline{d}_{l,t}\underline{s}_{l,t}c_{l,t}^{o} + \sum_{k: l\in R_{k}} y_{k,l,t}c_{l,t}^{o} - \underline{d}_{l,t}\alpha_{l,t} - \underline{s}_{l,t}\beta_{l,t},&&\quad \forall l\in[L], t\in[T],\\
        &\quad \eta_{l,t} \geq -\overline{d}_{l,t}\overline{s}_{l,t}c_{l,t}^{o} + \sum_{k: l\in R_{k}} y_{k,l,t}c_{l,t}^{o} - \overline{d}_{l,t}\alpha_{l,t} - \overline{s}_{l,t}\beta_{l,t},&&\quad \forall l\in[L], t\in[T],\\
        &\quad \eta_{l,t} \geq \underline{d}_{l,t}\overline{s}_{l,t}c_{l,t}^{u} - \sum_{k: l\in R_{k}} y_{k,l,t}c_{l,t}^{u} - \underline{d}_{l,t}\alpha_{l,t} - \overline{s}_{l,t}\beta_{l,t},&&\quad \forall l\in[L], t\in[T],\\
        &\quad \eta_{l,t} \geq \overline{d}_{l,t}\underline{s}_{l,t}c_{l,t}^{u} - \sum_{k: l\in R_{k}} y_{k,l,t}c_{l,t}^{u} - \overline{d}_{l,t}\alpha_{l,t} - \underline{s}_{l,t}\beta_{l,t},&&\quad \forall l\in[L], t\in[T],\\
        &\quad \eta_{l,t} \geq \underline{d}_{l,t}\underline{s}_{l,t}c_{l,t}^{u} - \sum_{k: l\in R_{k}} y_{k,l,t}c_{l,t}^{u} - \underline{d}_{l,t}\alpha_{l,t} - \underline{s}_{l,t}\beta_{l,t},&&\quad \forall l\in[L], t\in[T],\\
        &\quad \eta_{l,t} \geq -\overline{d}_{l,t}\underline{s}_{l,t}c_{l,t}^{o} + \sum_{k: l\in R_{k}} y_{k,l,t}c_{l,t}^{o} - \overline{d}_{l,t}\alpha_{l,t} - \underline{s}_{l,t}\beta_{l,t},&&\quad \forall l\in[L], t\in[T],\\
        &\quad \eta_{l,t} \geq -\underline{d}_{l,t}\overline{s}_{l,t}c_{l,t}^{o} + \sum_{k: l\in R_{k}} y_{k,l,t}c_{l,t}^{o} - \underline{d}_{l,t}\alpha_{l,t} - \overline{s}_{l,t}\beta_{l,t},&&\quad \forall l\in[L], t\in[T]. \label{inner_max_appendix_last}
    \end{alignat}
\end{subequations}
}
\end{spacing}

\end{proof}

\clearpage
\newpage

\section{Proof of Proposition \ref{LinearH} }\label{Appx:H_Reform}
\setcounter{figure}{0}
\setcounter{table}{0}
\begin{proof}
Not that problem $Q^A(\textbf{\emph{x}}, \bm{\xi})$ in \eqref{A_dual} is separable by each $t$. Therefore, w.l.o., we can rewrite $Q^A(\textbf{\emph{x}}, \bm{\xi})= \sum_{t} Q_{t}^A(\textbf{\emph{x}}, \bm{d}_t, \bm{s}_t)$, where for fixed $\pmb{x}$, $t$, $\bm{d}_t \in [\bm{\underline{d}}_t, \bm{\overline{d}}_t]$, and $\bm{s}_t \in[\bm{\underline{s}}_t, \bm{\overline{s}}_t]$
{
\setlength{\abovedisplayskip}{4pt}
\setlength{\belowdisplayskip}{4pt}
\begin{subequations}
\label{A_dual_a}
\begin{alignat}{3}
    Q_t^{A}(\textbf{\emph{x}}, \bm{\xi})=\quad \max_{\bm{\rho}_t}  &\Bigg \{ \sum_{l=1}^{L} d_{l,t}s_{l,t}\rho_{l,t} + \sum_{k\in K} \lambda_{k,t}x_{k}h_{k} \Bigg \}\quad \label{A_dual_a_1} \\
    \mbox{s.t.} &\quad \rho_{l,t} + \lambda_{k,t} \leq c_{k,l,t}, && \forall k \in [K], l\in R_{k} \label{A_dual_a_2},\\
    &\quad \rho_{l,t} \leq c_{l,t}^{u}, && \forall l \in [L] \label{A_dual_a_3},\\
    &\quad \lambda_{k,t} \leq c_{k,t}^{o}, && \forall k \in [K] \label{A_dual_a_4}.
\end{alignat}
\end{subequations}}
Accordingly,
{
\small
\setlength{\abovedisplayskip}{4pt}
\setlength{\belowdisplayskip}{4pt}
\begin{align}
    \max \limits_{(\bm{d,s})\in \mathcal{U}} \left\{Q^A(\bm{x}, \bm{\xi}) - \sum \limits_{l=1}^L \sum \limits_{t =1}^T(d_{l,t}\alpha_{l,t} + s_{l,t}\beta_{l,t})\right\}& 
    \equiv \sum \limits_{t=1}^{T} \max\limits_{\substack{\bm{d}_t \in [\bm{\underline{d}}_t, \bm{\overline{d}}_t]\\ \bm{s}_t \in[\bm{\underline{s}}_t, \bm{\overline{s}}_t]}} \left\{Q_{t}^A(\textbf{\emph{x}}, \bm{d}_t, \bm{s}_t)- \sum \limits_{l=1}^L (d_{l,t}\alpha_{l,t} + s_{l,t}\beta_{l,t}) \right\}, \nonumber\\
     &\equiv \sum \limits_{t=1}^{T}\max\limits_{\substack{\bm{d}_t \in [\bm{\underline{d}}_t, \bm{\overline{d}}_t]\\ \bm{s}_t \in[\bm{\underline{s}}_t, \bm{\overline{s}}_t]}} h_{t}(\bm{x},\bm{d}_t, \bm{s}_t, \bm{\alpha}_t, \bm{\beta}_t)\label{maxH_f},
\end{align}}
where
{
\setlength{\abovedisplayskip}{4pt}
\setlength{\belowdisplayskip}{4pt}
\begin{equation}
    h_{t}(\bm{x}, \bm{d}_t, \bm{s}_t, \bm{\alpha}_t, \bm{\beta}_t) := \underset{\bm{\rho}_{t}\in P_{t}}{\mbox{max}} \left\{\sum_{l=1}^{L} d_{l,t}s_{l,t}\rho_{l,t} + \sum_{k\in K} \lambda_{k,t}x_{k}h_{k}  - \sum \limits_{l=1}^L (d_{l,t}\alpha_{l,t} + s_{l,t}\beta_{l,t}) \right\}\label{func_f}.
\end{equation}}
and $P_t := \{\eqref{A_dual_a_2}-\eqref{A_dual_a_4}\}$. Problem \eqref{func_f} is not solvable directly due to the trilinear terms between the dual variables $\bm{\rho}_t$, $\bm{d}_t$, and $\bm{s}_t$. Note that we could first perform maximization over $\bm{d}_t$ and $\bm{s}_t$, i.e.,
{
\setlength{\abovedisplayskip}{4pt}
\setlength{\belowdisplayskip}{4pt}
\begin{subequations}
\label{H_2}
\begin{alignat}{3}
    h_{t}(\bm{x}, \bm{\rho}_t,\bm{\lambda}_t, \bm{\alpha}_t, \bm{\beta}_t):=  & \max_{\bm{d}_t,\bm{s}_t}\Bigg \{ \sum_{l=1}^{L} d_{l,t}s_{l,t}\rho_{l,t} + \sum_{k\in K} \lambda_{k,t}x_{k}h_{k} - \sum_{l=1}^{L}(d_{l,t}\alpha_{l,t} && + s_{l,t}\beta_{l,t}) \Bigg \} \label{H_2_1} \\
    \mbox{s.t.} &\quad \underline{d}_{l,t} \leq d_{l,t} \leq   \overline{d}_{l,t}, && \forall l \in [L]\label{H_2_2},\\
    &\quad \underline{s}_{l,t} \leq s_{l,t} \leq  \overline{s}_{l,t}, && \forall l \in [L] \label{H_2_3}.
\end{alignat}
\end{subequations}}
\indent The above problem is separable in $l$. Thus, it is easy to verify that the optimal $d_{l,t}$ is attained at either $d_{l,t}= \underline{d}_{l,t}$ or $d_{l,t}=  \overline{d}_{l,t}$, and the optimal $s_{l,t}$ is attained at either $s_{l,t} =\underline{s}_{l,t}$ or $s_{l,t}=\overline{s}_{l,t}$. Accordingly,  we define binary variables $g_{l,t} = 1$ if $d_{l,t} = \overline{d}_{l,t}$ and $g_{l,t} = 0$ if $d_{l,t} = \underline{d}_{l,t}$. And for all $\in [L]$, we define binary variables $z_{l,t} = 1$ if $s_{l,t} = \overline{s}_{l,t}$ and $z_{l,t} = 0$ if $s_{l,t} = \underline{s}_{l,t}$. Then, we have:
\allowdisplaybreaks
{
\setlength{\abovedisplayskip}{4pt}
\setlength{\belowdisplayskip}{4pt}
\begin{align}
    d_{l,t}^{*} = \underline{d}_{l,t} + ( \overline{d}_{l,t} - \underline{d}_{l,t})g_{l,t},  \label{solution_d}\\
    s_{l,t}^{*} = \underline{s}_{l,t} + ( \overline{s}_{l,t} - \underline{s}_{l,t})z_{l,t},  \label{solution_s}
\end{align}}
\indent With \eqref{solution_d} and \eqref{solution_s}, we drive the following equivalent reformulation of $h_{t}(\bm{x}, \bm{d}_t, \bm{s}_t, \bm{\alpha}_t, \bm{\beta}_t)$ as \eqref{h_reform}.
\allowdisplaybreaks
{
\setlength{\abovedisplayskip}{4pt}
\setlength{\belowdisplayskip}{4pt}
\begin{subequations}
    \label{h_reform}
    \begin{alignat}{3}
       \underset{\substack{\bm{\rho}_{t}, \bm{\lambda}_t,\\ \bm{g}_t, \bm{z}_t}}{\mbox{max}}  &\  \sum_{k=1}^{K} x_{k}h_{k}\lambda_{k,t} - \sum_{l=1}^{L}\Big\{ \Delta s_{l,t}z_{l,t}\beta_{l,t}+ \underline{d}_{l,t}\alpha_{l,t} + \Delta d_{l,t}g_{l,t}\alpha_{l,t} + \underline{s}_{l,t}\beta_{l,t} \Big\} &&\quad \notag\\
        &\quad + \sum_{l=1}^{L}\Big\{ \underline{d}_{l,t}\underline{s}_{l,t}\rho_{l,t} +  \Delta s_{l,t}\underline{d}_{l,t}z_{l,t}\rho_{l,t} + \Delta d_{l,t}\underline{s}_{l,t}g_{l,t}\rho_{l,t}+ \Delta d_{l,t}\Delta s_{l,t}g_{l,t} && z_{l,t}\rho_{l,t}\Big\} \quad \label{h_reform_1}\\
        \mbox{s.t.}&\quad \bm{\rho}_t\in P_t, &&\quad \label{h_reform_2}\\
        &\quad g_{l,t}\in \{0,1\}, &&\quad \forall l \in [L],\\
        &\quad z_{l,t}\in \{0,1\}, &&\quad \forall l \in [L] \label{h_reform_3},
    \end{alignat}
\end{subequations}}
where $\Delta d_{l,t} = \overline{d}_{l,t} - \underline{d}_{l,t}$ and $\Delta s_{l,t} = \overline{s}_{l,t} - \underline{s}_{l,t}$.

Then, we define variables $v_{l,t} = g_{l,t}z_{l}$. And we introduce the following McCormick inequalities for variables $v_{l,t}$.
\allowdisplaybreaks
\begin{subequations}
\setlength{\abovedisplayskip}{1.5pt}
    \begin{alignat}{2}
        &v_{l,t} - g_{l,t} \leq 0, &&\quad v_{l,t}  - z_{l,t} \leq 0\label{mc_1},\\
        &v_{l,t} \geq 0, &&\quad v_{l,t}- g_{l,t} - z_{l,t} \geq -1\label{mc_2}.
        \end{alignat}
\end{subequations}

Accordingly, we can derive the following  mixed integer quadratic programming (MIQP) reformulation of \eqref{func_f}:
{
\setlength{\abovedisplayskip}{4pt}
\setlength{\belowdisplayskip}{4pt}
\begin{subequations}
\setlength{\abovedisplayskip}{1.5pt}
    \label{h_temp_a}
    \begin{alignat}{3}
        \underset{\substack{\bm{\rho}_t, \bm{\lambda}_t, \bm{g}_t,\\ \bm{z}_t, \bm{v}_t}}{\mbox{max}} &\quad \sum_{k=1}^{K} x_{k}h_{k}\lambda_{k,t} -\sum_{l=1}^{L}\left\{\underline{d}_{l,t}\alpha_{l,t} + \Delta d_{l,t}g_{l,t}\alpha_{l,t} + \underline{s}_{l,t}\beta_{l,t} + \Delta s_{l,t}z_{l,t}\beta_{l,t}\right\} &&\quad \notag\\
        &\quad + \sum_{l=1}^{L} \left\{\underline{d}_{l,t}\underline{s}_{l,t}\rho_{l,t} +  \Delta s_{l,t}\underline{d}_{l,t}z_{l,t}\rho_{l,t} + \Delta d_{l,t}\underline{s}_{l,t}g_{l,t}\rho_{l,t} + \Delta d_{l,t}\Delta s_{l,t}v_{l,t}\rho_{l,t} \right\} &&\quad \label{h_temp_a1}\\
        \mbox{s.t.}&\quad \eqref{A_dual_a_2}-\eqref{A_dual_a_4}, \eqref{h_reform_2}-\eqref{h_reform_3}, &&\quad \label{htemp_a2}\\
        &\quad \eqref{mc_1}-\eqref{mc_2}, &&\quad \label{h_temp_a3}.
    \end{alignat}
\end{subequations}}

In formulation \ref{h_temp_a}, there are three bilinear terms in objective function \eqref{h_temp_a1}. To linearize the problem, we introduce three decision variables $e_{l,t} = z_{l,t} \rho_{l,t}$, $q_{l,t} = g_{l,t} \rho_{l,t}$, $r_{l,t}=v_{l,t} \rho_{l,t}$ and the following McCormick inequalities for them.
{
\setlength{\abovedisplayskip}{4pt}
\setlength{\belowdisplayskip}{4pt}
\begin{subequations}
\setlength{\abovedisplayskip}{1.5pt}
    \begin{alignat}{2}
        &e_{l,t} - \underline{\rho}_{l,t} v_{l,t}\geq 0, &&\quad e_{l,t}  - c_{l,t}^{u}v_{l,t} \leq 0\label{mc_3}\\
        &e_{l,t} - \rho_{l,t} + \underline{\rho}_{l,t}(1-v_{l,t})\leq 0, &&\quad e_{l,t}- \rho_{l,t} +c_{l,t}^u (1-v_{l,t}) \geq 0\label{mc_4},\\
        &q_{l,t} - \underline{\rho}_{l,t} g_{l,t}\geq 0, &&\quad q_{l,t}  - c_{l,t}^{u}g_{l,t} \leq 0\label{mc_5},\\
        &q_{l,t} - \rho_{l,t} + \underline{\rho}_{l,t}(1-g_{l,t})\leq 0, &&\quad q_{l,t}- \rho_{l,t} +c_{l,t}^u (1-g_{l,t}) \geq 0\label{mc_6},\\
        &r_{l,t} - \underline{\rho}_{l,t} v_{l,t}\geq 0, &&\quad r_{l,t}  - c_{l,t}^{u}v_{l,t} \leq 0\label{mc_7},\\
        &r_{l,t} - \rho_{l,t} + \underline{\rho}_{l,t}(1-v_{l,t})\leq 0, &&\quad r_{l,t}- \rho_{l,t} +c_{l,t}^u (1-v_{l,t}) \geq 0\label{mc_8},
        \end{alignat}
\end{subequations}}
where $\underline{\rho}_{l,t}$ is a big-M coefficient. Hence, $\max\limits_{\substack{\bm{d}_t \in [\bm{\underline{d}}_t, \bm{\overline{d}}_t]\\ \bm{s}_t \in[\bm{\underline{s}}_t, \bm{\overline{s}}_t]}} h_{t}(\bm{x}, \bm{d}_t, \bm{s}_t, \bm{\alpha}_t, \bm{\beta}_t)$ is equivalent to the following MILP problem.
\begin{subequations}
\setlength{\abovedisplayskip}{1.5pt}
    \label{h_final_a}
    \begin{alignat}{3}
        \underset{\substack{\bm{\rho}_t, \bm{\lambda}_t, \bm{g}_t, \bm{z}_t,\\ \bm{v}_t, \bm{r}_t, \bm{q}_t, \bm{e}_t}}{\mbox{max}} &\quad \sum_{k=1}^{K} x_{k}h_{k}\lambda_{k,t} -\sum_{l=1}^{L}\left\{ \right.\underline{d}_{l,t}\alpha_{l,t} + \Delta d_{l,t}g_{l,t}\alpha_{l,t} + \underline{s}_{l,t}\beta_{l,t} + \Delta s_{l,t}z_{l,t} && \beta_{l,t} \left.\right\} \quad \notag\\
        &\quad + \sum_{l=1}^{L} \left\{\underline{d}_{l,t}\underline{s}_{l,t}\rho_{l,t} +  \Delta s_{l,t}\underline{d}_{l,t}r_{l,t} + \Delta d_{l,t}\underline{s}_{l,t}q_{l,t} + \Delta d_{l,t}\Delta s_{l,t}e_{l,t} \right\} &&\quad \label{h_final_a1}\\
        \mbox{s.t.}&\quad \rho_{l,t} \geq \underline{\rho}_{l,t},&&\quad \forall l\in [L],\\
        &\quad \eqref{A_dual_a_2}-\eqref{A_dual_a_4}, \eqref{h_reform_2}-\eqref{h_reform_3}, &&\quad \label{h_final_a2}\\
        &\quad \eqref{mc_1}-\eqref{mc_2}, \eqref{mc_3}-\eqref{mc_8}. &&\quad\label{h_final_a3}
    \end{alignat}
\end{subequations}

Replacing $\max h(\cdot)$ in \ref{maxH_f} with its equivalent reformulation in \eqref{h_final_a} and summing over all $t \in[T]$, we derive the following equivalent reformulation.
\begin{subequations}
\setlength{\abovedisplayskip}{6pt}
    \label{whole_h_final_a}
    \begin{alignat}{3}
        \underset{\substack{\bm{\rho}, \bm{\lambda}, \bm{g}, \bm{z},\\ \bm{v}, \bm{r}, \bm{q}, \bm{e}}}{\mbox{max}} &\quad \sum_{t=1}^{T} \Bigg\{ \sum_{k=1}^{K} x_{k}h_{k}\lambda_{k,t} -\sum_{l=1}^{L}\Big[\underline{d}_{l,t}\alpha_{l,t} + \Delta d_{l,t}g_{l,t}\alpha_{l,t} + \underline{s}_{l,t}\beta_{l,t} + \Delta s_{l,t} &&  z_{l,t}\beta_{l,t}\Big] \quad \notag\\
        &\quad + \sum_{l=1}^{L} \left[\underline{d}_{l,t}\underline{s}_{l,t}\rho_{l,t} +  \Delta s_{l,t}\underline{d}_{l,t}r_{l,t} + \Delta d_{l,t}\underline{s}_{l,t}q_{l,t} + \Delta d_{l,t}\Delta s_{l,t}e_{l,t} \right] \Bigg \} &&\quad \label{whole_h_final_a1}\\
        \mbox{s.t.}&\quad \rho_{l,t} \geq \underline{\rho}_{l,t},&&\quad \forall l\in [L], t\in [T],\\
        &\quad \eqref{A_dual_a_2}-\eqref{A_dual_a_4}, \eqref{h_reform_2}-\eqref{h_reform_3}, &&\quad \forall t\in [T], \label{whole_h_final_a2}\\
        &\quad \eqref{mc_1}-\eqref{mc_2}, \eqref{mc_3}-\eqref{mc_8} &&\quad \forall t\in [T] \label{whole_h_final_a3}.
    \end{alignat}
\end{subequations}
\end{proof}
\clearpage
\newpage

\section{Proof of Proposition \ref{pro_rho}}\label{pro_rho_proof}
\setcounter{figure}{0}
\setcounter{table}{0}
\begin{proof}
According to \eqref{A_dual_2}, we have $\rho_{l,t} \leq c_{k,l,t} - \lambda_{k,t}, \forall k\in [K], l\in R_{k}$. Hence, we have \eqref{lower_1}.
{
\setlength{\abovedisplayskip}{4pt}
\setlength{\belowdisplayskip}{4pt}
\begin{equation}\label{lower_1}
\setlength{\abovedisplayskip}{6pt}
    \rho_{l,t} \leq \underset{k:l\in R_k}{\mbox{min}} (c_{k,l,t} - \lambda_{k,t})
\end{equation}}
As we consider a maximization problem in \eqref{A_dual} and $d_{l,t} s_{l,t} > 0$ always holds, it is always optimal to increase $\rho_{l,t}$ in the feasible region. Hence, we have $\rho_{l,t}^{*} = \underset{k:l\in R_k}{\mbox{min}} (c_{k,l,t} - \lambda_{k,t})$ or $\rho_{l,t}^{*} = c_{l,t}^u$.

Then, according to constraints \eqref{A_dual_4}, we have $c_{k,l,t} - \lambda_{k,t} \geq c_{k,l,t} - c_{k,t}^o$. This leads to $\rho_{l,t}^{*} = \underset{k:l\in R_k}{\mbox{min}} (c_{k,l,t} - \lambda_{k,t}) \geq \underset{k:l\in R_k}{\mbox{min}} (c_{k,l,t} - c_{k,t}^o)$. As we assume $c_{l,t}^u + c_{k,t}^o > c_{k,l,t}$, we have $\underset{k:l\in R_k}{\mbox{min}} (c_{k,l,t} - c_{k,t}^o) \leq c_{l,t}^u$.

Finally. we have that $\rho_{l,t} \geq \underset{k:l\in R_k}{\mbox{min}} (c_{k,l,t} - c_{k,t}^o)$ in the optimal solution.

\end{proof}

\newpage

\section{Proof of Proposition \ref{pos_alpha}}\label{proof_alpha}
\setcounter{figure}{0}
\setcounter{table}{0}
\begin{proof}
Observe from the objective of problem \eqref{A_inner} that variable $\alpha_{l,t}$ are multiplied by parameters $\mu_{l,t}^{d}$ and variables $d_{l,t}$, for all $l\in [L], t\in [T]$. And so, for fixed $\bm{\beta}$ and $\bm{x}$, the joint contribution $\bm{\alpha}$ and $\bm{d}$ to the objective of problem \eqref{A_inner} equals:
{
\setlength{\abovedisplayskip}{4pt}
\setlength{\belowdisplayskip}{4pt}
\begin{equation}
\setlength{\abovedisplayskip}{6pt}
    \sum_{t\in [T]}\sum_{l\in [L]} \mu_{l,t}^{d}\alpha_{l,t} + \sum_{t\in [T]} \underset{d_{l,t} \in [\underline{d}_{l,t} ,  \overline{d}_{l,t}]}{max}\sum_{l\in [L]} d_{l,t}(s_{l,t}\rho_{l,t} - \alpha_{l,t}) +K,
\end{equation}}
where $K$ is a non-negative term that equal to:
{
\setlength{\abovedisplayskip}{4pt}
\setlength{\belowdisplayskip}{4pt}
\begin{equation}
\setlength{\abovedisplayskip}{6pt}
    \sum_{t\in [T]}\sum_{l\in [L]} \mu_{l,t}^{s}\beta_{l,t} + \sum_{t\in [T]} \underset{s_{l,t} \in [\underline{s}_{l,t} ,  \overline{s}_{l,t}]}{\mbox{max}}\sum_{l\in [L]} - s_{l,t}\beta_{l,t}.
\end{equation}}

First, we suppose that $\alpha_{l,t} > 0$. If $\alpha_{l,t} > \overline{s}_{l,t} c_{l,t}^u$ where $c_{l,t}^u$ is the upper bound of $\rho_{l,t}$. The optimal solution of $d_{l,t}^{*}$ will be $\underline{d}_{l,t}$. In this case, $\alpha_{l,t}$ contributes to the objective value of problem \eqref{A_inner} by $(\mu_{l,t}^{s} - \underline{d}_{l,t})\alpha_{l,t}$. Let $\alpha_{l,t}'  = \alpha_{l,t} - \epsilon$ with $\epsilon > 0$. Since $(\mu_{l,t}^{s} - \underline{d}_{l,t}) > 0$, then $(\mu_{l,t}^{s} - \underline{d}_{l,t}) \alpha_{l,t}'  < (\mu_{l,t}^{s} - \underline{d}_{l,t}) \alpha_{l,t}$, i.e., $\alpha_{l,t}' $ improves the objective value of problem \eqref{A_inner}. It follows that without loss of optimality $\overline{\alpha}_{l,t} =  \overline{s}_{l,t} c_{l,t}^u$ is valid upper bound on variables $\alpha_{l,t}$.

Second, we suppose that $\alpha_{l,t} < 0$. If $\alpha_{l,t} < -  \overline{s}_{l,t} |\underline{\rho}_{l,t}|$ which leads to $s_{l,t}\rho_{l,t} - \alpha_{l,t} \geq 0$, the optimal solution of $d_{l,t}^{*}$ always is $ \overline{d}_{l,t}$. In this case, $\alpha_{l,t}$ contributes to the objective value of problem \eqref{A_inner} by $(\mu_{l,t}^{s} -  \overline{d}_{l,t})\alpha_{l,t}$. Let $\alpha_{l,t}'  = \alpha_{l,t} + \epsilon$ with $\epsilon > 0$. Since $(\mu_{l,t}^{s} -  \overline{d}_{l,t}) \leq 0$, then $(\mu_{l,t}^{s} -  \overline{d}_{l,t}) \alpha_{l,t}'  < (\mu_{l,t}^{s} - \underline{d}_{l,t}) \alpha_l^{t}$, i.e., $\alpha_{l,t}' $ improves the objective value of problem \eqref{A_inner}. It follows that without loss of optimality $\underline{\alpha}_{l,t} = - \overline{s}_{l,t} |\underline{\rho}_{l,t}|$ is valid lower bound on variables $\alpha_{l,t}$.

Follow the same proof process, we can prove the $\overline{\beta}_{l,t}$ and $\underline{\beta}_{l,t}$, defined in \eqref{a_b_lower}, are respectively valid upper and lower bounds on variable $\beta_{l,t}$.

\end{proof}

\clearpage
\newpage

\section{Sample average approximation}\label{SAA}
\setcounter{figure}{0}
\setcounter{table}{0}
First, we solve the sample average approximation (SAA) formulations of \eqref{eq:optim2} as formulation \eqref{eq:optim9} and \eqref{FA} as formulation \eqref{FA_SAA} which are mixed-integer linear programs (MILP).
\begin{subequations}
\setlength{\abovedisplayskip}{4pt}
    \label{eq:optim9}
    \begin{alignat}{3}
        {\text{min}} &\quad \sum_{k=1}^{K} c_{k}x_{k} + \sum_{k=1}^{K} \sum_{l=1}^{L} \sum_{t=1}^{T} c_{k,l,t}y_{k,l,t} + \frac{1}{N} \sum_{n=1}^{N} Q^{n}(\textbf{\emph{x}}, \bm{y}, \bm{\xi}^{n})  &&\quad \label{eq:const91} \\
        \mbox{s.t.} &\quad \eqref{eq:cont22}-\eqref{eq:cont24}, &&\quad\\
        &\quad \sum_{k: l\in R_{k}} y_{k,l,t}^n + u_{l,t}^n - o_{l,t}^n = d_{l,t}^n s_{l,t}^n, &&\quad \forall l \in [L], t \in [T], \label{eq:const92}\\
        &\quad o_{l,t}^n, u_{l,t}^n \in \mathbb{R}_{+}, &&\quad \forall l\in [L], t\in [T], n\in [N]. \label{eq:const93}
    \end{alignat}
\end{subequations}
\begin{subequations}
\setlength{\abovedisplayskip}{4pt}
    \label{FA_SAA}
    \begin{alignat}{3}
        {\text{min}} &\quad \sum_{k=1}^{K} c_{k}x_{k} + \frac{1}{N} \sum_{n=1}^{N}Q^{A,n}(\textbf{\emph{x}}, \bm{\xi}^n)&&\quad \label{FA_SAA1} \\
        \mbox{s.t.} &\quad \eqref{FA_1}-\eqref{FA_3}, &&\quad \label{FA_SAA2}\\
        &\quad \sum_{l\in R_{k}} y_{k,l,t}^n + o_{k,t}^n = x_{k}h_{k}, && \forall k \in [K], t\in [T], n\in [N] \label{FA_SAA3},\\
        &\quad \sum_{k: l\in R_{k}} y_{k,l,t}^n + u_{l,t}^n = d_{l,t}^n s_{l,t}^n, && \forall l \in [L], t\in [T], n\in [N],\\
        &\quad y_{k,l,t}^n, o_{k,t}^n, u_{l,t}^n \in \mathbb{R}_{+},  &&\forall k\in [K], l\in [L], t\in [T], n\in [N]. \label{FA_SAA4}
    \end{alignat}
\end{subequations}
\indent In formulation \eqref{eq:optim9} and \eqref{FA_SAA}, we replace all scenario-dependent parameters, variables and constraints with scenario index $n\in[N]$. Parameters $d_{l,t}$ and $s_{l,t}$ are replaced by $d_{l,t}^{n}$ and $s_{l,t}^{n}$ to represent the realized services demands and duration in scenario $n\in [N]$. And we use variables $y_{k,l,t}^{n}$, $o_{l,t}^{n}$, and $u_{l,t}^{n}$ to represent the second stage decision variables under scenarios $n\in [N]$.

According to Law of Large Number, the optimal objective values of formulation \eqref{eq:optim9} will converge to the optimal values of formulation \eqref{eq:optim2} with probability one as $N\rightarrow \infty$ \citep{homem2014monte, shehadeh2021using}. However, the computational effort and solution time of the SAA formulations will increase when we increase the sample size. Hence, we use Algorithm \ref{algorithm1} to determine an appropriate sample size N and obtain near-optimal solutions of our model via SAA formulation.

In Algorithm \ref{algorithm1}, we start with a small initial sample size. In first step, we will solve the SAA formulation with designated sample size. In step 1.1, we generate $N$ random samples of demands and service durations. Then, in step 1.2, we will solve the SAA formulation with samples generated in step 1.1 and record its optimal objective value $v_{N}^{k}$ and solution $x_{N}^{k}$. In step 1.3, we generate $N'$ random samples of demands and service duration. We evaluate the solution $x_{N}^{k}$ by Monte Carlo simulation with the generated $N'$ random samples and record the objective value $v_{N'}^{k}$. We will run step 1 with $K$ replicates. All results are recorded. In second step, we calculate the average of $v_{N}^{k}$ and $v_{N'}^{k}$ as $\overline{v}_{N}$ and $\overline{v}_{N'}$. According to \cite{shehadeh2021using}, $\overline{v}_{N}$ and $\overline{v}_{N'}$ are statistical lower and upper bound of the optimal value of original formulation. In step 3, we calculate the approximate optimality index as $AOI_{N} = \frac{\overline{v}_{N'} - \overline{v}_{N}}{\overline{v}_{N'}}$. $AOI_{N}$ approximately estimate the optimality gap of lower bound and upper bound of original formulation. Therefore, if $AOI_{N}$ is smaller than the predetermined termination tolerance $\epsilon$, we believe the optimality gap is small enough. The algorithm will terminate and output the termination sample size $N$ , objective value $\overline{v}_{N}$, $\overline{v}_{N'}$ and approximate optimality index $AOI_{N}$. Otherwise, we will increase the sample size then go to step 1.
\begin{algorithm}
\caption{Monte Carlo optimization (MCO) Method}
\label{algorithm1}
\begin{algorithmic}[1]
    \STATE Initial sample size $N_{0}(10)$, the number of replicates $K(10)$, the number of scenarios in the monte carlo simulation step, $N'(1000)$, and a termination tolerance $\epsilon(0.01)$.
	\STATE The sample size of SAA $N$, lower bound $\overline{v}_{N}$ and upper bound $\overline{v}_{N'}$ and the approximate optimality index $AOI_{N}$.
	\STATE $AOI_{0} \longleftarrow 1$;
	\STATE $N\longleftarrow N_{0} / 2$;
	\WHILE{$\left| AOI_{N} \right| > \epsilon$}
	    \STATE $N\longleftarrow 2N$;
	    \STATE \textbf{Step 1: MCO Procedure};
		\FOR{$k=1,...,K$}
		    \STATE \textbf{Step 1.1} Scenario Generation
		    \STATE \quad Generate $N$ i.i.d scenarios of $\bm{d, s}$ for all $t\in [T], k\in[K], l\in[L]$.
		    \STATE \textbf{Step 1.2} Solving the SAA formulation
		    \STATE \quad Solve the SAA formulation with the scenarios generated in Step 1.1 and record the corresponding optimal objective value $\overline{v}_{N}^{k}$ and optimal staffing number $(\bm{x})_{N}^{k}$.
		    \STATE \textbf{Step 1.3} Cost Evaluation using Monte Carlo Simulation
		    \STATE \quad  Generate $N'$ i.i.d scenarios of $\bm{d, s}$ for all $t\in [T], k\in[K], l\in[L]$.
		    \STATE \quad  Use $(\bm{x})_{N}^{k}$ and $\bm{d, s}$ to solve $N'$ second stage optimization problems and evaluate the objective value $\overline{v}_{N'}^{k}$ as:
		    \begin{equation*}
		        \overline{v}_{N'}^{k} = \frac{1}{N'} \sum_{n'=1}^{N'} \overline{v}_{n'}^{k}
		    \end{equation*}
		   \ENDFOR
		\STATE \textbf{Step 2: Compute $\overline{v}_{N}$ and $\overline{v}_{N'}$}
		\begin{equation*}
		\setlength{\abovedisplayskip}{4pt}
		\setlength{\belowdisplayskip}{4pt}
		    \overline{v}_{N} = \frac{1}{K}\sum_{k}^{K} \overline{v}_{N}^{k}\quad \quad \overline{v}_{N'} = \frac{1}{K}\sum_{k}^{K} \overline{v}_{N'}^{k}
		\end{equation*}
		\STATE\textbf{Step 3: Compute the Approximate Optimality}
		\begin{equation*}
		    AOI_{N} = \frac{\overline{v}_{N'} - \overline{v}_{N}}{\overline{v}_{N'}}
		\end{equation*}
	\ENDWHILE
\end{algorithmic}
\end{algorithm}

In \ref{CR_1} and \ref{CR_2}, we present the result of Approximate Optimality Index and 95\% Confidence Interval of the statistical lower bound and upper bound on the objective value, termination sample size and CPU time for 15 random generated instances with lognormal, normal and uniform distributions.

\newpage
\section{HSCP Instances}\label{instance_a}
\setcounter{figure}{0}
\setcounter{table}{0}
\begin{table}[htbp]
  \centering
  \caption{HSCP Instances}
    \begin{tabular}{cccc}
    \toprule
    Instance & L & K & T \\
    \midrule
    1 & 4 & 4 & 30 \\
    2 & 4 & 4 & 90 \\
    3 & 4 & 4 & 180 \\
    4 & 4 & 6 & 30 \\
    5 & 4 & 6 & 90 \\
    6 & 4 & 6 & 180 \\
    7 & 4 & 8 & 30 \\
    8 & 4 & 8 & 90 \\
    9 & 4 & 8 & 180 \\
    10 & 6 & 6 & 30 \\
    11 & 6 & 6 & 90 \\
    12 & 6 & 6 & 180 \\
    13 & 6 & 8 & 30 \\
    14 & 6 & 8 & 90 \\
    15 & 6 & 8 & 180 \\
    \bottomrule
    \end{tabular}%
  \label{instance}%
\end{table}%

\newpage
\section{Convergence Results of model for EA decision maker at the Termination Sample Size}\label{CR_1}
\setcounter{figure}{0}
\setcounter{table}{0}
In this section, we show the convergence results of SAA formulation \eqref{eq:optim9}. without spacial mention, we will use the termination MIP Gap 1 \% (We allowed all instances to run for one day, however, the MIP Gaps of some instances remained at around 1 \%). And if there are more than 50000 iterations with no MIP Gap improvement, the Gurobi will also terminate.
{
\renewcommand{\arraystretch}{0.9}
\small
\begin{longtable}{cm{2.cm}<{\centering}ccc}
\caption{The Approximate Optimality Index between the statistical lower bound and upper bound on the objective value and 95 \% Confidence interval (95\% CI) at the terminated sample size. The distribution of demand and service duration are both truncated lognormal distribution.}
\label{MCO_result}\\
\hline
\textbf{Instance} & \textbf{Termination sample size} & \textbf{\textbf{95 \%CI } \bm{$\overline{v}_{N}^{E}$}} & \textbf{\textbf{95 \%CI }  \bm{$\overline{v}_{N'}^{E}$}} & $\left| \textbf{AOI} \right|$\\
\endfirsthead
\multicolumn{5}{r}%
{{\bfseries -- continued from previous page}}\\
	\hline 
	\endhead
	\hline
	\endfoot
	\hline
	\endlastfoot
	\hline
    1 & 160 & [1305223, 1310901] & [1314541, 1318775] & 0.0065  \\
    2 & 160 & [4025209, 4042119] & [4051562, 4062731] & 0.0058  \\
    3 & 160 & [8145422, 8181920] & [8193853, 8227162] & 0.0057  \\
    4 & 80 & [1332869, 1341384] & [1346114, 1352009] & 0.0088  \\
    5 & 160 & [3878104, 3892090] & [3889652, 3894498] & 0.0018  \\
    6 & 160 & [8144640, 8177446] & [8178239, 8211080] & 0.0041  \\
    7 & 160 & [1318371, 1326795] & [1323659, 1327991] & 0.0024  \\
    8 & 160 & [4024278, 4031950] & [4049433, 4060713] & 0.0066  \\
    9 & 160 & [8133427, 8148144] & [8139783, 8153705] & 0.0007  \\
    10 & 160 & [1930063, 1938342] & [1951127, 1953218] & 0.0092  \\
    11 &160 & [5996573, 6022297] & [6032443, 6059512] & 0.0060  \\
    12 & 160 & [12161825, 12190077] & [12236514, 12275592] & 0.0065  \\
    13 & 160 & [1994892, 2006005] & [2011921, 2020811] & 0.0079  \\
    14 & 160 & [5997995, 6027363] & [6028977, 6055946] & 0.0049  \\
    15 & 160 & [12209479, 12249137] & [12274856, 12314217] & 0.0053  \\
    \hline
\end{longtable}%
}

\newpage
{
\renewcommand{\arraystretch}{0.9}
\small
\begin{longtable}{cm{2.5cm}<{\centering}cc}
\caption{The range of CPU time at the Termination Sample Size. The distribution of demand and service duration are both truncated lognormal distribution.}
\label{CPU_time}\\
\hline
\textbf{Instance} & \textbf{Termination sample size} & \textbf{CPU time mean (s)} & \textbf{CPU time range (s)} \\
\endfirsthead
\multicolumn{4}{r}%
{{\bfseries -- continued from previous page}}\\
	\hline 
	\endhead
	\hline
	\endfoot
	\hline
	\endlastfoot
	\hline
    1 & 160 & 0.7 & [0.64, 0.77] \\
    2 & 160 & 2.18 & [2.08, 2.26] \\
    3 & 160 & 4.74 & [4.51, 4.98] \\
    4 & 80 & 0.36 & [0.33, 0.39] \\
    5 & 160 & 2.39 & [2.29, 2.57] \\
    6 & 160 & 4.77 & [4.54, 4.91] \\
    7 & 160 & 0.8 & [0.75, 0.83] \\
    8 & 160 & 2.53 & [2.4, 2.77] \\
    9 & 160 & 5.4 & [5.1, 6.06] \\
    10 & 160 & 1.13 & [1.05, 1.3] \\
    11 & 160 & 3.69 & [3.45, 3.89] \\
    12 & 160 & 7.97 & [7.67, 8.46] \\
    13 & 160 & 1.2 & [1.1, 1.34] \\
    14 & 160 & 3.6 & [3.49, 3.89] \\
    15 & 160 & 7.64 & [7.24, 8.16] \\
    \bottomrule
\end{longtable}%
}

{
\renewcommand{\arraystretch}{0.9}
\small
\begin{longtable}{ccccc}
\caption{CPU time of SAA formulation with sample size = 500 and $L=6$. The distribution of demand and service duration are both truncated lognormal distribution.}
\label{CPU_500} \\
\hline
\textbf{K} & \textbf{T} & \textbf{CPU time mean (s)} & \textbf{CPU time range (s)}& $\left| \textbf{AOI} \right|$\\
\endfirsthead
\multicolumn{5}{r}%
{{\bfseries -- continued from previous page}}\\
	\hline 
	\endhead
	\hline
	\endfoot
	\hline
	\endlastfoot
	\hline
    \multirow{3}[0]{*}{6} & 30 & 3.89 & [3.65, 4.11] &0.0021\\
      & 90 & 11.53 & [11.2, 11.98] & 0.0025\\
      & 180 & 24.38 & [23.17, 26.97] & 0.0025\\
      \hline
    \multirow{3}[0]{*}{8} & 30 & 4.46 & [4.11, 4.67] & 0.0018\\
      & 90 & 12.71 & [11.91, 13.48] & 0.0020\\
      & 180 & 25.46 & [24.57, 27.09] & 0.0006\\
    \hline
\end{longtable}%
}

{
\renewcommand{\arraystretch}{0.9}
\small
\begin{longtable}{cm{2.5cm}<{\centering}ccc}
\caption{The Approximate Optimality Index between the statistical lower bound and upper bound on the objective value and 95 \% Confidence interval (95\% CI) at the terminated sample size. The distribution of demand is uniform distribution and service duration is normal distribution.}
\label{MCO_result_normal}\\
\hline
\textbf{Instance} & \textbf{Termination sample size} & \textbf{\textbf{95 \%CI } \bm{$\overline{v}_{N}^{E}$}} & \textbf{\textbf{95 \%CI }  \bm{$\overline{v}_{N'}^{E}$}} & $\left| \textbf{AOI} \right|$\\
\endfirsthead
\multicolumn{5}{r}%
{{\bfseries -- continued from previous page}}\\
	\hline 
	\endhead
	\hline
	\endfoot
	\hline
	\endlastfoot
	\hline
    1 & 80 & [1381556, 1388663] & [1394684, 1401987] & 0.0095  \\
    2 & 160 & [3966419, 3982868] & [3990798, 3998740] & 0.0050  \\
    3 & 160 & [8348880, 8374187] & [8389175, 8412060] & 0.0047  \\
    4 & 80 & [1348343, 1354935] & [1362006, 1364680] & 0.0086  \\
    5 & 160 & [3976175, 3985779] & [3996330, 3999179] & 0.0042  \\
    6 & 80 & [8342528, 8393709] & [8432410, 8468796] & 0.0098  \\
    7 & 80 & [1380734, 1389435] & [1394946, 1399992] & 0.0089  \\
    8 & 80 & [3945443, 3971844] & [3987209, 4004585] & 0.0093  \\
    9 & 160 & [8373784, 8385241] & [8408210, 8420127] & 0.0041  \\
    10 & 160 & [1998803, 2007622] & [2015306, 2019747] & 0.0071  \\
    11 & 160 & [6140490, 6161328] & [6173546, 6189661] & 0.0050  \\
    12 & 160 & [12514482, 12542256] & [12581618, 12606323] & 0.0052  \\
    13 & 160 & [2056244, 2071048] & [2072046, 2081958] & 0.0064  \\
    14 & 80 & [6123108, 6147422] & [6179068, 6208607] & 0.0095  \\
    15 & 80 & [12521793, 12572312] & [12629120, 12682200] & 0.0086  \\
    \hline
\end{longtable}%
}

{
\renewcommand{\arraystretch}{0.9}
\small
\begin{longtable}{cm{2cm}<{\centering}cc}
\caption{The range of CPU time at the Termination Sample Size. The distribution of demand and service duration are both truncated normal distribution.}
\label{CPU_time_n}\\
\hline
\textbf{Instance} & \textbf{Termination sample size} & \textbf{CPU time mean (s)} & \textbf{CPU time range (s)} \\
\endfirsthead
\multicolumn{4}{r}%
{{\bfseries -- continued from previous page}}\\
	\hline 
	\endhead
	\hline
	\endfoot
	\hline
	\endlastfoot
	\hline
    1 & 80 & 0.35 & [0.33, 0.41] \\
    2 & 160 & 2.14 & [2.08, 2.19] \\
    3 & 160 & 4.57 & [4.43, 4.72] \\
    4 & 80 & 0.36 & [0.33, 0.39] \\
    5 & 160 & 2.29 & [2.19, 2.42] \\
    6 & 80 & 2.45 & [2.35, 2.58] \\
    7 & 80 & 0.41 & [0.38, 0.44] \\
    8 & 80 & 1.24 & [1.21, 1.31] \\
    9 & 160 & 6.7 & [6.12, 7.04] \\
    10 & 160 & 1.13 & [1.05, 1.19] \\
    11 & 160 & 3.53 & [3.46, 3.64] \\
    12 & 160 & 7.49 & [7.13, 7.93] \\
    13 & 160 & 1.22 & [1.11, 1.38] \\
    14 & 80 & 1.87 & [1.79, 2.0] \\
    15 & 80 & 4.09 & [3.82, 5.17] \\
    \bottomrule
\end{longtable}%
}

{
\renewcommand{\arraystretch}{0.9}
\small
\begin{longtable}{ccccc}
\caption{CPU time of SAA formulation with sample size = 500 and $L = 6$. The distribution of demand and service duration are both truncated normal distribution.}
\label{CPU_normal_500} \\
\hline
\textbf{K} & \textbf{T} & \textbf{CPU time mean (s)} & \textbf{CPU time range (s)}& $\left| \textbf{AOI} \right|$\\
\endfirsthead
\multicolumn{5}{r}%
{{\bfseries -- continued from previous page}}\\
	\hline 
	\endhead
	\hline
	\endfoot
	\hline
	\endlastfoot
	\hline
    \multirow{3}[0]{*}{6} & 30 & 3.73 & [3.59, 3.91]&0.0017 \\
      & 90 & 11.69 & [11.4, 11.94] & 0.0021 \\
      & 180 & 22.98 & [22.46, 23.59] &0.0021 \\
    \hline
    \multirow{3}[0]{*}{8} & 30 & 3.95 & [3.83, 4.17]&0.0014 \\
      & 90 & 11.61 & [11.38, 11.85] & 0.0018\\
      & 180 & 23.45 & [23.01, 24.09] &0.0018\\
    \hline
\end{longtable}%
}

{
\renewcommand{\arraystretch}{0.9}
\small
\begin{longtable}{cm{2.5cm}<{\centering}ccc}
\caption{The Approximate Optimality Index between the statistical lower bound and upper bound on the objective value and 95 \% Confidence interval (95\% CI) at the terminated sample size. The distribution of demand and service duration are both uniform distribution.}
\label{MCO_result_uniform}\\
\hline
\textbf{Instance} & \textbf{Termination sample size} & \textbf{\textbf{95 \%CI } \bm{$\overline{v}_{N}^{E}$}} & \textbf{\textbf{95 \%CI }  \bm{$\overline{v}_{N'}^{E}$}} & $\left| \textbf{AOI} \right|$\\
\endfirsthead
\multicolumn{5}{r}%
{{\bfseries -- continued from previous page}}\\
	\hline 
	\endhead
	\hline
	\endfoot
	\hline
	\endlastfoot
	\hline
    1 & 80 & [977155, 983918] & [984647, 989363] & 0.0066  \\
    2 & 80 & [3068198, 3073211] & [3093554, 3095096] & 0.0076  \\
    3 & 80 & [6177607, 6185352] & [6221865, 6225124] & 0.0068  \\
    4 & 80 & [962606, 967277] & [969887, 971768] & 0.0061  \\
    5 & 80 & [3050475, 3055750] & [3071758, 3072860] & 0.0062  \\
    6 & 80 & [6115238, 6138472] & [6148787, 6177390] & 0.0059  \\
    7 & 80 & [1003086, 1008106] & [1009733, 1011627] & 0.0050  \\
    8 & 80 & [3014248, 3026031] & [3034057, 3046354] & 0.0066  \\
    9 & 80 & [6136129, 6161040] & [6171792, 6193458] & 0.0055  \\
    10 & 80 & [1469014, 1477824] & [1482836, 1490789] & 0.0090  \\
    11 & 80 & [4546500, 4560592] & [4584809, 4597373] & 0.0082  \\
    12 & 80 & [8915235, 8926078] & [8985911, 8998265] & 0.0079  \\
    13 & 80 & [1538918, 1545766] & [1551418, 1557400] & 0.0078  \\
    14 & 80 & [4427922, 4432286] & [4462700, 4465527] & 0.0076  \\
    15 & 80 & [8876987, 8885097] & [8941807, 8945098] & 0.0070  \\
    \hline
\end{longtable}
}
\newpage

{
\renewcommand{\arraystretch}{0.8}
\small
\begin{longtable}{cm{2.5cm}<{\centering}cc}
\caption{The range of CPU time at the Termination Sample Size. The distribution of demand and service duration are both uniform distribution.}
\label{CPU_time_uniform}\\
\hline
\textbf{Instance} & \textbf{Termination sample size} & \textbf{CPU time mean (s)} & \textbf{CPU time range (s)} \\
\endfirsthead
\multicolumn{4}{r}%
{{\bfseries -- continued from previous page}}\\
	\hline 
	\endhead
	\hline
	\endfoot
	\hline
	\endlastfoot
	\hline
    1 & 80 & 0.44 & [0.39, 0.56] \\
    2 & 80 & 1.33 & [1.28, 1.44] \\
    3 & 80 & 2.75 & [2.57, 3.36] \\
    4 & 80 & 0.46 & [0.41, 0.63] \\
    5 & 80 & 1.42 & [1.33, 1.68] \\
    6 & 80 & 3.08 & [2.76, 3.56] \\
    7 & 80 & 0.45 & [0.44, 0.49] \\
    8 & 80 & 1.49 & [1.39, 1.77] \\
    9 & 80 & 3.44 & [2.96, 4.14] \\
    10 & 80 & 0.73 & [0.61, 1.07] \\
    11 & 80 & 2.15 & [2.01, 2.54] \\
    12 & 80 & 4.38 & [4.2, 4.48] \\
    13 & 80 & 0.82 & [0.66, 1.38] \\
    14 & 80 & 2.23 & [2.05, 2.7] \\
    15 & 80 & 4.67 & [4.31, 5.92] \\
    \bottomrule
\end{longtable}%
}

{
\renewcommand{\arraystretch}{0.8}
\small
\begin{longtable}{ccccc}
\caption{CPU time of SAA formulation with sample size = 500 and $L=6$. The distribution of demand and service duration are both uniform distribution.}
\label{CPU_500_U} \\
\hline
\textbf{K} & \textbf{T} & \textbf{CPU time mean (s)} & \textbf{CPU time range (s)}& $\left| \textbf{AOI} \right|$\\
\endfirsthead
\multicolumn{5}{r}%
{{\bfseries -- continued from previous page}}\\
	\hline 
	\endhead
	\hline
	\endfoot
	\hline
	\endlastfoot
	\hline
    \multirow{3}[0]{*}{6} & 30 & 5.85 & [5.62, 6.38]&0.0013 \\
      & 90 & 15.88 & [14.25, 18.09]&0.0012 \\
      & 180 & 31.78 & [26.83, 34.83]&0.0010 \\
    \multirow{3}[0]{*}{8} & 30 & 4.57 & [4.41, 4.92]&0.0010 \\
      & 90 & 13.43 & [13.03, 14.23]&0.0009 \\
      & 180 & 27.48 & [25.97, 31.86]&0.0009 \\
    \hline
\end{longtable}%
}

\newpage
\section{Convergence Results of model for FA decision maker at the Termination Sample Size}\label{CR_2}
\setcounter{figure}{0}
\setcounter{table}{0}
In this section, we show the convergence results of SAA formulation \eqref{FA_SAA}. Without spacial mention, we will use the termination MIP Gap 0.1 \% (We allowed all instances to run for one day, however, the MIP Gaps of some instances remained at around 0.1 \%). And if there are more than 50000 iterations with no MIP Gap improvement, the Gurobi will also terminate.

{
\renewcommand{\arraystretch}{0.9}
\small
\begin{longtable}{cm{2.5cm}<{\centering}ccc}
\caption{The Approximate Optimality Index between the statistical lower bound and upper bound on the objective value and 95 \% Confidence interval (95\% CI) at the terminated sample size. The distribution of demand and service duration are both truncated lognormal distribution.}
\label{MCO_result_2}\\
\hline
\textbf{Instance} &  \textbf{Termination sample size} & \textbf{\textbf{95 \%CI } \bm{$\overline{v}_{N}^{F}$}} & \textbf{\textbf{95 \%CI }  \bm{$\overline{v}_{N'}^{F}$}} & $\left| \textbf{AOI} \right|$\\
\endfirsthead
\multicolumn{5}{r}%
{{\bfseries -- continued from previous page}}\\
	\hline 
	\endhead
	\hline
	\endfoot
	\hline
	\endlastfoot
	\hline
    1 & 10 & [797757, 819500] & [808593, 810567] & 0.0012 \\
    3 & 10 & [2476877, 2515034] & [2497294, 2501413] & 0.0014 \\
    5 & 10 & [4802916, 4845726] & [4826203, 4829671] & 0.0007 \\
    7 & 10 & [902648, 930749] & [909915, 912329] & 0.0061 \\
    9 & 10 & [2719304, 2759871] & [2731831, 2743086] & 0.0008 \\
    11 & 10 & [5454537, 5504424] & [5489562, 5506422] & 0.0034 \\
    13 & 10 & [825863, 847786] & [837489, 842430] & 0.0037 \\
    15 & 10 & [2495168, 2533598] & [2516358, 2520467] & 0.0016 \\
    17 & 10 & [4832324, 4876883] & [4851456, 4856927] & 0.0001 \\
    19 & 20 & [1136585, 1156616] & [1142789, 1147654] & 0.0012 \\
    21 & 10 & [3329970, 3376483] & [3357607, 3367852] & 0.0028 \\
    23 & 10 & [6778415, 6833764] & [6799658, 6841951] & 0.0022 \\
    25 & 10 & [1305861, 1336127] & [1325533, 1331182] & 0.0055 \\
    27 & 10 & [3674874, 3730706] & [3705872, 3709887] & 0.0014 \\
    29 & 10 & [7434423, 7497220] & [7470434, 7477952] & 0.0011 \\
    \hline
\end{longtable}%
}
\newpage

{
\renewcommand{\arraystretch}{0.9}
\small
\begin{longtable}{cm{2.5cm}<{\centering}cc}
\caption{The range of CPU time at the Termination Sample Size. The distribution of demand and service duration are both truncated lognormal distribution.}
\label{CPU_time_2}\\
\hline
\textbf{Instance}  & \textbf{Termination sample size} & \textbf{CPU time mean (s)} & \textbf{CPU time range (s)} \\
\endfirsthead
\multicolumn{4}{r}%
{{\bfseries -- continued from previous page}}\\
	\hline 
	\endhead
	\hline
	\endfoot
	\hline
	\endlastfoot
	\hline
    1 & 10 & 0.39 & [0.16, 0.67] \\
    2 & 10 & 0.96 & [0.5, 1.41] \\
    3 & 10 & 2.77 & [1.94, 3.56] \\
    4 & 10 & 0.23 & [0.19, 0.29] \\
    5 & 10 & 0.68 & [0.63, 0.71] \\
    6 & 10 & 1.53 & [1.43, 1.81] \\
    7 & 10 & 0.37 & [0.33, 0.44] \\
    8 & 10 & 1.19 & [1.11, 1.36] \\
    9 & 10 & 2.92 & [2.43, 3.55] \\
    10 & 20 & 0.69 & [0.63, 0.75] \\
    11 & 10 & 1.08 & [1.05, 1.19] \\
    12 & 10 & 2.69 & [2.51, 3.06] \\
    13 & 10 & 0.39 & [0.36, 0.47] \\
    14 & 10 & 1.34 & [1.27, 1.53] \\
    15 & 10 & 3.02 & [2.87, 3.17] \\
    \bottomrule
\end{longtable}%
}

{
\renewcommand{\arraystretch}{0.9}
\small
\begin{longtable}{ccccc}
\caption{CPU time of SAA formulation with sample size = 500 and $L = 6$. The distribution of demand and service duration are both truncated lognormal distribution.}
\label{CPU_500_2} \\
\hline
\textbf{K} & \textbf{T} & \textbf{CPU time mean (s)} & \textbf{CPU time range (s)}& $\left| \textbf{AOI} \right|$\\
\endfirsthead
\multicolumn{4}{r}%
{{\bfseries -- continued from previous page}}\\
	\hline 
	\endhead
	\hline
	\endfoot
	\hline
	\endlastfoot
	\hline
    \multirow{3}[0]{*}{6} & 30 & 100.7 & [76.8, 122.77] & 0.0004 \\
      & 90 & 1201.48 & [933.6, 1520.07] & 0.0062 \\
      & 180 & 1021.57 & [877.25, 1239.94] & 0.0010 \\
    \hline
    \multirow{3}[0]{*}{8} & 30 & 78.22 & [69.67, 87.12] & 0.0005 \\
      & 90 & 487.06 & [472.21, 518.52] & 0.0004\\
      & 180 & 353.69 & [301.54, 395.47] & 0.0006 \\
    \hline
\end{longtable}%
}

{
\renewcommand{\arraystretch}{0.9}
\small
\begin{longtable}{cm{2.5cm}<{\centering}ccc}
\caption{The Approximate Optimality Index between the statistical lower bound and upper bound on the objective value and 95 \% Confidence interval (95\% CI) at the terminated sample size. The distribution of demand and service duration are both truncated normal distribution.}
\label{MCO_result_2_N}\\
\hline
\textbf{Instance} & \textbf{Termination sample size} & \textbf{\textbf{95 \%CI } \bm{$\overline{v}_{N}^{F}$}} & \textbf{\textbf{95 \%CI }  \bm{$\overline{v}_{N'}^{F}$}} & $\left| \textbf{AOI} \right|$\\
\endfirsthead
\multicolumn{5}{r}%
{{\bfseries -- continued from previous page}}\\
	\hline 
	\endhead
	\hline
	\endfoot
	\hline
	\endlastfoot
	\hline
    1 & 10 & [893882, 913589] & [903220, 905948] & 0.0009  \\
    2 & 10 & [2766925, 2800500] & [2783988, 2786829] & 0.0006  \\
    3 & 10 & [5282073, 5325175] & [5304111, 5308345] & 0.0005  \\
    4 & 10 & [989657, 1003856] & [995193, 997365] & 0.0005  \\
    5 & 10 & [3075881, 3105107] & [3089870, 3093619] & 0.0004  \\
    6 & 10 & [5844939, 5888534] & [5859421, 5874568] & 0.0000  \\
    7 & 10 & [884214, 904160] & [893870, 896854] & 0.0013  \\
    8 & 10 & [2798140, 2829687] & [2812171, 2828840] & 0.0023  \\
    9 & 10 & [5309894, 5351091] & [5329231, 5335519] & 0.0004  \\
    10 & 20 & [1280817, 1300915] & [1284966, 1288345] & 0.0033  \\
    11 & 10 & [3686181, 3722289] & [3707673, 3729439] & 0.0039  \\
    12 & 10 & [7436500, 7510222] & [7471978, 7499459] & 0.0017  \\
    13 & 10 & [1402538, 1426072] & [1418141, 1422015] & 0.0041  \\
    14 & 10 & [4016029, 4066268] & [4041444, 4046779] & 0.0007  \\
    15 & 10 & [8440223, 8499236] & [8444406, 8458554] & 0.0022  \\
    \hline
\end{longtable}%
}

{
\renewcommand{\arraystretch}{0.9}
\small
\begin{longtable}{cm{2.5cm}<{\centering}cc}
\caption{The range of CPU time at the Termination Sample Size. The distribution of demand and service duration are both truncated normal distribution.}
\label{CPU_time_2_N}\\
\hline
\textbf{Instance} &  \textbf{Termination sample size} & \textbf{CPU time mean (s)} & \textbf{CPU time range (s)} \\
\endfirsthead
\multicolumn{4}{r}%
{{\bfseries -- continued from previous page}}\\
	\hline 
	\endhead
	\hline
	\endfoot
	\hline
	\endlastfoot
	\hline
    1 & 10 & 0.34 & [0.17, 0.89] \\
    2 & 10 & 0.58 & [0.55, 0.61] \\
    3 & 10 & 1.87 & [1.28, 5.33] \\
    4 & 10 & 0.25 & [0.22, 0.38] \\
    5 & 10 & 0.71 & [0.66, 0.77] \\
    6 & 10 & 1.51 & [1.41, 1.58] \\
    7 & 10 & 0.4 & [0.36, 0.47] \\
    8 & 10 & 1.32 & [1.03, 2.26] \\
    9 & 10 & 4.24 & [3.0, 6.91] \\
    10 & 20 & 0.8 & [0.72, 0.86] \\
    11 & 10 & 1.34 & [1.12, 1.64] \\
    12 & 10 & 2.89 & [2.64, 3.17] \\
    13 & 10 & 0.41 & [0.36, 0.44] \\
    14 & 10 & 1.37 & [1.3, 1.45] \\
    15 & 10 & 2.83 & [2.69, 3.19] \\
    \bottomrule
\end{longtable}%
}

{
\renewcommand{\arraystretch}{0.9}
\small
\begin{longtable}{ccccc}
\caption{CPU time of SAA formulation with sample size = 500 and $L = 6$. The distribution of demand and service duration are both truncated normal distribution.}
\label{CPU_500_n_2} \\
\hline
\textbf{K} & \textbf{T} & \textbf{CPU time mean (s)} & \textbf{CPU time range (s)}& $\left| \textbf{AOI} \right|$\\
\endfirsthead
\multicolumn{4}{r}%
{{\bfseries -- continued from previous page}}\\
	\hline 
	\endhead
	\hline
	\endfoot
	\hline
	\endlastfoot
	\hline
    \multirow{3}[0]{*}{6} & 30 & 103.3 & [99.91, 108.59]&0.0008 \\
      & 90 & 763.9 & [672.18, 898.69] & 0.0009\\
      & 180 & 539.09 & [348.71, 1040.81] & 0.0010 \\
    \hline
    \multirow{3}[0]{*}{8} & 30 & 64.86 & [58.01, 73.57] & 0.0028\\
      & 90 & 472.58 & [426.11, 531.31] & 0.0005 \\
      & 180 & 361.29 & [301.67, 424.21] & 0.0005 \\
    \hline
\end{longtable}%
}

{
\renewcommand{\arraystretch}{0.9}
\small
\begin{longtable}{cm{2.5cm}<{\centering}ccc}
\caption{The Approximate Optimality Index between the statistical lower bound and upper bound on the objective value and 95 \% Confidence interval (95\% CI) at the terminated sample size. The distribution of demand and service duration are both uniform distribution.}
\label{MCO_result_2_U}\\
\hline
\textbf{Instance} & \textbf{Termination sample size} & \textbf{\textbf{95 \%CI } \bm{$\overline{v}_{N}^{F}$}} & \textbf{\textbf{95 \%CI }  \bm{$\overline{v}_{N'}^{F}$}} & $\left| \textbf{AOI} \right|$\\
\endfirsthead
\multicolumn{5}{r}%
{{\bfseries -- continued from previous page}}\\
	\hline 
	\endhead
	\hline
	\endfoot
	\hline
	\endlastfoot
	\hline
    1 & 10 & [878450, 897084] & [888420, 894034] & 0.0039  \\
    2 & 10 & [2623308, 2648105] & [2636686, 2646287] & 0.0022  \\
    3 & 10 & [5574262, 5624156] & [5573696, 5578969] & 0.0041  \\
    4 & 10 & [978637, 998742] & [995634, 998667] & 0.0085  \\
    5 & 10 & [2914379, 2940861] & [2934246, 2942536] & 0.0037  \\
    6 & 10 & [6151851, 6188425] & [6169999, 6174902] & 0.0004  \\
    7 & 10 & [866115, 884066] & [877110, 882576] & 0.0054  \\
    8 & 10 & [2640531, 2665741] & [2655277, 2664965] & 0.0026  \\
    9 & 10 & [5587834, 5637298] & [5587666, 5594447] & 0.0038  \\
    10 & 10 & [1320896, 1343043] & [1338610, 1340305] & 0.0056  \\
    11 & 10 & [3847830, 3896153] & [3881127, 3887101] & 0.0031  \\
    12 & 10 & [7399004, 7453099] & [7431801, 7440199] & 0.0013  \\
    13 & 10 & [1372670, 1398602] & [1390656, 1392612] & 0.0043  \\
    14 & 10 & [4189521, 4244430] & [4235522, 4242783] & 0.0052  \\
    15 & 10 & [8099272, 8149726] & [8148273, 8156607] & 0.0034  \\
    \hline
\end{longtable}%
}

{
\renewcommand{\arraystretch}{0.9}
\small
\begin{longtable}{cm{2.5cm}<{\centering}cc}
\caption{The range of CPU time at the Termination Sample Size. The distribution of demand and service duration are both uniform distribution.}
\label{CPU_time_2_u}\\
\hline
\textbf{Instance} & \textbf{Termination sample size} & \textbf{CPU time mean (s)} & \textbf{CPU time range (s)} \\
\endfirsthead
\multicolumn{4}{r}%
{{\bfseries -- continued from previous page}}\\
	\hline 
	\endhead
	\hline
	\endfoot
	\hline
	\endlastfoot
	\hline
    1 & 10 & 0.34 & [0.17, 0.72] \\
    2 & 10 & 0.82 & [0.52, 1.95] \\
    3 & 10 & 1.32 & [1.2, 1.39] \\
    4 & 10 & 0.25 & [0.23, 0.33] \\
    5 & 10 & 1.89 & [1.09, 3.9] \\
    6 & 10 & 8.91 & [3.72, 14.43] \\
    7 & 10 & 2.35 & [0.39, 8.84] \\
    8 & 10 & 13.95 & [1.17, 93.71] \\
    9 & 10 & 3.57 & [2.95, 4.72] \\
    10 & 10 & 0.5 & [0.34, 1.23] \\
    11 & 10 & 1.35 & [1.08, 2.05] \\
    12 & 10 & 3.41 & [2.65, 5.89] \\
    13 & 10 & 1.14 & [0.37, 7.72] \\
    14 & 10 & 3.64 & [1.23, 8.98] \\
    15 & 10 & 12.94 & [3.05, 30.12] \\
    \bottomrule
\end{longtable}%
}

{
\renewcommand{\arraystretch}{0.9}
\small
\begin{longtable}{ccccc}
\caption{CPU time of SAA formulation with sample size = 500 and $L=6$. The distribution of demand and service duration are both uniform distribution.}
\label{CPU_500_2_U} \\
\hline
\textbf{K} & \textbf{T} & \textbf{CPU time mean (s)} & \textbf{CPU time range (s)}& $\left| \textbf{AOI} \right|$\\
\endfirsthead
\multicolumn{5}{r}%
{{\bfseries -- continued from previous page}}\\
	\hline 
	\endhead
	\hline
	\endfoot
	\hline
	\endlastfoot
	\hline
    \multirow{3}[0]{*}{6} & 30 & 88.04 & [83.63, 91.42] & 0.0005\\
      & 90 & 671.75 & [650.71, 688.06]& 0.0005 \\
      & 180 & 588.94 & [310.25, 1189.11] &0.0002\\
    \hline
    \multirow{3}[0]{*}{8} & 30 & 579.79 & [287.39, 830.13] & 0\\
      & 90 & 733.06 & [599.18, 1170.73]&0\\
      & 180 & 3157.05 & [2001.41, 4666.1]&0 \\
    \hline
\end{longtable}%
}

\newpage

\section{Results of optimal staffing patterns for EA models}\label{OSP_A}
\setcounter{figure}{0}
\setcounter{table}{0}
\begin{figure}[!ht]
    \centering
    \subcaptionbox{$c_{l,t}^o = 1$\label{EA-10-4-8-30}}{
        \includegraphics[scale=0.35]{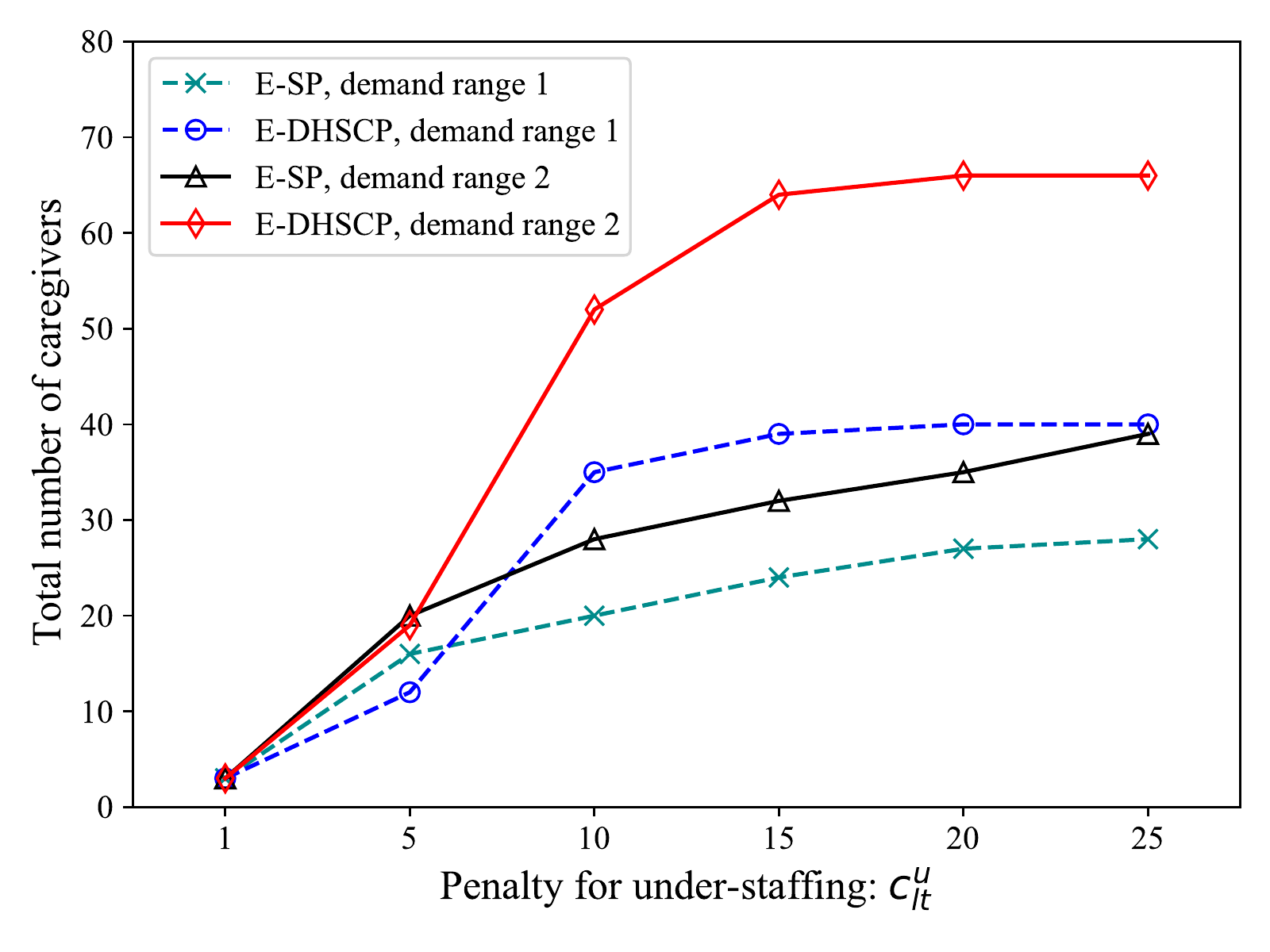}}
    \subcaptionbox{$c_{l,t}^u = 10$\label{EA-10-4-8-30-o}}{
        \includegraphics[scale=0.35]{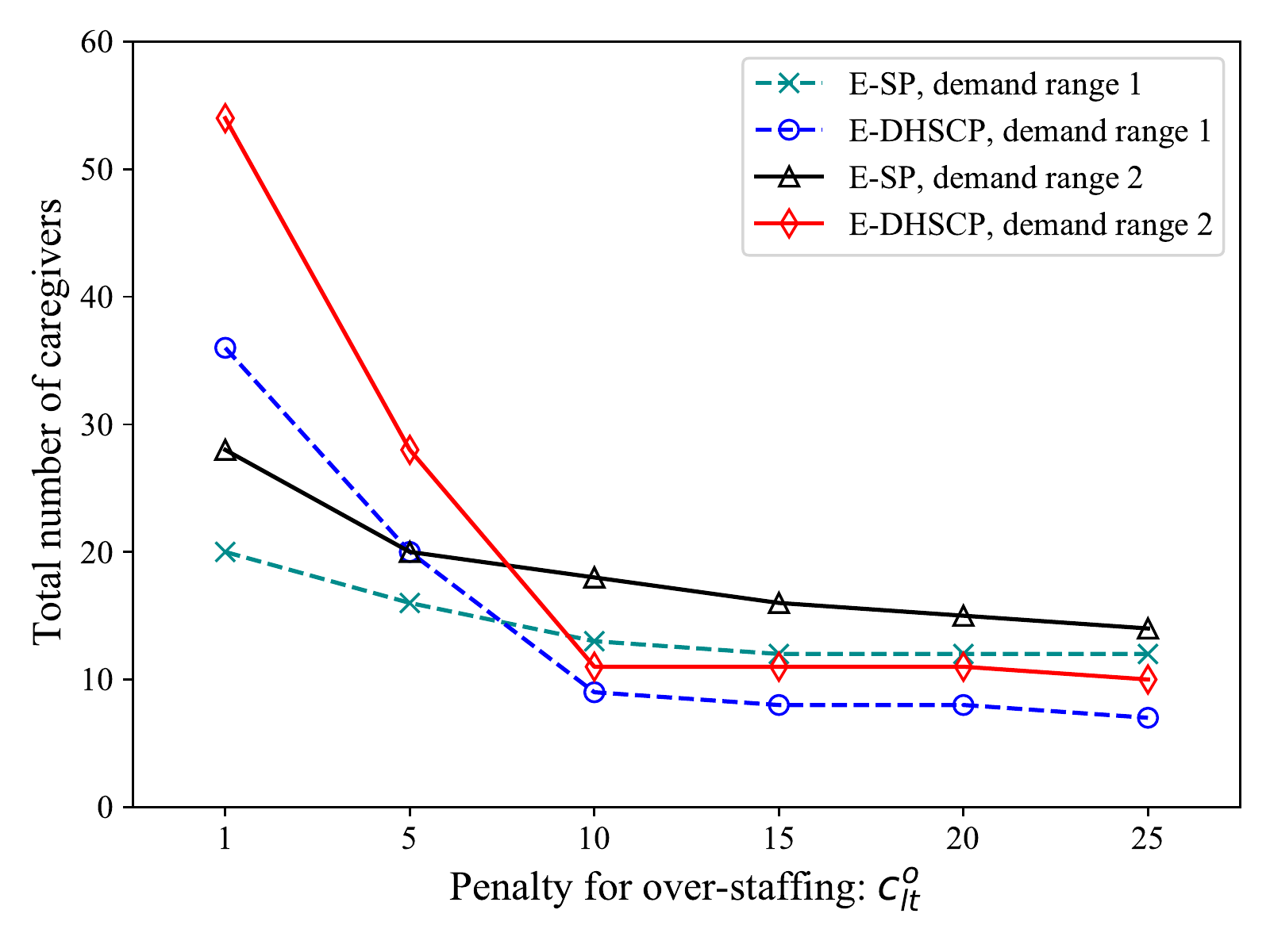}}
    \caption{Total number of caregivers hired by E-SP and E-DHSCP models for Instance 7.}
    \label{total_num_fig_7}
\end{figure}

\begin{figure}[!ht]
    \centering
    \subcaptionbox{Instance 10, demand range 1\label{EA-10-2-6-6-30-4060}}{
        \includegraphics[scale=0.35]{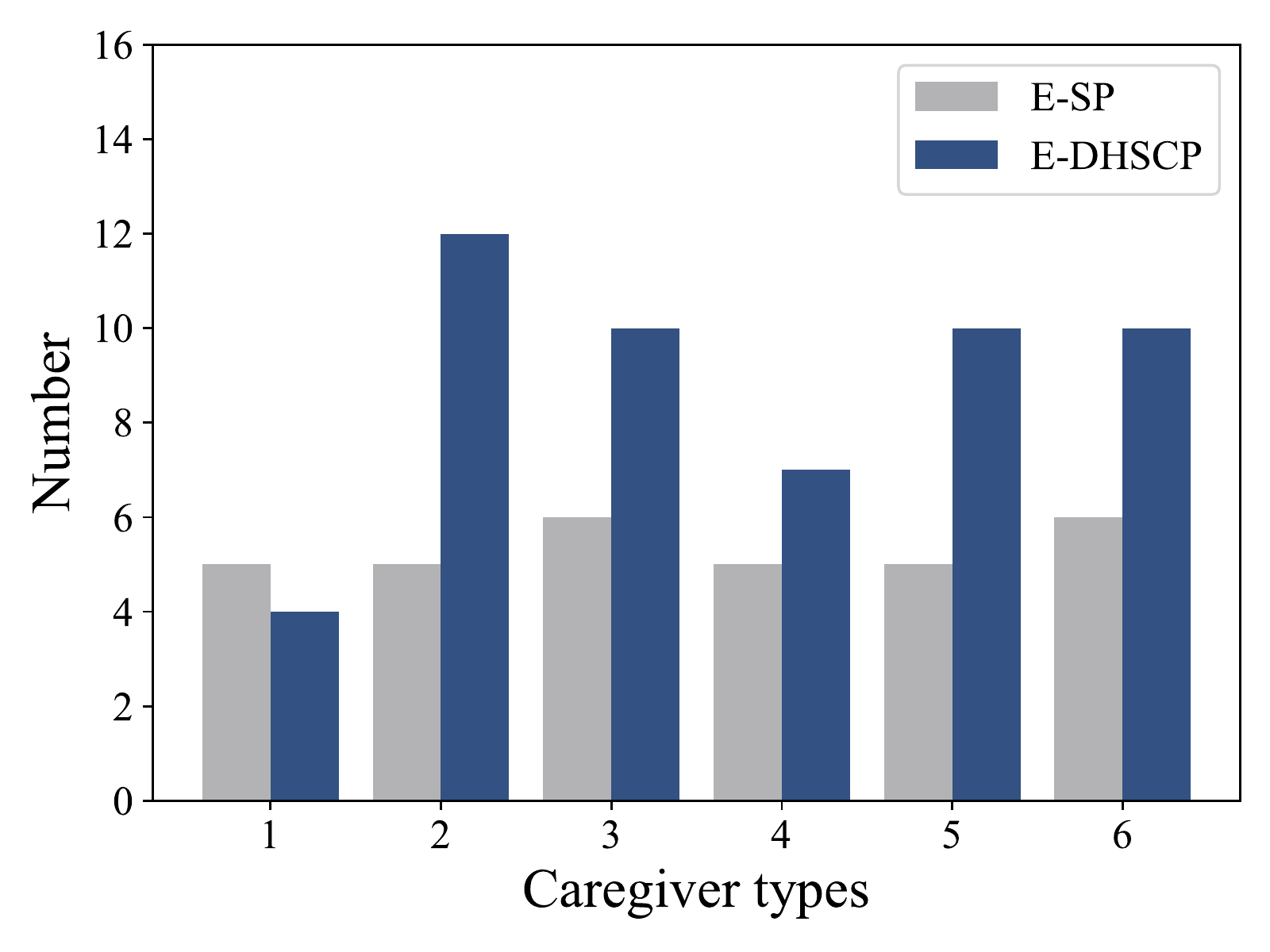}}
    \subcaptionbox{Instance 10, demand range 2\label{EA-10-2-6-6-30-1080}}{
        \includegraphics[scale=0.35]{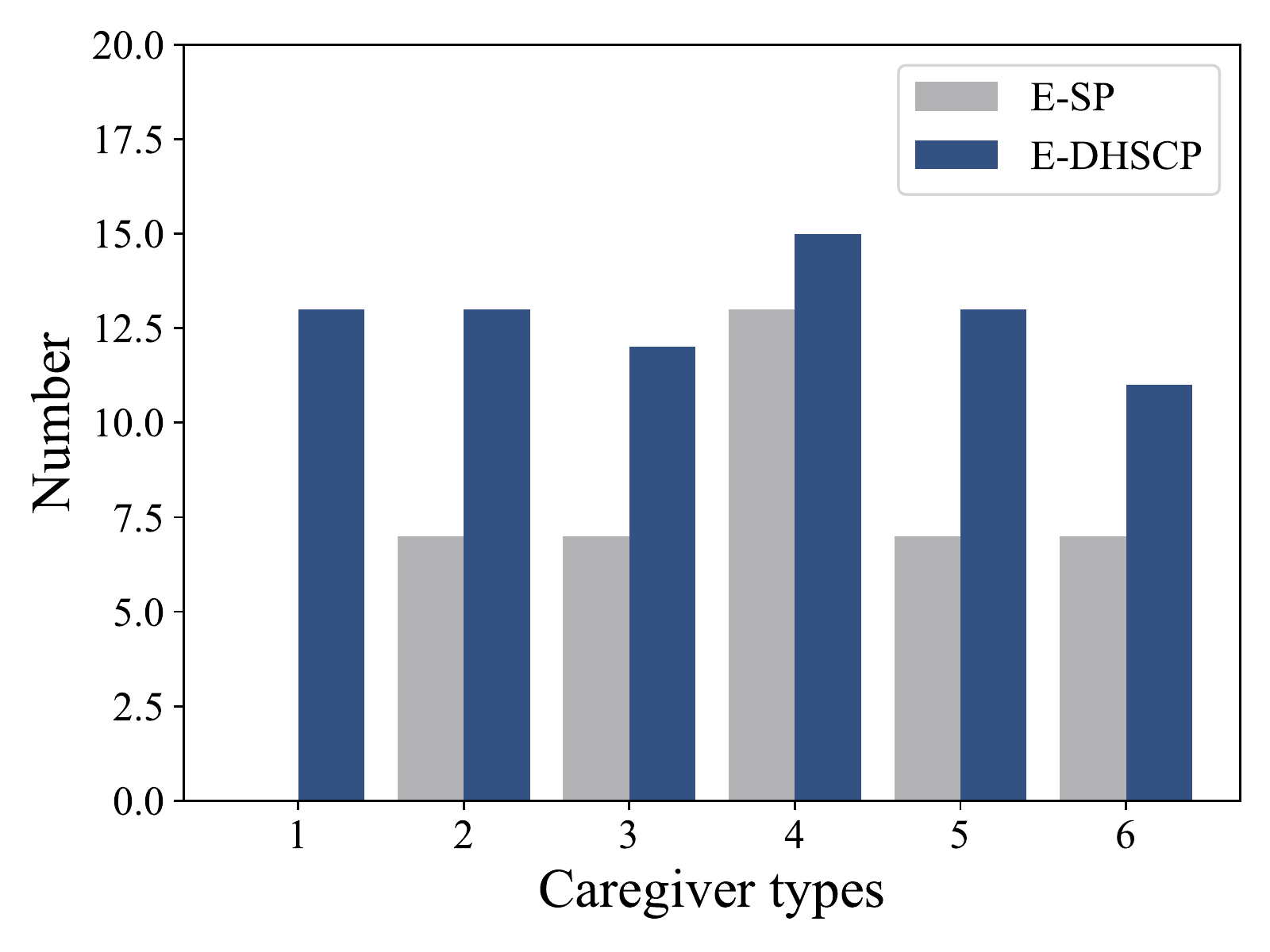}}
    \subcaptionbox{Instance 7, demand range 1\label{EA-10-2-4-8-30-4060}}{
        \includegraphics[scale=0.35]{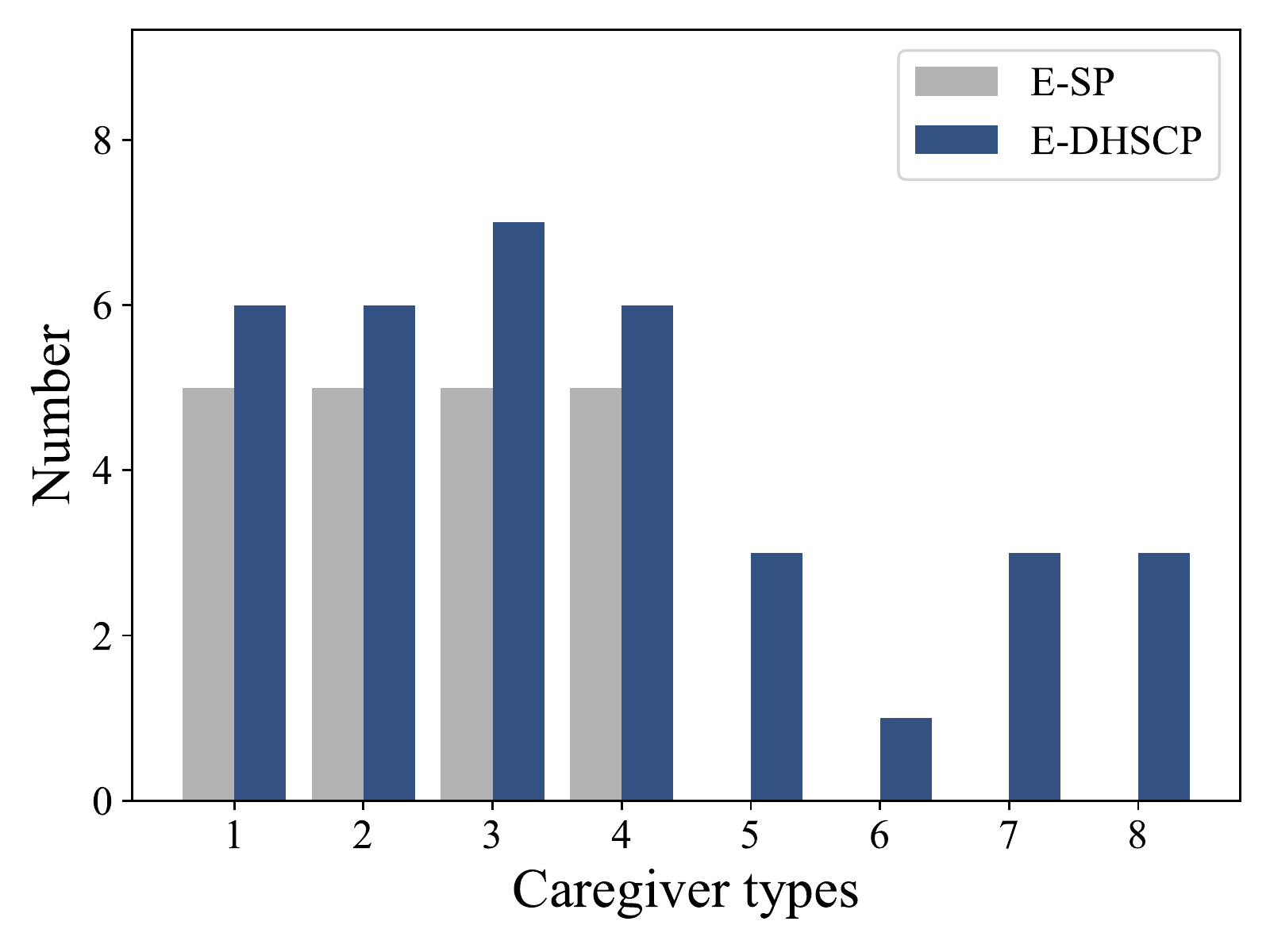}}
    \subcaptionbox{Instance 7, demand range 2\label{EA-10-2-4-8-30-1080}}{
        \includegraphics[scale=0.35]{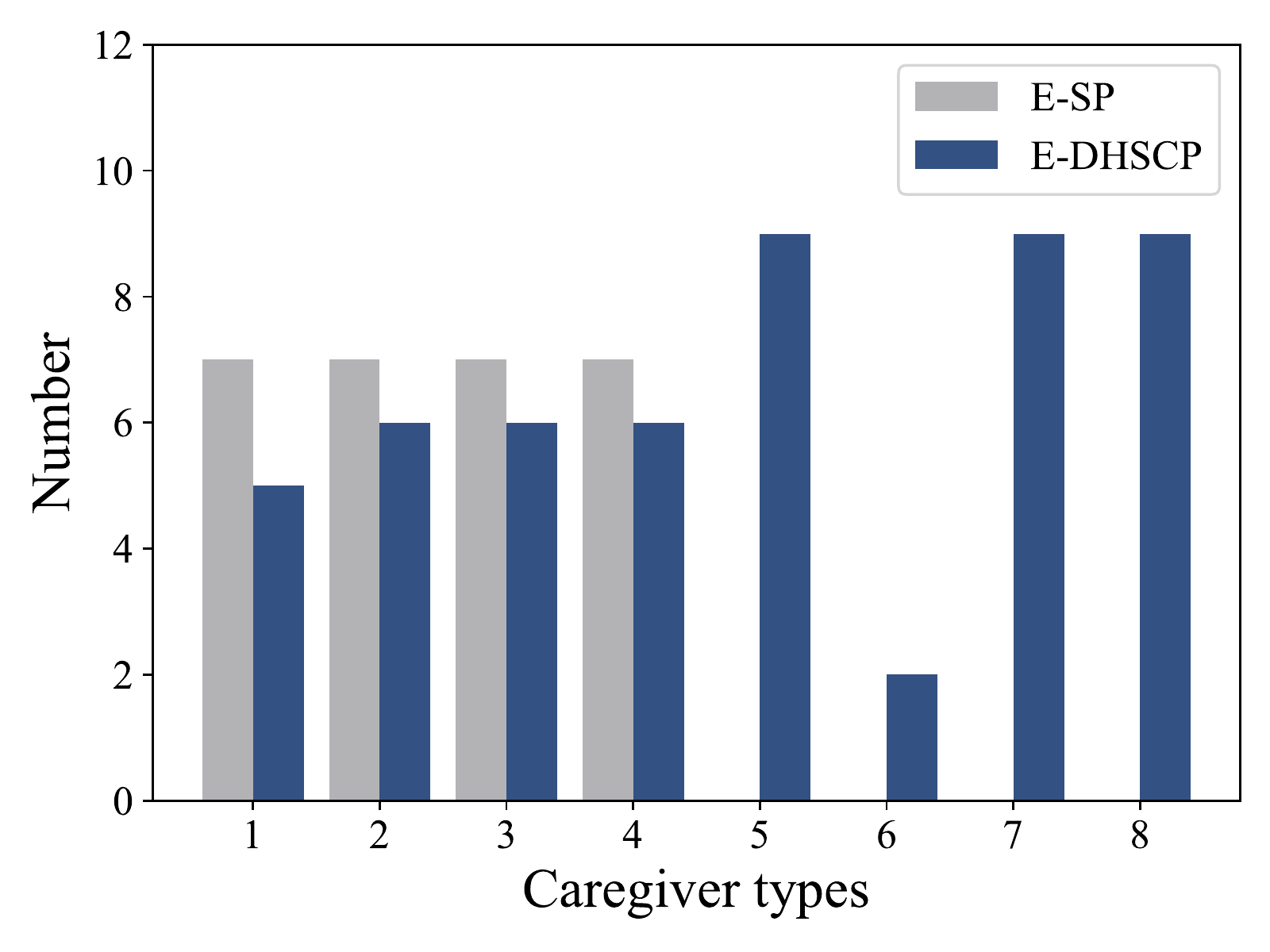}}
    \caption{Staffing patterns of EA decision-makers under cost structure $(c_{l,t}^{u}, c_{l,t}^{o}) = (10, 1)$.}
    \label{EA-solution}
\end{figure}
\clearpage
\newpage

\section{Out-of-sample performance of EA models}\label{oop-EA}
\subsection{Out-of-sample performance of EA models for Instance 10}\label{oop-EA-10}
\setcounter{figure}{0}
\setcounter{table}{0}
\begin{figure}[!ht]
    \centering
    \subcaptionbox{\label{EA-inc-10-2-6-6-30-4060}}{
        \includegraphics[scale=0.31]{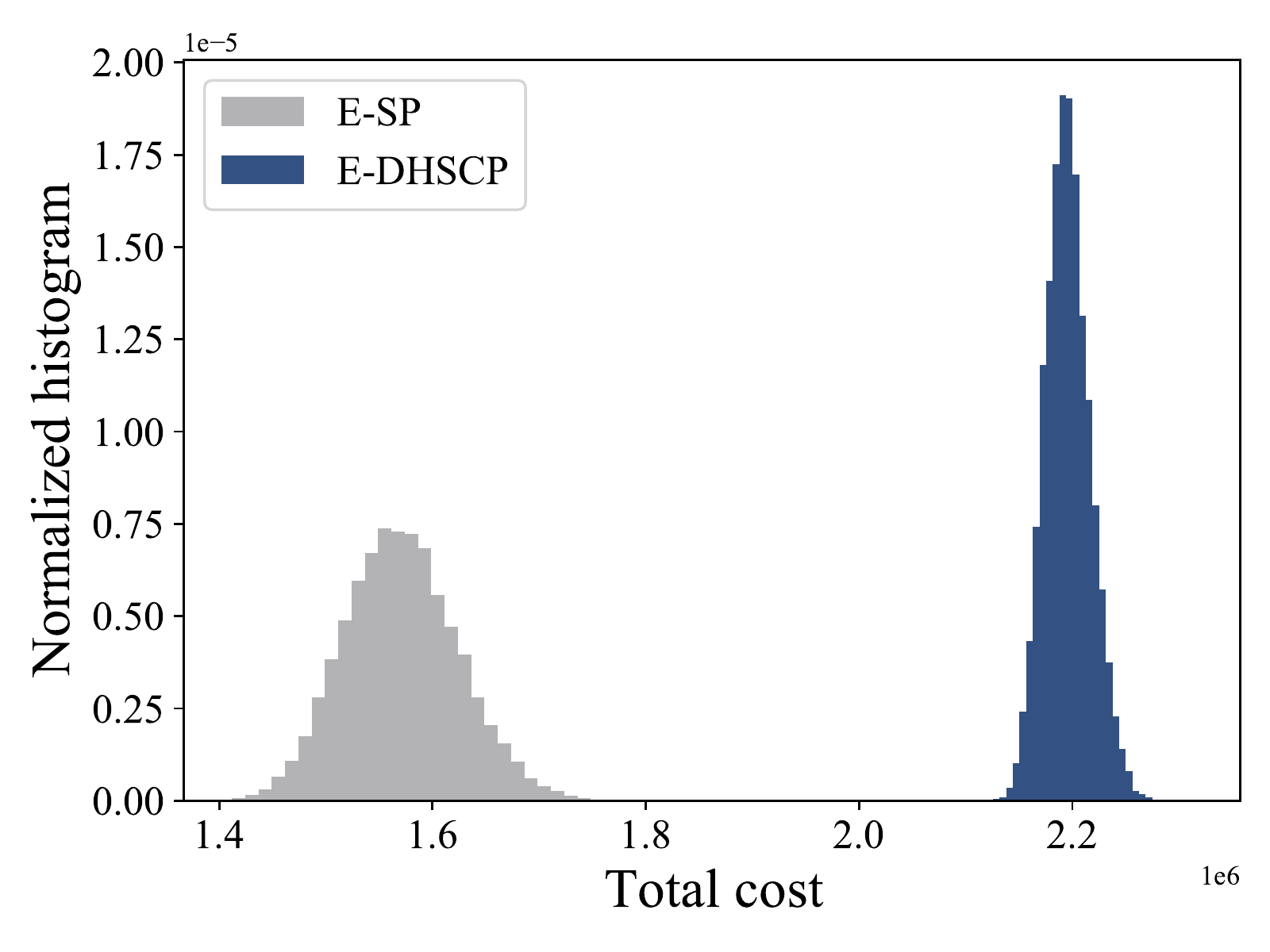}}
    \subcaptionbox{\label{EA-inc2-10-2-6-6-30-4060}}{
        \includegraphics[scale=0.31]{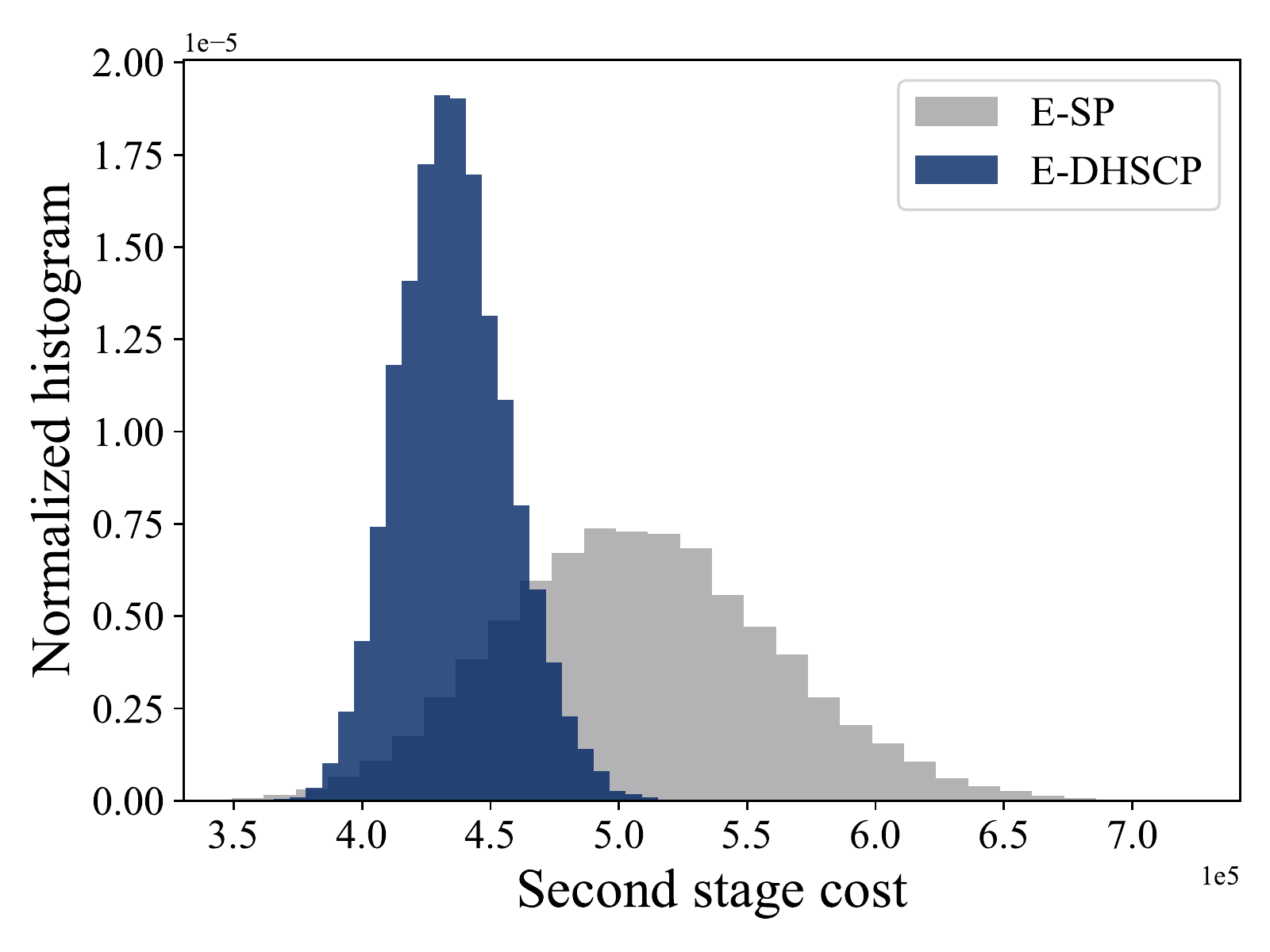}}
    \subcaptionbox{\label{EA-inu-10-2-6-6-30-4060}}{
        \includegraphics[scale=0.31]{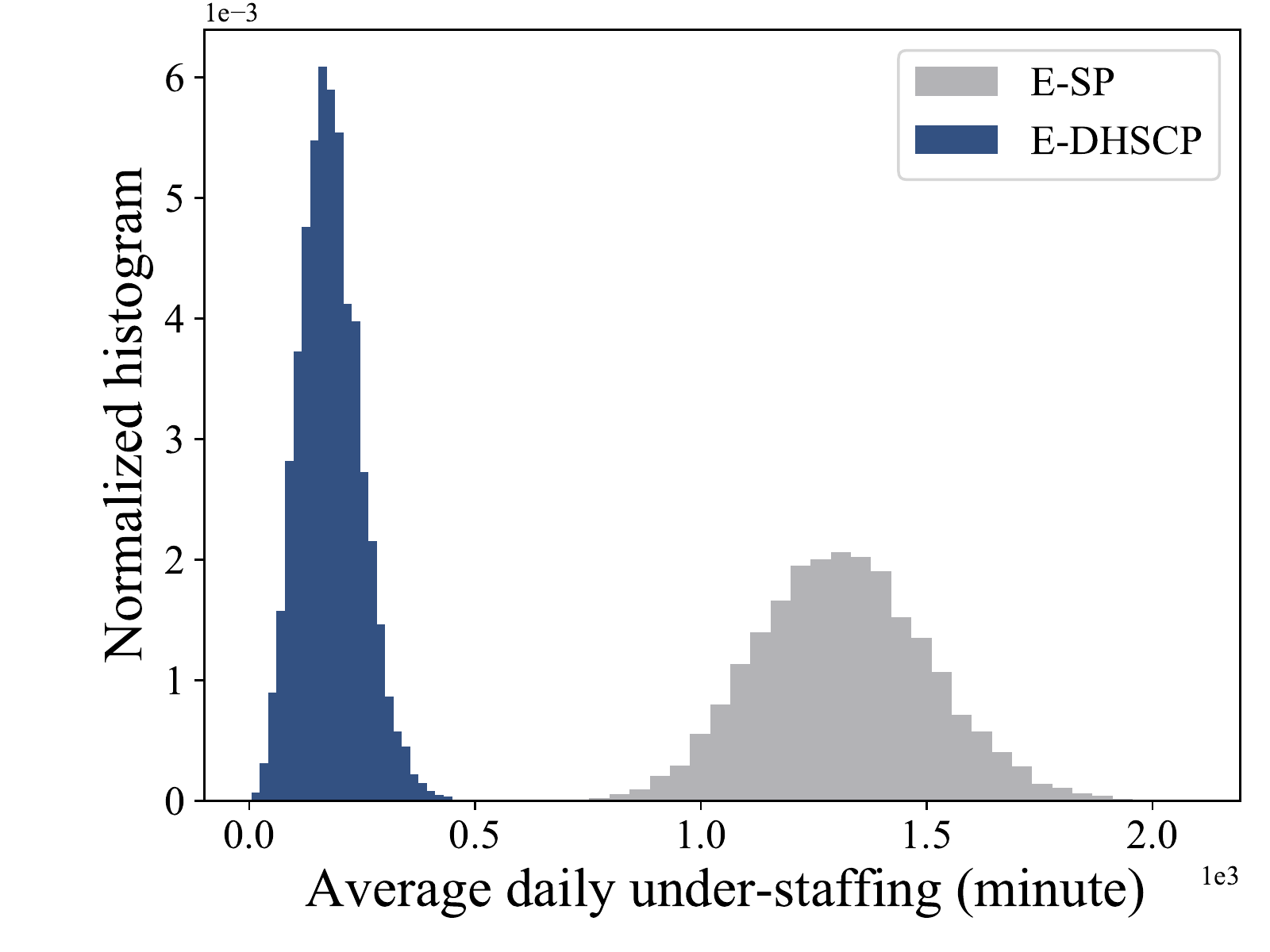}}
    \caption{Out-of-sample performance for Instance 10, demand range 1 under Set 1}
    \label{EA-oop-in-6}
\end{figure}

\begin{figure}[!ht]
    \centering
    \subcaptionbox{$\Delta = 0$\label{EA-outc-10-2-6-6-30-1080}}{
        \includegraphics[scale=0.32]{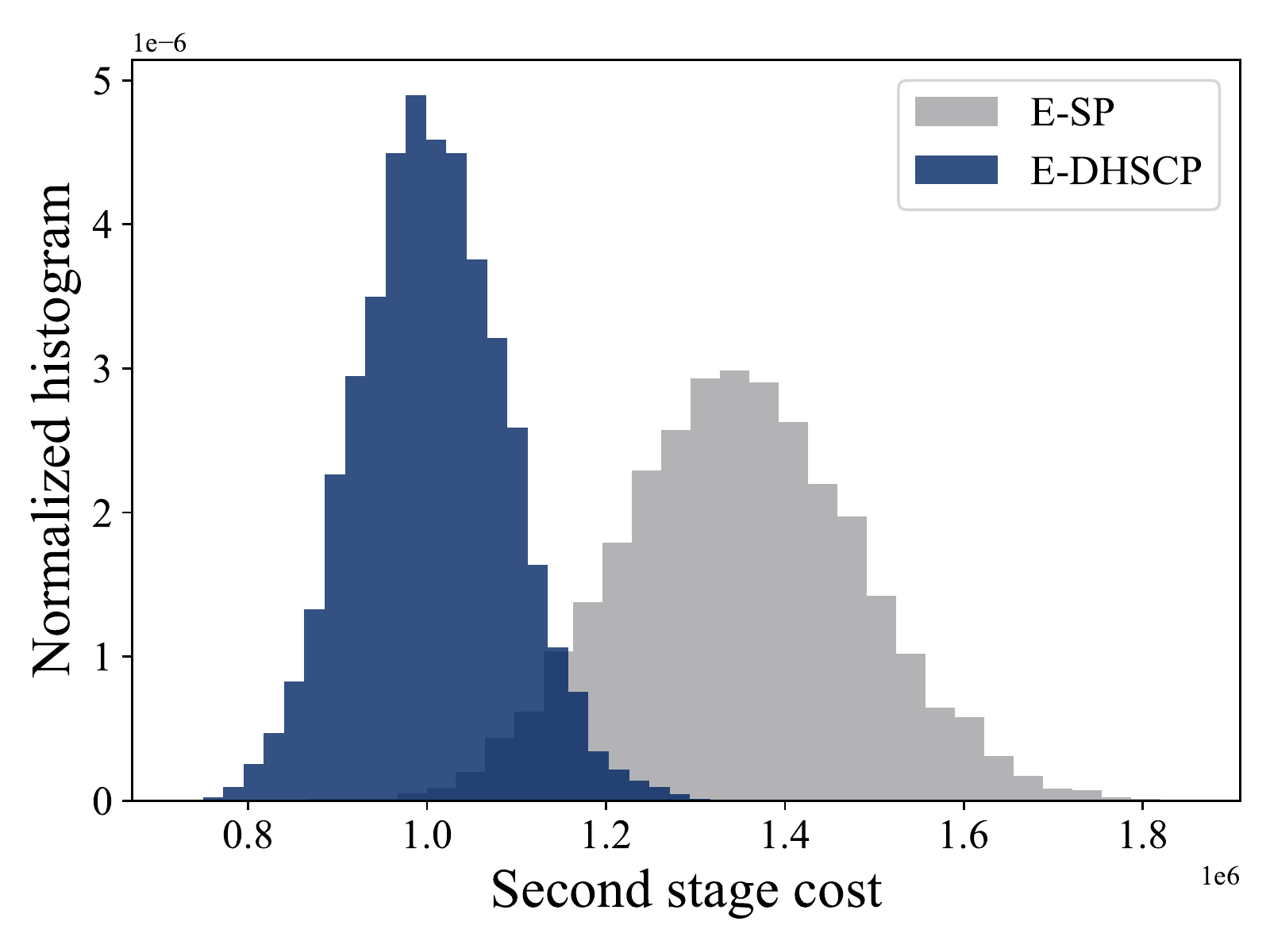}}
    \subcaptionbox{$\Delta = 0.1$\label{EA-out01c2-10-2-6-6-30-1080}}{
        \includegraphics[scale=0.32]{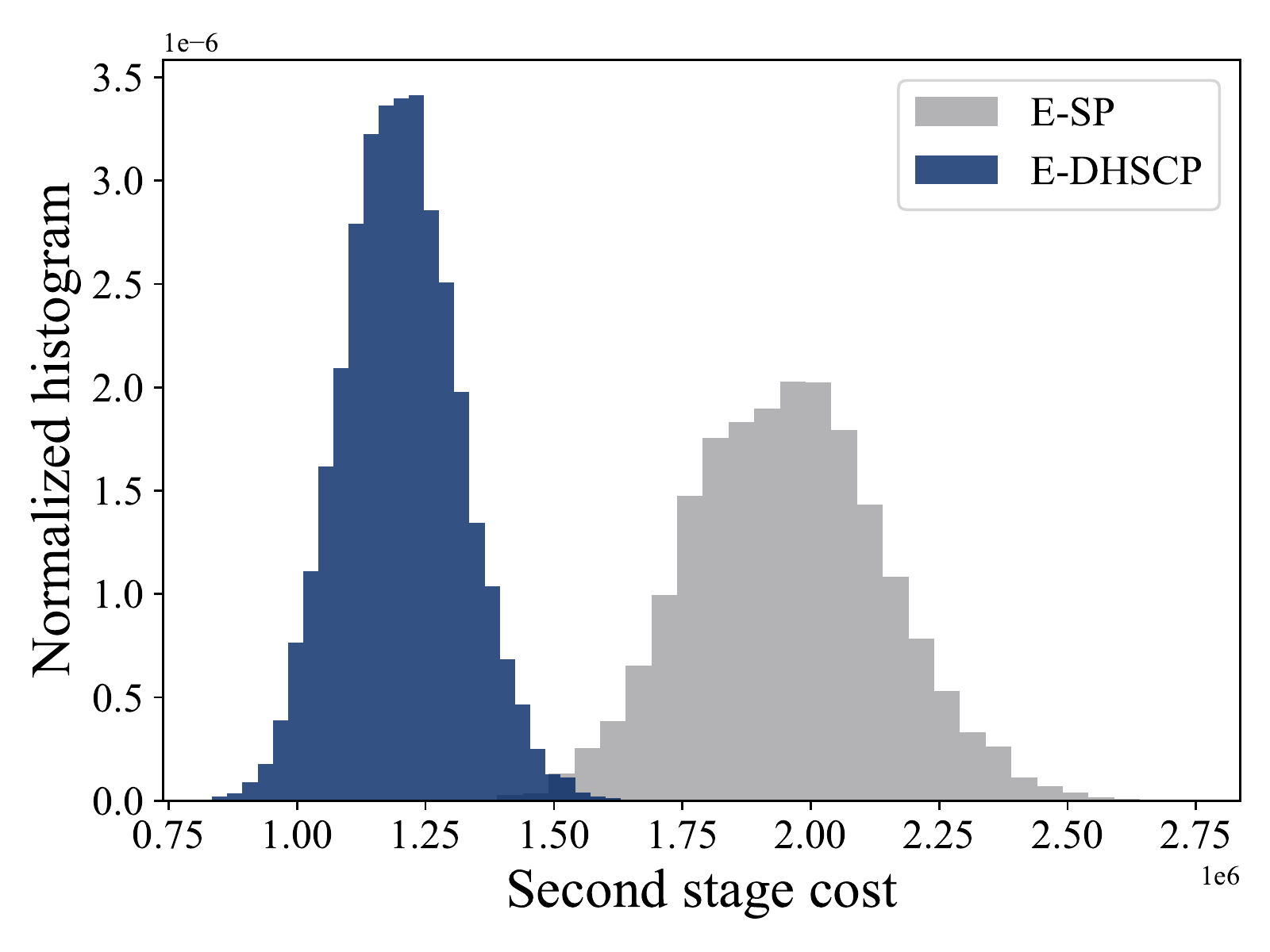}}
    \subcaptionbox{$\Delta = 0.25$\label{EA-out025c2-10-2-6-6-30-1080}}{
        \includegraphics[scale=0.32]{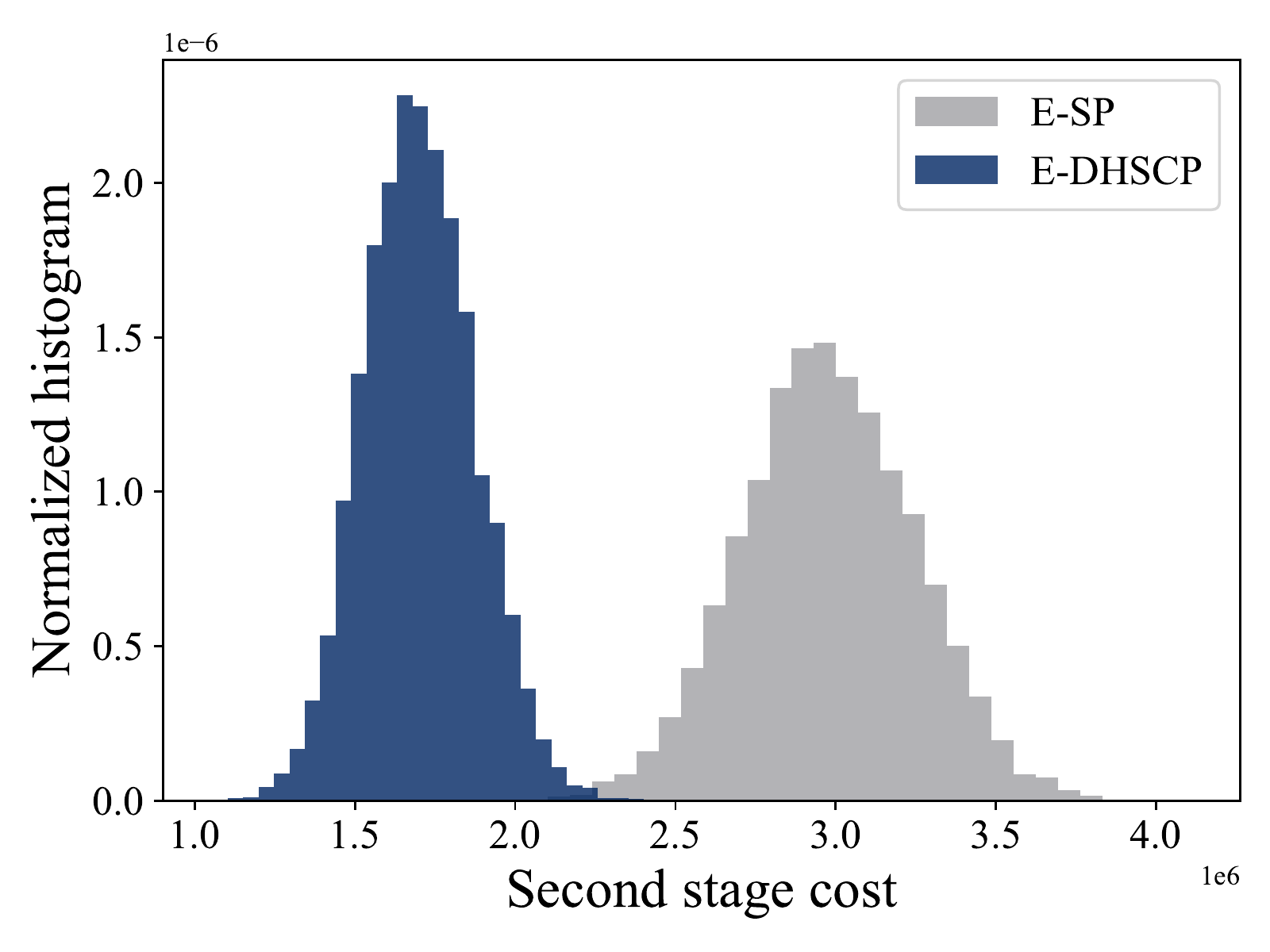}}
    \subcaptionbox{$\Delta = 0.5$\label{EA-out05c2-10-2-6-6-30-1080}}{
        \includegraphics[scale=0.32]{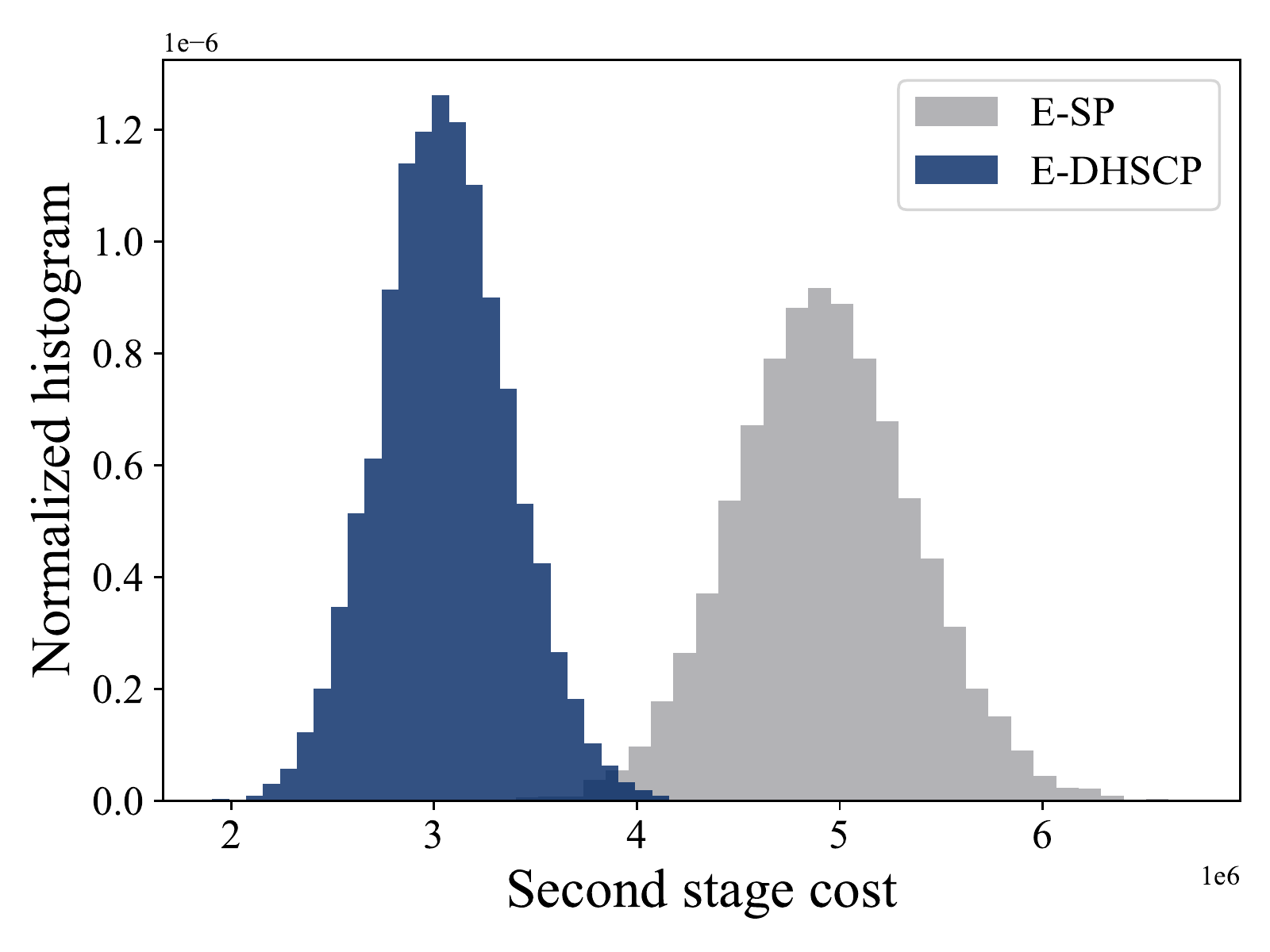}}
    \caption{Out-of-sample second stage cost for Instance 10, demand range 2 under Set 2}
    \label{EA-oop-6-1080}
\end{figure}

\begin{figure}[!ht]
    \centering
    \subcaptionbox{$\Delta = 0$\label{EA-outu-10-2-6-6-30-4060}}{
        \includegraphics[scale=0.35]{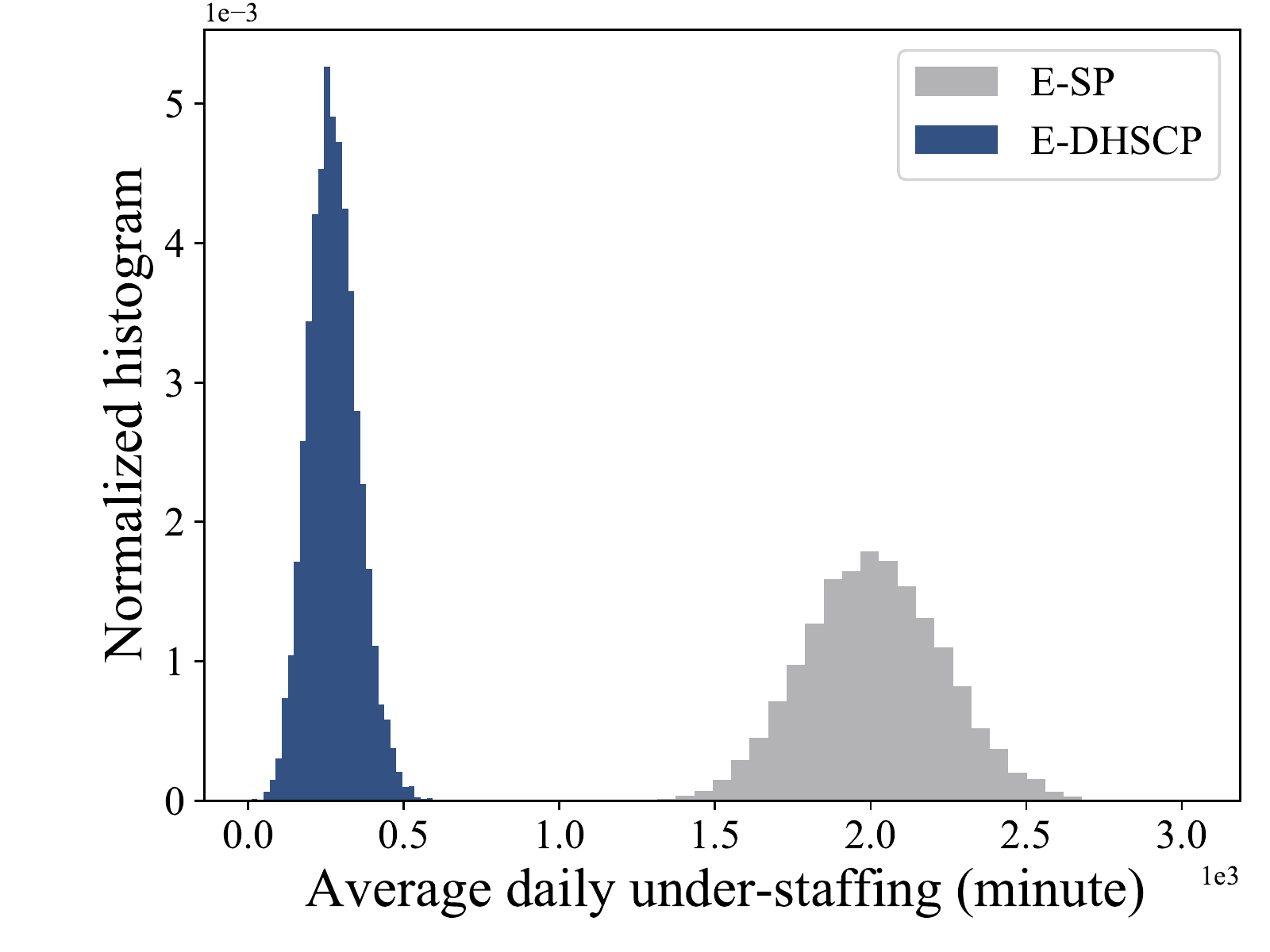}}
    \subcaptionbox{$\Delta = 0.1$\label{EA-out01u-10-2-6-6-30-4060}}{
        \includegraphics[scale=0.35]{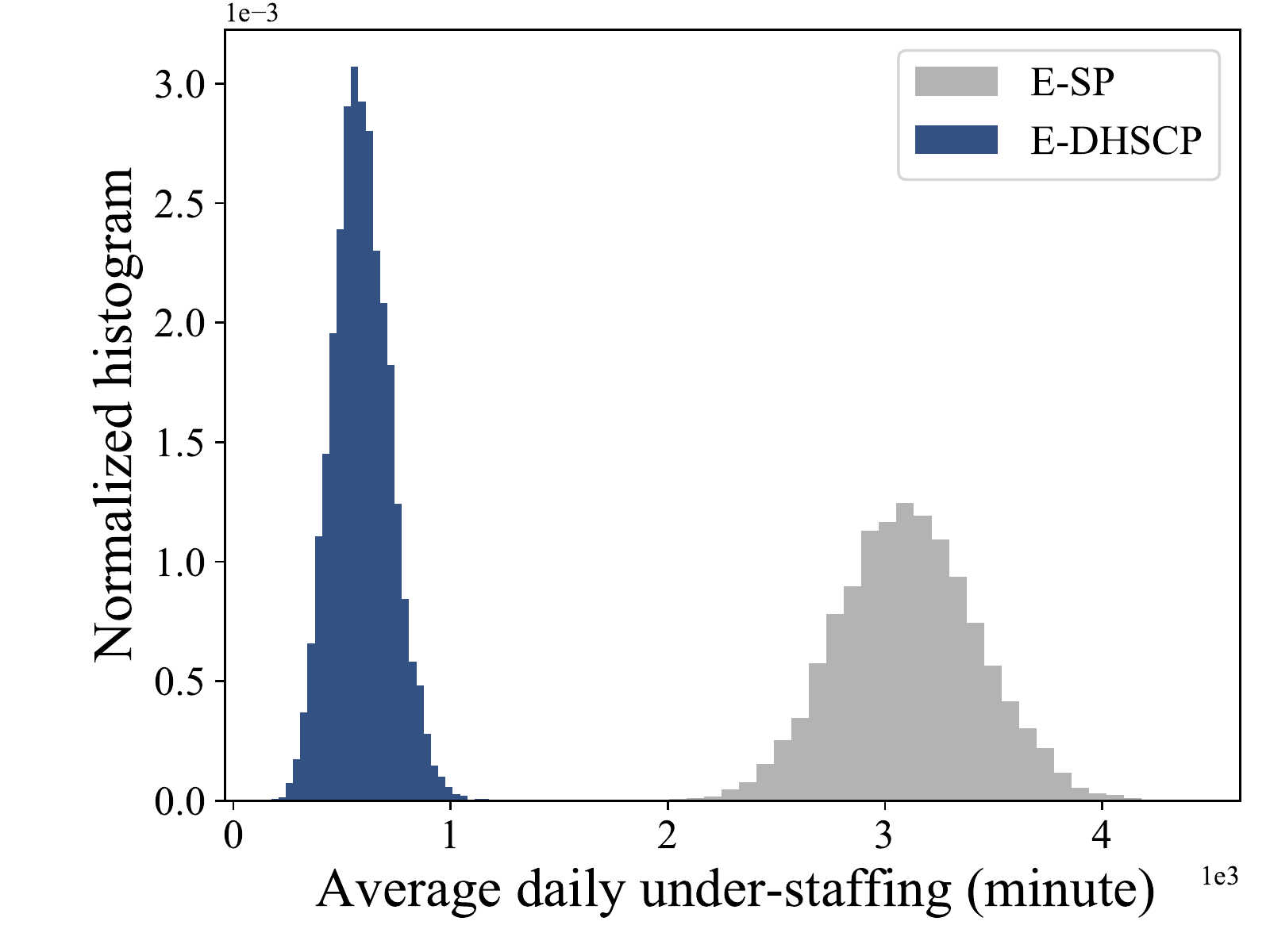}}
    \subcaptionbox{$\Delta = 0.25$\label{EA-out025u-10-2-6-6-30-4060}}{
        \includegraphics[scale=0.35]{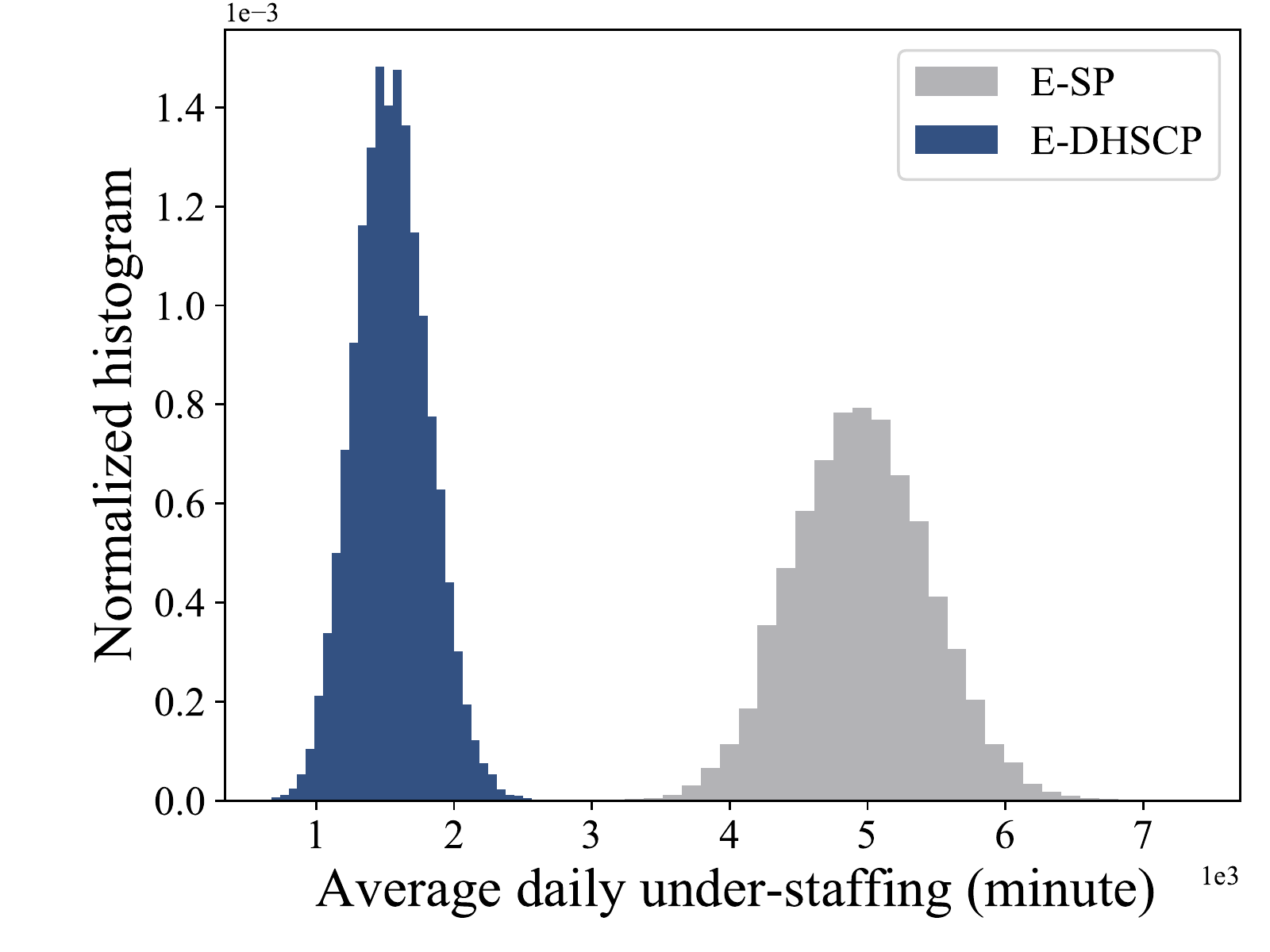}}
    \subcaptionbox{$\Delta = 0.5$\label{EA-out05u-10-2-6-6-30-4060}}{
        \includegraphics[scale=0.35]{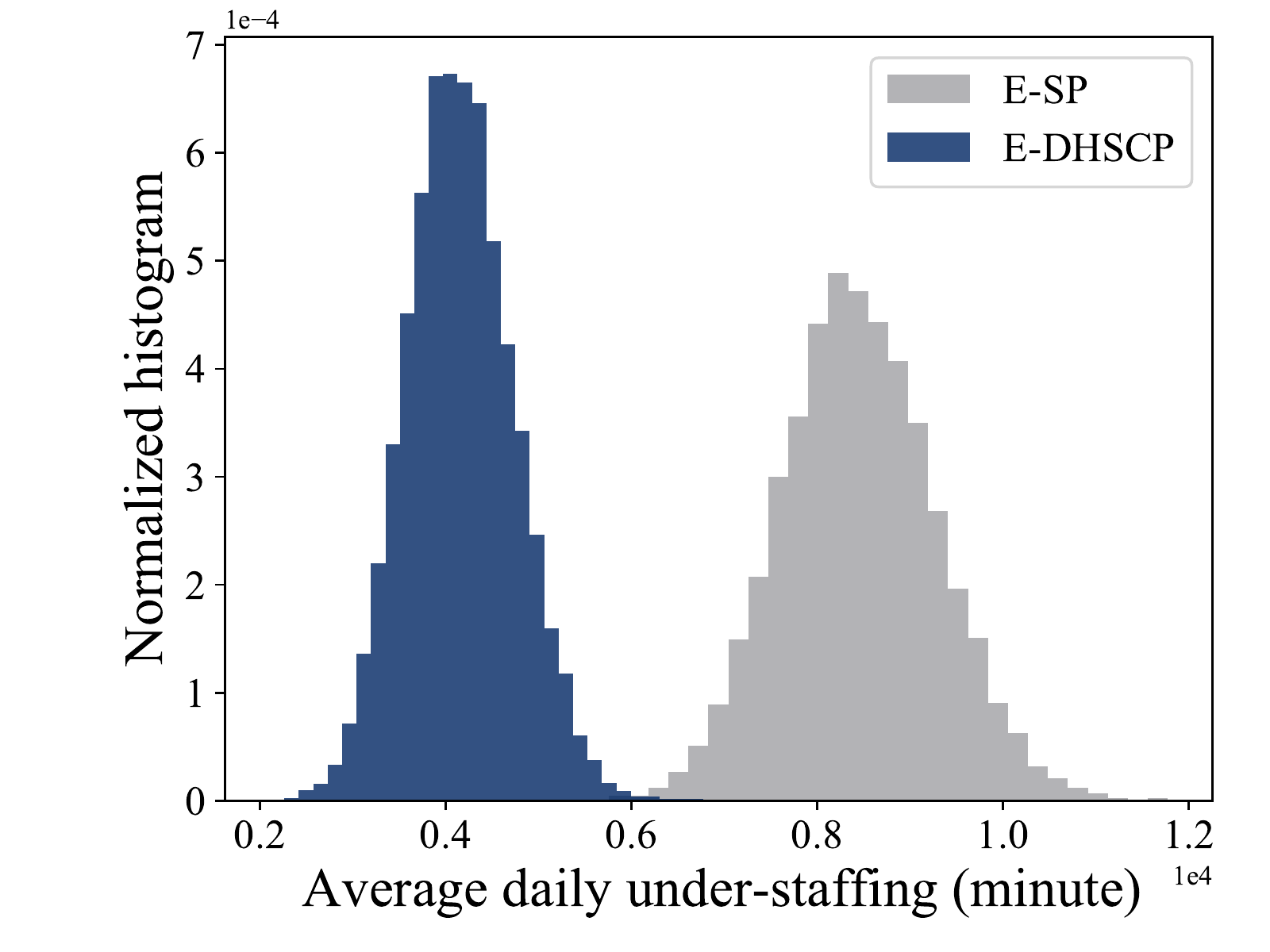}}
    \caption{Out-of-sample under-staffing for Instance 10, demand range 1 under Set 2}
    \label{EA-oop-u-6-4060}
\end{figure}

\newpage
\subsection{Out-of-sample performance of EA models for Instance 7}\label{oop-EA-7}
\begin{figure}[!htbp]
    \centering
    \subcaptionbox{Total cost\label{EA-inc-10-2-4-8-30-4060}}{
        \includegraphics[scale=0.31]{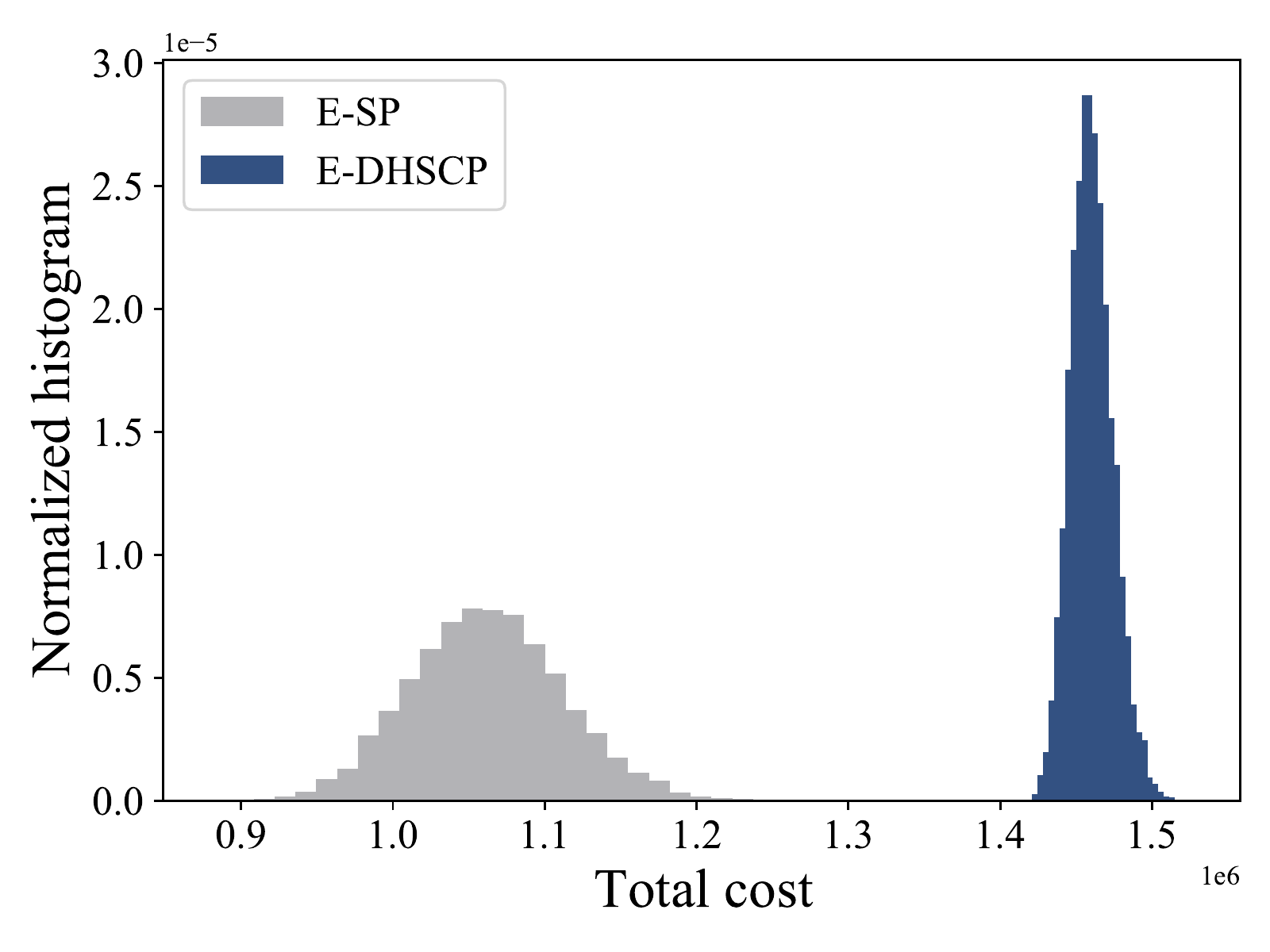}}
    \subcaptionbox{Second stage cost\label{EA-inc2-10-2-4-8-30-4060}}{
        \includegraphics[scale=0.31]{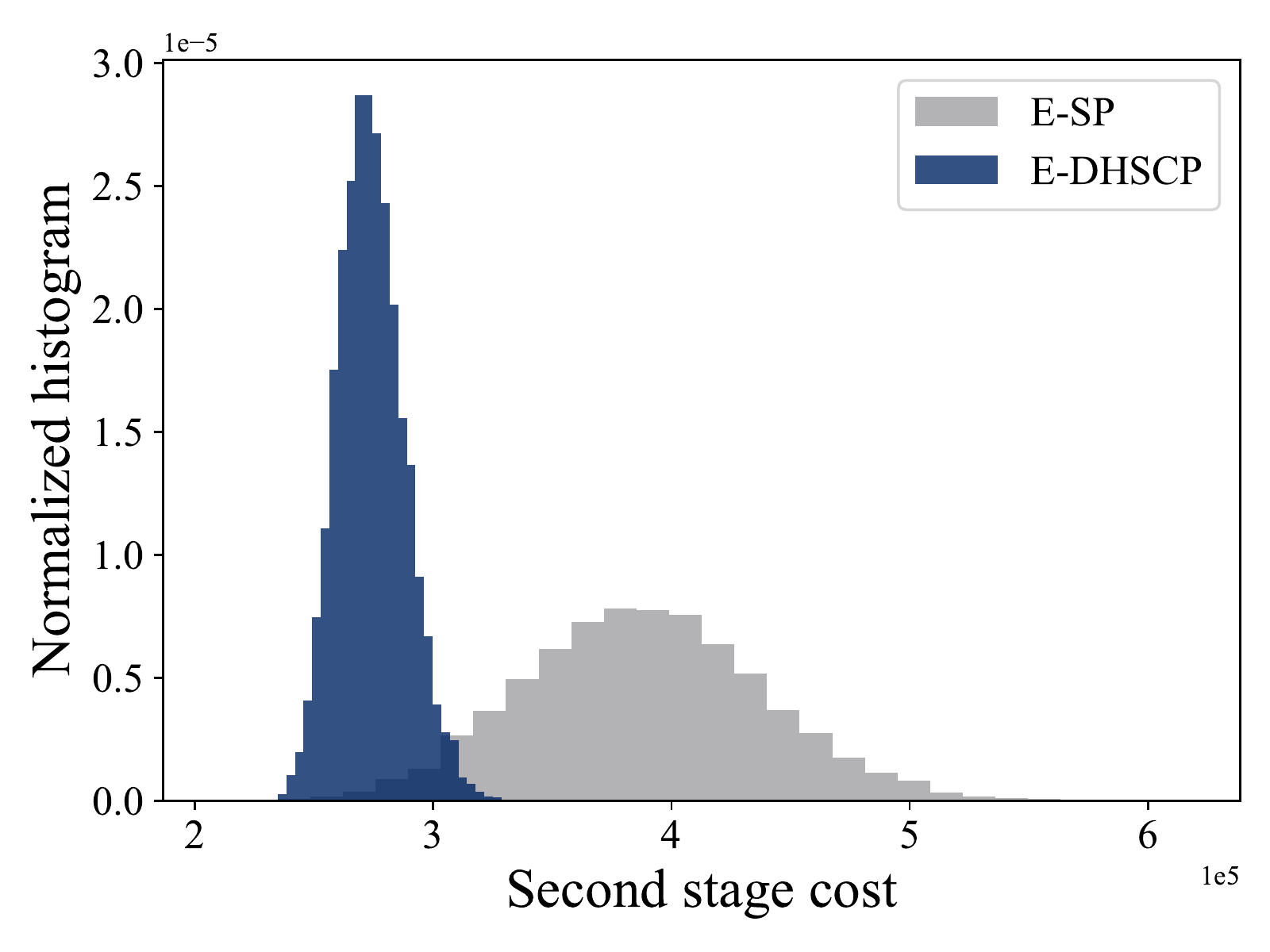}}
    \subcaptionbox{Average daily under-staffing\label{EA-inu-10-2-4-8-30-4060}}{
        \includegraphics[scale=0.31]{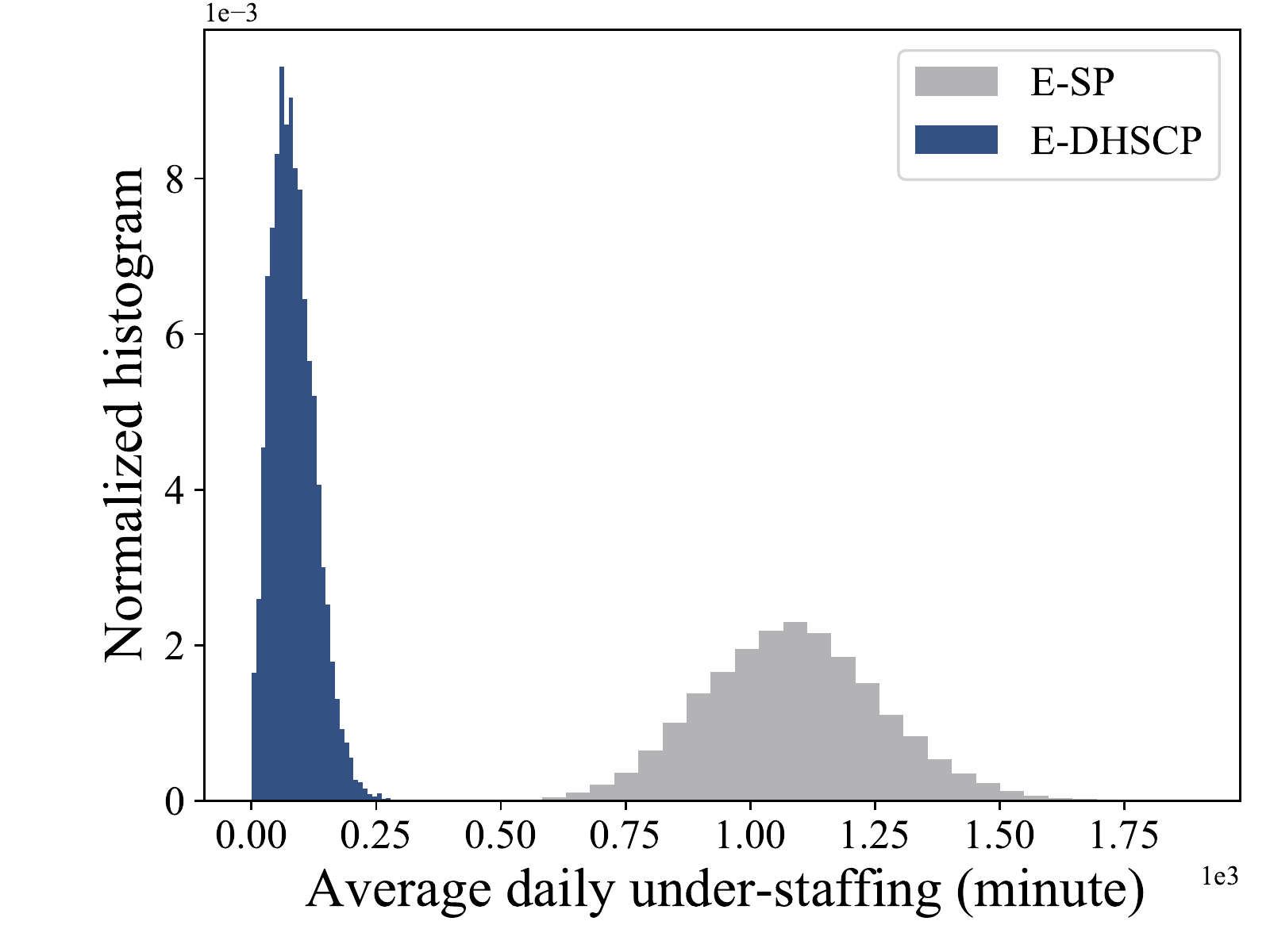}}
    \caption{Out-of-sample performance for Instance 7, demand range 1 under Set 1}
    \label{EA-oop-in-4}
\end{figure}

\begin{figure}[!htbp]
    \centering
    \subcaptionbox{$\Delta = 0$\label{EA-outc2-10-2-4-8-30-4060}}{
        \includegraphics[scale=0.35]{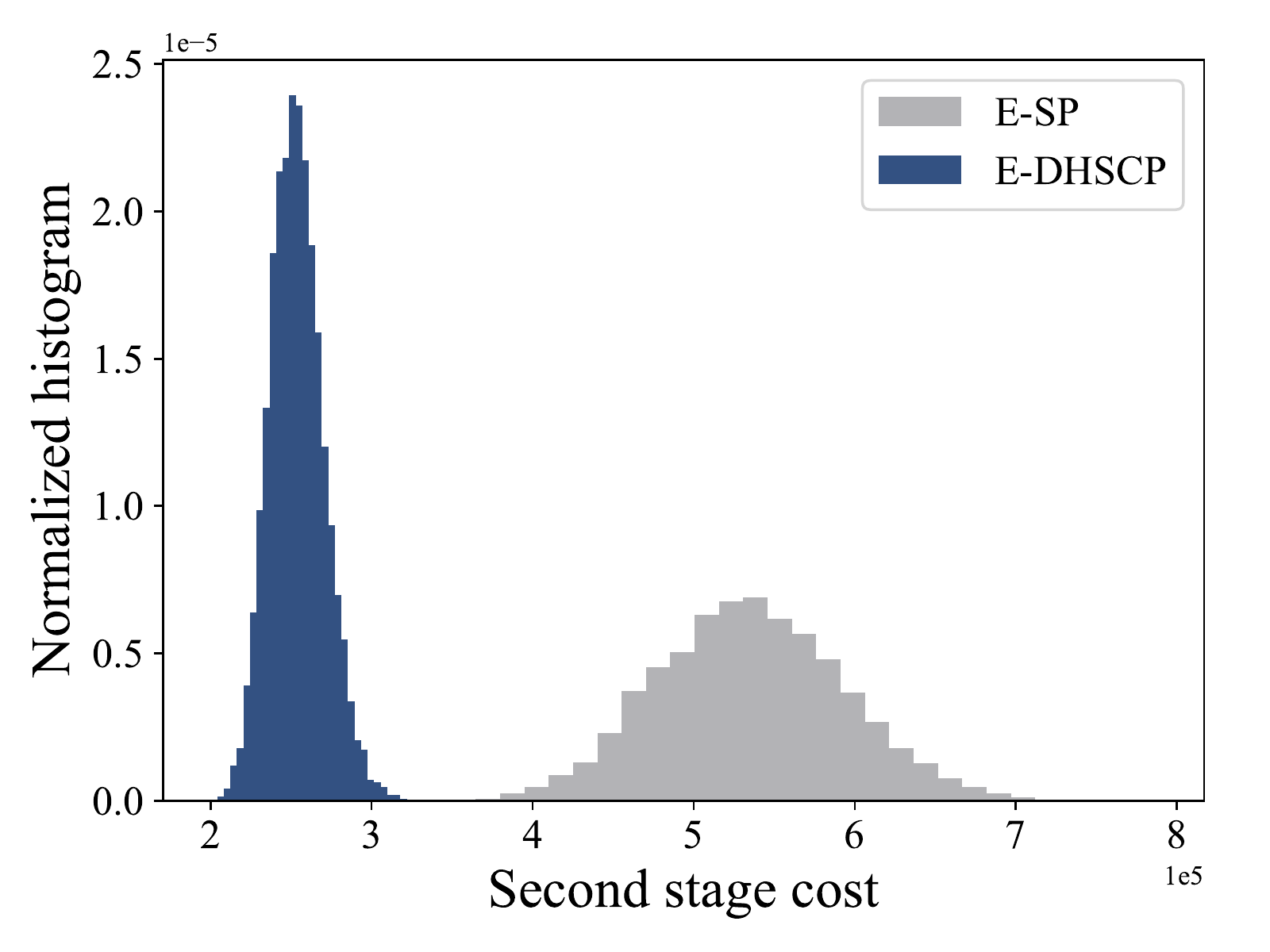}}
    \subcaptionbox{$\Delta = 0.1$\label{EA-out01c2-10-2-4-8-30-4060}}{
        \includegraphics[scale=0.35]{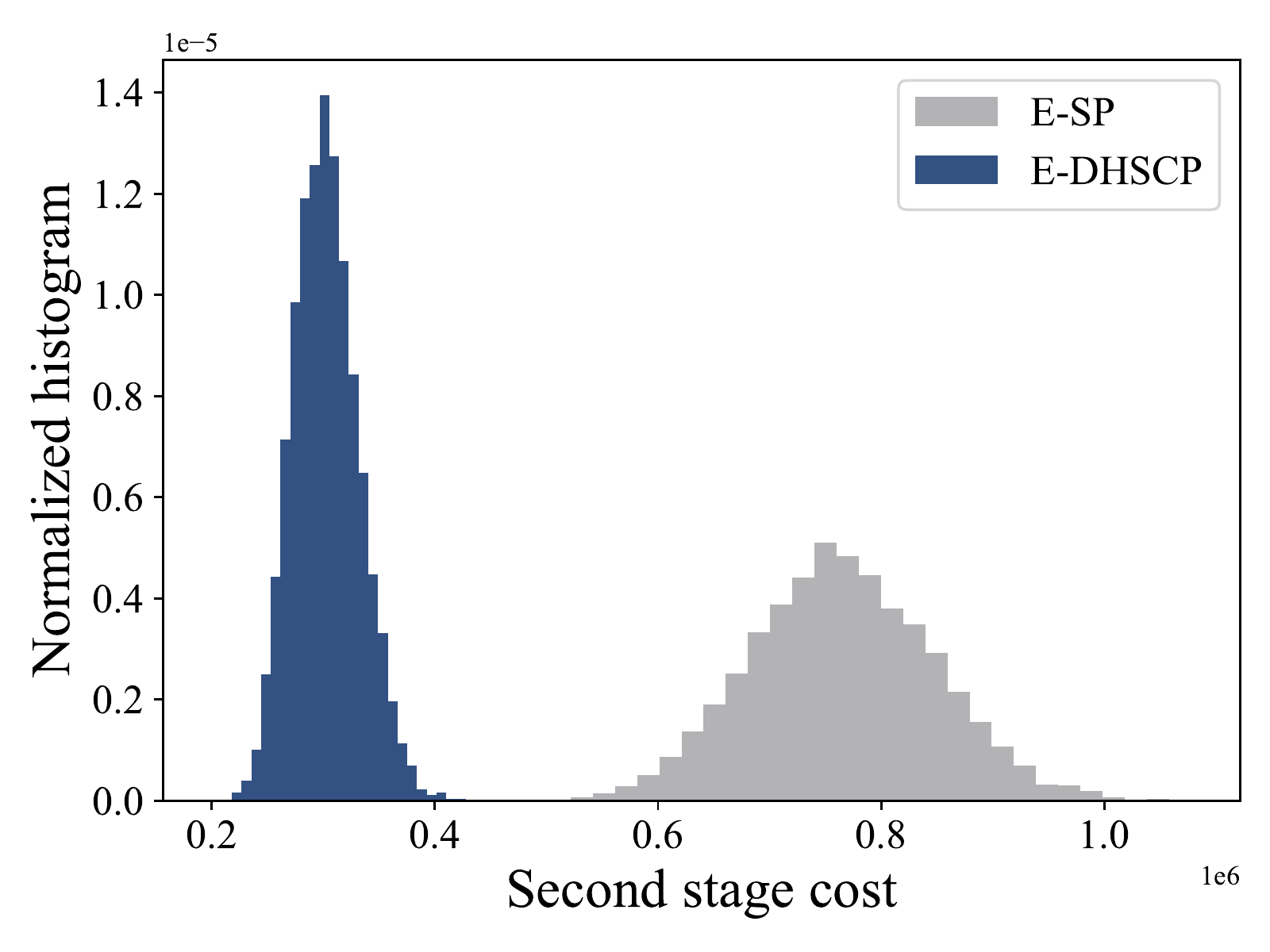}}
    \subcaptionbox{$\Delta = 0.25$\label{EA-out025c2-10-2-4-8-30-4060}}{
        \includegraphics[scale=0.35]{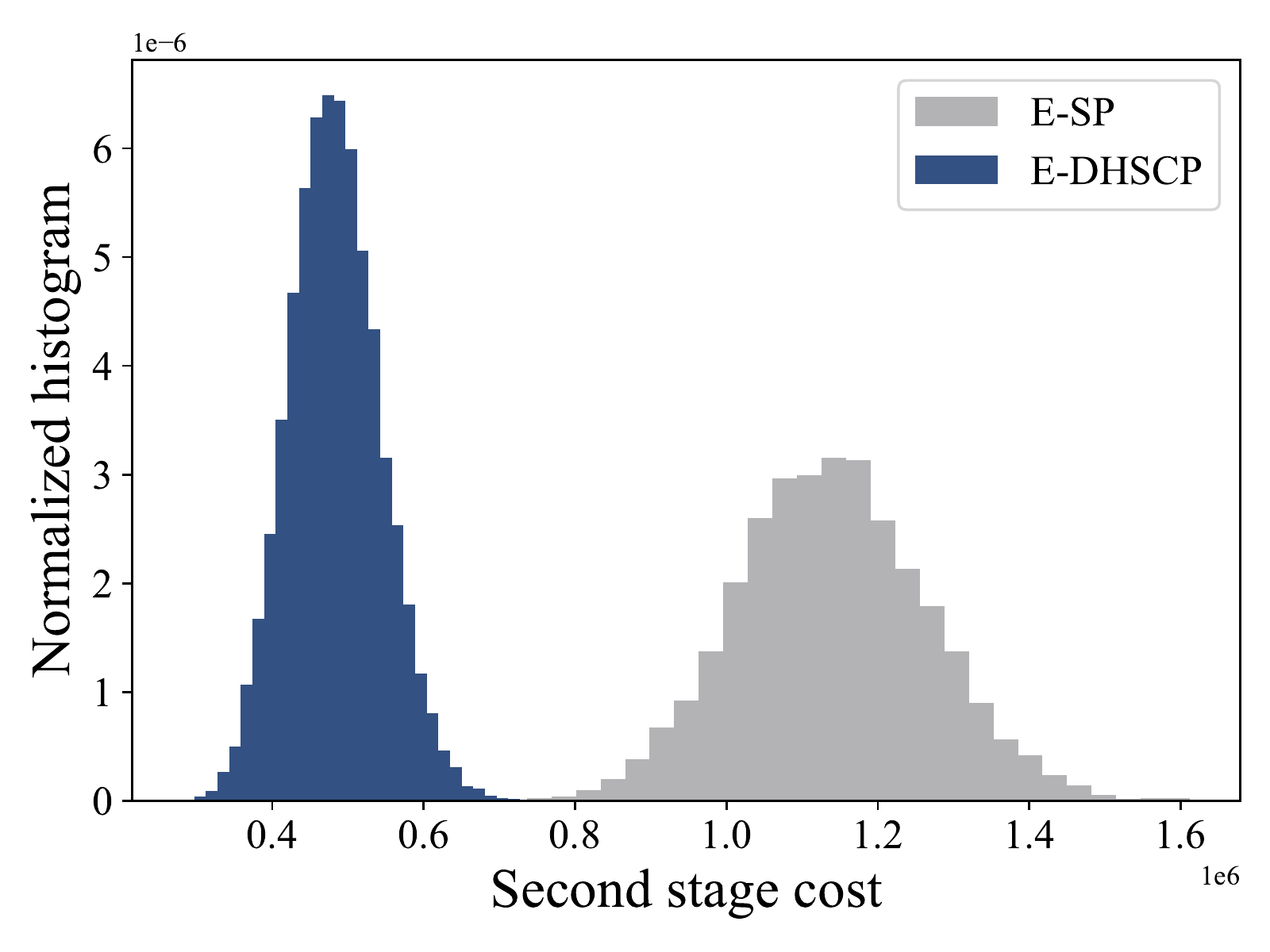}}
    \subcaptionbox{$\Delta = 0.5$\label{EA-out05c2-10-2-4-8-30-4060}}{
        \includegraphics[scale=0.35]{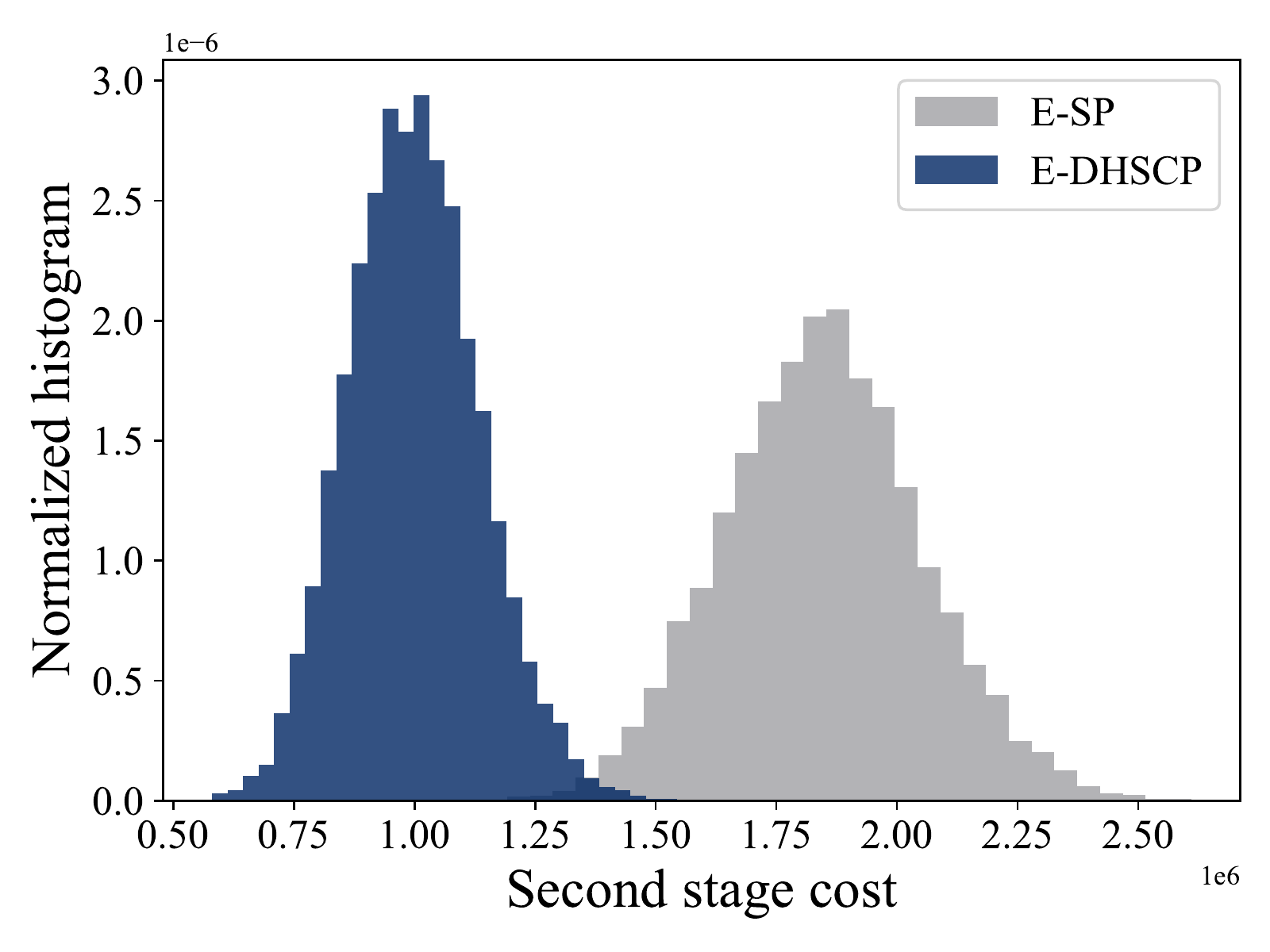}}
    \caption{Out-of-sample second stage cost for Instance 7, demand range 1 under Set 2}
    \label{EA-oop-4}
\end{figure}

\begin{figure}[!ht]
    \centering
    \subcaptionbox{$\Delta = 0$\label{EA-outu-10-2-4-8-30-4060}}{
        \includegraphics[scale=0.35]{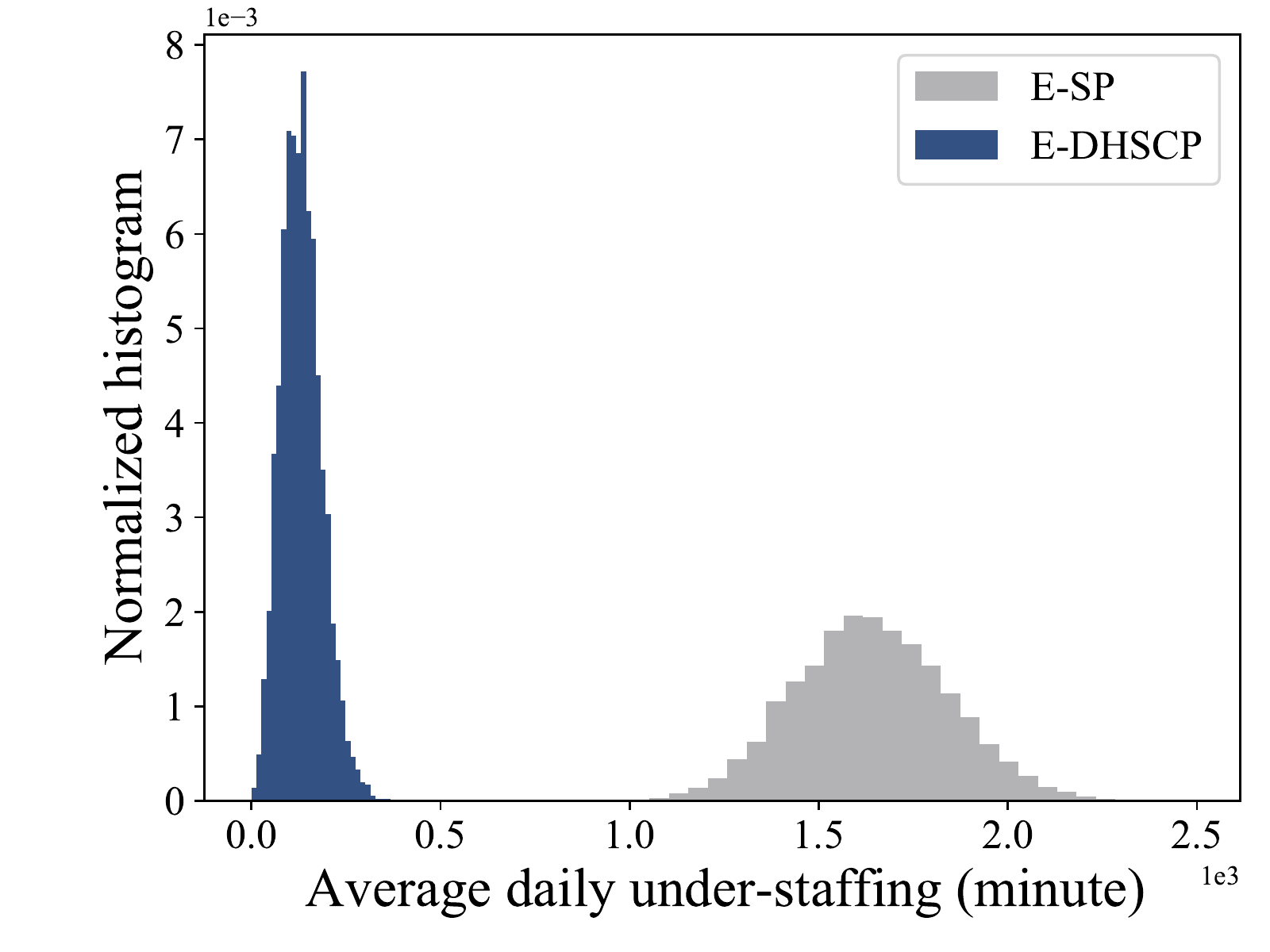}}
    \subcaptionbox{$\Delta = 0.1$\label{EA-out01u-10-2-4-8-30-4060}}{
        \includegraphics[scale=0.35]{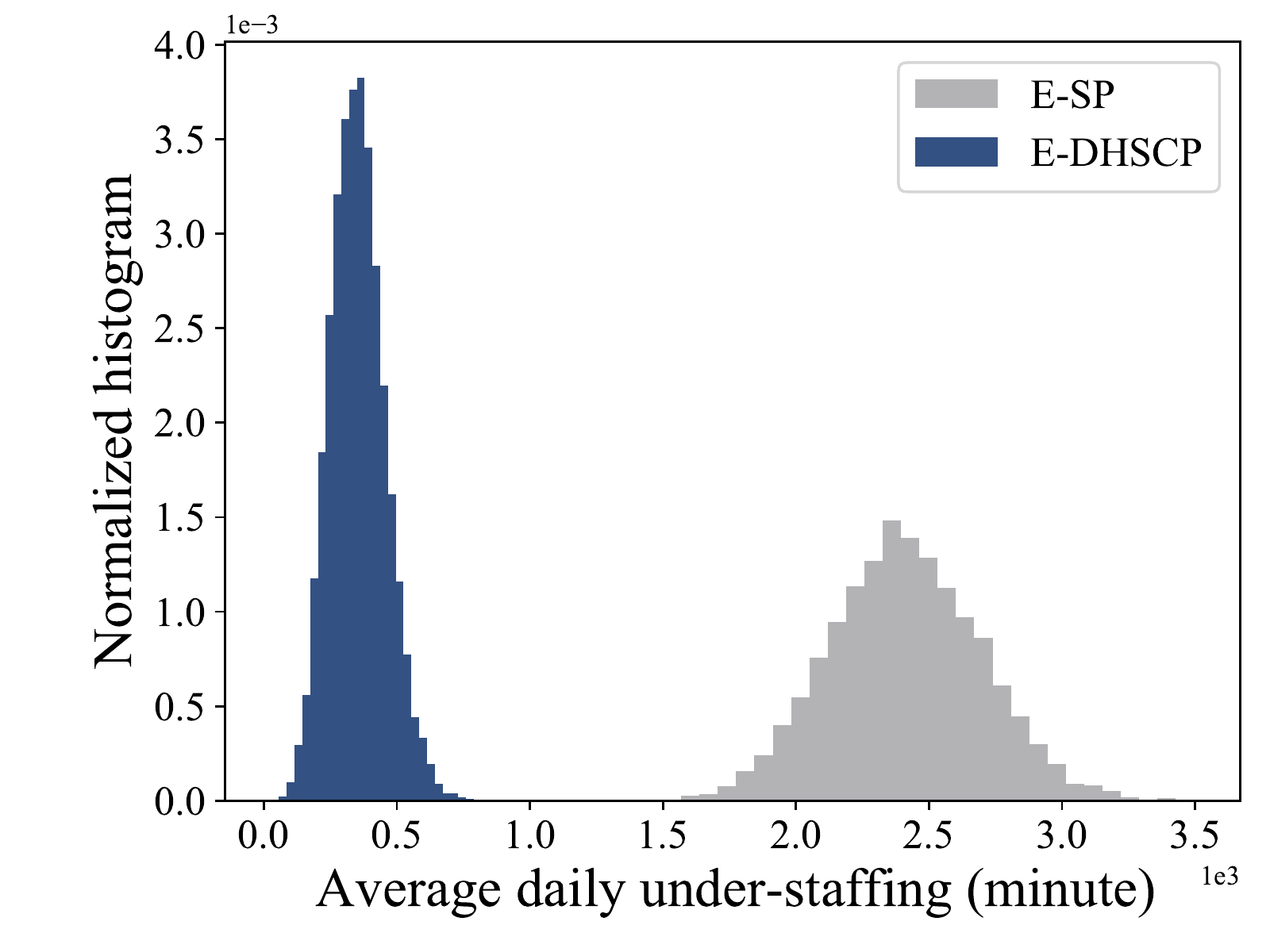}}
    \subcaptionbox{$\Delta = 0.25$\label{EA-out025u-10-2-4-8-30-4060}}{
        \includegraphics[scale=0.35]{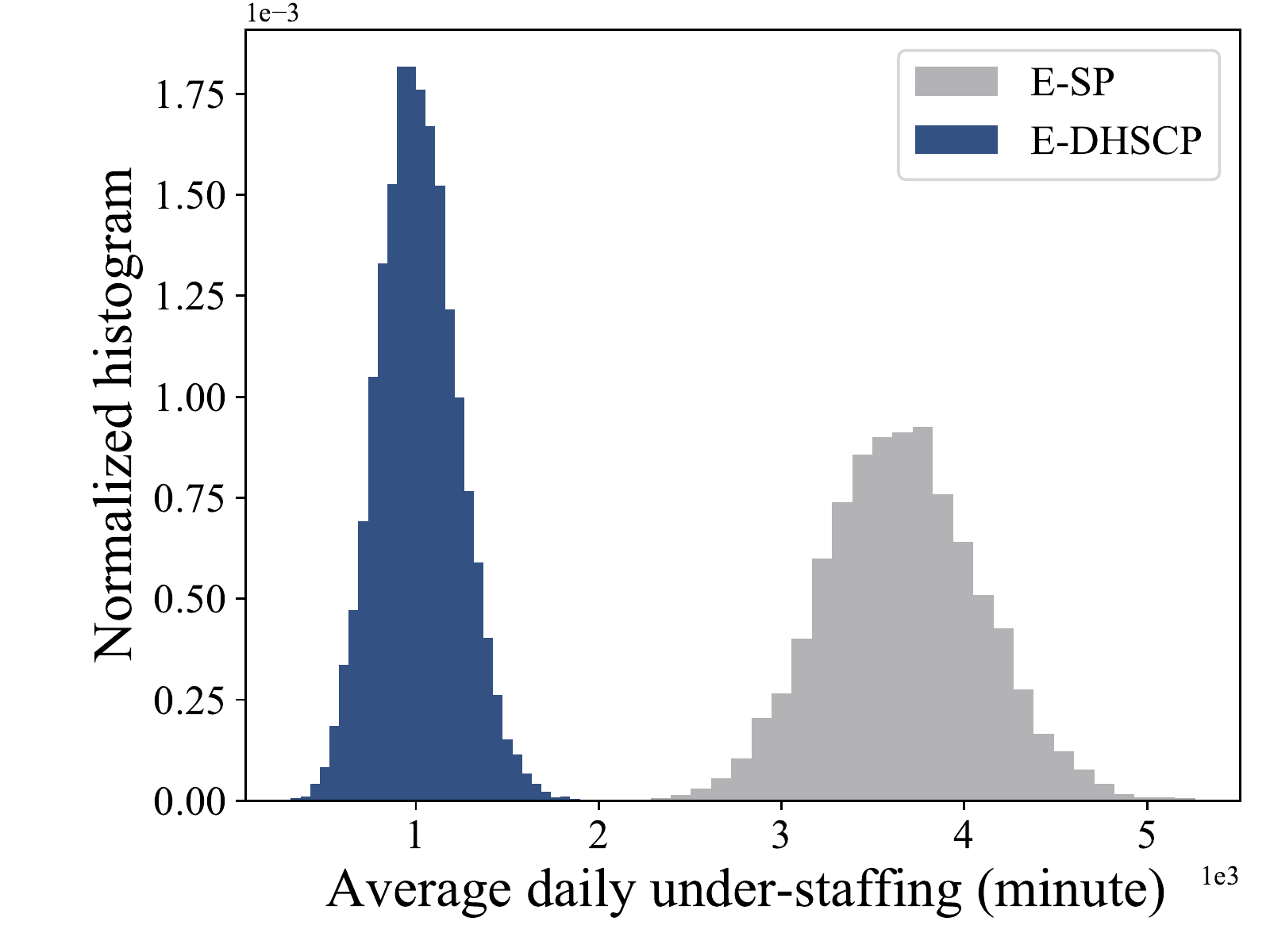}}
    \subcaptionbox{$\Delta = 0.5$\label{EA-out05u-10-2-4-8-30-4060}}{
        \includegraphics[scale=0.35]{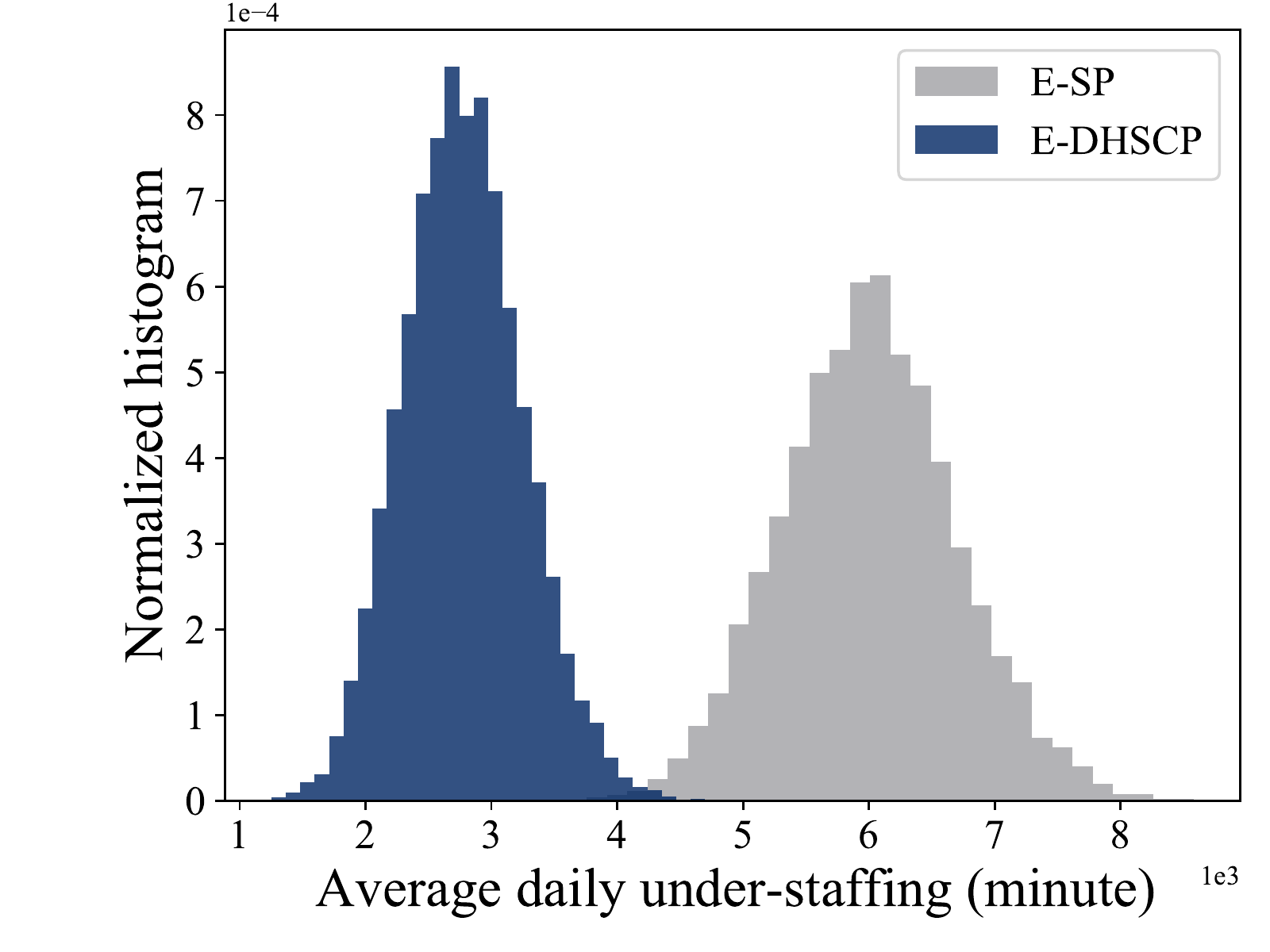}}
    \caption{Out-of-sample under-staffing for Instance 7, demand range 1 under Set 2}
    \label{EA-oop-u-4-4060}
\end{figure}
\begin{figure}[!ht]
    \centering
    \subcaptionbox{$\Delta = 0$\label{EA-outdis-10-2-4-8-30-4060}}{
        \includegraphics[scale=0.35]{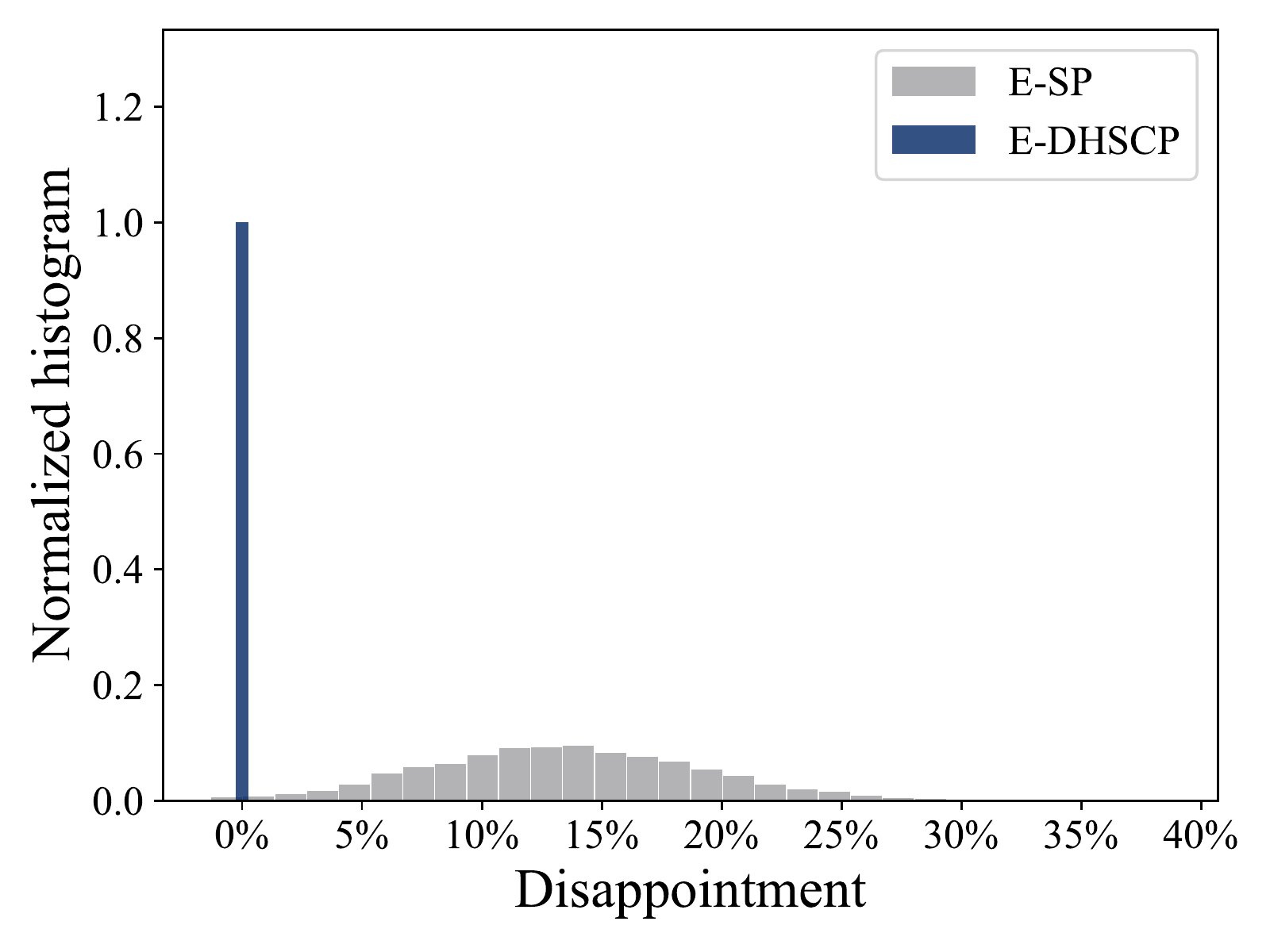}}
    \subcaptionbox{$\Delta = 0.1$\label{EA-out01dis-10-2-4-8-30-4060}}{
        \includegraphics[scale=0.35]{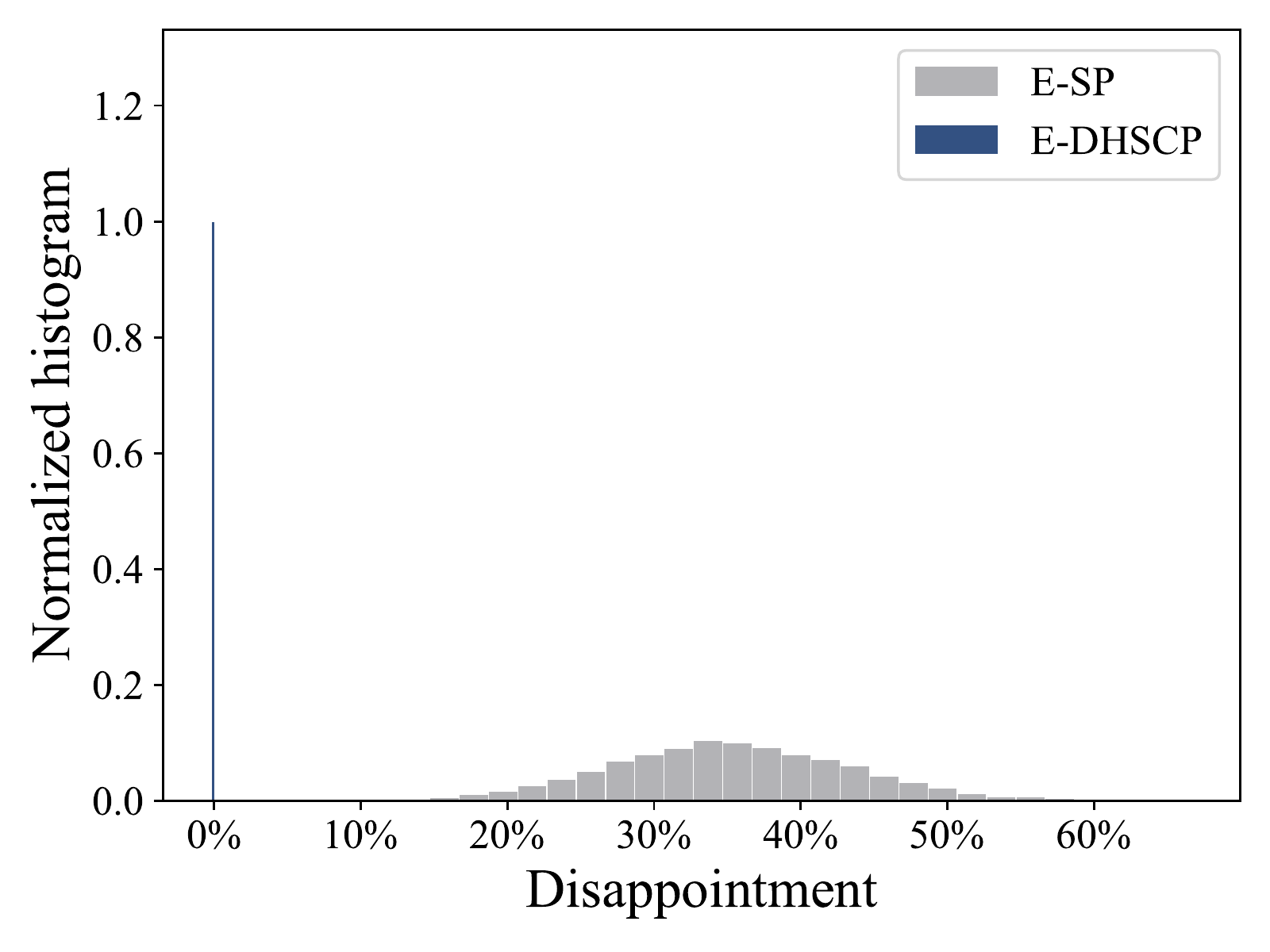}}
    \subcaptionbox{$\Delta = 0.25$\label{EA-out025dis-10-2-4-8-30-4060}}{
        \includegraphics[scale=0.35]{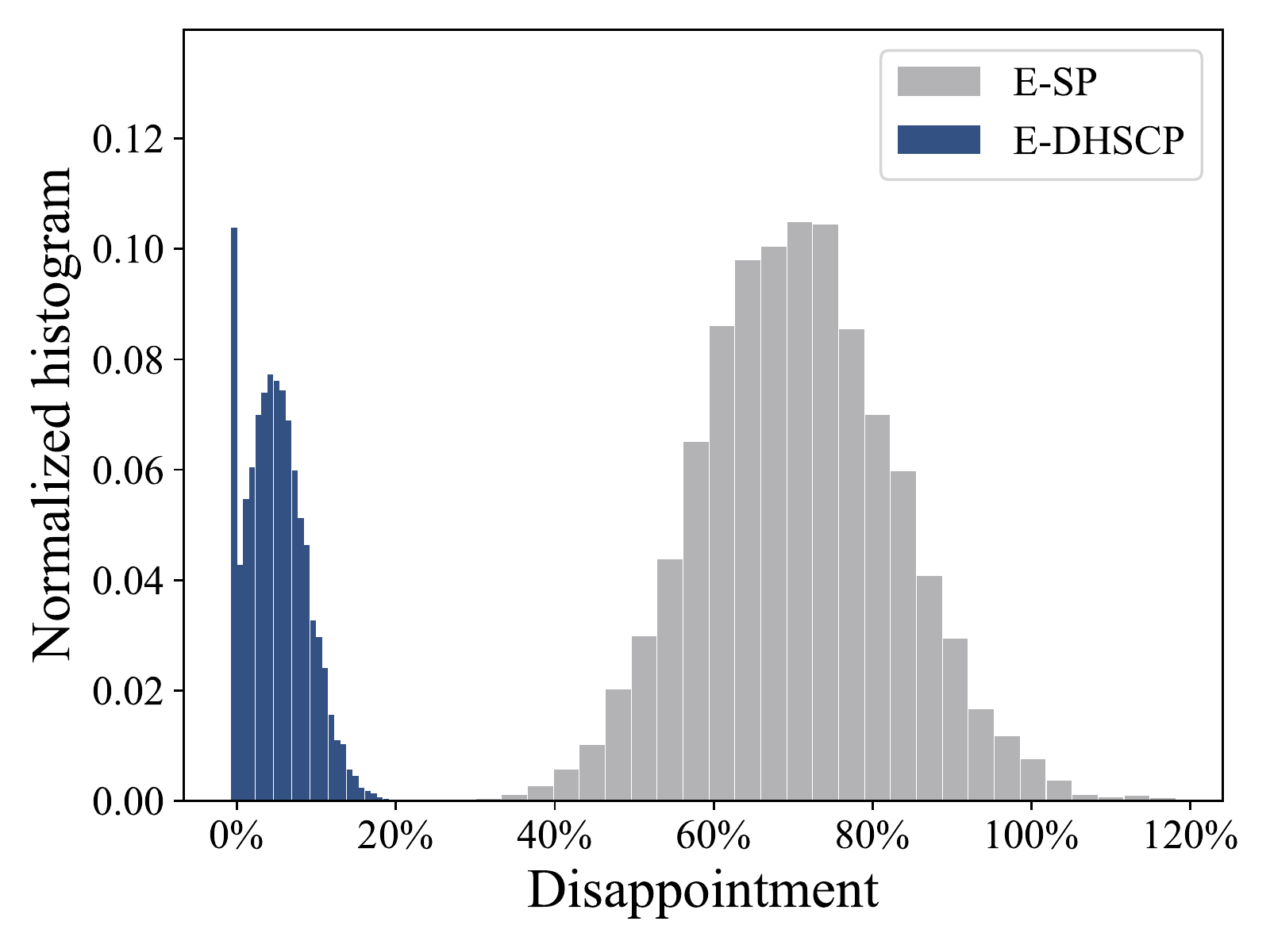}}
    \subcaptionbox{$\Delta = 0.5$\label{EA-out05dis-10-2-4-8-30-4060}}{
        \includegraphics[scale=0.35]{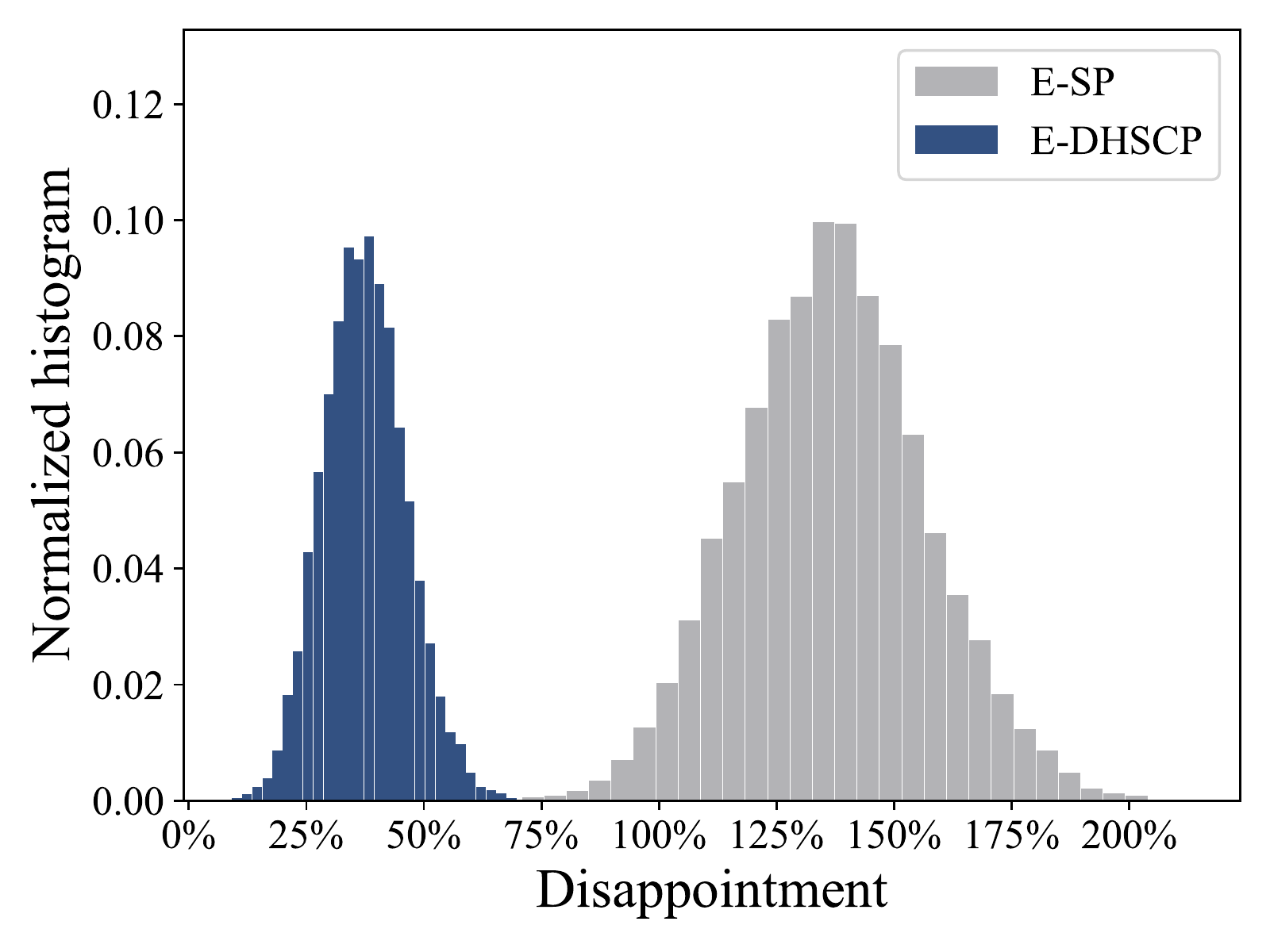}}
    \caption{Out-of-sample performance for Instance 7 demand range 1 under Set 2}
    \label{EA-oop-dis-4}
\end{figure}

\clearpage
\newpage

\section{Results of optimal staffing patterns for FA models}\label{OSP_F}
\setcounter{figure}{0}
\setcounter{table}{0}
\begin{figure}[!ht]
    \centering
    \subcaptionbox{$c_{l,t}^o = 1$\label{EA-10-4-8-30-7}}{
        \includegraphics[scale=0.35]{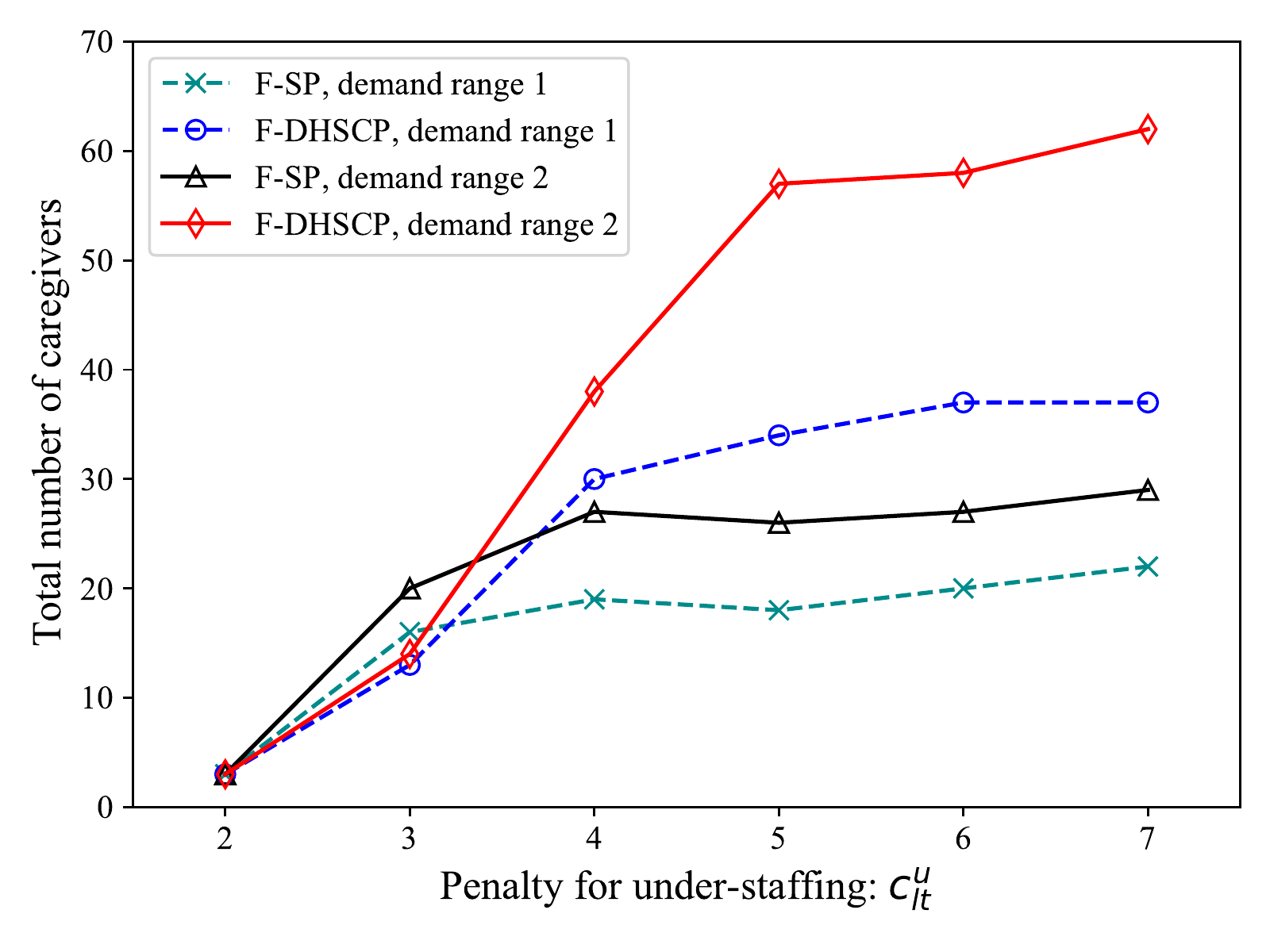}}
    \subcaptionbox{$c_{l,t}^u = 10$\label{EA-10-4-8-30-7-o}}{
        \includegraphics[scale=0.35]{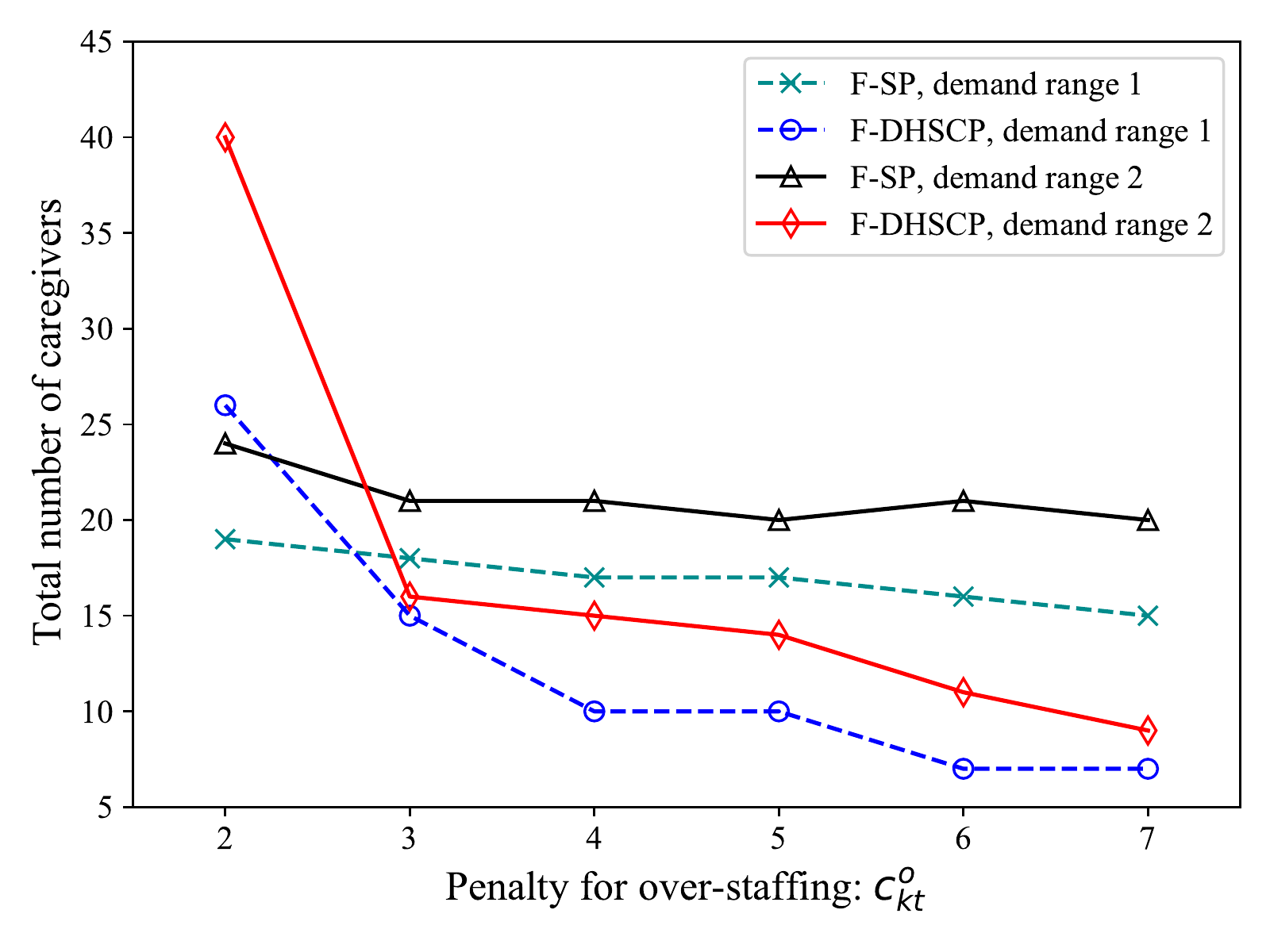}}
    \caption{Total number of caregivers hired by F-SP and F-DHSCP models for Instance 7.}
    \label{total_num_fig_f_a}
\end{figure}

\begin{figure}[!ht]
    \centering
    \subcaptionbox{Instance 10, demand range 1\label{FA-10-2-6-6-30-4060}}{
        \includegraphics[scale=0.35]{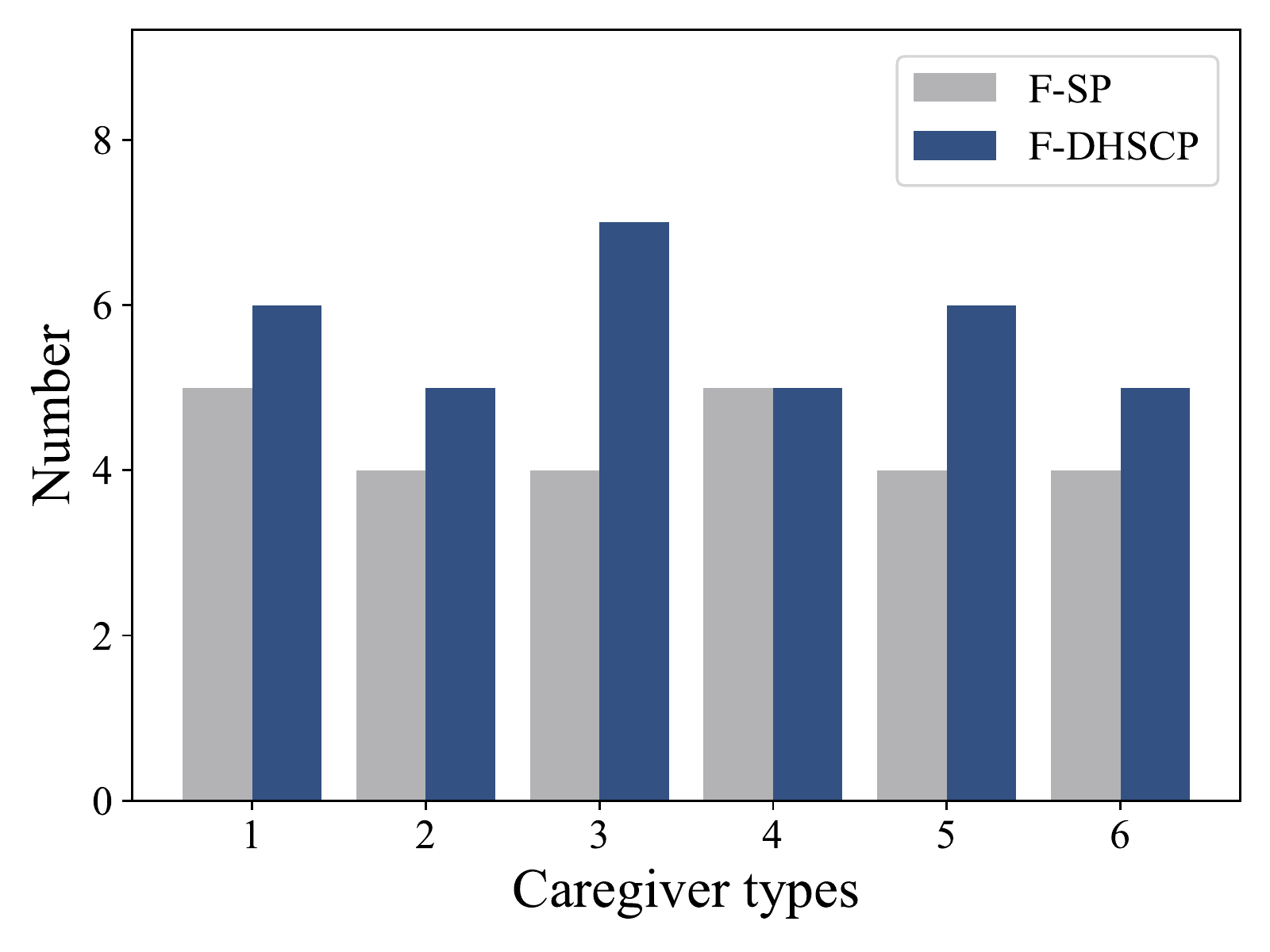}}
    \subcaptionbox{Instance 10, demand range 2\label{FA-10-2-6-6-30-1080}}{
        \includegraphics[scale=0.35]{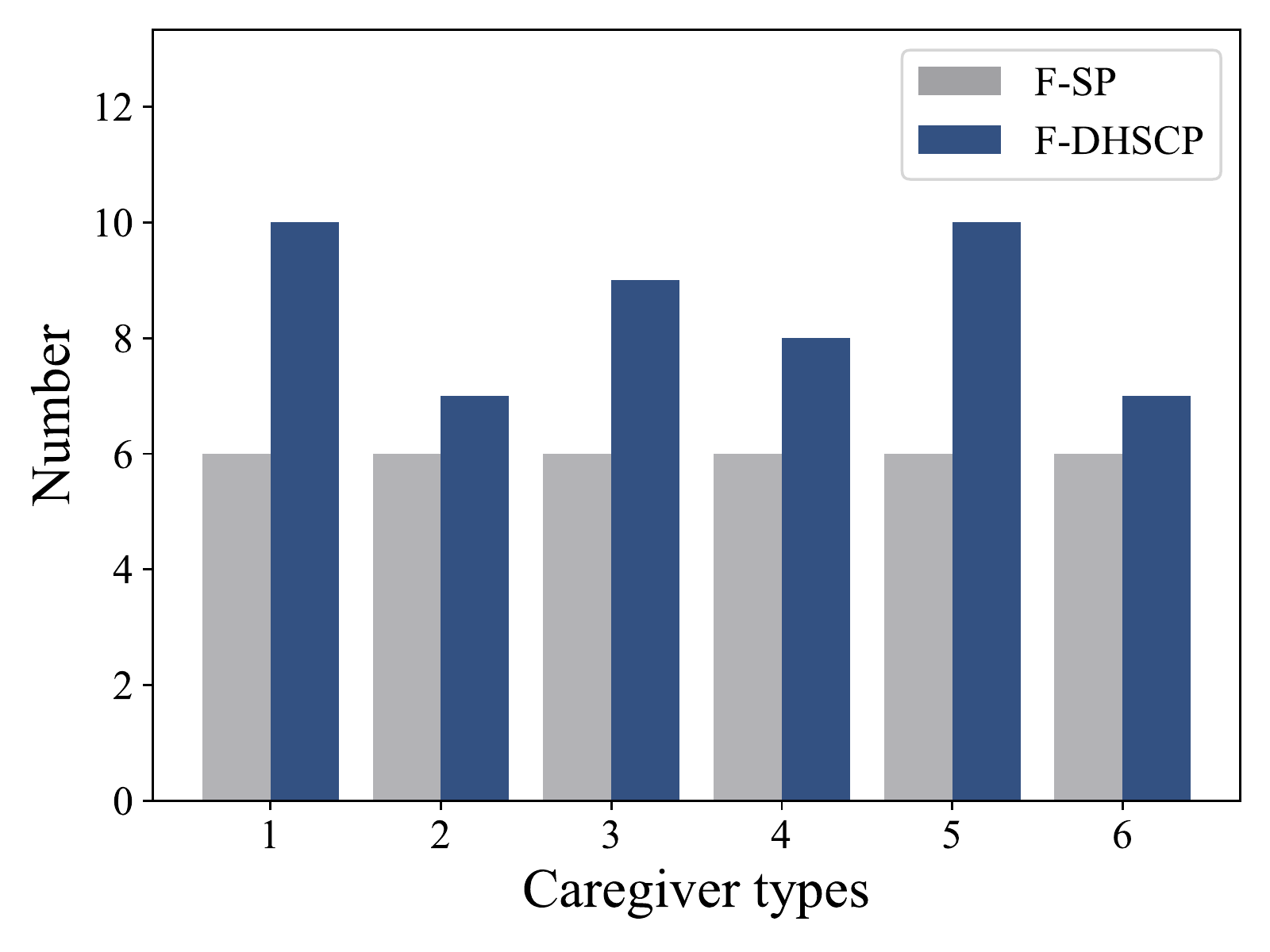}}
    \subcaptionbox{Instance 7, demand range 1\label{FA-10-2-4-8-30-4060}}{
        \includegraphics[scale=0.35]{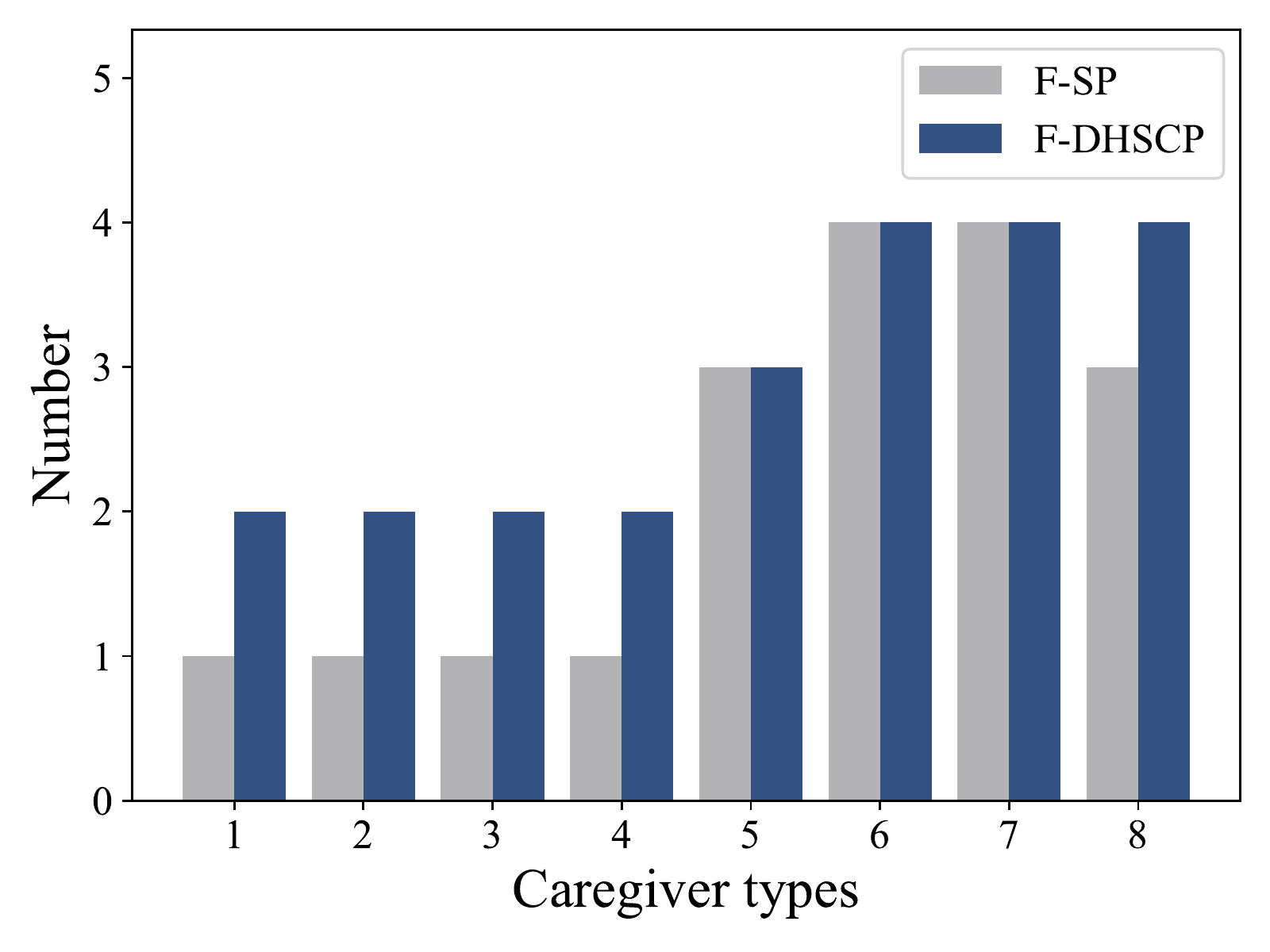}}
    \subcaptionbox{Instance 7, demand range 2\label{FA-10-2-4-8-30-1080}}{
        \includegraphics[scale=0.35]{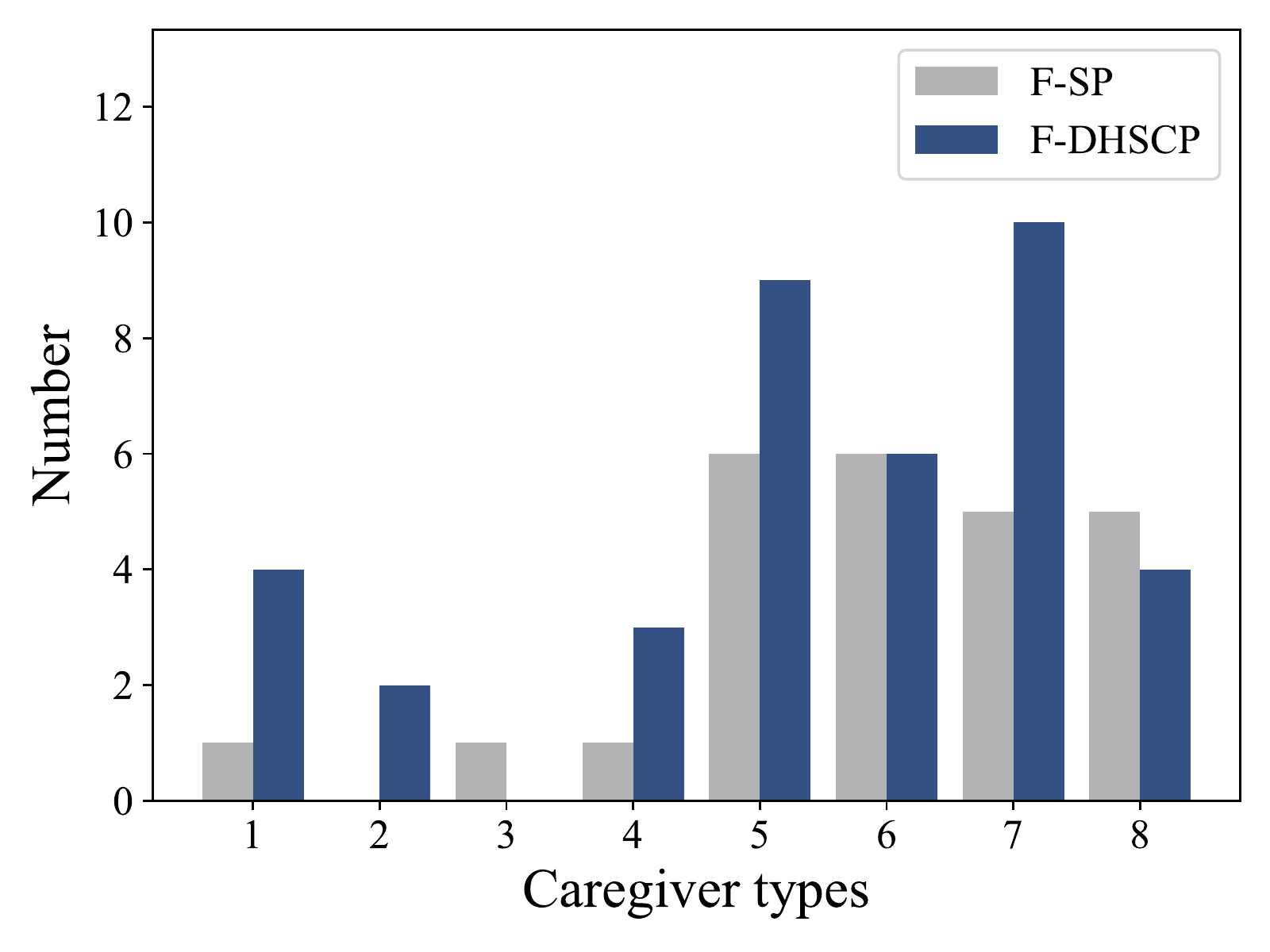}}
    \caption{Staffing patterns of FA decision-makers under cost structure $(c_{l,t}^{u}, c_{k,t}^{o}) = (5, 2)$.}
    \label{FA-solution}
\end{figure}

\clearpage
\newpage

\section{Out-of-sample performance of FA models}\label{oop-FA}
\subsection{Out-of-sample performance of FA models for Instance 10}\label{oop-FA-10}
\setcounter{figure}{0}
\setcounter{table}{0}
\begin{figure}[!ht]
    \centering
    \subcaptionbox{\label{FA-1outc-10-2-6-6-30-1080}}{
        \includegraphics[scale=0.3]{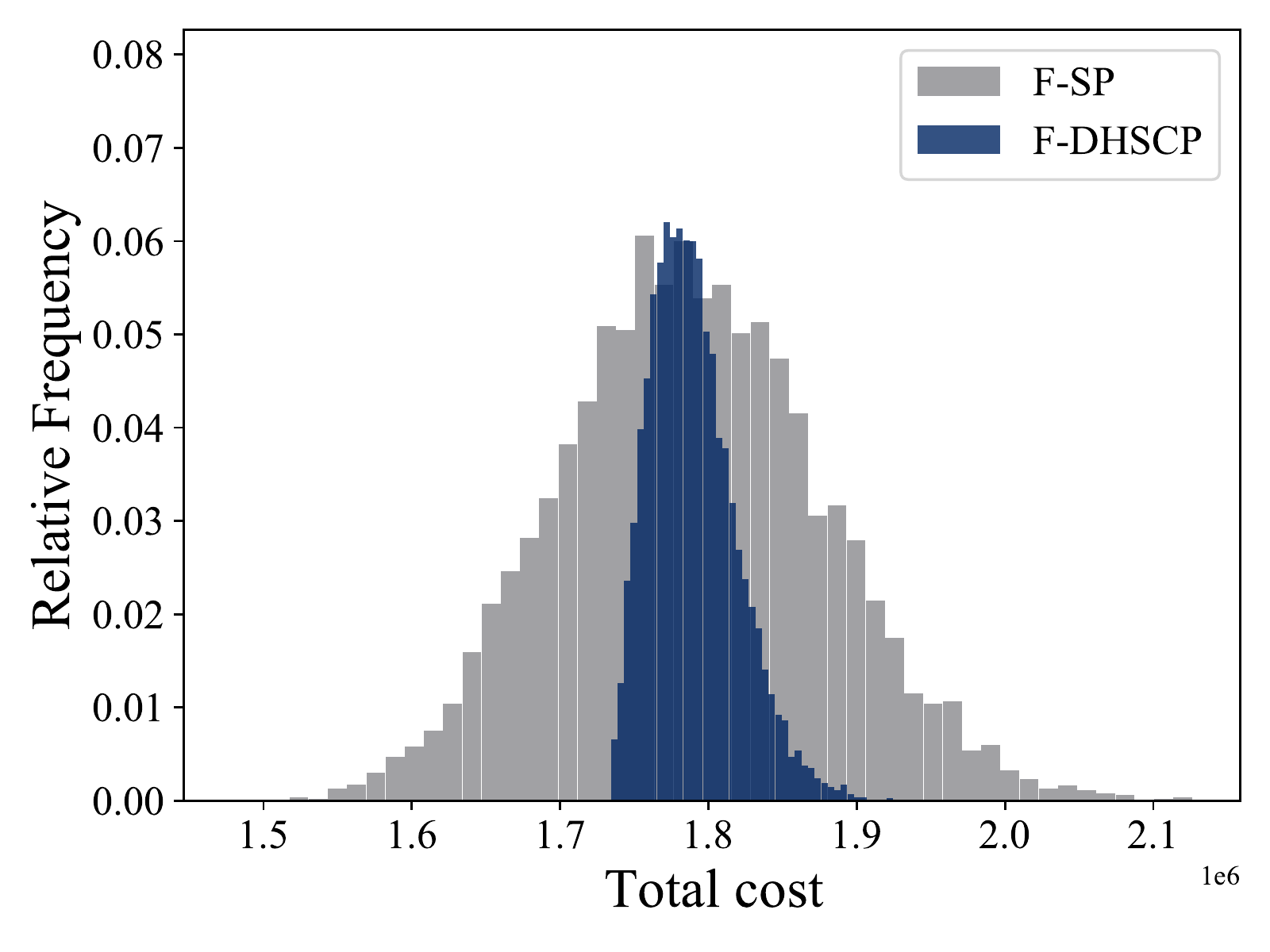}}
    \subcaptionbox{\label{FA-1outc2-10-2-6-6-30-1080}}{
        \includegraphics[scale=0.3]{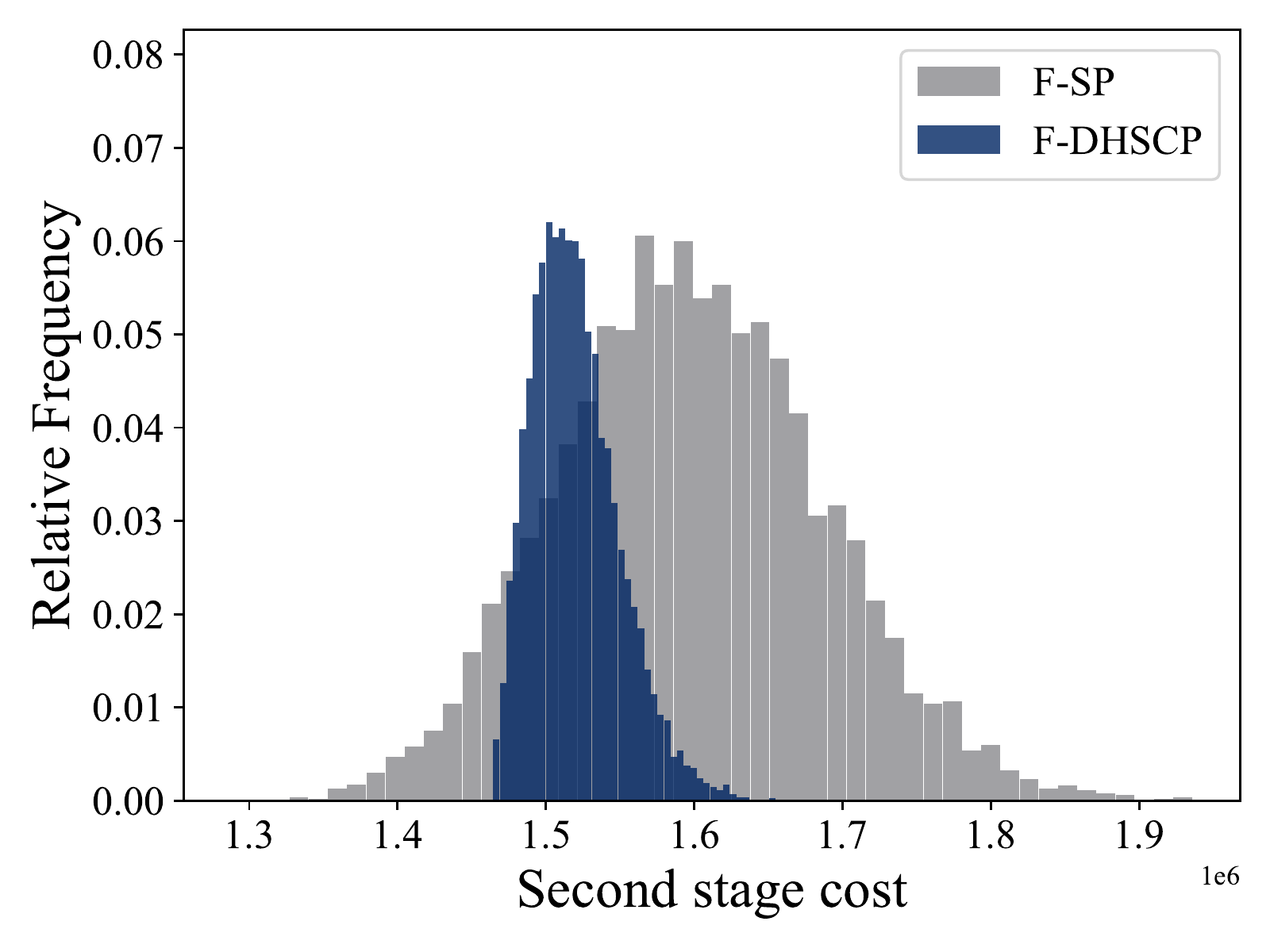}}
    \subcaptionbox{\label{FA-1outu-10-2-6-6-30-1080}}{
        \includegraphics[scale=0.3]{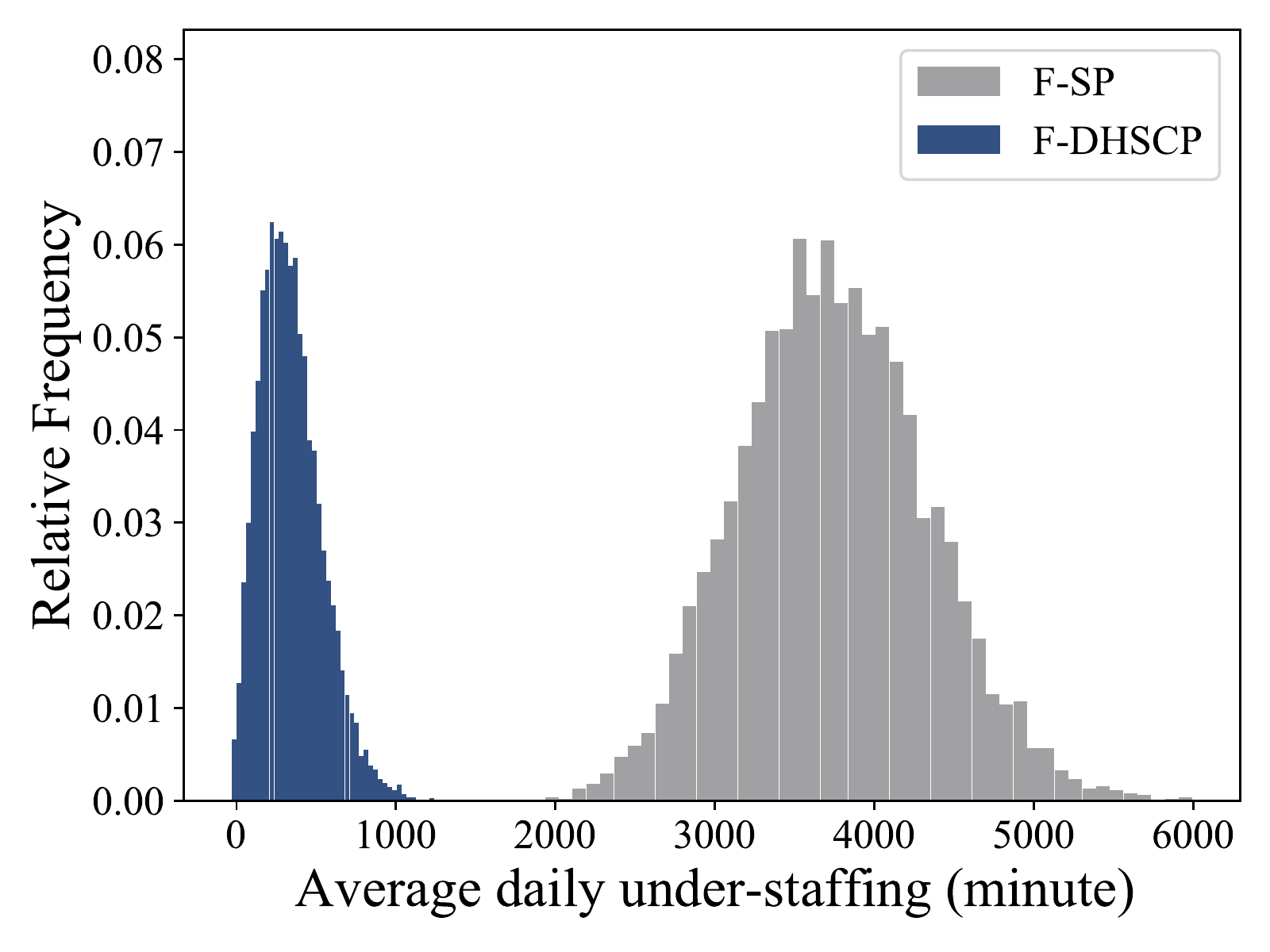}}
    \caption{Out-of-sample performance of FA models for Instance 10, demand range 2 under Set 2, $\Delta = 0$}
    \label{FA-oop-o-6-2}
\end{figure}
\begin{figure}[!ht]
    \centering
    \subcaptionbox{$\Delta = 0$\label{FA-outc2-10-2-6-6-30-4060}}{
        \includegraphics[scale=0.33]{Foutsamplecost2.pdf}}
    \subcaptionbox{$\Delta = 0.1$\label{FA-out01c2-10-2-6-6-30-4060}}{
        \includegraphics[scale=0.33]{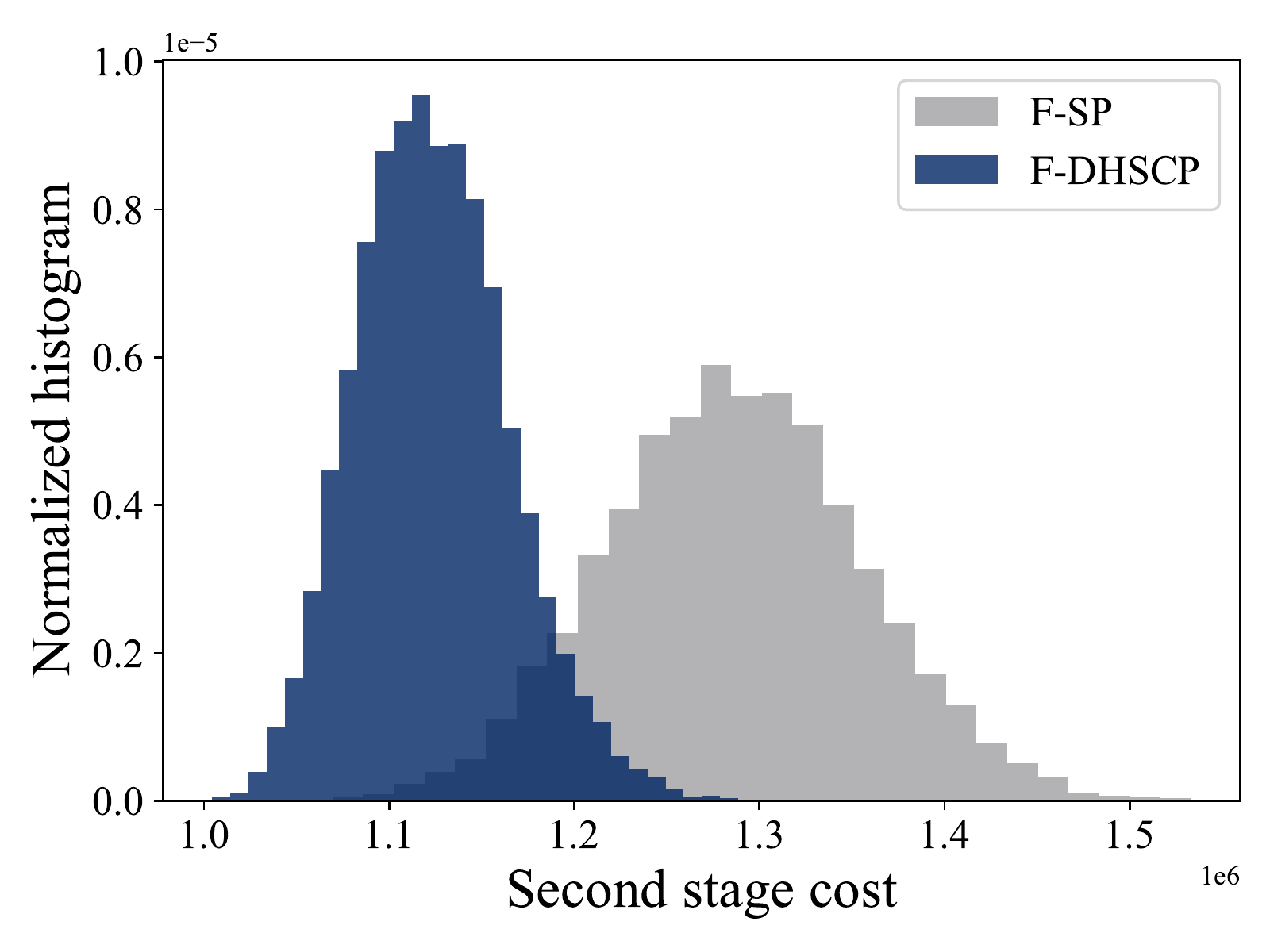}}
    \subcaptionbox{$\Delta = 0.25$\label{FA-out025c2-10-2-6-6-30-4060}}{
        \includegraphics[scale=0.33]{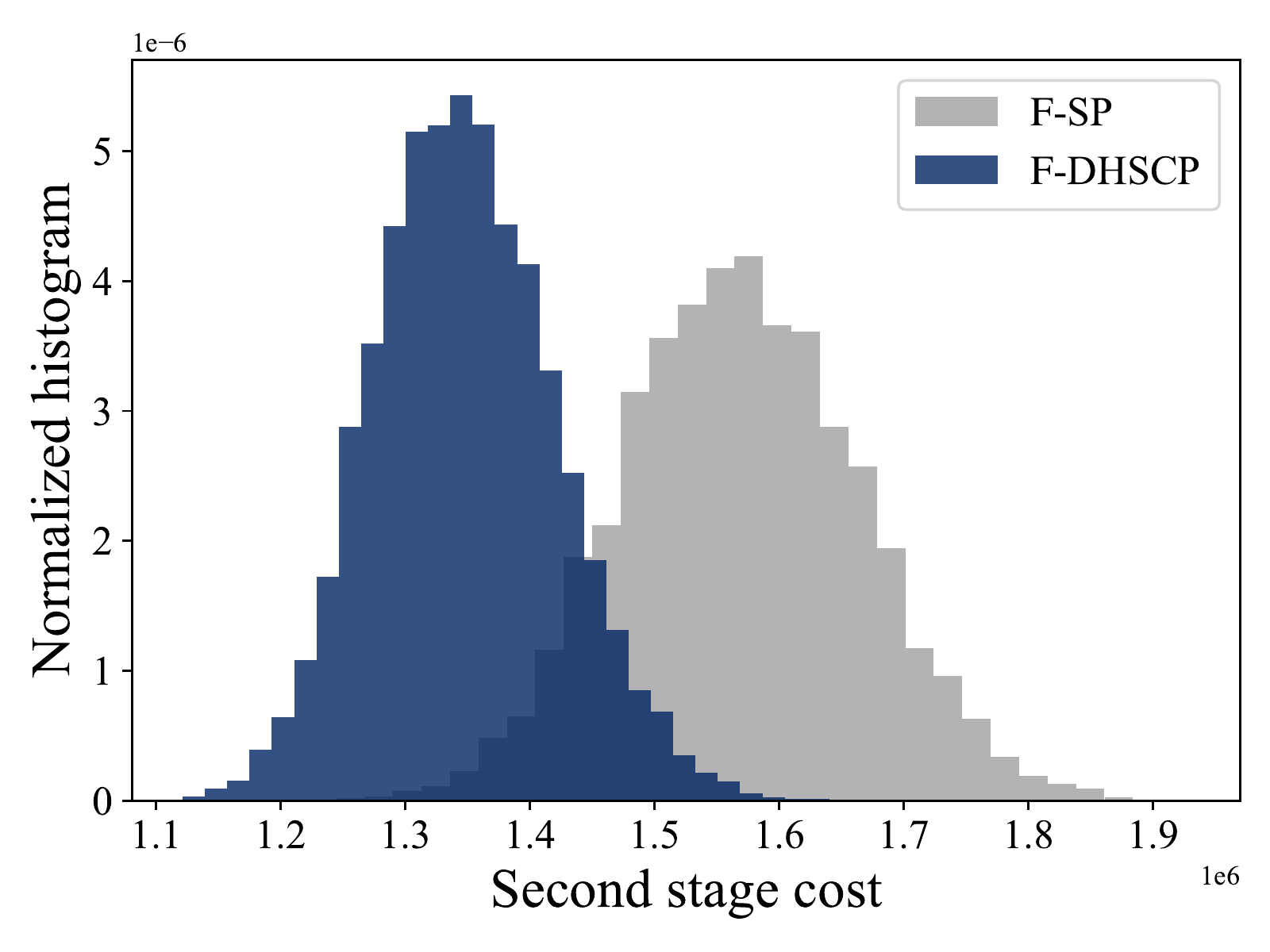}}
    \subcaptionbox{$\Delta = 0.5$\label{FA-out05c2-10-2-6-6-30-4060}}{
        \includegraphics[scale=0.33]{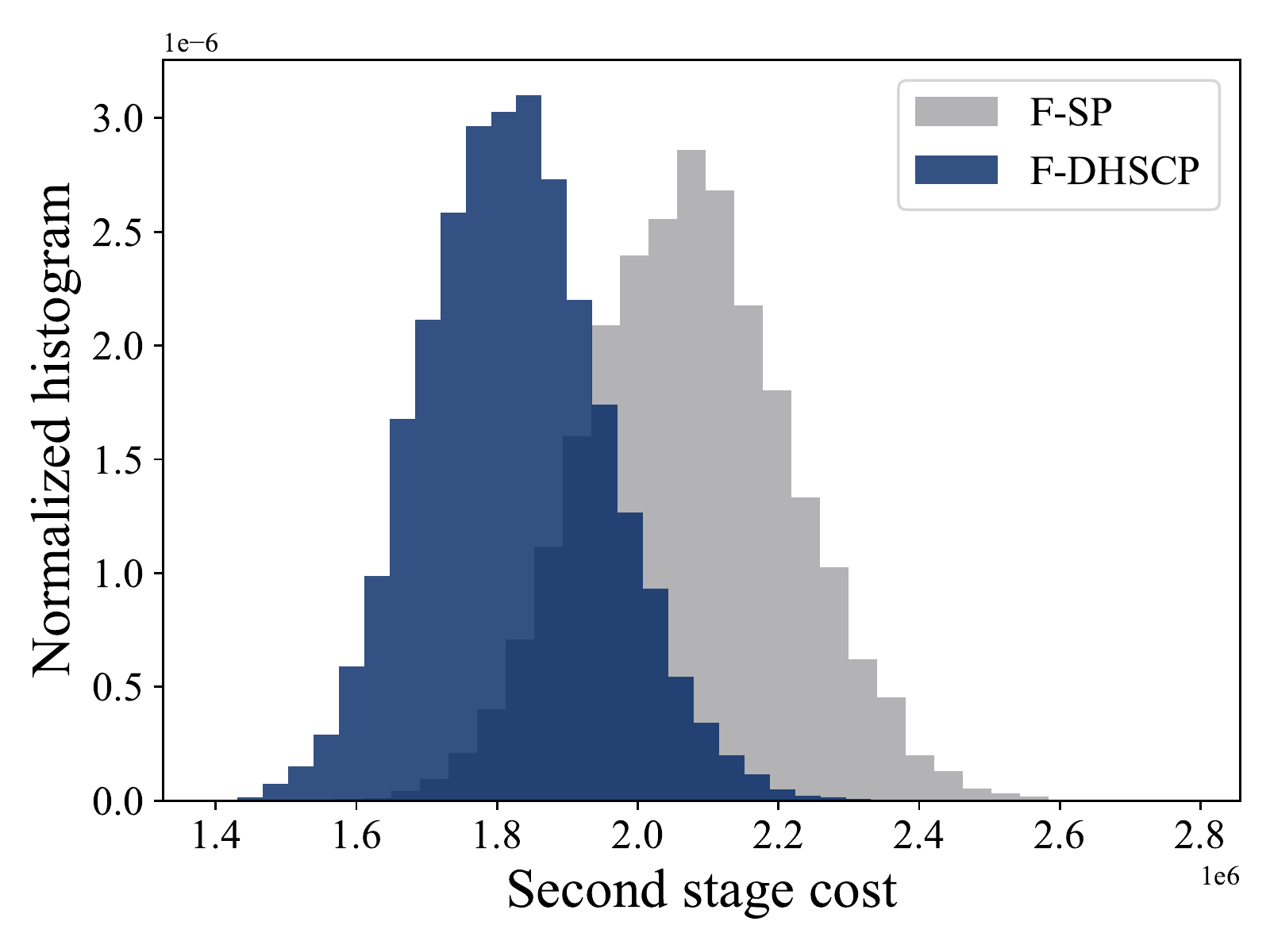}}
    \caption{Out-of-sample second-stage cost of FA models for Instance 10, demand range 1 under Set 2}
    \label{FA-oop-6}
\end{figure}

\begin{figure}[!ht]
    \centering
    \subcaptionbox{$\Delta = 0$\label{FA-outu-10-2-6-6-30-4060}}{
        \includegraphics[scale=0.35]{Foutsampleunder.pdf}}
    \subcaptionbox{$\Delta = 0.1$\label{FA-out01u-10-2-6-6-30-4060}}{
        \includegraphics[scale=0.35]{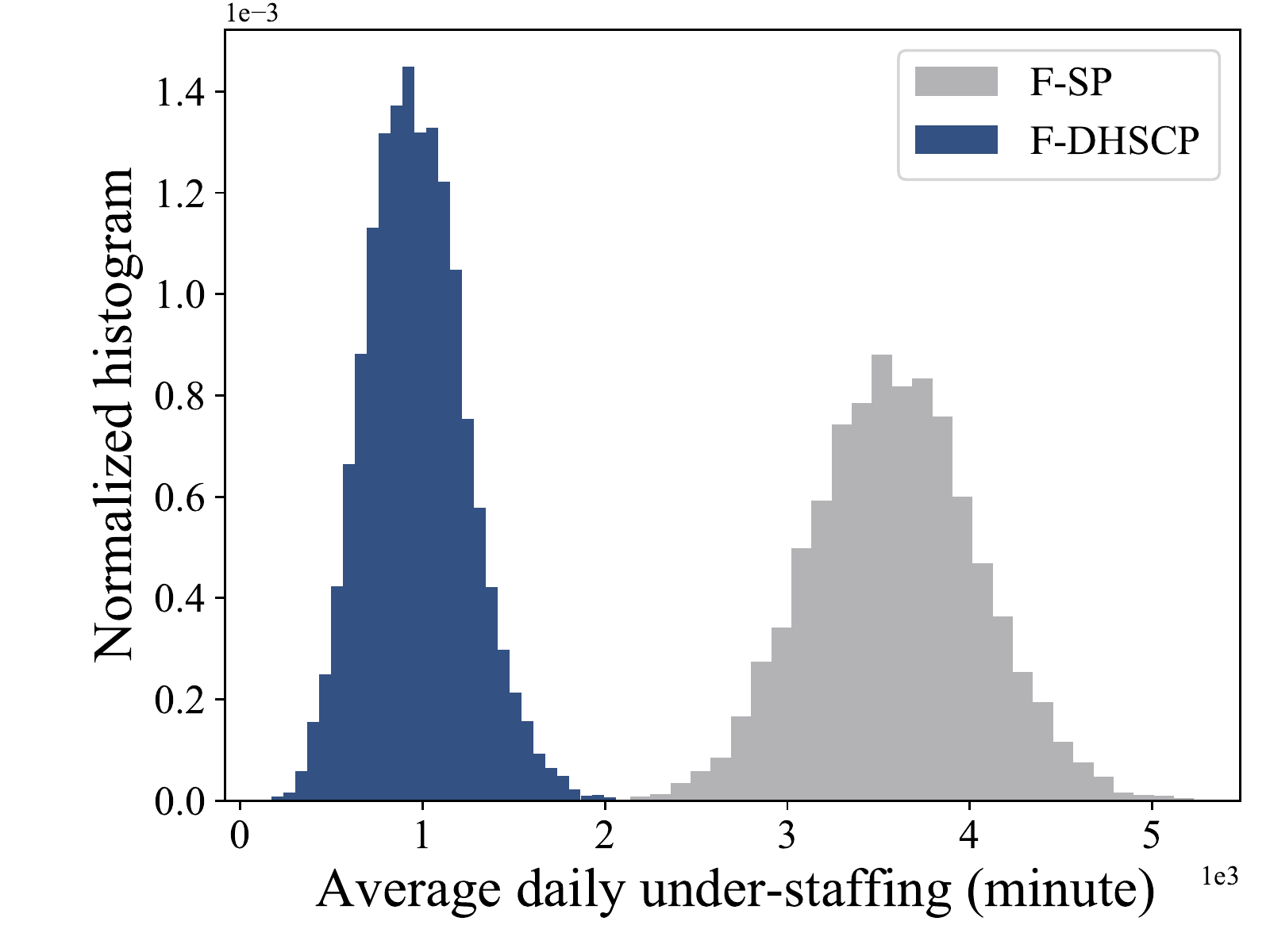}}
    \subcaptionbox{$\Delta = 0.25$\label{FA-out025u-10-2-6-6-30-4060}}{
        \includegraphics[scale=0.35]{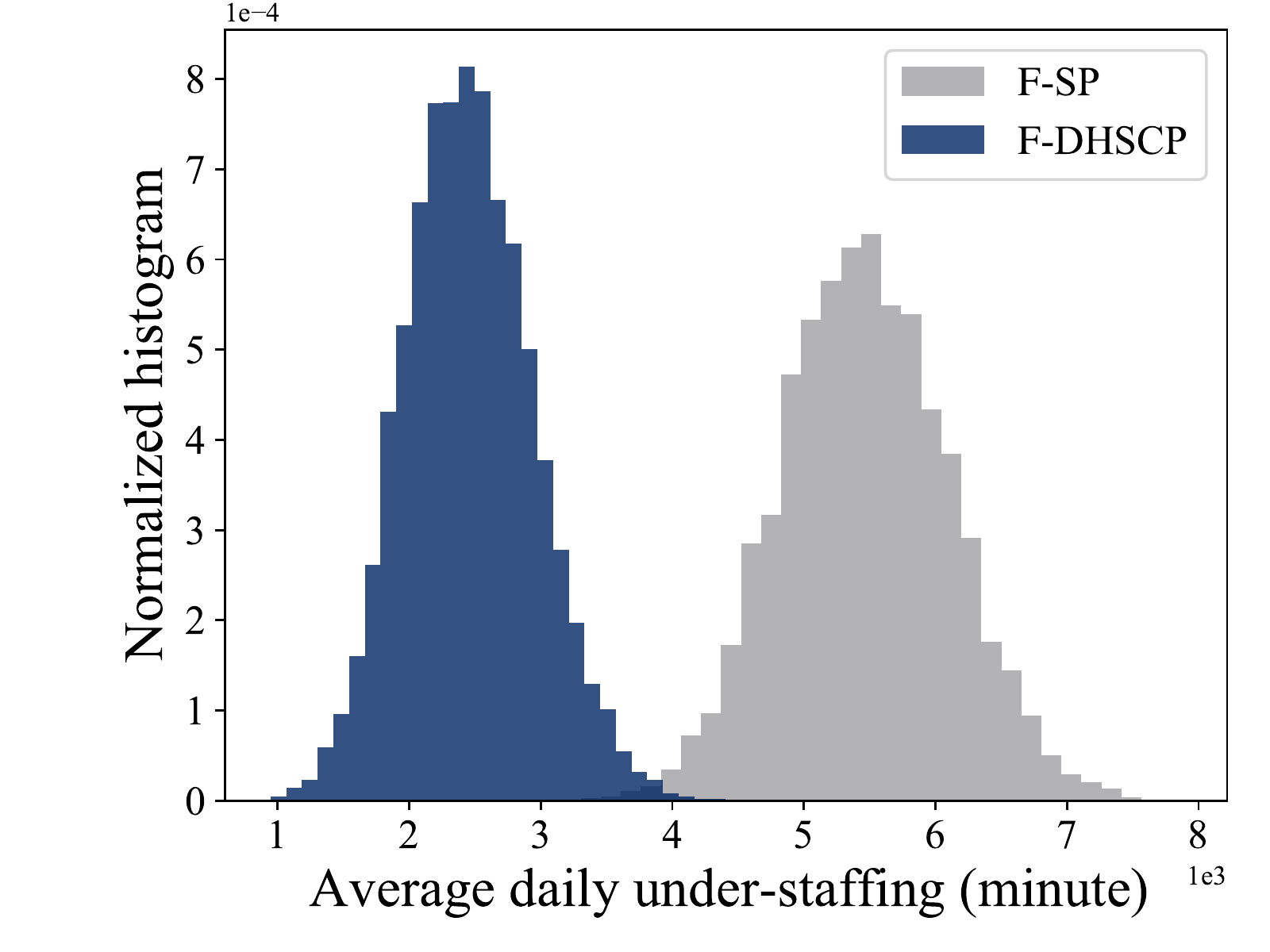}}
    \subcaptionbox{$\Delta = 0.5$\label{FA-out05u-10-2-6-6-30-4060}}{
        \includegraphics[scale=0.35]{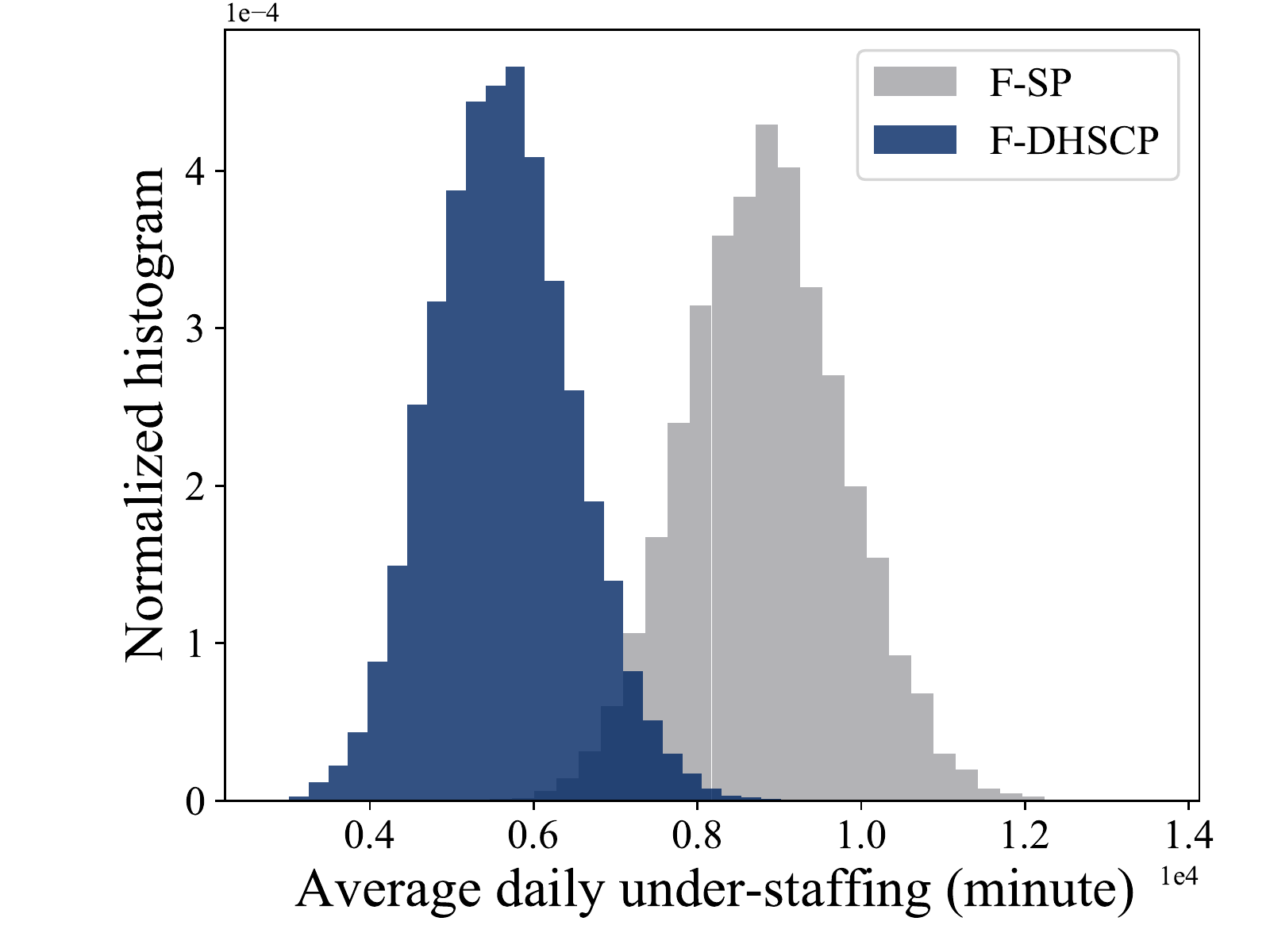}}
    \caption{Out-of-sample under-staffing for Instance 10, demand range 1 under Set 2}
    \label{FA-oop-u-6-4060}
\end{figure}

\newpage

\subsection{Out-of-sample performance of FA models for Instance 7}\label{oop-FA-7}
\begin{figure}[!ht]
    \centering
    \subcaptionbox{\label{FA-1outc-7-2-4-8-30-4060}}{
        \includegraphics[scale=0.3]{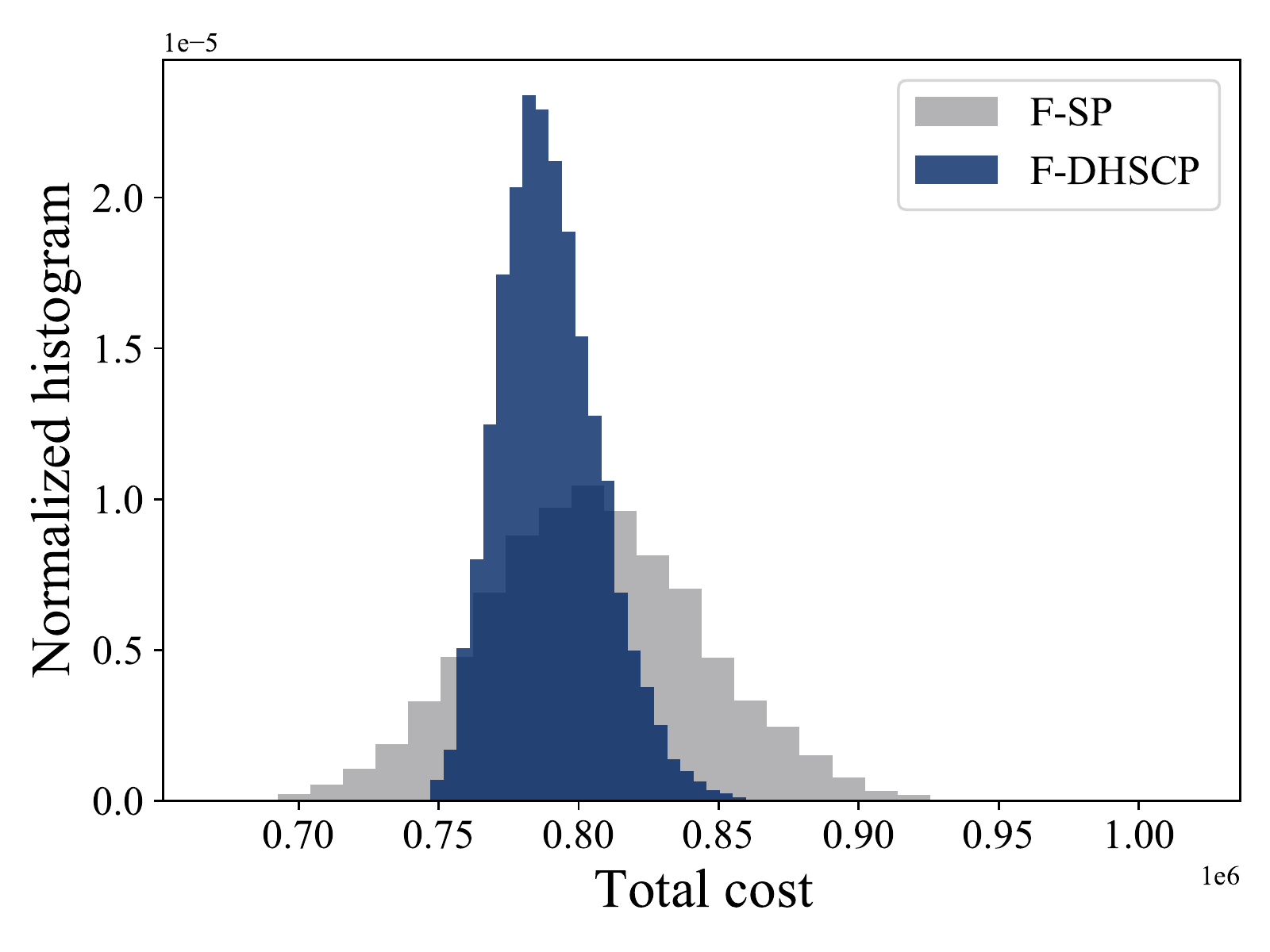}}
    \subcaptionbox{\label{FA-1outc2-7-2-4-8-30-4060}}{
        \includegraphics[scale=0.3]{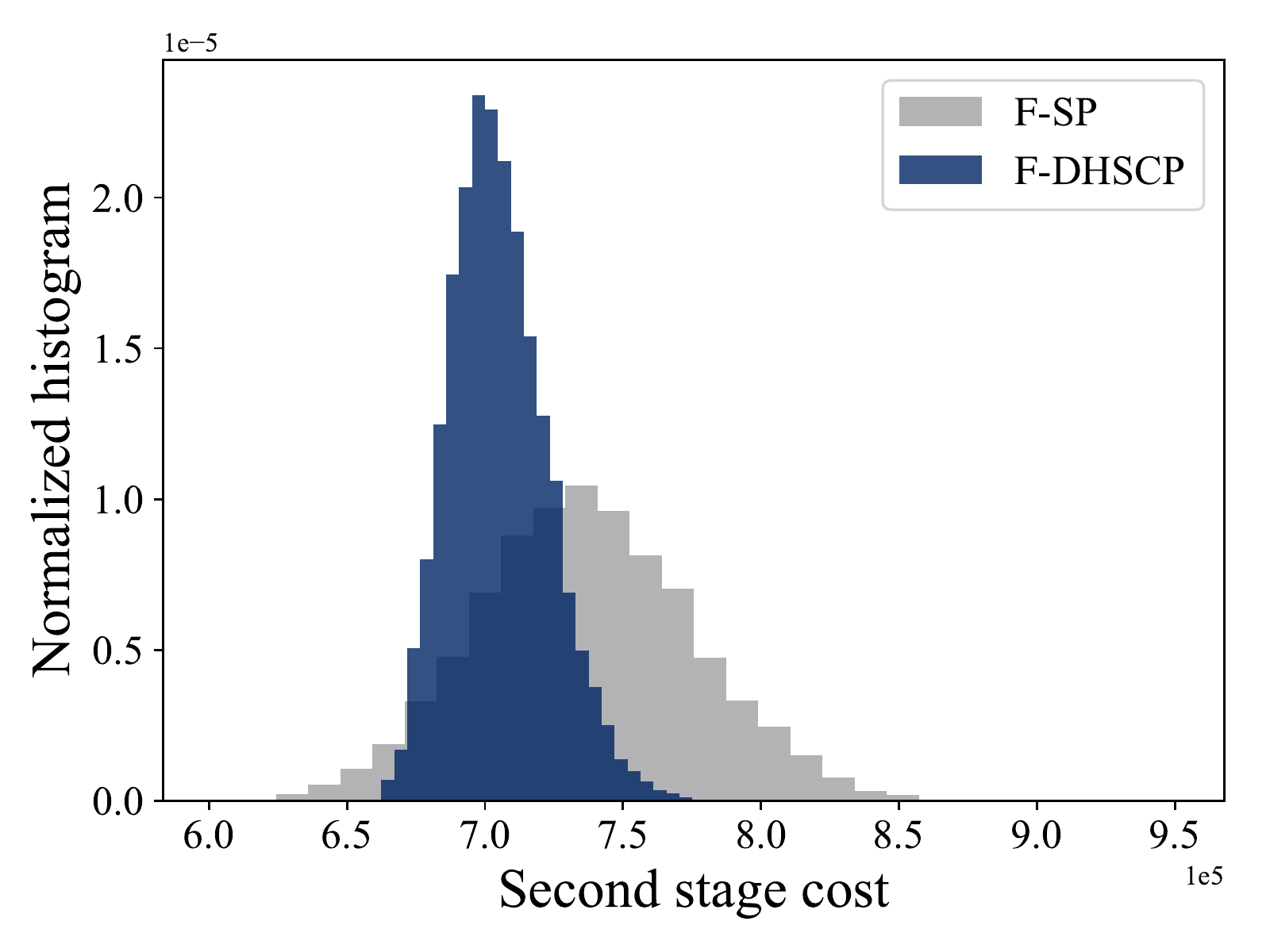}}
    \subcaptionbox{\label{FA-1outu-7-2-4-8-30-4060}}{
        \includegraphics[scale=0.3]{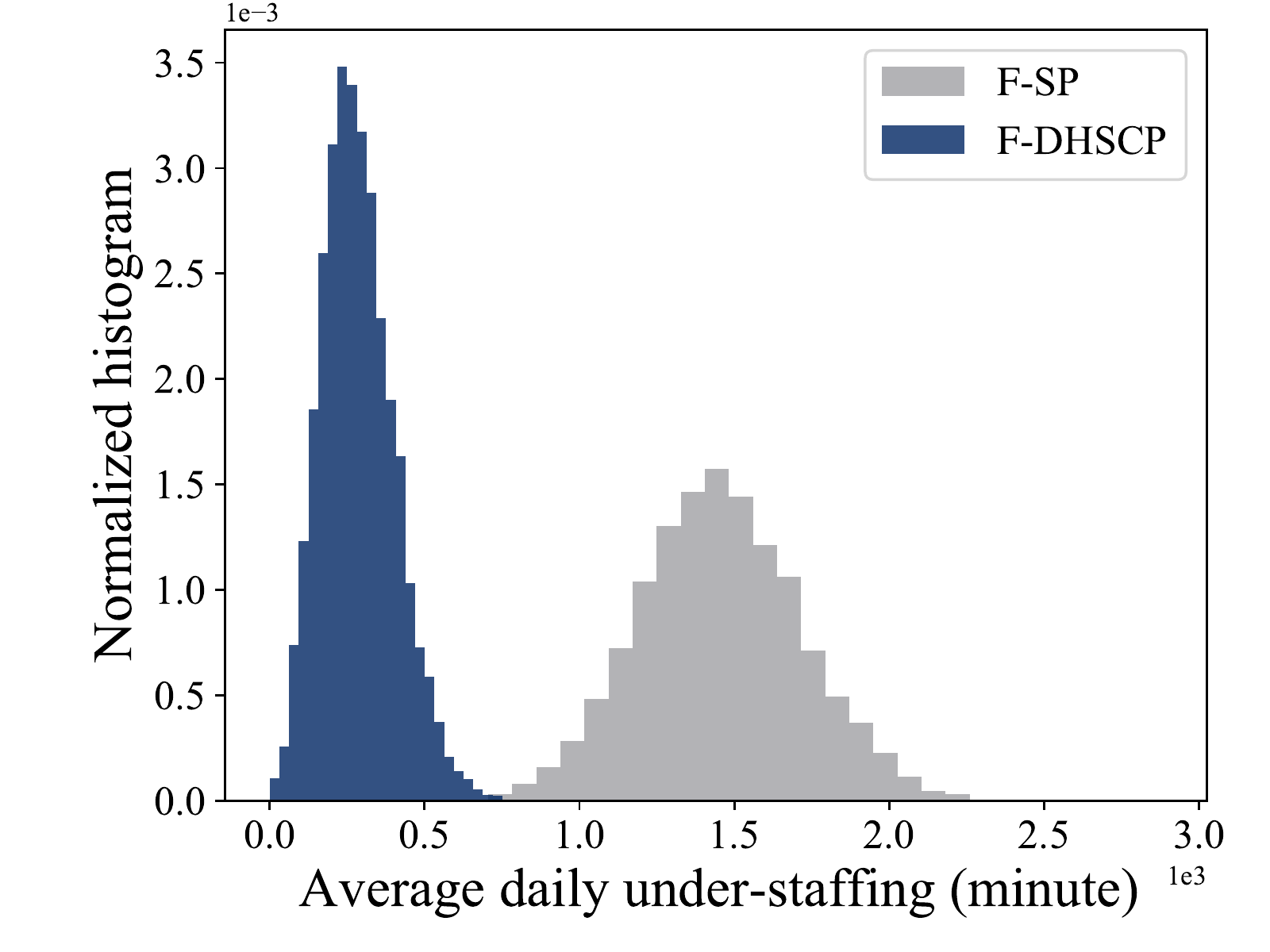}}
    \caption{Out-of-sample performance of FA models for Instance 7, demand range 1 under Set 2, $\Delta = 0$}
    \label{FA-oop-o-4}
\end{figure}
\begin{figure}[!ht]
    \centering
    \subcaptionbox{$\Delta=0$\label{FA-outdis-7-2-4-8-30-4060}}{
        \includegraphics[scale=0.3]{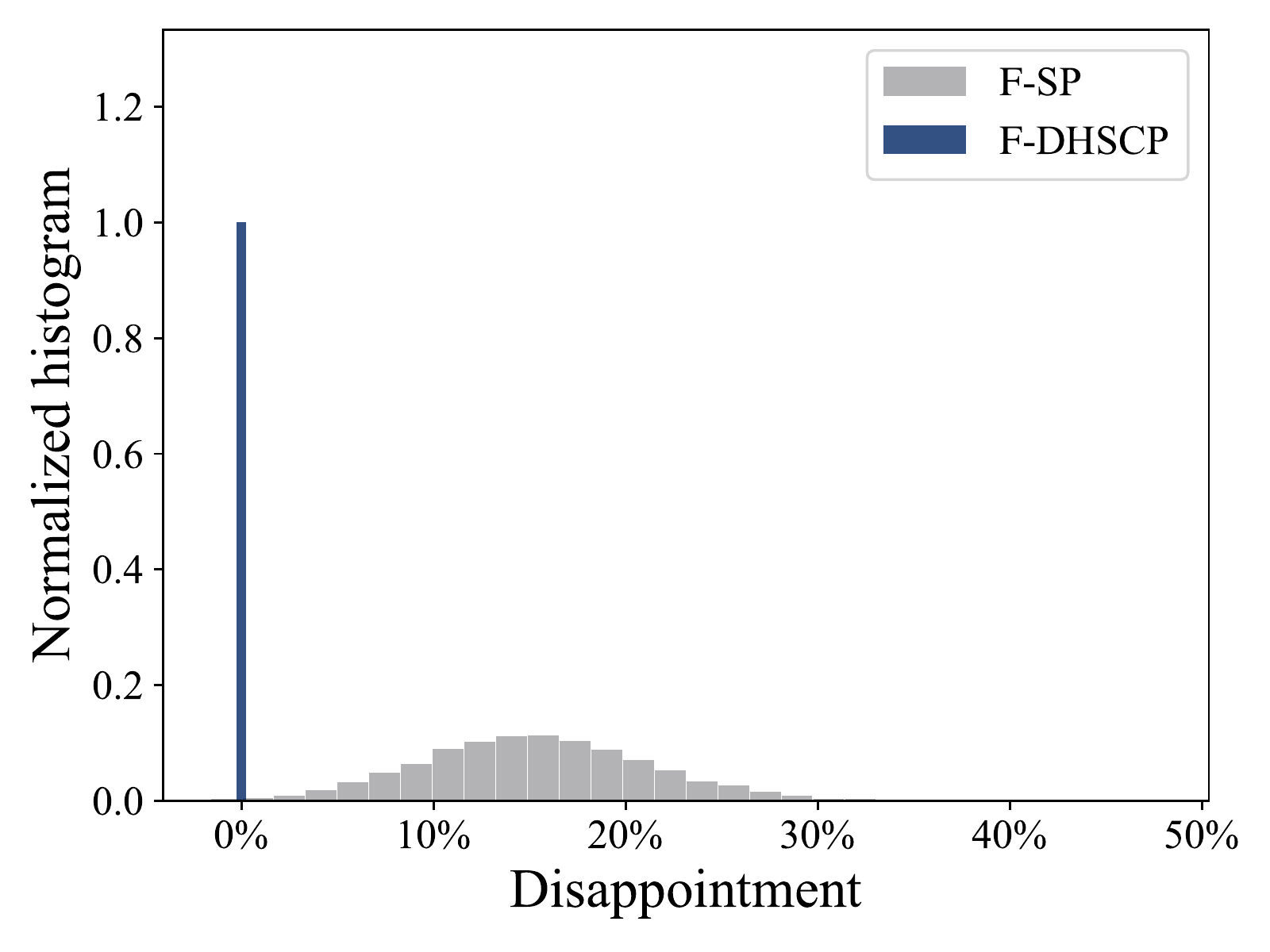}}
    \subcaptionbox{$\Delta=0.25$\label{FA-out025dis-7-2-4-8-30-4060}}{
        \includegraphics[scale=0.3]{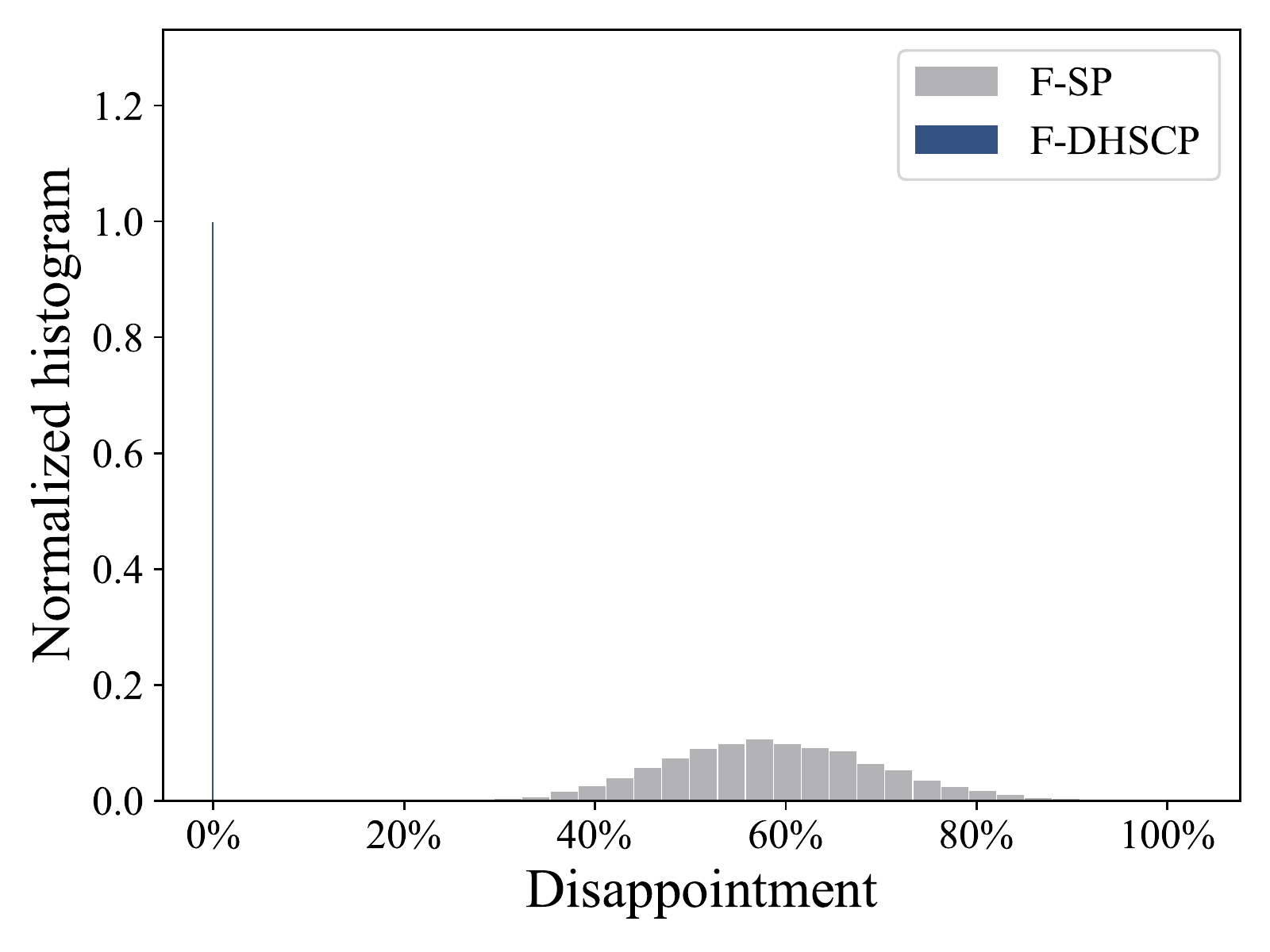}}
    \subcaptionbox{$\Delta=0.5$\label{FA-out05dis-7-2-4-8-30-4060}}{
        \includegraphics[scale=0.3]{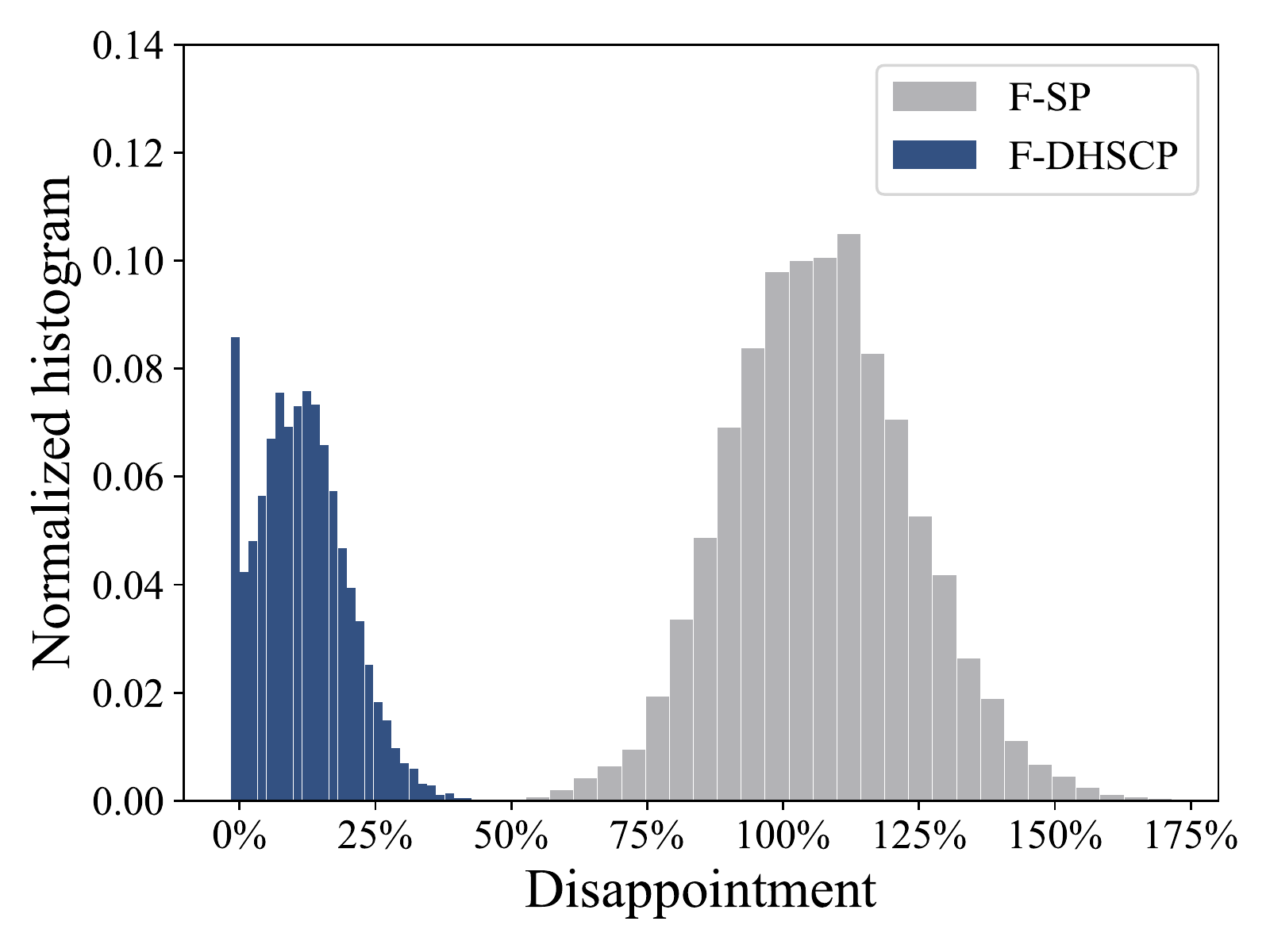}}
    \caption{Out-of-sample disappointment for Instance 7, demand range 1 under Set 2}
    \label{FA-oop-dis-4}
\end{figure}

\begin{figure}[!ht]
    \centering
    \subcaptionbox{$\Delta = 0$\label{FA-outc2-7-2-4-8-30-4060}}{
        \includegraphics[scale=0.33]{4Foutsamplecost2.pdf}}
    \subcaptionbox{$\Delta = 0.1$\label{FA-out01c2-7-2-4-8-30-4060}}{
        \includegraphics[scale=0.33]{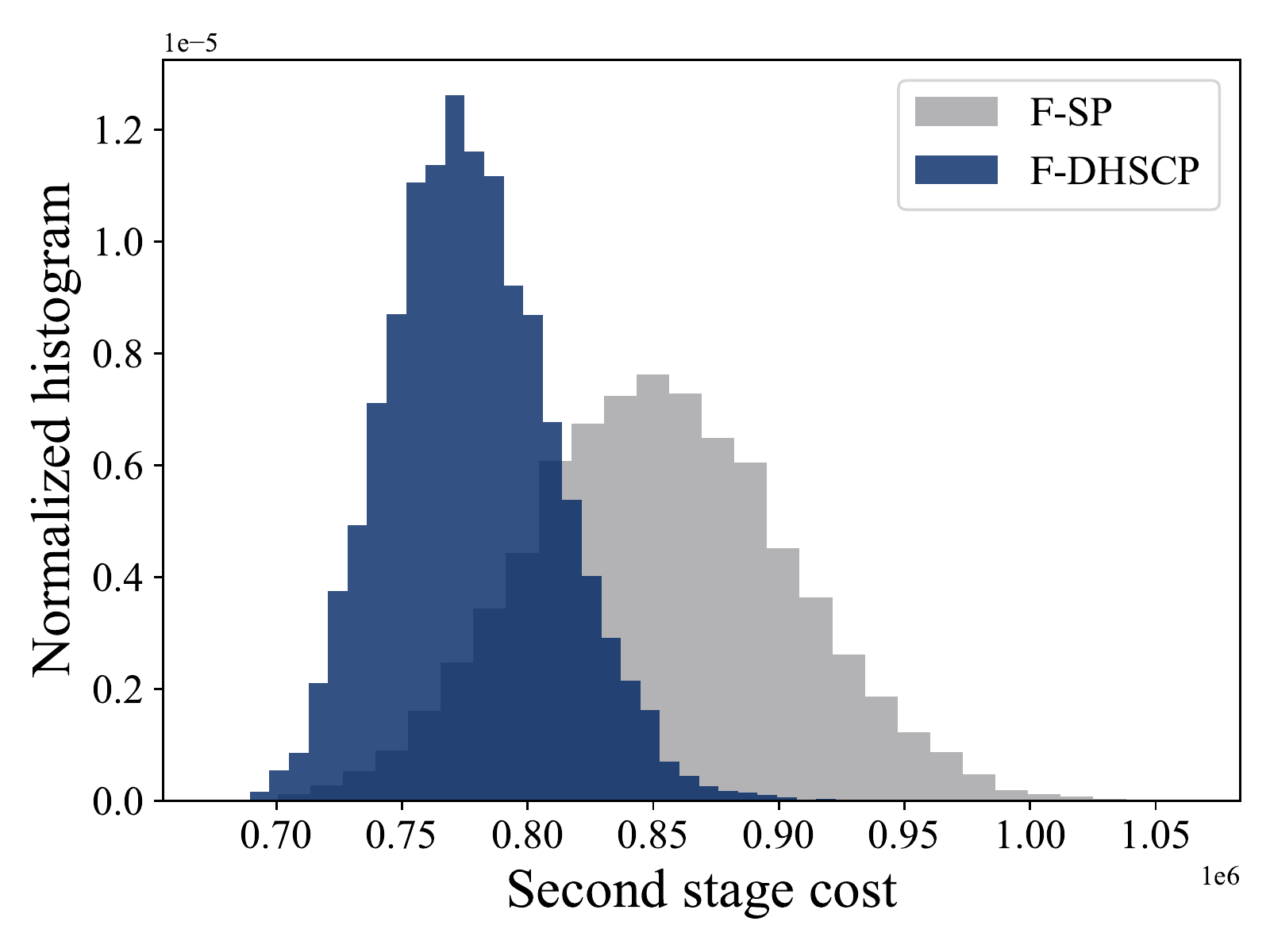}}
    \subcaptionbox{$\Delta = 0.25$\label{FA-out025c2-7-2-4-8-30-4060}}{
        \includegraphics[scale=0.33]{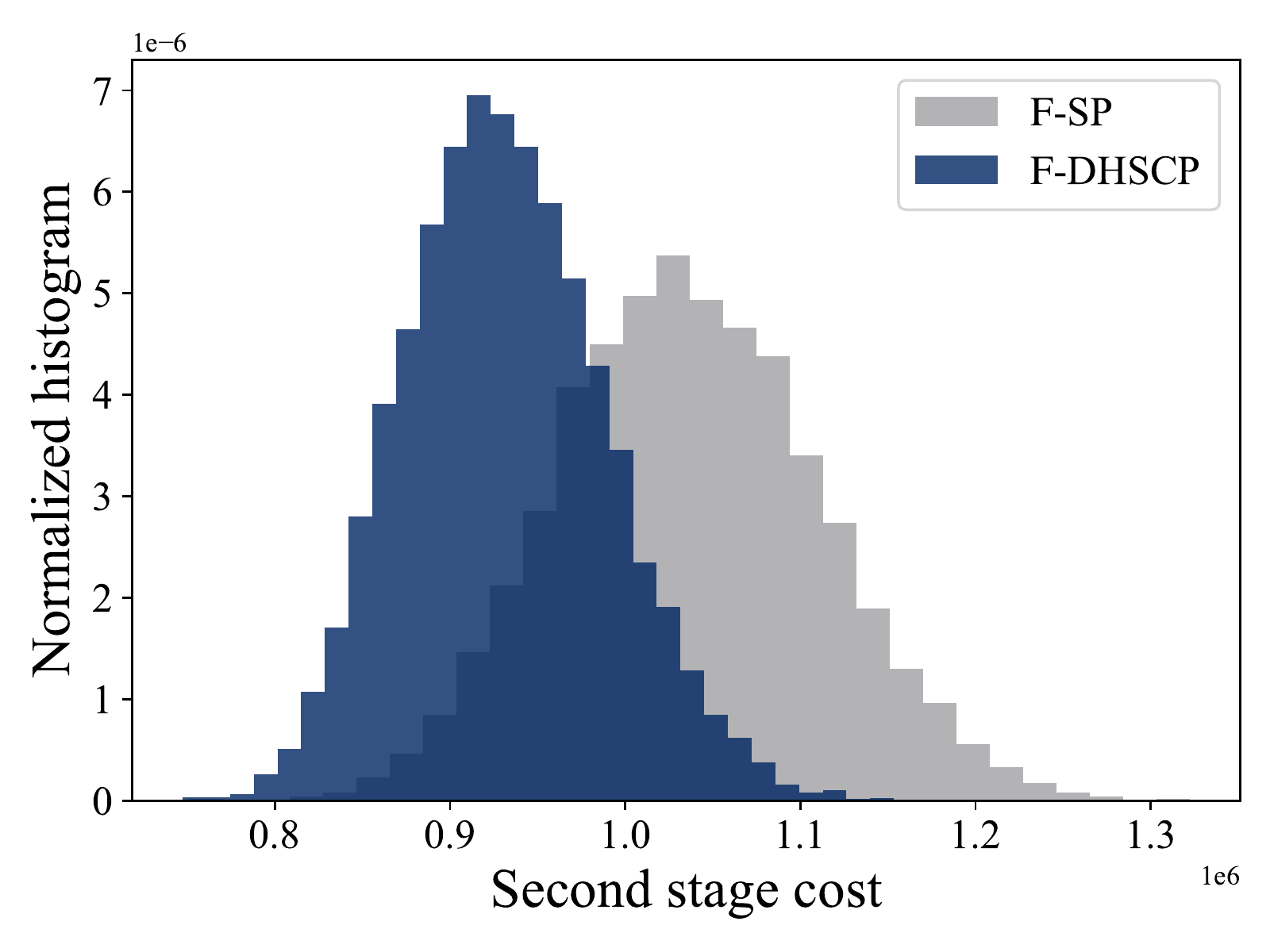}}
    \subcaptionbox{$\Delta = 0.5$\label{FA-out05c2-7-2-4-8-30-4060}}{
        \includegraphics[scale=0.33]{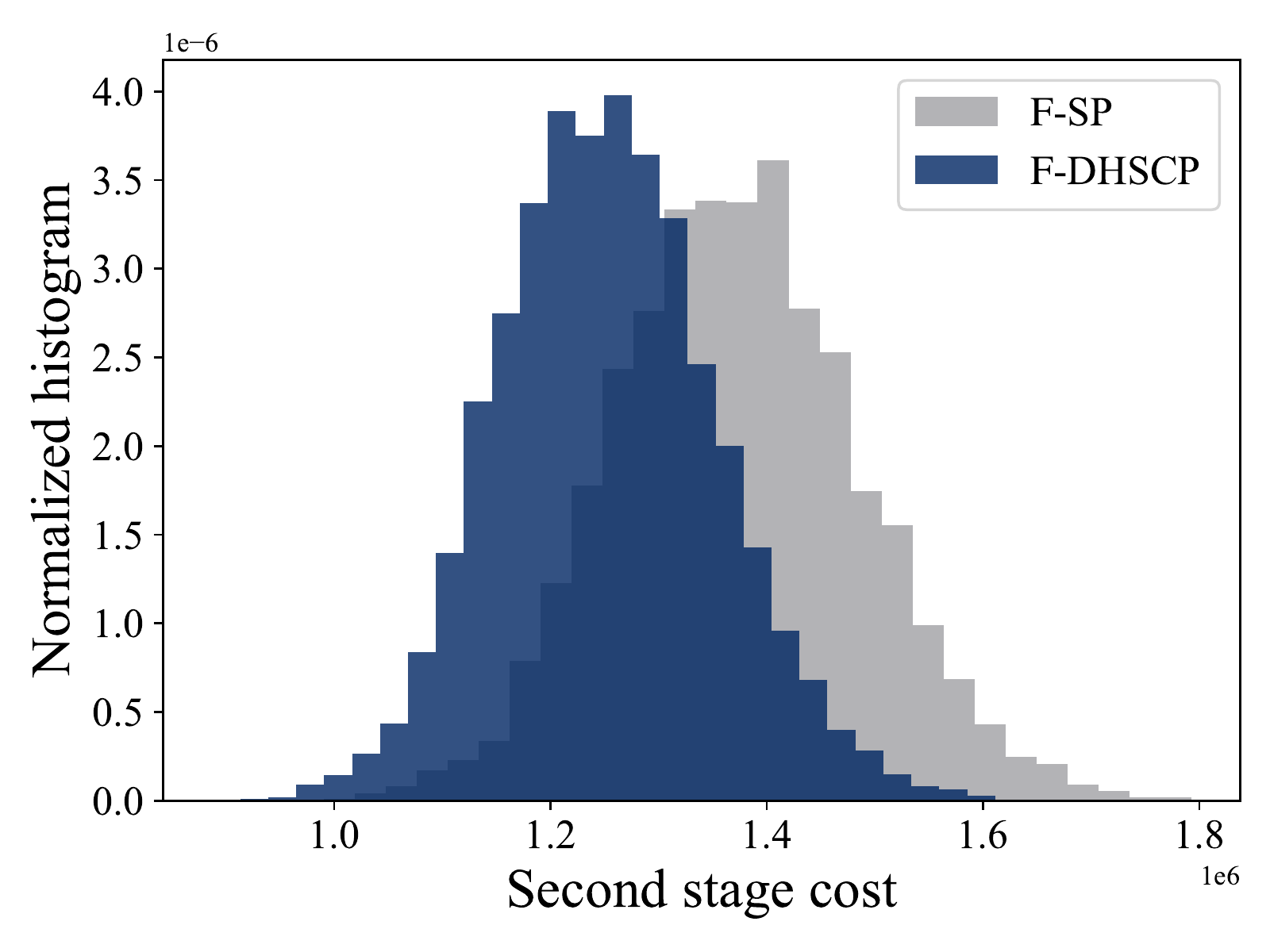}}
    \caption{Out-of-sample performance of FA models for Instance 7, demand range 1 under Set 2}
    \label{FA-oop-4}
\end{figure}

\begin{figure}[!ht]
    \centering
    \subcaptionbox{$\Delta = 0$\label{FA-outu-7-2-4-8-30-4060}}{
        \includegraphics[scale=0.35]{4Foutsampleunder.pdf}}
    \subcaptionbox{$\Delta = 0.1$\label{FA-out01u-7-2-4-8-30-4060}}{
        \includegraphics[scale=0.35]{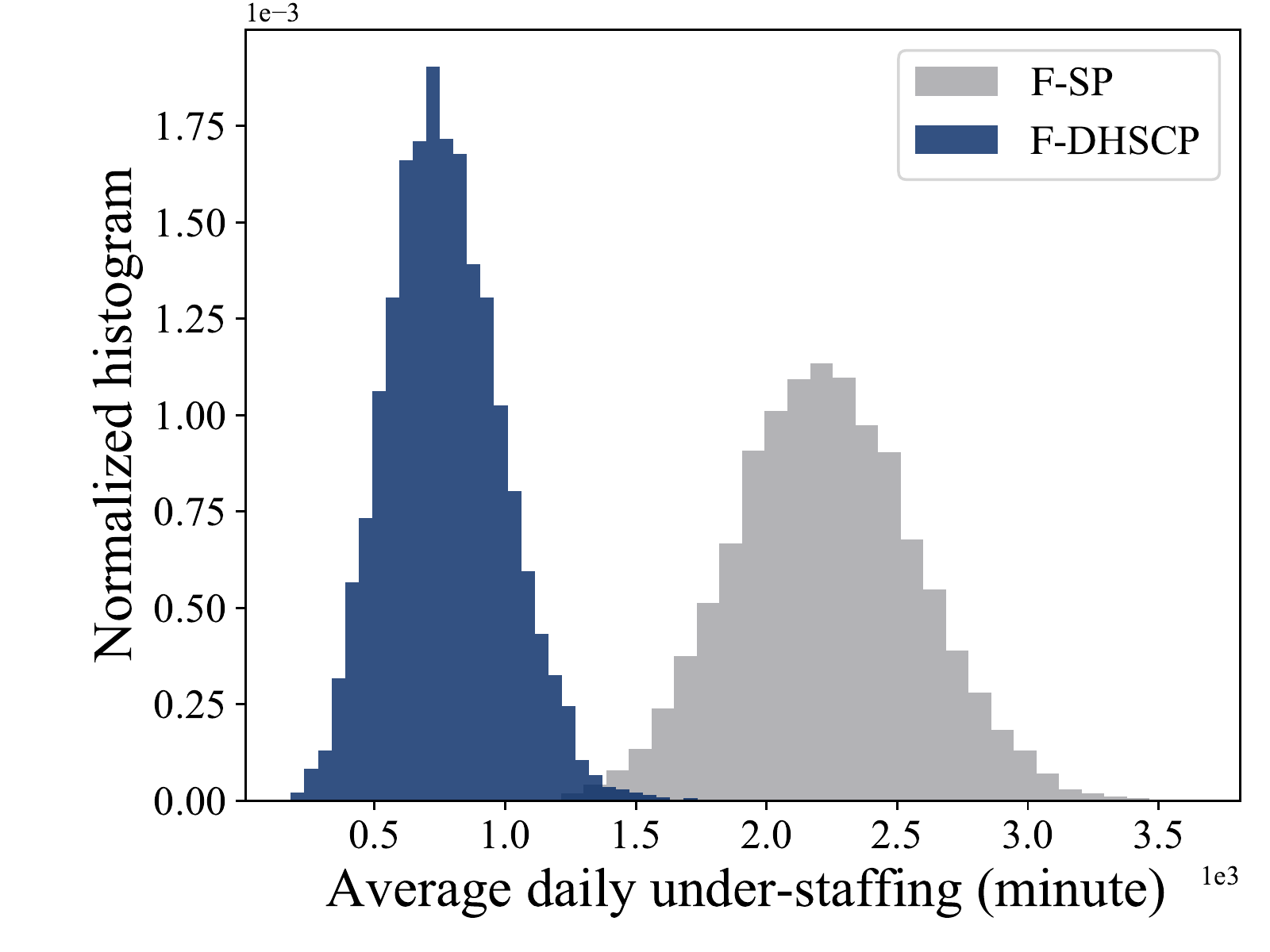}}
    \subcaptionbox{$\Delta = 0.25$\label{FA-out025u-7-2-4-8-30-4060}}{
        \includegraphics[scale=0.35]{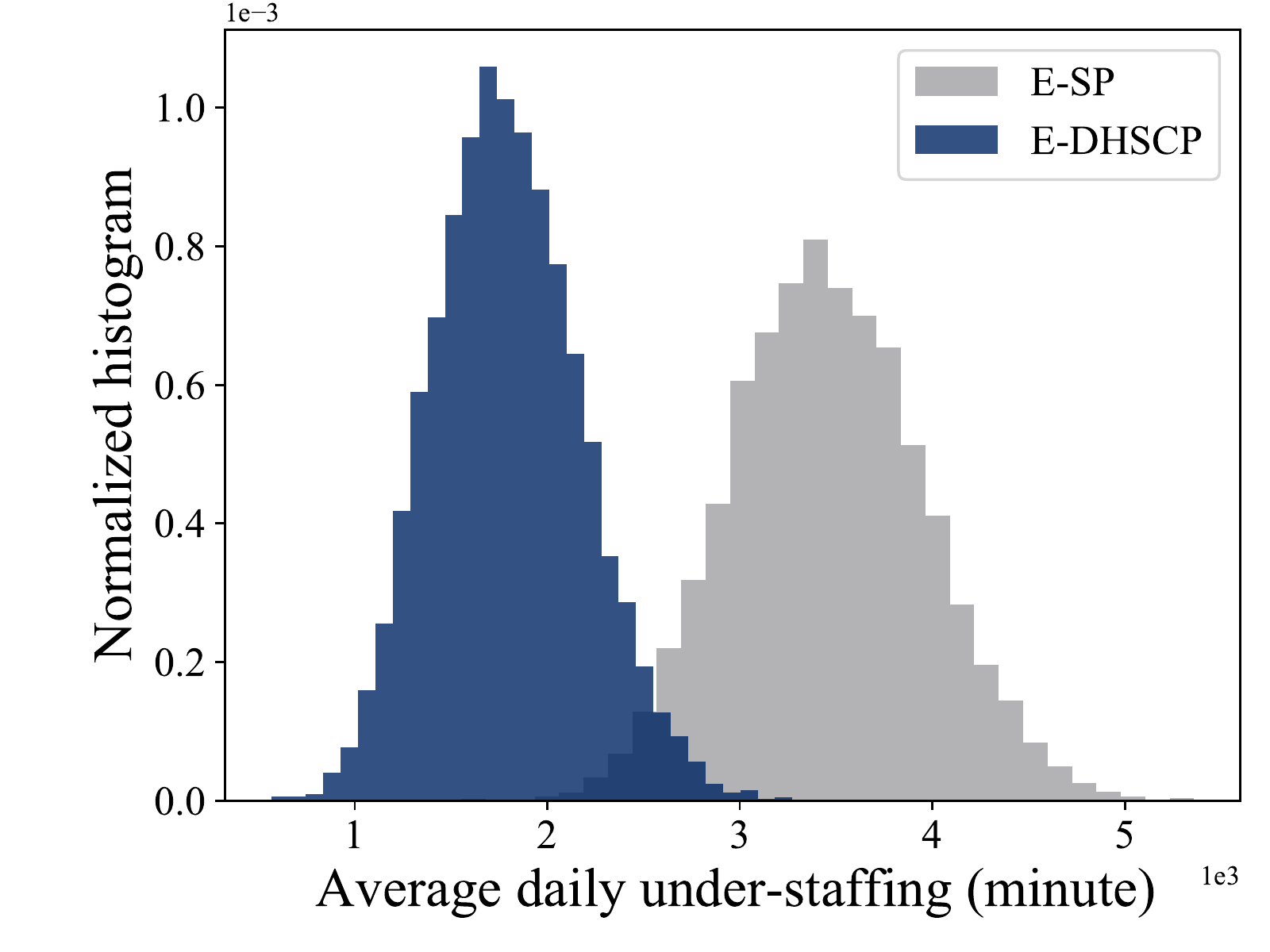}}
    \subcaptionbox{$\Delta = 0.5$\label{FA-out05u-7-2-4-8-30-4060}}{
        \includegraphics[scale=0.35]{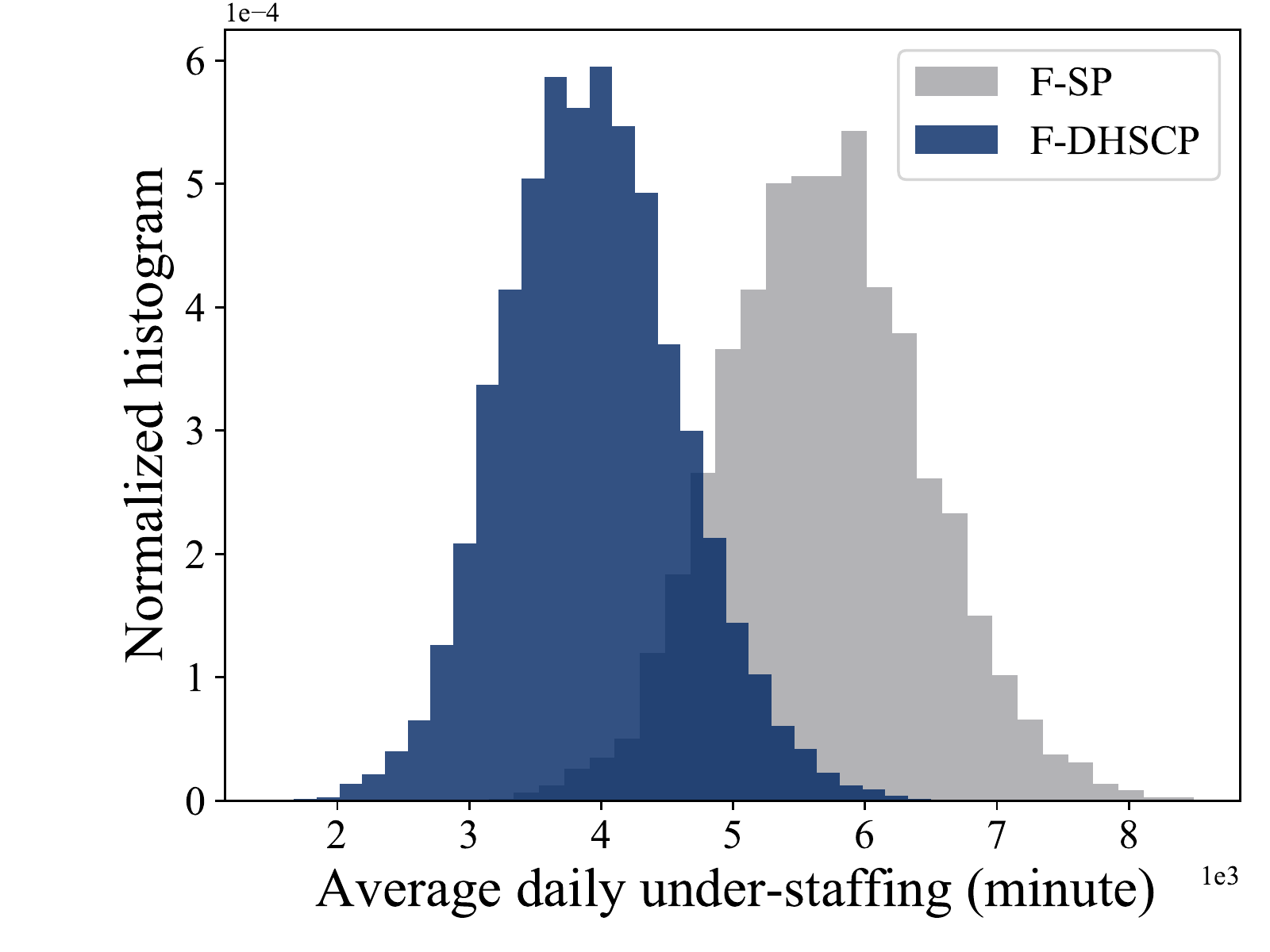}}
    \caption{Out-of-sample under-staffing for Instance 7, demand range 1 under Set 2}
    \label{FA-oop-u-4-4060}
\end{figure}

\end{document}